\input ./style/arxiv-general.cfg
\documentclass[aop,MSNbibl,seceqn,dvips]{arximspdf}
\makeatletter
   \@ifpackageloaded{graphicx}{}{\usepackage{graphicx}}
\makeatother
\usepackage{mathrsfs,mathbh}
\usepackage{url,breakurl}

\doi{10.1214/14-AOP962} % kopijuoti is laisko
\volume{43}
\issue{6}
\pubyear{2015}
\firstpage{3359}
\lastpage{3467}
\docsubty{FLA}

\makeatletter
\newcommand{\complementtt}{\complement}
\newcommand{\BES}{\mathrm{BES}}
\newcommand{\SBM}{\mathrm{SBM}}
\newcommand{\lleft}{\left}
\newcommand{\rrvert}{\vert}
\newcommand{\rright}{\right}
\newcommand{\rrVert}{\Vert}
\newcommand{\llvert}{\vert}
\newcommand{\llVert}{\Vert}
\renewcommand{\mid}{|}
\newcommand{\xrightarrow}[1]{\mathop{-\!\!\!-\!\!\!-\!\!\!\longrightarrow}_{#1}}
\newcommand{\xrightarroww}[2]{\mathop{-\!\!\!-\!\!\!-\!\!\!\longrightarrow}_{#1}^{#2}}
\newcommand{\1}{\mathbh{1}}
\newcommand{\G}{\mathscr G}
\newcommand{\C}{{\mathscr C}}
\newcommand{\F}{\mathscr F}
\newcommand{\B}{\mathscr B}
\newcommand{\R}{{\mathbb R}}
\newcommand{\lra}{\longrightarrow}
\newcommand{\loc}{\mathrm{loc}}
\newcommand{\lmt}{\longmapsto}
\newcommand{\Plim}{{\mathbb P\mbox{-}}}
\newcommand{\rest}{\upharpoonright}
\newcommand{\ind}{\perp\!\!\!\perp}
\newcommand{\supp}{\operatorname{supp}}
\newcommand{\bs}{\mathbf}
\newcommand{\ms}{\mathscr}
\renewcommand{\P}{{\mathbb P}}
\newcommand{\Q}{{\mathbb Q}}
\newcommand{\E}{{\mathbb E}}
\newcommand{\mc}{\mathcal}
\newcommand{\rap}{\mathrm{rap}}
\newcommand{\deq}{\stackrel{(\mathrm{d})}{=}}
\newcommand{\llang}{\langle\!\langle}
\newcommand{\rrang}{\rangle\!\rangle}
\newcommand{\llanga}{\bigl\langle\bigr\langle}
\newcommand{\rranga}{\bigl\rangle\bigr\rangle}

\newcommand{\propleq}{
\mathrel{\vcenter{\offinterlineskip
\hbox{$$}\vskip.8ex \hbox{$<$}\vskip-.6ex\hbox{$\smallfrown$}}}}
\newtheorem{teo}{Theorem}[section]
\newtheorem{teoi}{Theorem}
\newtheorem{prop}[teo]{Proposition}
\newtheorem{lem}[teo]{Lemma}
\newtheorem{cor}[teo]{Corollary}
\newproclaim{defi}[teo]{Definition}
\newproclaim{ass}[teo]{Assumption}
\newproclaim{rmk}[teo]{Remark}
\newproclaim{nota}[teo]{Notation}
\newproclaim{cas}{Choice of approximating solutions}
\newproclaim{ob}{Observation}
\makeatother

\begin{document}
\begin{frontmatter}

\title{Pathwise nonuniqueness for the SPDEs of some super-Brownian
motions with immigration\thanksref{F1}}
\runtitle{SPDEs of super-Brownian motions with immigration}

\begin{aug}
\author{\fnms{Yu-Ting} \snm{Chen}\corref{}\ead[label=e1]{yuting\_chen@cmsa.fas.harvard.edu}}
%\and
%\author{\fnms{}~\snm{}}
\runauthor{Y.-T. Chen}
\affiliation{Harvard University}
%\dedicated{}
\address{Center of Mathematical Sciences and Applications\\
Harvard University\\
Cambridge, Massachusetts 02138\\
USA\\
\printead{e1}}
\end{aug}
\thankstext{F1}{Supported in part by UBC Four Year Doctoral Fellowship
and CRM-ISM Postdoctoral Fellowship.}

% HISTORY:
\received{\smonth{8} \syear{2013}}
\revised{\smonth{7} \syear{2014}}
%\accepted{\smonth{} \syear{}}

% ABSTRACT
%
\begin{abstract}
We prove pathwise nonuniqueness in the stochastic partial differential
equations (SPDEs) for
some one-dimensional super-Brownian motions with immigration. In
contrast to a closely related
case investigated by Mueller, Mytnik and Perkins
[\textit{Ann. Probab.} \textbf{42} (2014) 2032--2112], the
solutions of the present
SPDEs are assumed to be nonnegative and have very different properties
including uniqueness in law. In proving possible separation of solutions,
we derive delicate properties
of certain correlated approximating solutions, which is based on a
novel coupling method called
continuous decomposition. In general, this method may be of independent
interest in furnishing solutions of
SPDEs with intrinsic adapted structure.
\end{abstract}

% KEYWORDS
\end{frontmatter}

%s1 #&#
\section{Introduction}\label{secintro}
In this work, we consider some one-dimensional super-Brownian motions
with (continuous) immigration, and
construct pairs of \mbox{distinct} \emph{nonnegative} solutions to the
associated stochastic partial differential equations (SPDEs). Hence, we
resolve in the negative the long-standing open problem concerning the
pathwise uniqueness in the SPDEs for one-dimensional super-Brownian
motions, when additional immigration is present (cf. page~217 in Perkins~\cite{PDW}).

We start with some informal descriptions
for the class of super-Brownian motions with immigration which are
considered throughout this work. See Dawson~\cite{DMMP}, Dynkin \cite
{DBMP}, Le Gall \cite{LGSBM}, Perkins \cite{PDW} and several
others for
super-Brownian motions as well as their connections with branching processes.
Imagine that, in the barren territory $\R$,
clouds of independent immigrants with infinitesimal initial mass land
randomly in space and throughout time. The underlying immigration
mechanism is time-homogeneous and gives a high intensity of arrivals of
immigrants so that the inter-landing times are infinitesimal. After
landing, each of the immigrant processes evolves independently of each
other as a super-Brownian motion, obeying the SPDE
%
%
%e1.1 #&#
\begin{eqnarray}
\label{eqmainSPDE0} \frac{\partial X}{\partial t} (x,t)&=& \frac{\Delta
X}{2}(x,t)+X(x,t)^{1/2}
\dot{W}(x,t),\qquad X\geq0,
\end{eqnarray}
subject to infinitesimal initial mass, where $W$ is (two-parameter)
space--time white noise on $\R\times\R_+$.
Superposition of their masses defines a super-Brownian motion with
immigration and zero initial value.
See Section~1.2 of Dawson \cite{DMMP}, Konno and Shiga \cite
{KSSPDE}, Section~III.4 of Perkins \cite{PDW} and Reimers \cite
{RSPDE} for the connection between solutions to the SPDE (\ref
{eqmainSPDE0}) and super-Brownian motions. See Sections~1.2~and~3.2 in Chen \cite{C} for some heuristic interpretations of the
terms of the SPDE (\ref{eqmainSPDE0}). Note that
super-Brownian motions with immigration can also be constructed by
Poisson point processes (see \cite{CD}).

We study the particular super-Brownian motions with immigration which
have densities, and the density processes obey the SPDEs:
%
%
%e1.2 #&#
\begin{eqnarray}
\label{eqmainSPDE} \frac{\partial X}{\partial t} (x,t)&=& \frac{\Delta
X}{2}(x,t)+\psi
(x)+X(x,t)^{1/2}\dot{W}(x,t),\qquad X\geq0,
\nonumber\\[-8pt]\\[-8pt]\nonumber
X(x,0)&=& 0.
\end{eqnarray}
Here, $\C_c^+(\R)$ being the function space of nonnegative continuous
functions on $\R$ with compact support, the \textit{immigration
functions} $\psi$ satisfy
%
%
%e1.3 #&#
\begin{equation}
\label{eqimmfunc} \psi\in\C_{c}^+(\R)\qquad\mbox{with }\psi\neq0
\end{equation}
and can be thought informally as the density of immigrants landing
within an infinitesimal amount of time.

To fix ideas, we give the precise definition of the pair $(X,W)$ in the
SPDE (\ref{eqmainSPDE}) before further discussions.
We need a filtration $(\G_t)$ which satisfies the usual conditions,
and it facilitates the following definitions of $W$ and $X$. We require
that $W$ be a $(\G_t)$-space--time white noise in the sense that it is
a family of $(\G_t)$-Brownian motions indexed by $L^2(\R)$, and the
Brownian motions satisfy the following properties: for any $d\in\mathbb
N$, $\phi_1,\ldots,\phi_d\in L^2(\R)$ and $a_1,\ldots,a_d\in\R$,
%
%
%e1.4 #&#
\begin{equation}
\label{WNlinear} W \Biggl(\sum_{j=1}^da_j
\phi_j \Biggr)=\sum_{j=1}^da_jW(
\phi_j)\qquad\mbox{a.s.}
\end{equation}
and $ (W(\phi_1),\ldots,W(\phi_d) )$ is a $d$-dimensional
$(\G_t)$-Brownian motion starting at zero with zero initial value and
covariance matrix $ [\langle\phi_i,\phi_j\rangle_{L^2(\R
)} ]_{1\leq i,j\leq d}$ (cf. Section~3 of Khoshnevisan \cite
{KSPDE} or Chapter~1 of Walsh \cite{WSPDE} for the standard
definition of space--time white noise). Since the immigration function
under consideration has compact support, it can be shown that the
density process of the corresponding super-Brownian motion with
immigration takes values in $\C_c^+(\R)$ (cf. Section~III.4 of
Perkins \cite{PDW}). Let
$\C_\rap(\R)$ denote the function space of rapidly decreasing
functions $f$:
%
%
%e1.5 #&#
\begin{equation}
\label{eqnorm-lambda} \llvert f\rrvert_\lambda\triangleq\sup
_{x\in\R} \bigl\llvert f(x)\bigr\rrvert e^{\lambda\llvert x\rrvert
}<\infty
\qquad\forall\lambda\in(0,\infty).
\end{equation}
Equip $\C_\rap(\R)$ with the complete separable metric
%
%
%e1.6 #&#
\begin{equation}
\label{eqCrap-norm} \llVert f\rrVert_\rap\triangleq\sum
_{\lambda=1}^\infty\frac{\llvert f\rrvert _\lambda
\wedge1}{2^\lambda}.
\end{equation}
For convenience, we follow the convention in Shiga \cite{SSPDE} and
use $\C_\rap(\R)$ as the underlying state space. Then by saying that
$X=(X_t)$ is a \textit{solution} to the SPDE~(\ref{eqmainSPDE}), we
require that $X$ be a nonnegative $(\G_t)$-adapted continuous process
with state space $\C_\rap(\R)$ and satisfy the following weak
formulation of (\ref{eqmainSPDE}):
%
%
%e1.7 #&#
\begin{eqnarray}
\label{distribution} \qquad X_t(\phi)&=&\int_0^t
X_s \biggl(\frac{\Delta\phi}{2} \biggr)\,ds+t\langle\psi,\phi\rangle+
\int_0^t\!\int_\R
X(x,s)^{1/2}\phi(x)\,dW(x,s)
\end{eqnarray}
for any test function $\phi\in\C_c^\infty(\R)$.
Here, we identify any locally integrable function $f$ on $\R$ as a
signed measure on $\B(\R)$ in the natural way and write
%
%
%e1.8 #&#
\begin{equation}
\label{signmeas} f(\phi)=\langle f,\phi\rangle\equiv\int_\R
f(x)\phi(x)\,dx,
\end{equation}
whenever there is no risk of confusion. For the last term in (\ref
{distribution}) and other two-parameter stochastic integrals in the
sequel, see Section~5 of Khoshnevisan \cite{KSPDE} or Chapter~2 of
Walsh \cite{WSPDE} for the construction.

A fundamental question for the SPDE (\ref{eqmainSPDE}) concerns its
uniqueness theory, and the major difficulty arises from the presence of
a non-Lipschitz diffusion coefficient. Uniqueness in law for the SPDE
(\ref{eqmainSPDE}) holds and can be proved by the duality method via
Laplace transforms (cf. Section~1.6 of Etheridge \cite{ESP} or the
proof of Lemma~\ref{lemred2}). In fact, it holds even if we impose
general nonnegative initial conditions for the SPDE (\ref
{eqmainSPDE0}) for super-Brownian motion and the SPDEs (\ref
{eqmainSPDE}) under consideration.
Nonetheless, duality methods for more general SPDEs of the form
%
%
%e1.9 #&#
\begin{equation}
\label{eqgenSPDE} \frac{\partial X}{\partial t} (x,t) = \frac{\Delta
X}{2}(x,t)+b \bigl(X(x,t)
\bigr)+\sigma\bigl(X(x,t) \bigr)\dot{W}(x,t)
\end{equation}
up to now seem only available when $b$ and $\sigma$ are of rather
special forms, and
hence are nonrobust. (See\vspace*{1pt} Mytnik \cite{MUL} for the duality method
for the case $b=0$ and $\sigma(x)=x^p$, where $p\in(\frac{1}{2},1)$
and nonnegative solutions are assumed.) After all, duality is based on
exactness and may become difficult to obtain by even slight changes of
coefficients in the context of SPDEs.

Under the classical theory of stochastic differential equations (SDEs),
uniqueness in law in an SDE is a consequence of
pathwise uniqueness of its solutions (cf. Theorem IX.1.7 of Revuz and
Yor \cite{RYCMB}). The strength of this point of view is that it has
emphasis on the \emph{values} of the H\"older exponents of
coefficients, instead of on the particular forms of coefficients. Then
a natural question is whether the duality method can be circumvented by
proving pathwise uniqueness in the SPDEs (\ref{eqgenSPDE}) instead.
Here, \textit{pathwise uniqueness} in an SPDE ensures that any two
solutions subject to the same space--time white noise and initial value
always coincide almost surely. Our objective in the present work is to
settle the question of pathwise uniqueness in the particular SPDEs
(\ref{eqmainSPDE}).

Let us discuss some results on pathwise uniqueness in various SDEs and
SPDEs which are closely related to the SPDEs (\ref{eqmainSPDE}).
We focus on the role of non-Lipschitz diffusion coefficients in
determining pathwise uniqueness.

For one-dimensional SDEs with H\"older-$p$ diffusion coefficients, the
famous Girsanov example (see Section V.26 of Rogers and Williams \cite
{RWDMM}) shows
the necessity of the condition $p\geq\frac{1}{2}$ for pathwise
uniqueness of solutions. The sufficiency was later confirmed
in the seminal work Yamada and Watanabe \cite{YWSDE} as far as the
cases with sufficiently regular drift coefficients are concerned.
In fact, the work \cite{YWSDE} showed that a finite-dimensional SDE
defined by
%
%
%e1.10 #&#
\begin{equation}
\label{eqYWSDE} dX^i_t=b_i(X_t)
\,dt+\sigma_i\bigl(X^i_t\bigr)\,dB^i_t,
\qquad1\leq i\leq d,
\end{equation}
enjoys pathwise uniqueness
as long as all $b_i$'s are Lipschitz continuous and each $\sigma_i$ is
H\"older $p$-continuous, for any $p\geq\frac{1}{2}$.

In view of the complete results for SDEs and the strong parallels
between (\ref{eqgenSPDE}) and (\ref{eqYWSDE}), it had been hoped
for decades that pathwise uniqueness would also hold in (\ref
{eqgenSPDE}) if the diffusion coefficient $\sigma$ is H\"older-$p$
continuous whenever $p\geq\frac{1}{2}$.
It was shown in Mytnik and Perkins \cite{MP11} that this is the case
if $\sigma$ is H\"older-$p$ for $p>\frac{3}{4}$, but in Burdzy,
Mueller and Perkins \cite{BMP} and Mueller, Mytnik and Perkins \cite
{MMP} that pathwise uniqueness in
%
%
%e1.11 #&#
\begin{eqnarray}
\label{eqmainSPDE+1} \frac{\partial X}{\partial t} (x,t)&=& \frac{\Delta
X}{2}(x,t)+\bigl\llvert
X(x,t)\bigr\rrvert^{p}\dot{W}(x,t),
\nonumber\\[-8pt]\\[-8pt]\nonumber
X(x,0)&=& 0
\end{eqnarray}
fails for any $p\in(0,\frac{3}{4})$. Here, a nonzero solution to
(\ref{eqmainSPDE+1}) exists and, as $0$ is obviously another
solution, both pathwise uniqueness and uniqueness in law fail. All
these results point to the general conclusion that pathwise uniqueness
of solutions holds for H\"older-$p$ diffusion coefficients $\sigma$
for $p>\frac{3}{4}$ but can fail for $p\in(0,\frac{3}{4})$. See also
Mytnik, Perkins and Sturm \cite{MPS} for the case of colored noises.

In this work, we confirm pathwise \emph{nonuniqueness} in the SPDEs
(\ref{eqmainSPDE}). We stress that by definition, only \emph
{nonnegative} solutions are considered in this regard and hence are
unique in law by the duality argument mentioned above.
Our main result is given by the following theorem.

%
%
%th1 #&#
\begin{teoi}[(Pathwise nonuniqueness)]\label{teoSPDE}
For any nonzero immigration function $\psi\in\C_c^+(\R)$, there
exists a filtered probability space $(\Omega,\F,(\G_t),\P)$ which carries
a $(\G_t)$-space--time white noise $W$ and two solutions $X$ and $Y$ of
the SPDE~(\ref{eqmainSPDE}) with respect to $(\G_t)$
such that $\P(X\neq Y )>0$. Hence, there is pathwise
nonuniqueness in the SPDE (\ref{eqmainSPDE}) for $\psi$ as above.
\end{teoi}

A comparison of diffusion coefficients may suggest that
the construction in Mueller, Mytnik and Perkins \cite{MMP} of a
nonzero signed solution to the particular case
%
%
%e1.12 #&#
\begin{eqnarray}
\label{eqmainSPDEpm} \frac{\partial X}{\partial t} (x,t)&=& \frac{\Delta
X}{2}(x,t)+\bigl\llvert
X(x,t)\bigr\rrvert^{1/2}\dot{W}(x,t),
\nonumber\\[-8pt]\\[-8pt]\nonumber
X(x,0)&=& 0
\end{eqnarray}
for (\ref{eqmainSPDE+1})
should be closely related to our case (\ref{eqmainSPDE}).
Nonetheless, solutions to (\ref{eqmainSPDE}) are subject to
the assumed nonnegativity, and uniqueness in law in the SPDEs (\ref
{eqmainSPDE}) \emph{does hold}. These facts mean that
the goal will be to find two nonzero solutions which have the same law and
are nontrivially correlated through the
shared white noise. Although many features of our arguments will follow
their counterparts in~\cite{MMP}, a number of new problems, including
the choice of solutions to work with, arise in dealing with these
distinct properties.

For a fixed nonzero immigration function $\psi\in\C_c^+(\R)$, we construct
the pair of distinct solutions to the corresponding SPDE (\ref
{eqmainSPDE}) by approximation. Basic properties of the approximating
solutions are as follows.
An ${\varepsilon}$-approximating pair, still denoted by
$(X,Y)$ but under $\P_{\varepsilon}$, consists of super-Brownian
motions with \textit{intermittent immigration} and subject to the \emph
{same} space--time white noise.
Here, a~super-Brownian motion with intermittent immigration is defined
as a discrete sum of certain immigrant processes. The immigrants land
after intervals of deterministic and equal length and at i.i.d.
targets, and then, along with their offspring, evolve independently as
true super-Brownian motions.
In more detail, the pairs $(X,Y)$ satisfy the following properties.
The initial masses of the immigrant processes are of the form $\psi(\1
)J^a_{\varepsilon}( \cdot)$ with $a$ denoting the target, where
%
%
%e1.13 #&#
\begin{equation}
\label{eqJxvep} J^a_{\varepsilon}(x)\equiv{\varepsilon}^{1/2}J
\bigl((a-x ){\varepsilon}^{-1/2} \bigr),\qquad x\in\R,
\end{equation}
for a fixed even $\C_c^+(\R)$-function $J$ which is bounded by $1$,
has topological support contained in $[-1,1]$, and satisfies $\int_\R
J(x)\,dx=1$. In addition,
the landing times of the immigrants are interlaced as
%
%
%e1.14 #&#
\begin{eqnarray}
\label{sttime} s_i&=& \bigl(i-\tfrac{1}{2} \bigr){\varepsilon}
\quad\mbox{and}\quad t_i=i{\varepsilon}\qquad\mbox{for }i\in\mathbb N,
\end{eqnarray}
and the targets associated with the immigrants of $X$ and $Y$ are given
by i.i.d. spatial variables $x_i$ and $y_i$ at $s_i$ and $t_i$,
respectively, where
%
%
%e1.15 #&#
\begin{equation}
\label{eqxiyilaw} \P_{\varepsilon}(x_i\in dx) =\P_{\varepsilon}(y_i
\in dx)\equiv\frac{\psi(x)\,dx}{\psi(\1)}.
\end{equation}
Details of these approximating solutions and their convergence to true
solutions of the SPDE (\ref{eqmainSPDE}) can be found in Sections~\ref{secISBM}~and~\ref{secpropwc}.

At this point, we only describe two ${\varepsilon}$-approximating
solutions which share the same space--time white noise, and what can be
deduced from this relation in understanding their interactions seems
limited. The perspective of the present work
is to emphasize the role of
immigrant processes, and the readers will see that they lead to a
detailed comparison of local masses for the approximating pairs in
particular. On the other hand, by adopting this point of view,
we are faced with an issue of defining approximating solutions by appropriate
immigrant processes, as will be discussed in more detail later on.

We notice that similar ${\varepsilon}$-approximating solutions appear
in Mueller, Mytnik and Perkins \cite{MMP} for the construction of a
nonzero solution to the SPDE (\ref{eqmainSPDEpm}). In this case,
each approximating solution
is obtained by specifying its ``positive part'' and ``negative part''
as two super-Brownian motions with intermittent immigration, but now
subject to pairwise annihilation upon collision. Both parts are in turn
defined by sums of their own immigrant processes undergoing
annihilation, and all of the summands can be seen as
i.i.d. super-Brownian immigrant processes taken off annihilated
individuals and their possible offspring (cf. equations~(2.1), (2.4) and
(2.6) of \cite{MMP} and Lemma~\ref{lemred2}). The latter property
implies fairly explicit stochastic calculus for the immigrant
processes, and is the key to
make further analysis possible in \cite{MMP}.

For our case, while a super-Brownian motion with intermittent
immigration can be defined as a sum of independent immigrant processes,
the question for the same construction of two, with interlacing
immigrating times and subject to the \emph{same} space--time white
noise, lies in the interactions between immigrants through space--time
white noises.
The major difficulty here is in specifying a family of correlated
immigrants so that the corresponding approximating solutions not only conform
to the same space--time white noise but also generate two distinct
solutions to the SPDE (\ref{eqmainSPDE}).
After all, in contrast to the counterexample in Mueller, Mytnik and
Perkins \cite{MMP} which stems from annihilation of colliding
individuals in two independent population processes, it is still not
known whether a similar interpretation applies to the SPDE (\ref
{eqmainSPDE}) under consideration, since in our case the existence of
two different solutions means a special kind of coexistence of two
population processes.
We need a different point of view to choose approximating solutions.

The correlated immigrant processes which meet our needs are chosen
through a reverse analysis for the ${\varepsilon}$-approximating solutions.
The aim is to find immigrant processes subject a ``tractable''
correlation structure for every coupled pair of approximating solutions
(readers interested in more details about the motivation may see
Section~3.2 of \cite{C} for a nonrigorous proof of our main result).
Our main machinery is a novel coupling method called \textit{continuous
decomposition}. By this method, essentially we can
\emph{elicit} certain immigrant processes from any pair of
${\varepsilon}$-approximating solutions so that the integrals of $\C
_c^\infty(\R)$-functions against them
define continuous semimartingales starting at the respective landing
times, all with respect to the \emph{same} filtration (see
Theorem~\ref{propequiallocation}). Here, $\C_c^\infty(\R)$ denotes
the space of infinitely differentiable functions on $\R$ with compact
support. We remark that this semimartingale property of immigrant
processes does not follow directly from the general theory of coupling
(see Section~\ref{seccont} for a discussion).

The immigrant processes from continuous decomposition satisfy natural
distributional properties including the one that both families of
immigrants, say $\{X^i\}$ and $\{Y^i\}$ for $X$ and $Y$, respectively,
consist of independent processes. Moreover, they can be chosen such
that for all $\phi,\varphi\in\C_c^\infty(\R)$, the ``coarse''
\mbox{(predictable)} covariations
$\langle X^i(\phi),Y(\varphi)\rangle$ and $\langle X(\phi
),Y^j(\varphi)\rangle$, rather than the covariations $\langle
X^i(\phi),Y^j(\varphi)\rangle$ between immigrants,
admit explicit expressions (see Proposition~\ref{propcovar}). The
expressions are simple enough for one to conjecture that the immigrant
processes should satisfy the \textit{coexistence condition}:
%
%
%e1.16 #&#
%e1.17 #&#
\begin{eqnarray}\label{keycovar}
\bigl\langle X^i(\phi),Y^j(\varphi)\bigr
\rangle_t=\int_0^t\!\int
_\R\frac
{X^i(x,s)Y^j(x,s)}{X(x,s)^{1/2}Y(x,s)^{1/2}}\phi(x)\varphi
(x)\,dx\,ds,
\nonumber\\[-10pt]\\[-10pt]
\eqntext{i,j\geq1,}
\end{eqnarray}
where $0/0$ is read as $0$. Note that the covariations in (\ref
{keycovar}) are given by \emph{two-fold} integrals. See Section~\ref
{seccovar-imm} for other possibilities of covariations for stochastic
integrals with respect to space--time white noises, and also Remark~\ref
{rmkmistake} on related issues. By classical arguments, it can be
verified that immigrant processes subject to the coexistence condition
(\ref{keycovar}) do exist. See Theorem~\ref{teoapprox-v2} for the
precise result. We will restrict our attention to the corresponding
${\varepsilon}$-approximating solutions in the remaining of
Section~\ref{secintro}.

Let us explain why the ${\varepsilon}$-approximating solutions remain
separated if we pass ${\varepsilon}$ to zero. We
switch to the conditional probability measure under which the total
mass process of a generic immigrant, say $X^i$, hits $1$.
Let us call such a conditional probability measure $\Q^i_{\varepsilon
}$ from now on.
The motivation to invoke these conditional probabilities is that with
high $\P_{\varepsilon}$-probability, there is at least one immigrant
from $X$ whose total mass will hit $1$ by the independence of the
immigrants for $X$, so whenever $X$ and $Y$ are separated with
sufficiently high probability under every~$\mathbb Q^i_{\varepsilon}$, we
should be able to integrate these immigrant-wise phenomena of \textit
{conditional separation} of $X$ and $Y$ into a kind of separation under
$\P_{\varepsilon}$.

The readers may notice that the above argument to obtain separation
under $\P_{\varepsilon}$ is reminiscent of the use of excursion
theory in studying
pathwise uniqueness in SDEs and SPDEs (cf. Bass, Burdzy and Chen \cite
{BBC} and Burdzy, Mueller and Perkins \cite{BMP}). The major
difference, however, is that in the present case, the immigrant
processes can overlap
in time without waiting until the earlier ones die out. In order to use
conditional separation of the approximating pairs, we resort to an
inclusion--exclusion argument as in Mueller, Mytnik and Perkins \cite
{MMP}. The result is \textit{uniform separation} of the approximating
pairs under $\P_{\varepsilon}$. It states that for some constants
$T,\Delta\in(0,\infty)$ independent of ${\varepsilon}$, $\sup_{0\leq
s\leq T}\llVert X_s-Y_s\rrVert _\rap$ under $\P_{\varepsilon}$ are
uniformly bounded below by $\Delta$ with uniformly positive
probability for all small ${\varepsilon}\in(0,1)$, or more precisely
\[
\liminf_{{\varepsilon}\searrow0}\P_{\varepsilon} \Bigl(\sup
_{0\leq s\leq T}\llVert X_s-Y_s\rrVert
_\rap\geq\Delta\Bigr)>0.
\]
Then it is not difficult to argue that any two true solutions to (\ref
{eqmainSPDE}) as a limit of our approximation pairs separate with
strictly positive probability.
See Section~\ref{secsepSOL} for the details.

The conditional separation under $\mathbb Q^i_{\varepsilon}$ of the two
approximating solutions concerns a comparison of their local masses over
a growing space--time region.
We envelop the support processes
of $X^i$ and $Y^j$ by \emph{approximating} parabolas of the form
%
%
%e1.18 #&#
\begin{equation}
\label{parab-0} \mc P^{(a,s)}_{\beta}(t)= \bigl\{(x,r)\in\R
\times[s,t];\llvert a-x\rrvert\leq{\varepsilon}^{1/2}+(r-s)^\beta
\bigr\}
\end{equation}
for $\beta$ near $1/2$
and consider the propagation of these parabolas
instead of that of the support processes. The known modulus of
continuity for the support of super-Brownian motion
implies that, for example, the support of $X^i$ satisfies
%
%
%e1.19 #&#
\begin{equation}
\label{suppXi} \operatorname{supp}\bigl(X^i\bigr)\cap\bigl(\R
\times[s_i,t] \bigr)\subseteq\mc P_\beta^{(x_i,s_i)}(t)
\qquad\mbox{for }t-s_i\mbox{ small},
\end{equation}
where $\operatorname{supp}(X^i)$ is the space--time support of the
random function
$(x,s)\lmt X^i(x,s)$,
and $x_i$ and $s_i$ denote the landing target and landing time of
$X^i$, respectively
(see Section~\ref{sec2ndsep-setup} and Proposition~\ref{propsupppre}).
As in Mueller, Mytnik and Perkins \cite{MMP}, the total mass process
$X^i(\1)$ under $\mathbb Q^i_{\varepsilon}$ can be shown to be a constant
multiple of a $4$-dimensional Bessel squared process near its landing
time and hence has a known growth rate. Thanks to (\ref{suppXi}), this
growth rate of the total mass is the same as the growth rate of the
local mass of $X^i$ over its support envelope, and then a lower bound
of the associated local mass of $X$ follows from
the nonnegativity of the immigrant processes.

We prove that the local mass of $Y$ over the envelope for $X^i$ grows
at a smaller rate.
This involves a subcollection of immigrants from $Y$ which we choose
now. The $\Q^i_{\varepsilon}$-probability that one of the $Y^j$
clusters preceding $X^i$ ever invades the ``territory'' of $X^i$ by
time $t\in(s_i,\infty)$ can be made relatively small as long as
$t-s_i$ small, which follows from an argument similar to the proof of
Lemma~8.4 of Mueller, Mytnik and Perkins \cite{MMP} (see
Proposition~\ref{propoutclus-si}).
These $Y^j$ clusters can henceforth be excluded from our consideration.
Then the simple geometry of the approximating parabolas (\ref
{parab-0}) yields the space--time rectangles
\[
\mc R^i(t)= \bigl[x_i-2 \bigl({\varepsilon}^{1/2}+(t-s_i)^\beta
\bigr),x_i+2 \bigl({\varepsilon}^{1/2}+(t-s_i)^\beta
\bigr) \bigr]\times[s_i,t]
\]
so that the immigrant processes $Y^j$ landing inside $\mc R^i(t)$ are
the only possible invaders of the support envelope for $X^i$ by time
$t$. This results in a family of clusters, say, $\{Y^j;j\in\mc J^i(t)\}
$ to the effect that the local mass of $Y$ over the growing envelope
for $X^i$ is dominated by the sum of total masses of these clusters. We
further classify them into \textit{critical clusters} and \textit{lateral
clusters}. In essence, the critical clusters land near the territory of
$X^i$ so the interactions between these clusters and $X^i$ are
significant. In contrast, the lateral clusters must evolve for
relatively larger amounts of time before they start to interfere with $X^i$.

Up to this point, the framework we set for investigating conditional
separation of approximating solutions is very similar to that in
Mueller, Mytnik and Perkins \cite{MMP}. The \emph{interactions}
between the approximating solutions considered in both cases are,
however, very different in nature. For example, bounding the finite
variation process of the semimartingale $Y^j(\1)$ under $\Q
^i_{\varepsilon}$ is the main source of difficulty in our case, but
this creates no difficulty in \cite{MMP}.
Hence, our case calls for a new analysis again. Our result for the
conditional separation can be captured quantitatively by saying that
for arbitrarily small $\delta>0$,
%
%
%e1.20 #&#
\begin{eqnarray}\label{statcondsep+}
&& \mbox{with high $\Q^i_{\varepsilon}$-probability, }
X^i_t(\1) \geq \operatorname{constant}\,\cdot\,(t-s_i)^{1+\delta}\mbox{ and}
\nonumber\\[-11pt]\\[-12pt]\nonumber
&&\textstyle\sum_{j\in\mathcal J^i(t)}Y^j_t(\1)\leq\operatorname{constant}\,\cdot\,
(t-s_i)^{{3}/{2}-\delta},
\mbox{ for $t$ close to $s_i+$.}
\end{eqnarray}
Here, the initial behavior of $X^i(\1)$ under $\Q^i_{\varepsilon}$ as
a constant multiple of a \mbox{$4$-}dimensional Bessel squared process
readily gives the first part of (\ref{statcondsep+}) (see
Section~\ref{sectauvep}).
On the other hand, the extra order, which is roughly $(t-s_i)^{1/2}$,
for the sum of the (potential) invaders $Y^j$ can be seen as the result
of spatial structure.

In fact, the above framework needs to be further modified in a critical
way due to a technical difficulty which arises in our setting (but not
in Mueller, Mytnik and Perkins \cite{MMP}). We must consider a
slightly looser definition for critical clusters, and a slightly more
stringent definition for lateral clusters. It will be convenient to
consider this modified classification for the $Y^j$ clusters, still
indexed by $j\in\mathcal J^i(t)$,
landing inside a slightly larger rectangle in place of $\mc R^i(t)$. Write
\[
\mc J^i(t)=\mc C^i(t)\cup\mc L^i(t),
\]
where $\mc C^i(t)$ and $\mc L^i(t)$ are the random index sets
associated with critical clusters and lateral clusters, respectively.
See Section~\ref{sec2ndsep-setup} for the precise classification.

Let us discuss the method to bound the sum of the total masses $Y^j_t(\1
)$, $j\in\mathcal J^i(t)$,
under $\Q^i_{\varepsilon}$ [recall (\ref{statcondsep+})].
As in Mueller, Mytnik and Perkins \cite{MMP}, this part plays a major
role in the present work besides the selection of approximating solutions.
The treatment of the sum
is through an analysis of its first moment, or more precisely an
analysis of the expected finite variation process of $Y^j(\1)$ under
$\mathbb Q^i_{\varepsilon}$ for $j\in\mc J^i(t)$.

For the critical clusters $Y^j$, the finite variation processes of
their total masses under $\mathbb Q^i_{\varepsilon}$ have bounds given by
%
%
%e1.21 #&#
\begin{equation}
\label{qv} \int_{t_j}^t \frac{[Y^j_s(\1)]^{1/2}}{[X^i_s(\1)]^{1/2}}\,ds
\end{equation}
for $t$ sufficiently close to $t_j+$ (cf. Lemma~\ref{lemsep1} below),
so only the \emph{total masses} of the clusters need to be handled.
In this direction, we use an \textit{improved modulus of continuity} of
the total mass processes $Y^j(\1)$ and
the lower bound of $X^i(\1)$ in~(\ref{statcondsep+}) to give
deterministic bounds for the integrands in (\ref{qv}).
The overall effect is a bound for the expected sum of the total masses
$Y^j_t(\1)$, $j\in\mc C^i(t)$, and this can be used to show that the
corresponding random sum
has growth similar to that in the second part of (\ref
{statcondsep+}). See Section~\ref{seccondsep-1}.

The lateral clusters pose an additional difficulty here which is not present
in Mueller, Mytnik and Perkins \cite{MMP} due to the possibly
nontrivial covariations between these clusters and $X^i$. The question
is whether or not conditioning on $X^i$ being significant can pull
along the nearby $Y^j$'s at a greater rate, even though any of these
$Y^j$ does not interfere with $X^i$ upon their landing. In order to
help bound the contributions of these clusters, we
argue that a lateral cluster $Y^j$ is independent of $X^i$ until they
collide (cf. Lemma~\ref{lemoc} and Proposition~\ref{propocXY}).
This allows us to adapt the arguments for the critical clusters and
furthermore bound the growth rate of the sums of the total masses
$Y^j_t(\1)$, $j\in\mc L^i(t)$, by the desired order.
See the discussion in Section~\ref{seccondsep-2} for more on
this issue.

We close our discussion in this section with an immediate corollary for
the SPDE (\ref{eqmainSPDE}) in which the immigration function has
small total mass $\psi(\1)$ and the initial value is replaced by a
nonzero nonnegative $\C_\rap(\R)$-function. In this case, pathwise
nonuniqueness remains true. This
follows from the Markov property of super-Brownian motions with
immigration and the recurrence of Bessel squared processes with small
dimensions (cf. page~442 in Revuz and Yor \cite{RYCMB}). In detail,
we can run a copy of such a super-Brownian motion with immigration
until its total mass first hits zero, and then the required distinct
solutions can be obtained by concatenating this piece with the
separating solutions in Theorem~\ref{teoSPDE}.

This paper is organized as follows. In Section~\ref{secISBM}, we give
the precise definition of the pairs of approximating solutions from
which we choose particular ones for the proof of our main result, and
discuss their basic properties.
In Section~\ref{seccont}, we explain the idea of continuous
decomposition of a super-Brownian motion with intermittent immigration
and then give the rigorous proof for the continuous decompositions of
the approximation solutions specified in Section~\ref{secISBM}.
Covariations of the resulting immigrant processes are studied in
Section~\ref{seccovar-imm}. By the results in Sections~\ref{secISBM}--\ref{seccovar-imm}, we identify a system of SPDEs for
immigrant processes and prove the existence of its solutions in
Section~\ref{seccappro}. Except in Section~\ref{secpropwc}, we
restrict our attention to the corresponding approximating solutions
from Section~\ref{sec2ndsep} on.

In Section~\ref{sec2ndsep}, we proceed to conditional separation of
the approximating solutions. Some basic results are stated in
Section~\ref{sec2ndsep-basic}, and the setup is given in Section~\ref{sec2ndsep-setup}.
Due to the complexity, the main two lemmas of Section~\ref{sec2ndsep}
are proved in Sections~\ref{seccondsep-1}~and~\ref{seccondsep-2}, respectively, with some preliminaries set in
Section~\ref{secpre}. In Section~\ref{secsepSOL}, we show the
uniform separation of approximating solutions under $\P_{\varepsilon
}$, which completes the proof of our main result.

In Sections~\ref{secpropwc}~and~\ref{sectauvep}, we prove
Propositions~\ref{propwc}~and~\ref{propseptime},
respectively, which are
two technical results. In Section~\ref{secsupp}, we discuss some
properties of the support processes for immigrants. Finally, in
Section~\ref{secimprovmod}, we study the improved modulus of
continuity for functions satisfying certain Gronwall-type integral inequalities.

%s2 #&#
\section{Approximating solutions}\label{secSBMIM}
%s2.1 #&#
\subsection{Interlacing pairs of approximating solutions}\label{secISBM}
In this section, we give details for the approximating solutions of the
SPDE (\ref{eqmainSPDE}) which are discussed in Section~\ref
{secintro}, and state their basic properties. Recall that we identify
every locally integrable function $f$ on $\R$ as a signed measure by
(\ref{signmeas}). We will further write
$f(\Gamma)=f(\1_\Gamma)$
for Borel sets $\Gamma\in\B(\R)$,
whenever the right-hand side makes sense.

For ${\varepsilon}\in(0,1]$, the ${\varepsilon}$-approximating
solutions $X$ and $Y$ in Section~\ref{secintro} obey the equations
given as follows.
The first solution $X$ is a nonnegative c\`adl\`ag $\C_\rap(\R
)$-valued process and is continuous within each time interval
$[s_i,s_{i+1})$ for $s_0=0$ and $s_1,s_2,\ldots$ defined by (\ref
{sttime}). Its time evolution is given by
%
%
%e2.1 #&#
\begin{eqnarray}\label{eqXvep}
X_t(\phi)&=&\int_0^t
X_s \biggl(\frac{\Delta}{2}\phi\biggr) \,ds+\int
_{(0,t]} \int_\R\phi(x)\,dA^{X}(x,s)
\nonumber\\[-8pt]\\[-8pt]\nonumber
&&{} +
\int_0^t\!\int_\R
X(x,s)^{1/2}\phi(x)\,dW(x,s)
\end{eqnarray}
for $\phi\in\C_c^\infty(\R)$.
In (\ref{eqXvep}), the nonnegative measure $A^{X}$ on $\R\times\R
_+$ is defined by
%
%
%e2.2 #&#
\begin{equation}
\label{eqAX} A^{X}\bigl(\Gamma\times[0,t]\bigr)\triangleq\sum
_{i\dvtx0<s_i\leq t} \psi(\1 )J^{x_i}_{\varepsilon}(
\Gamma),
\end{equation}
and $W$ is a space--time white noise. Here, in (\ref{eqAX}), recall
our notation $J^a_{\varepsilon}$ in (\ref{eqJxvep}) and the i.i.d.
spatial random points $\{x_i\}$ with law (\ref{eqxiyilaw}).
In terms of the interpretation in Section~\ref{secintro}, $A^X$ can be
thought of as
being contributed by the initial masses of the underlying immigrant
processes for $X$.

A similar characterization applies to the other approximating solution
$Y$. It is a nonnegative c\`adl\`ag $\C_\rap(\R)$-valued process satisfying
%
%
%e2.3 #&#
\begin{eqnarray}
\label{eqYvep} Y_t (\phi)&=&\int_0^t
Y_s \biggl(\frac{\Delta}{2}\phi\biggr) \,ds+\int
_{(0,t]} \int_\R\phi(x)\,dA^{Y}(x,s)
\nonumber\\[-8pt]\\[-8pt]\nonumber
&&{} +
\int_0^t\!\int_\R
Y(x,s)^{1/2}\phi(x)\,dW(x,s),
\end{eqnarray}
for $\phi\in\C_c^\infty(\R)$,
and is continuous over each $[t_i,t_{i+1})$ for $t_0=0$ and
$t_1,t_2,\ldots$ defined by (\ref{sttime}).
The nonnegative measure $A^{Y}$ on $\R\times\R_+$ is now defined by
\[
A^{Y}\bigl(\Gamma\times[0,t]\bigr)\triangleq\sum
_{j\dvtx0<t_i\leq t} \psi(\1 )J^{y_i}_{\varepsilon}(\Gamma).
\]

We observe that the equations (\ref{eqXvep}) and (\ref{eqYvep}) for
$X$ and $Y$ can be described completely in terms of the processes themselves.
For $X$, each random point $x_i$ in the definition (\ref{eqAX}) of
$A^X$ is a measurable function of the corresponding jump size $\Delta
X_{s_i}$ and conversely, where we write $\Delta Z_s=Z_s-Z_{s-}$ with
$Z_{0-}=0$ for a c\`adl\`ag process $Z$ taking values in a Polish space.
Indeed, we have
%
%
%e2.4 #&#
\begin{equation}
\label{impilc} \qquad x_i=\inf\biggl\{x\in\R;\Delta X_{s_i} \bigl((-
\infty,x] \bigr) > \frac
{\psi(\1){\varepsilon}}{2} \biggr\}\quad\mbox{and}\quad\Delta
X_{s_i}=\psi(\1)J^{x_i}_{\varepsilon},
\end{equation}
where the first equality follows since $\Delta X_{s_i}$ has total mass
$\psi(\1){\varepsilon}$ and defines a measure symmetric about the
center~$x_i$ of its topological support
[cf. (\ref{eqJxvep})]. For~$Y$, similar relations between the random
points $\{y_i\}$ and the jump sizes $\{\Delta Y_{t_i}\}$ hold.

As a summary of the above discussions,
we give in Definition~\ref{defstd} below a minimal description of
the approximating solutions considered throughout this paper. From
Section~\ref{sec2ndsep-basic} on, we will work with ${\varepsilon
}$-approximating pairs
subject to particular correlations.
Here and in the sequel, we use the notation ``$\G\ind\xi$'' to mean
that the $\sigma$-field $\G$ and the random element $\xi$ are
independent, and analogous notation applies to other pairs of objects
which allow probabilistic independence in the usual sense.

%
%
%de2.1 #&#
\begin{defi}\label{defstd}
Fix an immigration function $\psi\in\C_c^+(\R)\setminus\{0\}$. For
any ${\varepsilon}\in(0,1]$, an \textit{interlacing pair} of
${\varepsilon}$-approximating solutions is
a pair $(X,Y)$ defined on a filtered probability space $(\Omega,\F,(\F_t),\P_{\varepsilon})$, with $(\F_t)$ satisfying the usual
conditions, which carries an $(\F_t)$-space--time white noise $W$, and
such that:
\begin{longlist}[(iii)]
\item $X$ and $Y$ are two nonnegative $(\F_t)$-adapted $\C
_\rap(\R)$-valued processes satisfying (\ref{eqXvep}) and (\ref
{eqYvep}) with respect to $W$ for $x_i$ defined by the first equation
in (\ref{impilc}) and $y_i$ by the same equation with $\Delta
X_{s_i}$ replaced by $\Delta Y_{t_i}$, and have paths
which are c\`adl\`ag on $\R_+$ and continuous within each
$[s_i,s_{i+1})$ and $[t_i,t_{i+1})$, respectively,
\item the jumps $\{\Delta X_{s_i},\Delta Y_{t_i};i\in\mathbb
N\}$ are i.i.d. $\C_\rap(\R)$-valued random elements with law given
by (\ref{eqxiyilaw}) through the second equation in (\ref{impilc}), and
\item the random variables $x_i$ and $y_i$ take values in
the topological support of $\psi$ with
%
%
%e2.5 #&#
\begin{equation}
\label{ind-land} \forall i\in\mathbb N\qquad\sigma(X_{t},Y_{t}
; t< s_i) \ind x_i\quad\mbox{and}\quad\sigma(X_t,Y_t
; t<t_i) \ind y_i.
\end{equation}
\end{longlist}
\end{defi}

The existence of these pairs of approximation solutions can be obtained
by considering the so-called mild forms of solutions of SPDEs and then
resorting to the classical Peano's existence argument as in Theorem 2.6
of \cite{SSPDE}. We omit the details.

%
%
%no2.2 #&#
\begin{nota}\label{notavep}
The following convention will be in force throughout this paper unless
otherwise mentioned.
As before, we suppress the dependence on ${\varepsilon}$ for
quantities related to an interlacing pair of ${\varepsilon
}$-approximating solutions except the underlying probability measure
$\P_{\varepsilon}$.
The subscript ${\varepsilon}$ of $\P_{\varepsilon}$ is further
omitted in cases where there is no ambiguity, although in this context
we will remind the readers
of this practice.
%\qed
\end{nota}

The processes described in Definition~\ref{defstd} are genuine
approximating solutions to the SPDE (\ref{eqmainSPDE}) with respect
to the same white noise, as the following proposition states.

%
%
%pr2.3 #&#
\begin{prop}\label{propwc}
Equip $\C_\rap(\R)$ with the norm $\llVert \cdot\rrVert_\rap$ defined by
(\ref{eqCrap-norm}) and $D (\R_+,\C_\rap(\R) )$ with
Skorokhod's $J_1$-topology.
Let $({\varepsilon}_n)\subseteq(0, 1]$ be such that \mbox{${\varepsilon
}_n\searrow0$}, and $ ( (X,Y ),\P_{{\varepsilon
}_n} )$ be a sequence of interlacing pairs of ${\varepsilon
}_n$-approximating solutions. Then the sequence of laws of $ (
(X,Y ),\P_{{\varepsilon}_n} )$ is relatively compact in the
space of probability measures on the product space $D (\R_+,\C
_\rap(\R) )\times D (\R_+,\C_\rap(\R) )$ and every
subsequential limit defines the law of a pair of solutions to the SPDE
(\ref{eqmainSPDE}) subject to the same space--time white noise.
\end{prop}

The proof of Proposition~\ref{propwc} is given in Section~\ref{secpropwc}.
At this point, the readers should be convinced of the result upon
observing the limiting behavior of the random measures $A^X$: for any
$t\in(0,\infty)$,
%
%
%e2.6 #&#
\begin{equation}
\quad\mathbb P\mbox{-}\lim_{{\varepsilon}\searrow0}\int_{(0,t]}
\int_\R\phi(x)\,dA^{X}(x,s) =\mathbb P\mbox{-}
\lim_{{\varepsilon}\searrow0}\psi(\1){\varepsilon}\sum
_{i=1}^{\lfloor t{\varepsilon}^{-1}\rfloor}\phi(x_i) =t\langle\psi,\phi
\rangle\label{eqconvip}
\end{equation}
for any $\phi\in\C_c^\infty(\R)$,
by the law of large numbers. Here, $\Plim\lim$ denotes convergence in
probability, and
$\lfloor t\rfloor$ is the greatest integer less than or equal to $ t$.

We close this section with a property of the above approximating
solutions. Here and in the sequel, we use the following notation.
For any real-valued random function $Z\dvtx(x,s)\lmt Z(x,s)$, we write
%
%
%e2.7 #&#
\begin{eqnarray}
(Z\in\Gamma)&\triangleq& \bigl\{(x,s)\in\R\times\R_+;Z(x,s)\in\Gamma
\bigr\},
\qquad\Gamma\in\B(\R).\label{spacetime}
\end{eqnarray}
In addition, for a space--time white noise $W'$, we write $L^2_\loc
(W')$ for the set of functions $Z=Z(\omega,x,s)$, product measurable
in $(\omega,s)$ and $x$ with respect to the underlying predictable
$\sigma$-field and $\B(\R)$, so that
%
%
%e2.8 #&#
\begin{equation}
\label{eqL2W} \int_0^t\!\int
_{\R}Z(x,s)^2\,dx\,ds<\infty\qquad\forall t\in(0,
\infty)\mbox{ a.s.},
\end{equation}
and define processes of stochastic integrals as
%
%
%e2.9 #&#
\begin{equation}
\label{bullet} Z\bullet W'(\phi)\equiv\int_0^\cdot\!\int_\R Z(x,s)\phi(x) \,dW(x,s)
\end{equation}
for $Z\in L^2_\loc(W')$ and $\phi\in L^2(\R)$.

%
%
%pr2.4 #&#
\begin{prop}[(Cherny's substitution)]\label{propsbst}
For ${\varepsilon}\in(0,1]$, let $(X,Y)$ be an interlacing pair of
${\varepsilon}$-approximating solutions. By enlarging the underlying
filtered probability space $ (\Omega,\F,(\F_t),\P_{\varepsilon
} )$ if necessary, we can find random elements
\[
V^{X,i}= \bigl\{ \bigl(V^{X,i}_t(\phi)
\bigr)_{t\in[s_i,\infty)};\phi\in L^2(\R) \bigr\}
\]
and
\[
V^{Y,i}= \bigl\{ \bigl(V^{Y,i}_t(\phi)
\bigr)_{t\in[t_i,\infty)};\phi\in L^2(\R) \bigr\},
\]
for $i\in\mathbb N$, which satisfy the following properties:
\begin{longlist}[(iii)]
\item[(i)] Every $V^{X,i}$ is an $(\F_{t})_{t\in[s_i,\infty
)}$-space--time white noise and satisfies
%
%
%e2.10 #&#
\begin{equation}
\label{indsbst} V^{X,i}\ind\bigl\{\F_{s_i},
(X_t)_{t\in[s_i,s_{i+1})} \bigr\}.
\end{equation}
Here, an $(\F_t)_{t\in[s_i,\infty)}$-space--time white noise is
defined as an $(\F_t)$-space--time white noise except that
its components are $(\F_t)_{t\in[s_i,\infty)}$-Brownian motions
started at $s_i$ with zero initial value.
\item[(ii)] Every $V^{Y,i}$ satisfies the same properties in (i)
with $(X,s_i)$ replaced by $(Y,t_i)$.
\item[(iii)] The following substitution identities of space--time
white noises hold: for all $i\in\mathbb N$,
%
%
%e2.11 #&#
%e2.12 #&#
\begin{eqnarray}
\1_{[s_i,s_{i+1})}\1_{(X=0)}\bullet W&=&\1_{[s_i,s_{i+1})}\1
_{(X=0)}\bullet V^{X,i},\label{sbstX}
\\
\1_{[t_i,t_{i+1})}\1_{(Y=0)}\bullet W&=&\1_{[t_i,t_{i+1})}\1
_{(Y=0)}\bullet V^{Y,i}.\label{sbstY}
\end{eqnarray}
\end{longlist}
\end{prop}

Proposition~\ref{propsbst} will be used in Section~\ref{seccont} to
reinforce immigrant processes from continuous decomposition with
analogous properties
[condition (vi) of Theorem~\ref{propequiallocation}]. From these
properties, we will deduce some key equations for covariations of the
immigrant processes (see Proposition~\ref{propcovar}).

\begin{pf*}{Sketch of proof of Proposition~\ref{propsbst}}
The proof is a generalization of the proof of Theorem 3.1 in Cherny
\cite{CUIL} to the context of the SPDE (\ref{eqmainSPDE0}), and so
we only give a sketch. Below we consider the assertions for $V^{X,i}$
for $i\in\mathbb N$. The assertions for $V^{Y,i}$ follow similarly.

We define $V^{X,i}$ as a mixture of the original noise $W$ and another
space--time white noise, say $U^{X,i}$, which is independent of
$(X,Y,W)$ and adapted to the same filtration, by
%
%
%e2.13 #&#
\begin{equation}
\label{eqsbstV} V^{X,i}\triangleq\1_{[s_i,\infty)}\1_{(X=0)}
\bullet W+\1_{[s_i,\infty
)}\1_{(X>0)}\bullet U^{X,i}.
\end{equation}
Then $V^{X,i}$ is
an $(\F_t)_{t\in[s_i,\infty)}$-space--time white noise by L\'evy's
theorem (cf. Theorem IV.3.6 of \cite{RYCMB})
and gives the required substitution (\ref{sbstX}). We have proved the
first assertion in (i) and the assertion in (iii) for $V^{X,i}$.

It remains to prove the independence (\ref{indsbst}). Consider the
counterpart of $V^{X,i}$ (\ref{eqsbstV}):
\[
\widetilde{V}{}^{X,i}\triangleq\1_{[s_i,\infty)}\1_{(X>0)}
\bullet W+\1 _{[s_i,\infty)}\1_{(X=0)}\bullet U^{X,i}.
\]
By\vspace*{2pt} L\'evy's theorem again, $\widetilde{V}{}^{X,i}$ is an $(\F_t)_{t\in
[s_i,\infty)}$-space--time white noise and $V^{X,i}\ind\widetilde{V}{}^{X,i}$.
The latter property implies that $(X,\widetilde{V}{}^{X,i})$ over
$[s_i,s_{i+1})$ solves\vspace*{1pt} the SPDE (\ref{eqmainSPDE0}) of super-Brownian
motion with respect to $ (\F_t\vee\sigma(V^{X,i}) )_{t\in
[s_i,s_{i+1})}$. Recall that the martingale problem for super-Brownian
motion is well-posed (cf. Lemma~\ref{lemred2}), and\vspace*{2pt} note that
$V^{X,i}\ind\F_{s_i}$ by a standard property of Brownian motion. We
deduce that for any $\Gamma_0\in\sigma(V^{X,i})$ and $\Gamma_1\in
\F_{s_i}$ with $\P_{\varepsilon}(\Gamma_0\cap\Gamma_1)>0$,
%
%
%e2.14 #&#
\begin{equation}
\label{eqcondprob} \quad\frac{\P_{\varepsilon}(\Gamma_0\cap\Gamma_1 \cap\{
(X_t)_{t\in
[s_i,s_{i+1})} \in\cdot\})}{\P_{\varepsilon}(\Gamma_0\cap\Gamma
_1)}=\frac{\P_{\varepsilon}(\Gamma_1 \cap\{(X_t)_{t\in
[s_i,s_{i+1})} \in\cdot\})}{\P_{\varepsilon}(\Gamma_1)}.
\end{equation}
Indeed, under the conditional probabilities $\P_{\varepsilon}( \cdot
\mid\Gamma_0\cap\Gamma_1)$ and $\P_{\varepsilon}( \cdot\mid\Gamma_1)$,
the laws of $X_{s_i}$ are the same and $ ( (\F_t\vee\sigma
(V^{X,i}) )_{t\in[s_i,s_{i+1})},\P_{\varepsilon}
)$-martingales remain martingales.
See the proof of Theorem 4.4.2 of \cite{EK}. Using $V^{X,i}\ind\F
_{s_i}$ again, we can substitute $\P_{\varepsilon}(\Gamma_0\cap
\Gamma_1)$ on the left-hand side of (\ref{eqcondprob}) with $\P
_{\varepsilon}(\Gamma_0)\P_{\varepsilon}(\Gamma_1)$. The required
property (\ref{indsbst}) follows.
\end{pf*}

%s2.2 #&#
\subsection{Continuous decomposition}\label{seccont}
For every ${\varepsilon}\in(0,1]$, consider an interlacing pair
$(X,Y)$ of ${\varepsilon}$-approximating solutions to the SPDE (\ref
{eqmainSPDE}) (recall Definition~\ref{defstd}). From their informal
descriptions in Section~\ref{secintro},
it is reasonable to expect that they can be \emph{decomposed} into
%
%
%e2.15 #&#
\begin{equation}
X=\sum_{i=1}^\infty X^i\quad
\mbox{and}\quad Y= \sum_{i=1}^\infty
Y^i,\label{eqallocXY}
\end{equation}
where the summands $X^i$ and $Y^i$ are super-Brownian motions started
at $s_i$ and $t_i$ and with starting measures $\Delta X_{s_i}=\psi(\1
)J^{x_i}_{\varepsilon}$ and $\Delta Y_{t_i}=\psi(\1
)J^{y_i}_{\varepsilon}$, respectively, for each $i$, and each of the
families $\{X^i\}$ and $\{Y^i\}$ consists of independent random elements.

Let us give an elementary discussion on obtaining the decompositions in
(\ref{eqallocXY}). Later on, we will require additional properties of
the decompositions.
It follows from the uniqueness in law of super-Brownian motions and the
defining equation (\ref{eqXvep}) that $X$ is a (time-inhomogeneous)
Markov process and, for each $i\in\mathbb N$, $ (X_t )_{t\in
[s_i,s_{i+1})}$ defines a super-Brownian\vspace*{1pt} motion with initial
distribution $X_{s_{i}}$ (cf. the proof of Theorem 4.4.2 in \cite{EK}
and the martingale problem characterization of super-Brownian motion in
\cite{PDW}). Hence, each of the equalities in (\ref{eqallocXY})
holds in the sense of being \emph{identical in distribution}. Then we
recall the following general theorem (see Theorem 6.10 in \cite{KFMP}).

%
%
%th2.5 #&#
\begin{teo}\label{teose}
Fix any measurable space $E_1$ and Polish space $E_2$, and let $\xi
\stackrel{(\mathrm{d})}{=}\widetilde{\xi}$ and $\eta$ be random elements
taking values in $E_1$ and $E_2$, respectively. Here, we only assume
that $\xi$ and $\eta$ are defined on the same probability space.
Then there exists a measurable\vspace*{1pt} function $F\dvtx E_1\times[0,1]\lra E_2 $
such that for any random variable $\widetilde{U}$ uniformly
distributed over $[0,1]$ with $\widetilde{U}\ind\widetilde{\xi}$,
the random element $\widetilde{\eta}=F(\widetilde{\xi},\widetilde
{U})$ solves
$(\xi,\eta)\stackrel{(\mathrm{d})}{=}(\widetilde{\xi},\widetilde{\eta})$.
\end{teo}

By the preceding discussions and Theorem~\ref{teose}, we can
immediately construct the summands $X^i$ and $Y^i$ by introducing
additional independent uniform variables and validate the equalities
(\ref{eqallocXY}) as almost-sure equalities. Such decompositions,
however, are too crude because, for example, we are unable to say that
all the resulting random processes perform their defining properties
with respect to the \emph{same} filtration. This difficulty implies in
particular that we cannot do semimartingale calculations for them. A
finer decomposition method, however, does yield a solution to this
problem. The result is stated in Theorem~\ref{propequiallocation}
below. See also Figure~\ref{Figca} for a sketch of the decomposition
of $X$ along a particular value $x$.

%
%
%f1 #&#
\begin{figure}

\includegraphics{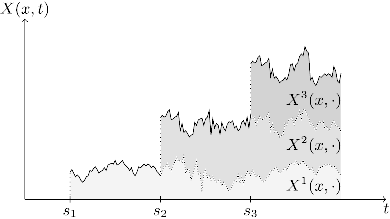}

\caption{Decomposition of $X$ along $x$.}\label{Figca}
\end{figure}

%
%
%th2.6 #&#
\begin{teo}[(Continuous decomposition)]\label{propequiallocation}
Fix ${\varepsilon}\in(0,1]$. Let $ (X,Y,W)$ be an interlacing
pair of ${\varepsilon}$-approximating solutions, and $\{V^{X,i}\},\{
V^{Y,i}\}$ be two families of space--time white noises chosen in
Theorem~\ref{propsbst}. By changing the underlying probability space
if necessary,\vspace*{1pt} we can find a filtration $(\G_t)$ satisfying the usual
conditions and two families $\{X^i\}$ and $\{Y^i\}$ of nonnegative $\C
_\rap(\R)$-valued processes, such that the following conditions are satisfied:
\begin{longlist}[(iii)]
\item[(i)] The processes $X^i$, $i\in\mathbb N$, are independent.
\item[(ii)] The equality in (\ref{eqallocXY}) for $X$ and $X^i$
holds almost surely.
\item[(iii)] Each $ (X^i_{t} )_{t\in[s_i,\infty)}$
has sample paths in $C ([s_i,\infty),\C_\rap(\R) )$ and is
a $(\G_t)_{t\geq s_i}$-super-Brownian motion started at time $s_i$
with starting measure $\psi(\1)J^{x_i}_{\varepsilon}$.
Also, $X^i_t\equiv0$ for every $t\in[0,s_i)$.
\item[(iv)] The processes $Y^i$, $i\in\mathbb N$, satisfy the same
properties as \textup{(i)--(iii)} with the roles of $X$ and $\{(X^i,x_i,s_i)\}$
replaced by $Y$ and $\{(Y^i,y_i,t_i)\}$, respectively.
\item[(v)] Conditions \textup{(i)--(iii)} of Definition~\ref{defstd} hold
with $(\F_t)$ replaced by $(\G_t)$, and~(\ref{ind-land}) in the same
definition is replaced by the stronger \textit{independent landing property}:
%
%
%e2.16 #&#
%e2.17 #&#
\begin{eqnarray}\label{coneilc}
\quad\sigma\bigl(X^j_t,Y^j_t; t<s_i, j\in\mathbb N \bigr)\ind x_i\quad\mbox{and}\quad
\sigma\bigl(X^j_t,Y^j_t;
t<t_i, j\in\mathbb N \bigr)\ind y_i
\nonumber\\[-8pt]\\[-8pt]
\eqntext{\forall i\in\mathbb N.}
\end{eqnarray}
\item[(vi)] Condition \textup{(iii)} of Proposition~\ref{propsbst} holds.
In addition,
%
%
%e2.18 #&#
\begin{eqnarray}
\label{Vind} \bigl\{\bigl(X^j_t\bigr)_{t\in[0,s_{i+1})};j
\in\mathbb N\mbox{ satisfying }s_j<s_{i+1} \bigr\}& \ind&
V^{X,i}\quad\mbox{and}
\nonumber\\[-8pt]\\[-8pt]\nonumber
\bigl\{\bigl(Y^j_t\bigr)_{t\in[0,t_{i+1})};j\in\mathbb N
\mbox{ satisfying }t_j<t_{i+1} \bigr\}& \ind&
V^{Y,i}\qquad\forall i\in\mathbb N.
\end{eqnarray}
\end{longlist}
\end{teo}

Due to the length of the proof of Theorem~\ref{propequiallocation},
we first outline its \emph{informal} idea for the convenience of
readers. Recall that the first immigration event for $X$ and $Y$ occurs
at $s_1=\frac{{\varepsilon}}{2}$. Take\vspace*{1pt} a grid of $[\frac
{{\varepsilon}}{2},\infty)$ containing all the points $s_i$ and $t_i$
for $i\in\mathbb N$ and with ``infinitesimal'' mesh size. Here, the mesh
size of a grid is the supremum of the distances between consecutive
grid points. The key observation in this construction is that,
over any subinterval $[t,t+\Delta t]\subseteq[s_i,s_{i+1})$ from this
grid, $(X_r;r\in[t,t+\Delta t])$ has the same distribution as the sum
of $i$ independent super-Brownian motions started at $t$ over
$[t,t+\Delta t]$, \emph{whenever} the sum of the initial conditions of
these independent super-Brownian motions has the same distribution as $X_{t}$.

This fact allows us to inductively decompose $X$ over the intervals of
infinitesimal lengths from this grid, such that the resulting
infinitesimal pieces of super-Brownian motions can be concatenated in
the natural way to obtain the desired immigrating super-Brownian motions.
More precisely, we apply Theorem~\ref{teose} by bringing in
independent uniform variables as ``allocators'' to obtain these
infinitesimal pieces.
A similar method applies to continuously decompose $Y$ into the desired
independent super-Brownian motions by another family of independent allocators.

Finally, because the path regularity of these concatenated processes
and $W$ allows us to characterize their laws over the entire time
horizon $\R_+$ by their laws over $[0,{\varepsilon}/2]$ and their
probabilistic transitions on this grid with infinitesimal mesh size,
the filtration obtained by sequentially adding the $\sigma$-fields of
the independent allocators will be the desired one. In particular, the
time evolutions of these stochastic processes are now consistent with
the ``progression'' of the enlarged filtration.

\begin{pf*}{Proof of Theorem~\ref{propequiallocation}}
Fix ${\varepsilon}\in(0,1]$ and we shall drop the subscript
${\varepsilon}$ of $\P_{\varepsilon}$.
Throughout the proof, we take, for each $m\in\mathbb N$, a countable
subset $D_m$ of $[\frac{{\varepsilon}}{2},\infty)$ which contains
$s_i$ and $t_i$ for any $i\in\mathbb N$ and satisfies $\# (D_m\cap
K )<\infty$ for any compact subset $K$ of $\R_+$.
We further assume that (1) $D_{m+1}\subseteq D_m$ for each $m$, (2)~between any two points $s_i$ and $t_i$ there is another point belonging
to~$D_1$, and hence to each $D_m$, and (3) the mesh size of $D_m$ goes
to zero as \mbox{$m\lra\infty$}.
In addition, we will write
$\{{\SBM}_t(\mu,d\nu);t\in\R_+\}$ for the semigroup of
super-Brownian motion on $\R$. When the \emph{density} of the
super-Brownian motion on $\R$
started at time $s$ and with starting measure $f(x)\,dx$ for a
nonnegative $\C_\rap(\R)$-function $f$ is concerned, we write ${\SBM}_{f,[s,t]}$
for the law of its $C ([s,\infty),\C_\rap(\R) )$-valued
density restricted to the time interval $[s,t]$.
\begin{longlist}[\textit{Step} 1.]
\item[\textit{Step} 1.] Fix $m\in\mathbb N$ and write $\frac{{\varepsilon
}}{2}=\tau_0<\tau_1<\cdots$ as the consecutive points of $D_m$.
Assume, by an enlargement of the underlying probability space where
$ (X,Y,W,\{V^{X,i}\},\{V^{Y,i}\} )$ lives if necessary, the
existence of i.i.d. variables $\{U_{j}^X,U^Y_j;j\in\mathbb N\}$ with
%
%
%e2.19 #&#
\begin{equation}
\label{eqUXUY} \quad U_1^X\mbox{ is uniformly distributed over
}[0,1]\mbox{ and }\bigl\{U_j^X,U^Y_j;j
\in\mathbb N \bigr\}\ind\F.
\end{equation}
In this step, we will decompose $X$ and $Y$ into the random elements
\begin{eqnarray*}
\mc X^{m}&=& \bigl(X^{m,1},X^{m,2},\ldots\bigr)
\quad\mbox{and}\quad\mc Y^{m}= \bigl(Y^{m,1},Y^{m,2},
\ldots\bigr),
\end{eqnarray*}
respectively, according to the grid $D_m$.
Here,
%
%
%e2.20 #&#
%e2.21 #&#
\begin{eqnarray}
\label{eqXmYmpath} X^{m,i}\in C \bigl([s_i,\infty),
\C_\rap(\R) \bigr)\quad\mbox{and}\quad Y^{m,i}\in C
\bigl([t_i,\infty),\C_\rap(\R) \bigr)
\nonumber\\[-8pt]\\[-8pt]
\eqntext{\mbox{with }X^{m,i}\equiv0 \mbox{ on } [0,s_i)
\mbox{ and } Y^{m,i}\equiv0 \mbox{ on } [0,t_i),}
\end{eqnarray}
so we need to specify $X^{m,i}$ over $[s_i,\infty)$ and $Y^{m,i}$ over
$[t_i,\infty)$.

Consider the construction of $\mc X^m$.
The decomposition of $X$ over $[s_1,s_2]$ should be evident. Over this
interval, set $X^{m,1}\equiv X $ on $[s_1,s_2)$ with
$X^{m,1}_{s_2}=X_{s_2-}$ and
\[
X^{m,2}_{s_2}= \psi(\1)J^{x_2}_{\varepsilon}=
\Delta X_{s_2}.
\]
We define $\mc X^m$ over $[s_2,\tau_j]$ by an induction on integers
$j\geq j_*^X$, where $j_*^X\in\mathbb N$ satisfies $s_2=\tau_{j_*^X}$,
such that:
\begin{longlist}[(3)]
\item[(1)] the following measurability condition holds:
%
%
%e2.22 #&#
%e2.23 #&#
\begin{eqnarray}\label{eqeqdistINF}
&& \bigl(X^{m,i}_s;s\in[0,
\tau_k], i\in\mathbb N \bigr)
\in\sigma\bigl(X_s;s\in[0, \tau_k]
\bigr)\vee\sigma\bigl(U_{i}^X,1\leq i\leq k \bigr)
\nonumber\\[-8pt]\\[-8pt]
\eqntext{\forall k\in\{ 0,\ldots,j\},}
\end{eqnarray}
with $\sigma(U_i^X,1\leq i\leq0 )$ understood to be the
trivial $\sigma$-field $\{\Omega,\varnothing\}$,
\item[(2)] the laws of $X^{m,i}$ obey
%
%
%e2.24 #&#
\begin{equation}
\label{eqextlaw} \qquad\cases{\mbox{(a)}\quad \ms L \bigl(X^{m,i}_{s};s
\in[s_i,\tau_j] \bigr)\sim{\SBM}_{\psi(\1)J^{x_i}_{\varepsilon},
[s_i,\tau_j ]}\mbox{ if $s_i\leq\tau_j$ and}
\cr
\mbox{(b)}\quad
\bigl(X^{m,i}_{s};s\in[s_i,\tau_j]
\bigr)\mbox{ for $i$ satisfying } s_i\leq\tau_j,\mbox{
are independent,}}
\end{equation}
and
\item[(3)] a preliminary decomposition of $X$ up to time $\tau_j$ holds:
%
%
%e2.25 #&#
\begin{equation}
X_s=\sum_{i=1}^\infty
X^{m,i}_s\qquad\forall s\in[0,\tau_j]\mbox{ a.s.}\label{eqeqdistMASS}
\end{equation}
\end{longlist}
By the foregoing definition of $\mc X^m$ over $[s_1,s_2]$ and (\ref{ind-land}),
we have handled the case that $j=j_*^X$, that is, the first step of our
inductive construction.

Assume that $\mc X^{m}$ has been defined up to time $\tau_j$ for some
integer $j\geq j_*^X$ such that (\ref{eqeqdistINF})--(\ref
{eqeqdistMASS}) are all satisfied. Then our goal is to extend $\mc
X^{m}$ over $[\tau_j,\tau_{j+1}]$ so that all of (\ref
{eqeqdistINF})--(\ref{eqeqdistMASS}) hold with $j$ replaced by
$j+1$. First, consider the case that
%
%
%e2.26 #&#
\begin{equation}
\label{eqtausk} [\tau_j,\tau_{j+1} ]\subseteq
[s_k,s_{k+1} )
\end{equation}
for some $k\geq2$.
In this case, we need to extend $X^{m,1},\ldots,X^{m,k}$ up to time
$\tau_{j+1}$.
Take an auxiliary nonnegative random element
\[
\xi= \bigl(\xi^1,\ldots,\xi^k \bigr)\in C \Biggl([
\tau_j,\infty),\prod_{i=1}^k
\C_\rap(\R) \Biggr)
\]
such that the coordinates $ (\xi^i_s;s\in[\tau_j,\infty) )$
are independent processes and each of them defines a super-Brownian
motion started at time $\tau_j$ with initial law
%
%
%e2.27 #&#
\begin{equation}
\label{eqinitialxiX} \ms L \bigl(\xi^i_{\tau_j} \bigr)=\ms L
\bigl(X^{m,i}_{\tau_j} \bigr)\qquad\forall i\in\{1,\ldots,k\}.
\end{equation}

We claim that it is possible to extend $X^{m,1},\ldots,X^{m,k}$
continuously over $[\tau_j,\tau_{j+1}]$ so that
%
%
%e2.28 #&#
\begin{eqnarray}
\label{eqeqindist}  && \bigl(X^{m,1}_{\tau_j},\ldots,X^{m,k}_{\tau
_j},(X_r)_{r\in
[\tau_j, \tau_{j+1} ]},
\bigl(X^{m,1}_r \bigr)_{r\in[\tau
_j,\tau_{j+1}]},\ldots,
\bigl(X^{m,k}_r \bigr)_{r\in[\tau_j,\tau_{j+1}]} \bigr)
\nonumber\\[-4pt]\\[-14pt]\nonumber
&&\qquad\deq\Biggl(\xi^1_{\tau_j},\ldots,
\xi^{k}_{\tau
_j}, \Biggl(\sum_{i=1}^k
\xi^i_r \Biggr)_{r\in[\tau_j,\tau
_{j+1} ]}, \bigl(
\xi^1_r \bigr)_{r\in[\tau_j,\tau
_{j+1}]},\ldots, \bigl(
\xi^k_r \bigr)_{r\in[\tau_j,\tau_{j+1}]} \Biggr).\hspace*{-15pt}
\end{eqnarray}
To prove our claim that (\ref{eqeqindist}) can be done, first we consider
\begin{eqnarray}\label{eqYwt3}
&&\P\bigl((X_r)_{r\in[\tau_j,\tau_{j+1} ]}\in\Gamma, X^{m,1}_{\tau_j}
\in A_1,\ldots,X^{m,k}_{\tau_j}\in A_k
\bigr)
\nonumber
\\
&&\qquad =\E\bigl[ \P\bigl((X_r)_{r\in[\tau_j,\tau_{j+1}
]}\in\Gamma\mid
X_{\tau_j} \bigr);X^{m,1}_{\tau_j}\in A_1,
\ldots,X^{m,k}_{\tau_j}\in A_k \bigr]
\nonumber\\[-8pt]\\[-8pt]\nonumber
&&\qquad =\E\bigl[{\SBM}_{X_{\tau_j},[\tau_j,\tau_{j+1}]} ( \Gamma
);X^{m,1}_{\tau_j}
\in A_1,\ldots,X^{m,k}_{\tau_j}\in A_k
\bigr]
\nonumber
\\
&&\qquad =\E\bigl[ {\SBM}_{\sum_{i=1}^kX^{m,i}_{\tau_j}, [\tau
_j,\tau_{j+1} ]} ( \Gamma);X^{m,1}_{\tau_j}
\in A_1,\ldots,X^{m,k}_{\tau_j}\in A_k
\bigr],\nonumber
\end{eqnarray}
where the first and the second equalities use the (time-inhomogeneous)
Markov property of $X$ and (\ref{eqeqdistINF}), and the last equality
follows from
the equality (\ref{eqeqdistMASS}) by induction. Second, by (\ref
{eqextlaw}) from induction and (\ref{eqinitialxiX}), we have
\[
\bigl(X^{m,1}_{\tau_j},\ldots,X^{m,k}_{\tau_j}
\bigr)\deq\bigl(\xi^1_{\tau_j},\ldots,\xi^k_{\tau_j}
\bigr).
\]
Hence, from (\ref{eqYwt3}), we get
\begin{eqnarray}\label{eqYwt4}
&&\P\bigl((X_r)_{r\in[\tau_j,\tau_{j+1} ]}\in\Gamma,X^{m,1}_{\tau_j}
\in A_1,\ldots,X^{m,k}_{\tau_j}\in A_k
\bigr)
\nonumber
\\
&&\qquad =\E\bigl[{\SBM}_{\sum_{i=1}^k\xi^i_{\tau_j},[\tau_j,\tau
_{j+1}]}(\Gamma);\xi^1_{\tau_j}
\in A_1,\ldots,\xi^k_{\tau_j}\in A_k
\bigr]
\nonumber\\[-6pt]\\[-10pt]\nonumber
&&\qquad =\E\Biggl[\P\Biggl( \Biggl(\sum_{i=1}^k
\xi^i_r \Biggr)_{r\in
[\tau_j, \tau_{j+1} ]}\in\Gamma\bigg|
\xi^1_{\tau_j},\ldots,\xi^k_{\tau_j}
\Biggr);\xi^1_{\tau_j}\in A_1,\ldots,
\xi^k_{\tau_j}\in A_k \Biggr]\hspace*{-20pt}
\nonumber
\\
&&\qquad =\P\Biggl( \Biggl(\sum_{i=1}^k
\xi^i_r \Biggr)_{r\in[\tau_j,
\tau_{j+1} ]}\in\Gamma,
\xi^1_{\tau_j}\in A_1,\ldots,\xi
^k_{\tau_j}\in A_k \Biggr).\nonumber
\end{eqnarray}
Here, the second equality follows from the convolution property of the
laws of super-Brownian motions:
\[
{\SBM}_{f_1,[s,t]}\star\cdots\star{\SBM}_{f_k,[s,t]}={\SBM}_{\sum
_{i=1}^kf_i,[s,t]}.
\]
Then (\ref{eqYwt4}) implies that
%
%
%e2.29 #&#
\begin{equation}
\label{eqeqindist1} \qquad\bigl(X^{m,1}_{\tau_j},\ldots,X^{m,k}_{\tau
_j},(X_r)_{r\in
[\tau_j, \tau_{j+1} ]}
\bigr) \deq\Biggl(\xi^1_{\tau_j},\ldots,
\xi^{k}_{\tau_j}, \Biggl(\sum_{i=1}^k
\xi^i_r \Biggr)_{r\in[\tau_j,\tau_{j+1}
]} \Biggr).
\end{equation}
Using the ``boundary condition'' (\ref{eqeqindist1}) and Theorem~\ref
{teose}, we can solve the stochastic equation on the left-hand side of
(\ref{eqeqindist}) by a Borel measurable function
\begin{eqnarray*}
&& F_{j}^{m}\dvtx\prod_{i=1}^k
\C_\rap(\R)\times C \bigl([\tau_j,\tau_{j+1}],
\C_\rap(\R) \bigr)\times[0,1]
\\
&&\qquad \lra\prod_{i=1}^kC
\bigl([\tau_j,\tau_{j+1}],\C_\rap(\R) \bigr)
\end{eqnarray*}
such that the desired extension of $\mc X^m$ over $[\tau_j,\tau
_{j+1}]$ can be defined by
%
%
%e2.30 #&#
\begin{eqnarray}
\label{eqfuncY}  && \bigl( \bigl(X^{m,1}_r
\bigr)_{r\in[\tau_j,\tau_{j+1}
]},\ldots, \bigl(X^{m,k}_r
\bigr)_{r\in[\tau_j,\tau_{j+1}
]} \bigr)
\nonumber\\[-8pt]\\[-8pt]\nonumber
&&\qquad=F_j^{m} \bigl(X^{m,1}_{\tau_j},
\ldots,X^{m,k}_{\tau
_j},(X_r)_{r\in[\tau_j,\tau_{j+1} ]},
U_{j+1}^X \bigr),
\end{eqnarray}
where the independent uniform variable $U_{j+1}^X$ now plays its role
to decompose $(X_r)_{r\in[\tau_j,\tau_{j+1}]}$. This proves
our claim on the continuous extension of $X^{m,1},\ldots,X^{m,k}$ over
$[\tau_j,\tau_{j+1}]$ satisfying (\ref{eqeqindist}). As a
consequence of
the equality~(\ref{eqeqindist}) in distribution, the following
equalities hold almost surely:
%
%
%e2.31 #&#
\begin{equation}
\label{eqYwt1} X^{m,1}_r+\cdots+X^{m,k}_r=X_r
\qquad\forall r\in[\tau_j, \tau_{j+1} ]
\end{equation}
and
%
%
%e2.32 #&#
\begin{eqnarray}
\label{eqYwt2} && \ms L \bigl( \bigl(X^{m,1}_r
\bigr)_{r\in[\tau_j,\tau
_{j+1} ]},\ldots, \bigl(X^{m,k}_r
\bigr)_{r\in[\tau_j,
\tau_{j+1} ]}\mid X^{m,1}_{\tau_j},\ldots,X^{m,k}_{\tau
_j},X_{\tau_j}
\bigr)\nonumber
\\
&&\qquad =\ms L \bigl( \bigl(X^{m,1}_r \bigr)_{r\in[\tau_j,\tau
_{j+1} ]},
\ldots, \bigl(X^{m,k}_r \bigr)_{r\in[\tau_j,\tau
_{j+1} ]}\mid
X^{m,1}_{\tau_j},\ldots,X^{m,k}_{\tau_j} \bigr)
\\
&&\qquad ={\SBM}_{X^{m,1}_{\tau_j}, [\tau_j,\tau_{j+1}
]}\otimes\cdots\otimes{\SBM}_{X^{m,k}_{\tau_j}, [\tau
_j,\tau_{j+1} ]},\nonumber
\end{eqnarray}
where the first equality of (\ref{eqYwt2}) follows from (\ref
{eqYwt1}), and the second from the definition of $\xi$.
By
induction and (\ref{eqfuncY}), the extension of $\mc X^m$ over $[\tau
_j,\tau_{j+1}]$ satisfies~(\ref{eqeqdistINF}) with $j$ replaced by $j+1$;
by induction and (\ref{eqYwt1}), it satisfies (\ref{eqeqdistMASS})
with $j$ replaced by $j+1$.

Let us verify that (\ref{eqextlaw}) is satisfied with $j$ replaced by
$j+1$. By (\ref{eqeqdistINF}), we can write
%
%
%e2.33 #&#
\begin{eqnarray}\label{eqsig1}
&&\P\bigl( \bigl(X^{m,i}_r \bigr)_{r\in[\tau_j,\tau_{j+1}]}\in
A_i, \bigl(X^{m,i}_r \bigr)_{r\in[s_i,\tau_{j}]}\in
B_i,\forall i\in\{1,\ldots,k\} \bigr)\hspace*{-15pt}
\nonumber
\\
&&\qquad =\E\bigl[\P\bigl( \bigl(X^{m,i}_r \bigr)_{r\in[\tau_j,\tau
_{j+1}]}
\in A_i, \forall i\in\{1,\ldots,k\}\mid\F_{\tau_j}\vee\sigma
\bigl(U^X_1,\ldots,U^X_j \bigr)
\bigr);\hspace*{-15pt}
\\
&&\hspace*{177pt}\bigl(X^{m,i}_r \bigr)_{r\in[s_i,\tau_{j}]}\in
B_i,\forall i\in\{1,\ldots,k\} \bigr].\hspace*{-15pt}\nonumber
\end{eqnarray}
To\vspace*{1pt} reduce the conditional probability on the right-hand side of (\ref
{eqsig1}) to a probability conditioned on $X^{m,1}_{\tau_j},\ldots,X^{m,k}_{\tau_j}$, we
review the defining equation (\ref{eqfuncY}) of~$\mc X^m$ over $[\tau
_j,\tau_{j+1}]$ and consider the calculation:
%
%
%e2.34 #&#
\begin{eqnarray}\label{eqsig2}
&&\E\bigl[g_1 \bigl(X^{m,1}_{\tau_j},
\ldots,X^{m,k}_{\tau_j} \bigr)g_2
\bigl((X_r)_{r\in[\tau_j,\tau_{j+1} ]} \bigr)g_3
\bigl(U^X_{j+1} \bigr)\mid\F_{\tau_j}\vee\sigma
\bigl(U^X_1,\ldots,U^X_j \bigr)
\bigr]\hspace*{-18pt}
\nonumber
\\
&&\qquad =g_1 \bigl(X^{m,1}_{\tau_j},\ldots,X^{m,k}_{\tau_j}
\bigr)\E\bigl[g_2 \bigl((X_r)_{r\in[\tau_j,\tau_{j+1} ]} \bigr)
\mid\F_{\tau_j}\vee\sigma\bigl(U^X_1,
\ldots,U^X_j \bigr) \bigr]\hspace*{-18pt}
\nonumber
\\
&&\qquad\quad{} \times\E\bigl[g_3 \bigl(U^X_{j+1}
\bigr) \bigr]\hspace*{-18pt}
\nonumber\\[-8pt]\\[-8pt]\nonumber
&&\qquad =g_1 \bigl(X^{m,1}_{\tau_j},\ldots,X^{m,k}_{\tau_j}
\bigr)\E\bigl[g_2 \bigl((X_r)_{r\in[\tau_j,\tau_{j+1} ]} \bigr)
\mid X^{m,1}_{\tau_j},\ldots,X^{m,k}_{\tau_j}
\bigr]\hspace*{-18pt}
\\
&&\quad\qquad{} \times\E\bigl[g_3 \bigl(U^X_{j+1} \bigr)\bigr]\hspace*{-18pt}\nonumber
\\
&&\qquad =\E\bigl[g_1 \bigl(X^{m,1}_{\tau_j},
\ldots,X^{m,k}_{\tau_j} \bigr)g_2
\bigl((X_r)_{r\in[\tau_j,\tau_{j+1} ]} \bigr)g_3
\bigl(U^X_{j+1} \bigr)\mid X^{m,1}_{\tau_j},
\ldots,X^{m,k}_{\tau_j} \bigr],\hspace*{-18pt}\nonumber
\end{eqnarray}
where the first equality follows again from (\ref{eqeqdistINF}) and
the second equality follows by using the $(\F_t)$-Markov property of
$X$ and considering the ``sandwich'' of \mbox{$\sigma$-}fields:
\[
\sigma(X_{\tau_j})\subseteq\sigma\bigl(X^{m,1}_{\tau_j},
\ldots,X^{m,k}_{\tau_j} \bigr)\vee\ms N\subseteq
\F_{\tau_j}\vee\sigma\bigl(U^X_1,
\ldots,U^X_j \bigr)
\]
with $\ms N$ being the collection of $\P$-null sets,
and the last equality (\ref{eqsig2}) follows since $U^X_{j+1}$ is not
yet used in the construction of $\mc X^m$ up to time $\tau_j$. Hence,
by (\ref{eqfuncY}) and (\ref{eqsig2}), we can continue our
calculation in (\ref{eqsig1}) as follows:
\begin{eqnarray*}
&&\P\bigl( \bigl(X^{m,i}_r \bigr)_{r\in[\tau_j,\tau_{j+1}]}\in
A_i, \bigl(X^{m,i}_r \bigr)_{r\in[s_i,\tau_{j}]}\in
B_i,\forall i\in\{1,\ldots,k\} \bigr)
\\
&&\qquad =\E\bigl[\P\bigl( \bigl(X^{m,i}_r \bigr)_{r\in[\tau_j,\tau
_{j+1}]}
\in A_i, \forall i\in\{1,\ldots,k\}\mid X^{m,1}_{\tau
_j},
\ldots,X^{m,k}_{\tau_j} \bigr);
\\
&&\hspace*{123pt}\qquad \bigl(X^{m,i}_r \bigr)_{r\in[s_i,\tau_{j}]}\in
B_i,\forall i\in\{1,\ldots,k\} \bigr]
\\
&&\qquad =\E\Biggl[\prod_{i=1}^k{\SBM}_{X^{m,i}_{\tau_j}, [\tau
_j,\tau_{j+1} ]}(A_i); \bigl(X^{m,i}_r
\bigr)_{r\in[s_i,\tau
_{j}]}\in B_i,\forall i\in\{1,\ldots,k\} \Biggr],
\end{eqnarray*}
where the second equality follows from (\ref{eqYwt2}). By (\ref
{eqextlaw}) and induction, the foregoing equality implies that (\ref
{eqextlaw}) with $j$ replaced by $j+1$ still holds. This completes our
inductive construction for the case (\ref{eqtausk}).

We also need to handle the case complementary to (\ref{eqtausk}) that
$[\tau_j,\tau_{j+1}]\subseteq(s_k,s_{k+1}]$ and $\tau
_{j+1}=s_{k+1}$ for some $k\geq2$. In this case, the construction of
$X^{m,1},\ldots,X^{m,k}$ over the time interval $[\tau_j,\tau
_{j+1}]$ is the same as before, but the extra coordinate $X^{m,k+1}$ is
defined to be $\psi(\1)J^{x_{k+1}}_{{\varepsilon}}$ at time $\tau
_{j+1}=s_{k+1}$. The properties (\ref{eqeqdistINF}) and (\ref
{eqeqdistMASS}) with $j$ replaced by $j+1$ remain true, by the
argument for the previous case. The property (\ref{eqextlaw}) with
$j$ replaced by $j+1$ follows too, if we notice that
the coordinate $X^{m,k+1}$ is independent of the others by time $\tau
_{j+1}$ by (iii) of Definition~\ref{defstd}. This completes our
inductive construction of $\mc X^{m}$.

The\vspace*{1pt} construction of $\mc Y^{m}$ is similar to that of $\mc X^{m}$. We
use $\{U^Y_j\}$ to validate decompositions, and
the immigration times $\{t_j;j\in\mathbb N\}$ are taken into
consideration for the construction instead. We omit other details.

We observe that $\mc X^m$ and $\mc Y^m$ satisfy properties analogous to
(\ref{coneilc}) and (\ref{Vind}).
First, from the constructions of $\mc X^m$ and $\mc Y^m$, (\ref
{eqUXUY}), and the property (iii) in Definition~\ref{defstd}, we see
that the following independent landing property is satisfied by $\mc
X^m$ and $\mc Y^m$:
%
%
%e2.35 #&#
\begin{eqnarray}
\label{coneilc-0} \sigma\bigl(X^{m,j}_s,Y^{m,j}_s
; s<s_i, j\in\mathbb N \bigr)&\ind& x_i\quad\mbox{and}
\nonumber\\[-8pt]\\[-8pt]\nonumber
\sigma\bigl(X^{m,j}_s,Y^{m,j}_s
; s<t_i, j\in\mathbb N \bigr)&\ind& y_i\qquad \forall i\in\mathbb N.
\end{eqnarray}
Second, since for all $j\in\mathbb N$ satisfying $s_j<s_{i+1}$ and $i\in
\mathbb N$, $ (X^{m,j}_r )_{r\in[0,s_{i+1})}$ is given by a
measurable function of the random elements $(X_r)_{r\in[0,s_{i+1})}$
and $\{U^{X}_k\}$ [cf.~(\ref{eqfuncY}) for the case (\ref{eqtausk})
and use the path regularity of $X$ in the complementary case],
we deduce from Theorem~\ref{propsbst}(i) and (\ref{eqUXUY}) that
%
%
%e2.36 #&#
\begin{eqnarray}
\label{Vindm} \qquad \bigl\{ \bigl(X_r^{m,j}
\bigr)_{r\in[0,s_{i+1})};j\in\mathbb N\mbox{ satisfying }s_j<s_{i+1}
\bigr\} &\ind& V^{X,i}\quad\mbox{and}
\nonumber\\[-8pt]\\[-8pt]\nonumber
\bigl\{ \bigl(Y_r^{m,j} \bigr)_{r\in[0,t_{i+1})};j\in\mathbb
N\mbox{ satisfying }t_j<t_{i+1} \bigr\} &\ind&
V^{Y,i}\qquad\forall i\in\mathbb N.
\end{eqnarray}
\end{longlist}

\begin{longlist}[\textit{Step} 2.]
\item[\textit{Step} 2.]
Let us define a filtration $(\G^{(m)}_t)$ with respect to which the
processes $X^{m,i}$, $Y^{m,i}$, and $W$ perform their defining
properties \emph{on the grid $D_m$}.
The filtration $(\G_t^{(m)})$ is larger than $(\F_t)$ and is defined by
\begin{eqnarray*}
\cases{ \G_t^{(m)}=\F_t, &\quad$t\in[0,
\tau_0]$,
\cr
\G_t^{(m)}=\F_{\tau_{j+1}}\vee
\sigma\bigl(U_k^X,U^Y_k;1\leq k
\leq j+1 \bigr), &\quad$t\in(\tau_j,\tau_{j+1}], j\in\mathbb
Z_+$.}
\end{eqnarray*}
In particular, it follows from (\ref{eqeqdistINF}) and the analogue
for $\mc Y^m$ that the processes $X^{m,i}$ and $Y^{m,i}$ are all $(\G
_t^{(m)})$-adapted.\vspace*{1pt}
Also, it is obvious that $X$, $Y$ and $W(\phi)$ for any $\phi\in
L^2(\R)$ are $(\G_t^{(m)})$-adapted.

We observe a key feature of $\mc X^m$:
%
%
%e2.37 #&#
%e2.38 #&#
\begin{eqnarray}
\label{eqgridMP}
\P\bigl(X^{m,i}_t\in\Gamma\mid
\G_{\tau_j}^{(m)} \bigr)&=&{\SBM}_{t-\tau_j}
\bigl(X^{m,i}_{\tau_j},\Gamma\bigr)
\nonumber\\[-8pt]\\[-8pt]
\eqntext{\forall t\in(
\tau_j,\tau_{j+1}]\mbox{ for } s_i\leq
\tau_j\mbox{ and } i\in\mathbb N,}
\end{eqnarray}
for any Borel measurable subset $\Gamma$ of the space of finite
measures on $\R$. To see (\ref{eqgridMP}), we consider a slight
generalization of the proof of (\ref{eqsig2}) by adding $\sigma
(U^Y_1,\ldots,U^Y_j)$ to the $\sigma$-field $\F_{\tau_j}\vee\sigma
(U^X_1,\ldots,U^X_j)$ in the first line therein and then apply (\ref
{eqextlaw}) to obtain
%
%e2.39 #&#
\begin{eqnarray}
\P\bigl(X_t^{m,i}\in\Gamma\mid\G_{\tau_j}^{(m)}
\bigr)&=&\P\bigl(X_t^{m,i}\in\Gamma\mid X^{m,1}_{\tau_j},X^{m,2}_{\tau
_j},\ldots\bigr)\nonumber
\\
&=&\P\bigl(X^{m,i}_t\in\Gamma\mid X^{m,i}_{\tau_j}
\bigr)
={\SBM}_{t-\tau_j} \bigl(X^{m,i}_{\tau_j},\Gamma\bigr)\nonumber
\\
\eqntext{\qquad\forall t\in(\tau_j,\tau_{j+1}].}
\end{eqnarray}
In particular, we deduce from iteration and the semigroup property of
$\{{\SBM}_t\}$ that
the following \textit{grid Markov property} is satisfied:
%
%
%e2.40 #&#
%e2.41 #&#
\begin{eqnarray}
\label{eqgridMPX} \P\bigl(X^{m,i}_t\in\Gamma\mid
\G^{(m)}_{\tau_j} \bigr)={\SBM}_{t-\tau_j}
\bigl(X^{m,i}_{\tau_j},\Gamma\bigr)
\nonumber\\[-8pt]\\[-8pt]
\eqntext{\forall t\in(
\tau_k,\tau_{k+1}]\mbox{ when }s_i\leq
\tau_j\leq\tau_k.}
\end{eqnarray}
We note that the foregoing display does \emph{not} say that $X^{m,i}$
is a $(\G_s^{(m)})_{s\geq s_i}$-super-Brownian motion because the
$\sigma$-fields which we can use in verifying the $(\G
_s^{(m)})_{s\geq s_i}$-Markov property are only $\G^{(m)}_{\tau_j}$,
rather than \emph{any} $\sigma$-field $\G_s^{(m)}$.
With a similar argument, we also have the grid Markov property of
$Y^{m,i}$ stated as
%
%
%e2.42 #&#
%e2.43 #&#
\begin{eqnarray}
\label{eqgridMPY} \P\bigl(Y^{m,i}_t\in\Gamma\mid
\G^{(m)}_{\tau_j} \bigr)={\SBM}_{t-\tau_j}
\bigl(Y^{m,i}_{\tau_j},\Gamma\bigr)
\nonumber\\[-8pt]\\[-8pt]
\eqntext{\forall t\in(
\tau_k,\tau_{k+1}]\mbox{ when }t_i\leq
\tau_j\leq\tau_k.}
\end{eqnarray}
With a much simpler argument, the space--time white noise $W$ has the
same grid Markov property:
%
%
%e2.44 #&#
%e2.45 #&#
\begin{eqnarray}
\label{eqgridMPW} \ms L \bigl(W_t(\phi)\mid\G^{(m)}_{\tau_j}
\bigr)=\mathcal N \bigl(W_{\tau_j}(\phi),(t-\tau_j)\llVert
\phi\rrVert^2_{L^2(\R)} \bigr)
\nonumber\\[-8pt]\\[-8pt]
\eqntext{\forall t\in(\tau_k,\tau_{k+1}]\mbox{ for }
\tau_j\leq\tau_k\mbox{ and }\phi\in L^2(\R),}
\end{eqnarray}
where $\mathcal N(\mu,\sigma^2)$ denotes the normal distribution with
mean $\mu$ and variance $\sigma^2$. Similar results hold for the
substituting space--time white noises $V^{X,i}$ and $V^{Y,i}$.
\end{longlist}

\begin{longlist}[\textit{Step} 3.]
\item[\textit{Step} 3.]
To facilitate our argument in the next step, we digress to a general
property of space--time white noises.

Let $W^1$ denote a space--time white noise, and suppose that $\{W^2(\phi
_n)\}$ is a family of Brownian motions indexed by a countable dense
subset $\{\phi_n\}$ of $L^2(\R)$ such that $\{W^1(\phi_n)\}$ and $\{
W^2(\phi_n)\}$ have the same law as random elements taking values in
$\prod_{n=1}^\infty C(\R_+,\R)$. Then,\vspace*{1.5pt}
whenever $(\phi_{n_k})$ is a subsequence converging to some $\phi$ in
$L^2(\R)$, the linearity of $W^1$ gives
%
%
%e2.46 #&#
%e2.47 #&#
\begin{eqnarray}\label{eqCauchyconv} \E\Bigl[\sup_{0\leq s\leq T}\bigl\llvert
W^2_s(\phi_{n_k})-W^2_s(
\phi_{n_\ell})\bigr\rrvert^2 \Bigr] &=&\E\Bigl[\sup
_{0\leq s\leq T}\bigl\llvert W^1_s(
\phi_{n_k}-\phi_{n_\ell
})\bigr\rrvert^2 \Bigr]
\nonumber
\\
&\leq&4T\llVert\phi_{n_k}-\phi_{n_\ell}\rrVert
_{L^2(\R)}^2\xrightarrow{k,\ell\to\infty} 0
\\
\eqntext{\forall T\in(0,\infty),}
\end{eqnarray}
where the inequality follows from Doob's $L^2$-inequality and the fact
that, for any $\phi\in L^2(\R)$, $W^1(\phi)$ is a Brownian motion with
$\ms L (W^1_1(\phi) )=\mathcal N (0,\llVert \phi\rrVert^2_{L^2(\R
)} )$.
The convergence in (\ref{eqCauchyconv}) implies that, for some
continuous process, say $W^2(\phi)$, we have
\[
W^2(\phi_{n_k})\lra W^2(\phi)\qquad
\mbox{uniformly on }[0,T]\mbox{ a.s., }\forall T\in(0,\infty).
\]
The\vspace*{1pt} same holds with $W^2$ replaced by $W^1$.
Hence, making comparisons with the reference space--time white noise
$W^1$, we obtain an extension of the map
$\phi\lmt W^2(\phi)$ to the entire space $L^2(\R)$ such that
$\{W^2(\phi);\phi\in L^2(\R)\}$ is a space--time white noise and, in
fact, is uniquely defined by $\{W^2(\phi_n)\}$.
\end{longlist}

\begin{longlist}[\textit{Step} 4.]
\item[\textit{Step} 4.]
In this step, we formalize the infinitesimal description outlined
before by shrinking the mesh size of $D_m$, that is, by passing $m\lra
\infty$, and then work with the limiting objects. To use our
observation in step~3, we work with a fixed countable dense subset $\{
\phi_n\}$ of $L^2(\R)$.

We have constructed in step~1 random elements $\mc X^m$ and $\mc
Y^m$, and hence determined the laws
%
%
%e2.48 #&#
\begin{equation}\label{disjlaw}
\ms L \bigl(X,Y,W,\bigl\{V^{X,i}\bigr\},\bigl\{V^{Y,i}\bigr\},
\mc X^m,\mc Y^m,\{x_i\},\{y_i\}
\bigr),\qquad m\in\mathbb N,
\end{equation}
as probability measures
on a countably infinite product of Polish spaces. More precisely, our
choice of the Polish spaces is through the following identifications of
state spaces.
We identify $X$ as a random element taking values in the closed subset
of $D (\R_+,\C_\rap(\R) )$ consisting of paths having
continuity over each interval $[s_i,s_{i+1})$ for $i\in\mathbb Z_+$
(recall $s_0=0$), with a similar identification applied to $Y$ (cf.
Proposition 5.3 and Remark 5.4 of \cite{EK}).
By step~3, we identify $W$ as the infinite-dimensional vectors $
(W(\phi_1),W(\phi_2),\ldots)$ whose coordinates\vspace*{1pt} are $C(\R
_+,\R)$-valued random elements. Similarly, $V^{X,i}$ and $V^{Y,i}$ are
infinite-dimensional vectors of $C([s_i,\infty),\R)$ and
$C([t_i,\infty),\R)$-valued random elements.
We identify each coordinate $X^{m,i}$ of $\mathcal X^m$ as a random
element taking values in $C ([s_i,\infty),\C_\rap(\R) )$,
with a similar identification applied to $\mathcal Y^m$. Finally, the
Polish spaces for the infinite sequences $\{x_i\}$ and $\{y_i\}$ are obvious.

We make an observation for the sequence of laws in (\ref{disjlaw}).
Note that $\ms L(\mc X^m)$ does not depend on $m$, because, by (\ref
{eqextlaw}), any of its $i$th coordinate $X^{m,i}$ is a super-Brownian
motion with initial measure $\psi(\1)J^{x_i}_{\varepsilon}$ and
started at $s_i$, and the coordinates are independent. Similarly, $\ms
L(\mc Y^m)$ does not depend on $m$. This implies that the sequence of
laws in (\ref{disjlaw})
is tight in the space of probability measures on the above infinite
product of Polish spaces. Hence, by taking a subsequence if necessary,
we may assume that this sequence converges in distribution. By
Skorokhod's representation, we may assume the existence of the vectors
of random elements in the following display as well as the almost-sure
convergence therein:
%
%
%e2.49 #&#
\begin{eqnarray}
\label{eqlimD}  && \bigl(\widetilde{X}{}^{(m)},\widetilde{Y}{}^{(m)},
\widetilde{W}{}^m,\bigl\{ \widetilde{V}{}^{X,i,m}\bigr\},\bigl\{
\widetilde{V}{}^{Y,i,m}\bigr\},\widetilde{\mc X}^m,\widetilde{
\mc Y}^m,\bigl\{\widetilde{x}_i^m\bigr\},
\bigl\{\widetilde{y}_i^m\bigr\} \bigr)
\nonumber\\[-8pt]\\[-8pt]\nonumber
&&\qquad \xrightarroww{m\to\infty} {\mathrm{a.s.}} \bigl(\widetilde
{X},\widetilde{Y},
\widetilde{W},\bigl\{\widetilde{V}{}^{X,i}\bigr\},\bigl\{\widetilde
{V}{}^{Y,i}\bigr\},\widetilde{\mc X},\widetilde{\mc Y},\{
\widetilde{x}_i\},\{\widetilde{y}_i\} \bigr).
\end{eqnarray}
Here, $\widetilde{x}_i$ and $\widetilde{y}_i$ take values in the
topological support of $\psi$ and
\begin{eqnarray*}
&& \ms L \bigl(\widetilde{X}{}^{(m)},\widetilde{Y}{}^{(m)},
\widetilde{W}{}^m,\bigl\{\widetilde{V}{}^{X,i,m}\bigr\},\bigl\{
\widetilde{V}{}^{Y,i,m}\bigr\},\widetilde{\mc X}^m,\widetilde{
\mc Y}^m,\bigl\{\widetilde{x}_i^m\bigr\},
\bigl\{\widetilde{y}_i^m\bigr\} \bigr)
\\
&&\qquad =\ms L \bigl(X,Y,W,\bigl\{V^{X,i}\bigr\},\bigl
\{V^{Y,i}\bigr\},\mc X^m,\mc Y^m,
\{x_i\},\{y_i\} \bigr)\qquad\forall m\in\mathbb N.
\end{eqnarray*}
\end{longlist}

\begin{longlist}[\textit{Step} 5.]
\item[\textit{Step} 5.]
We define\vspace*{1pt} $(\widetilde{\G}_t)$ to be the minimal
filtration satisfying the usual conditions to which the limiting objects
$\widetilde{X},\widetilde{Y},\widetilde{W},\{\widetilde{V}{}^{X,i}\},\{\widetilde{V}{}^{Y,i}\},\widetilde{\mc X},\widetilde{\mc Y}$
on the right-hand side of (\ref{eqlimD})
are adapted. We will complete the proof in this step by verifying that,
with an obvious adaptation of notation, all the limiting objects on the
right-hand side of (\ref{eqlimD}) along with the filtration
$(\widetilde{\G}_t)$ are the required objects satisfying conditions
(i)--(vi) of Theorem~\ref{propequiallocation}.

First, let us verify the easier properties (i) and (ii) for $\{
\widetilde{X}{}^i\}$ and the analogues for $\{\widetilde{Y}{}^i\}$. The
statement (i) and its analogue for $\{\widetilde{Y}{}^i\}$ obviously
hold, by the analogous properties of $\widetilde{\mc X}^m$ and
$\widetilde{\mc Y}^m$ [see (b) of (\ref{eqextlaw})]. To verify the
statement (ii), we use the property (\ref{eqeqdistMASS}) possessed by
$ (\widetilde{X}{}^{(m)},\widetilde{\mc X}^m )$ and then pass
limit, as is legitimate because the infinite series in (\ref
{eqeqdistMASS}) are always finite sums on compact time intervals.
Similarly, the analogue of (ii) holds for $(\widetilde{Y},\widetilde
{\mc Y})$.

Condition\vspace*{1pt} (iii) holds by the property (a) in (\ref{eqextlaw}) of
$\widetilde{\mc X}^m$, except that we still need to verify that each
$\widetilde{X}{}^{i}$ defines a $(\widetilde{\G}_t)_{t\geq s_i}$-super
Brownian motion, not just a super-Brownian motion in itself. From this
point on, we will use the continuity of the underlying objects and the
fact that $\bigcup_m D_m$ is dense in $[\frac{{\varepsilon
}}{2},\infty)$. Let $\frac{{\varepsilon}}{2}\leq s<t<\infty$ with
$s,t\in\bigcup_m D_m$. Then $s,t\in D_m$ from some large $m$ on by
the nesting property of the sequence $\{D_m\}$. For any bounded
continuous function $g$ on the path space of
\[
\bigl(\widetilde{X}{}^{(m)},\widetilde{Y}{}^{(m)},
\widetilde{W}{}^{m},\bigl\{ \widetilde{V}{}^{X,i,m}\bigr\},\bigl\{
\widetilde{V}{}^{Y,i,m}\bigr\},\widetilde{\mc X}^m,\widetilde{
\mc Y}^m \bigr)
\]
restricted to the time interval $[0,s]$, $\phi\in\C_c^+(\R)$, and
index $i$ such that $s_i\leq s$, the grid Markov property (\ref
{eqgridMPX}) entails that
%
%
%e2.50 #&#
\begin{eqnarray}
\label{eqLap} && \E\bigl[g \bigl(\widetilde{X}{}^{(m)},
\widetilde{Y}{}^{(m)},\widetilde{W}{}^{m},\bigl\{
\widetilde{V}{}^{X,i,m}\bigr\},\bigl\{\widetilde{V}{}^{Y,i,m}\bigr\},
\widetilde{\mc X}^m,\widetilde{\mc Y}^m
\bigr)e^{-\langle
\widetilde{X}{}^{(m),i}_t,\phi\rangle} \bigr]\nonumber
\\
&&\qquad =\E\biggl[g \bigl(\widetilde{X}{}^{(m)},\widetilde
{Y}{}^{(m)},\widetilde{W}{}^{m},\bigl\{\widetilde{V}{}^{X,i,m}
\bigr\},\bigl\{\widetilde{V}{}^{Y,i,m}\bigr\},\widetilde{\mc
X}^m,\widetilde{\mc Y}^m \bigr)
\\
&&\hspace*{84pt}\quad\qquad{} \times\int{\SBM}_{t-s} \bigl(\widetilde
{X}{}^{(m),i}_{s},d\nu\bigr)e^{-\langle\nu,\phi\rangle} \biggr].\nonumber
\end{eqnarray}
The formula of Laplace transforms of super-Brownian motion shows that
the map
\[
f\lmt\int{\SBM}_{t-s} (f,d\nu)e^{-\langle\nu,\phi
\rangle}
\]
has a natural extension to $\C_\rap(\R)$ which is continuous (cf.
Proposition~II.5.10 of~\cite{PDW}). Hence, passing $m\lra\infty$
for both sides of (\ref{eqLap}) leads to
%
%
%e2.51 #&#
\begin{eqnarray}
\label{eqLap1} && \E\bigl[g \bigl(\widetilde{X},\widetilde{Y},\widetilde
{W},\bigl
\{ \widetilde{V}{}^{X,i}\bigr\},\bigl\{\widetilde{V}{}^{Y,i}\bigr\},
\widetilde{\mc X},\widetilde{\mc Y} \bigr)e^{- \langle\widetilde
{X}{}^{i}_t,\phi\rangle} \bigr]\nonumber
\\
&&\qquad =\E\biggl[g \bigl(\widetilde{X},\widetilde{Y},\widetilde{W},\bigl
\{\widetilde{V}{}^{X,i}\bigr\},\bigl\{\widetilde{V}{}^{Y,i}\bigr\},
\widetilde{\mc X},\widetilde{\mc Y} \bigr)
\\
&&\hspace*{66pt}{}\times\int{\SBM}_{t-s} \bigl(
\widetilde{X}{}^{i}_{s},d\nu\bigr)e^{-\langle\nu,\phi\rangle
} \biggr].\nonumber
\end{eqnarray}
By the continuity of super-Brownian motion and the denseness of
$\bigcup_m D_m$ in $[\frac{{\varepsilon}}{2},\infty)$, the
foregoing display implies that
each coordinate $\widetilde{X}{}^{i}$ is truly a $(\widetilde{\G
}_t)_{t\geq s_i}$-super-Brownian motion. A similar argument shows that
each $\widetilde{Y}{}^{i}$ is\vspace*{2pt} a $(\widetilde{\G}_t)_{t\geq
t_i}$-super-Brownian motion. We have proved the statement (iii) and its
analogue for~$\widetilde{Y}{}^i$ in~(iv).

Next, we consider the assertions of (v) concerning conditions analogous
to (i) and (ii) of Definition~\ref{defstd}. By definition,
%
%
%e2.52 #&#
\begin{eqnarray}
\label{stationary} && \ms L \bigl(\widetilde{X}{}^m,\widetilde{Y}{}^m,
\widetilde{W}{}^m,\bigl\{ \widetilde{V}{}^{X,i,m}\bigr\},\bigl\{
\widetilde{V}{}^{Y,i,m}\bigr\},\bigl\{\widetilde{x}_i^m
\bigr\},\bigl\{\widetilde{y}_i^m\bigr\} \bigr)
\nonumber\\[-8pt]\\[-8pt]\nonumber
&&\qquad =\ms L \bigl(X,Y,W,\bigl\{V^{X,i}\bigr\},\bigl
\{V^{Y,i}\bigr\},\{x_i\},\{y_i\} \bigr)\qquad
\forall m\in\mathbb N,
\end{eqnarray}
and this stationarity gives
%
%
%e2.53 #&#
\begin{eqnarray}
\label{eqidlaw}
&& \ms L \bigl(\widetilde{X},\widetilde{Y},\widetilde
{W},\bigl\{
\widetilde{V}{}^{X,i}\bigr\},\bigl\{\widetilde{V}{}^{Y,i}\bigr\},\{
\widetilde{x_i}\},\{\widetilde{y}_i\} \bigr)
\nonumber\\[-8pt]\\[-8pt]\nonumber
&&\qquad =\ms L
\bigl(X,Y,W,\bigl\{V^{X,i}\bigr\},\bigl\{V^{Y,i}\bigr\},
\{x_i\},\{y_i\} \bigr).
\end{eqnarray}
Arguing\vspace*{1.5pt} as in the proof of (\ref{eqLap}) and using the grid Markov
property (\ref{eqgridMPW}) of~$\widetilde{W}{}^m$, we see that each
$\widetilde{W}(\phi_n)$ is a $(\widetilde{\G}_t)$-Brownian motion with
\[
\ms L \bigl(\widetilde{W}_1(\phi_n) \bigr)=\mathcal N
\bigl(0,\llVert\phi_n\rrVert_{L_2(\R)}^2 \bigr).
\]
It follows from (\ref{eqidlaw})
and our discussion in step~3 that $\widetilde{W}$ extends uniquely
to a $(\widetilde{\G}_t)$-space--time white noise.
In addition, one more application of (\ref{eqidlaw}) shows that the
defining\vspace*{1pt} equations (\ref{eqXvep}) and (\ref{eqYvep}) of $X$ and $Y$
by $\{(x_i,y_i)\}$ and $W$ carry over to the analogous equations for
$\widetilde{X}$ and $\widetilde{Y}$ by $\{(\widetilde
{x}_i,\widetilde{y}_i)\}$ and $\widetilde{W}$, respectively [recall
(\ref{impilc}) as well]. This proves that $ (\widetilde
{X},\widetilde{Y},\widetilde{W} )$ satisfies the analogous
property described in (i) and (ii) of Definition~\ref{defstd} with
$(\F_t)$ replaced by $(\widetilde{\G}_t)$.

By construction, $\widetilde{x}_i$ and $\widetilde{y}_i$ take values
in the topological support of $\psi$. Hence, to complete the proof of
(v), it remains to
obtain the independent landing property~(\ref{coneilc}). We\vspace*{2pt} recall
that an analogous property is satisfied by $ (\widetilde{\mc
X}^m,\widetilde{\mc Y}^m,\{\widetilde{x}_i^m\},\{\widetilde{y}_i^m\}
)$ in (\ref{coneilc-0}). Then arguing in the standard way as in
the proof of (\ref{eqLap}) with the use of bounded continuous
functions shows that the required independent landing property (\ref
{coneilc}) is satisfied by $ (\widetilde{\mc X},\widetilde{\mc
Y},\{\widetilde{x}_i\},\{\widetilde{y}_i\} )$. The proof of~(v)
is complete.

Finally, we explain the proof of the assertions in (vi). The proof of
(iii) of Proposition~\ref{propsbst} uses again the stationarity (\ref
{stationary}). The assertion that (\ref{Vind}) holds follows from its
discrete version (\ref{Vindm}) and a limiting argument as in the proof
of (iii). We have proved that all of the conditions (i)--(vi) of
Theorem~\ref{propequiallocation} hold. The proof is complete.\quad\qed
\end{longlist}\noqed
\end{pf*}

%s2.3 #&#
\subsection{Covariations of immigrant processes}\label{seccovar-imm}
In this section, we study the covariations
$\langle X^i(\phi_1),Y^j(\phi_2)\rangle$ of the immigrant processes
$\{X^i\}$ and $\{Y^i\}$ in Theorem~\ref{propequiallocation} for
${\varepsilon}\in(0,1]$. Here, the test functions $\phi_1,\phi_2$
belong to $\C_c^\infty(\R)$. Our goal is to understand how explicit
$\langle X^i(\phi_1),Y^j(\phi_2)\rangle$ can be in terms of the
immigrant processes.

For convenience, we attach space--time white noises to the immigrant
processes $\{X^i\}$ and $\{Y^i\}$.
Recall that by (iii) of Theorem~\ref{propequiallocation},
$(X^i_t )_{t\in[s_i,\infty)}$ is\vspace*{1pt} a $(\G_t)_{t\geq
s_i}$-super-Brownian motion for any $i\in\mathbb N$, and similarly, each
$ (Y^i_t )_{t\in[t_i,\infty)}$ is a $(\G_t)_{t\geq
t_i}$-super-Brownian motion.
By a classical argument,
we can find, by enlarging the filtered probability space if necessary,
two families of $(\G_t)$-white noises $ \{W^{X^i} \}$ and
$ \{W^{Y^i} \}$ such that $ (X^i,W^{X^i} )$ and $
(Y^i,W^{Y^i} )$ are solutions to the SPDE (\ref{eqmainSPDE0}) of
super-Brownian motion (up to appropriate translations of starting
time). See Theorem III.4.2 of \cite{PDW} for details. Moreover, by
(i) and (vi) of Theorem~\ref{propequiallocation}, we can assume that
each of the families $ \{W^{X^i} \}$ and $ \{W^{Y^i} \}
$ consists of independent adapted space--time white noises, and in
addition, the following independence holds:
%
%
%e2.54 #&#
\begin{eqnarray}
\label{WVind} \bigl\{X^i_t,W^{X^i}_{t}(
\phi);t\in[0,s_{j+1}),1\leq i\leq j,\phi\in L^2(\R) \bigr
\} &\ind& V^{X,j}\quad\mbox{and}
\nonumber\\[-8pt]\\[-8pt]\nonumber
\bigl\{Y^i_t,W^{Y^i}_{t}(\phi);t
\in[0,t_{j+1}),1\leq i\leq j,\phi\in L^2(\R) \bigr\} &\ind&
V^{Y,j}\qquad\forall j\in\mathbb N,\hspace*{-40pt}
\end{eqnarray}
where $V^{X,j}$ and $V^{Y,j}$ are the adapted space--time white noises
which substitute $W$ and satisfy (\ref{sbstX}) and (\ref{sbstY}).

%
%
%re2.7 #&#
\begin{rmk}\label{rmkmistake}
Let us point out an issue for the covariations $\langle X^i(\phi
_1),\break Y^j(\phi_2)\rangle$. An implication of the classical
Kunita--Watanabe inequality (cf. Proposition~IV.1.15 of \cite{RYCMB})
is that for any two $(\G_t)$-Brownian motions $B^1,B^2$ and
nonnegative locally bounded predictable processes $H^1$ and $H^2$, the
covariation of the (ordinary) stochastic integrals $H^i\bullet B^i$ of
$H^i$ with respect to $B^i$ satisfies
%
%
%e2.55 #&#
\begin{equation}
\label{BMcovar} \bigl\llvert\bigl\langle H^1\bullet
B^1,H^2\bullet B^2\bigr\rangle_t
\bigr\rrvert\leq\int_0^t H^1_sH^2_s\,ds
\end{equation}
(recall $d\langle B^i,B^i\rangle_s\equiv ds$).
On the other hand, let $W^1$ and $W^2$ be two $(\G_t)$-space--time
white noises, and $J^1,J^2\in L^2_\loc(W^1)=L^2_\loc(W^2)$ be
nonnegative. [Recall the notation (\ref{eqL2W}) and (\ref{bullet}).]
In this case, the measure $dx\,ds$ determines quadratic variations of
stochastic integrals with respect to a space--time white noise in the
sense that
$\langle J^i\bullet W^i(\1),J^i\bullet W^i(\1)\rangle_t=\int_0^t\int_\R
J^i(x,\break s)^2\,dx\,ds$.
In bounding the covariations of $J^1\bullet W^1(\1)$ and $J^2\bullet
W^2(\1)$, however, the following inequality, analogous to (\ref{BMcovar}), is
not always true:
%
%
%e2.56 #&#
\begin{equation}
\label{specialcovar} \bigl\llvert\bigl\langle J^1\bullet
W^1(\1),J^2\bullet W^2(\1)\bigr
\rangle_t\bigr\rrvert\leq\int_0^t\!\int_\R J^1(x,s)J^2(x,s)\,dx\,ds.
\end{equation}
A counterexample is given by taking $W^2$ to be a nonidentity spatial
translation of $W^1$.
Hence, the conclusion of Proposition 3.5 in \cite{C} is incorrect in
general, as pointed by the anonymous referee, and there it is used to
bound covariations of general adapted immigrant processes.
Nevertheless, we will show in Section~\ref{seccappro} that there do
exist \emph{some} immigrant processes whose covariations satisfy the
concluding inequality of Proposition 3.5 in \cite{C} (see Theorem~\ref
{teoapprox-v2} and Proposition~\ref{propcovarbdd}), and so the
arguments from Section~3.5 on in \cite{C} remain valid if these
particular immigrant processes are in force.
%\qed
\end{rmk}

To facilitate the forthcoming computation of covariations, we give some
terminology and notation. For a locally bounded signed measure $\mu$
on $\R\times\R_+$, we define a measure-valued process $J\bullet\mu
$ by
\[
J\bullet\mu(\phi)\triangleq\int_{(0, \cdot]}\int
_{\R}J(x,s)\phi(x)\,d\mu(x,s),\qquad\phi\in
\C_c^\infty(\R),
\]
whenever $J$ is a two-parameter random function satisfying $\int
_{(0,t]}\int_\R\hspace*{-0.5pt}\llvert J(x,s)\times\break \phi(x)\rrvert \llvert d\mu
(x,s)\rrvert <\infty$ a.s. for all
$t$ and $\phi\in\C_c^\infty(\R)$,
and we put $\mu(\phi)\equiv\1\bullet\mu(\phi)$. Let $U^i,V^i$ be
$(\G_t)$-space--time white noises and $J^i,K^i\in L^2_\loc
(U^i)=L^2_\loc(V^i)$ for $1\leq i\leq N$ and a natural number $N$.
Then the pair $(\sum_iJ^i\bullet U^i,\sum_iK^i\bullet V^i)$ of finite
sums of stochastic integrals is said to have a \textit{normal
covariation} if
%
%
%e2.57 #&#
%e2.58 #&#
\begin{eqnarray}
\label{eqcovarmeas} \biggl\langle\sum_iJ^i
\bullet U^i(\phi_1),\sum_iK^i
\bullet V^i(\phi_2) \biggr\rangle= H\bullet\lambda(
\phi_1\phi_2)
\nonumber\\[-10pt]\\[-10pt]
\eqntext{\mbox{a.s. }\forall
\phi_1,\phi_2\in\C_c^\infty(\R)}
\end{eqnarray}
for some $H\in L^2_\loc(U^i)$, where the measure $\lambda$ in (\ref
{eqcovarmeas}) is defined by
%
%
%e2.59 #&#
\begin{equation}
\label{deflambda} d\lambda(x,s)\triangleq dx\,ds.
\end{equation}
In this case, we write
$\llang\sum_iJ^i\bullet U^i,\sum_iK^i\bullet V^i\rrang$ for
$H\bullet\lambda$. If $Z=Z(x,t)$ and $Z'=Z'(x,t)$ are solutions to
SPDEs with stochastic integral terms characterized by $\sum_iJ^i\bullet
U^i$ and $\sum_i K^i\bullet V^i$,
respectively, then the pair $(Z,Z')$ is also said to have a normal
covariation and we write $\llang Z,Z'\rrang$ for $\llang\sum
_iJ^i\bullet U^i,\sum_i K^i\bullet V^i\rrang$.

%
%
%pr2.8 #&#
\begin{prop}\label{propcovar}
Fix ${\varepsilon}\in(0,1]$ and an interlacing pair $(X,Y)$ of
${\varepsilon}$-approxi\-mating solutions. Let $\{X^i\}$ and $\{Y^i\}$
be immigrant processes as in Theorem~\ref{propequiallocation} and let
$\{W^{X^i}\}$ and $\{W^{Y^i}\}$ be the auxiliary space--time white
noises chosen before Remark~\ref{rmkmistake}. Then for all $i,j\in
\mathbb N$, $(X^i,X)$ and $(X,Y^i)$ have normal covariations, and we have
%
%
%e2.60 #&#
\begin{eqnarray}
\label{eqcovardisc} \llanga X^i,Y\rranga&=&\1_{(X>0,Y>0)}
\bigl(X^i\bigr)^{1/2} \biggl(\frac
{X^i}{X}
\biggr)^{1/2}Y^{1/2}\bullet\lambda\quad\mbox{and}
\nonumber\\[-8pt]\\[-8pt]\nonumber
\llanga X,Y^j\rranga&=&\1_{(X>0,Y>0)}X^{1/2}
\bigl(Y^j\bigr)^{1/2} \biggl(\frac
{Y^j}{Y}
\biggr)^{1/2}\bullet\lambda,
\end{eqnarray}
where the measure $\lambda$ is defined by (\ref{deflambda}).
\end{prop}
\begin{pf}
We only show that $(X^i,Y)$ has a normal covariation and compute
$\llang X^i,Y\rrang$, as the covariations of $(X,Y^j)$ can be handled
similarly.
For $\phi_1,\phi_2\in\C_c^\infty(\R)$, write
%
%
%e2.61 #&#
\begin{eqnarray}
\label{dec}
&& \bigl\langle\bigl(X^i\bigr)^{1/2}\bullet
W^{X^i}(\phi_1),Y^{1/2}\bullet W(\phi
_2) \bigr\rangle\nonumber
\\
&&\qquad
= \bigl\langle\bigl(X^i
\bigr)^{1/2}\bullet W^{X^i}(\phi_1),\1
_{(X>0)}Y^{1/2}\bullet W(\phi_2) \bigr\rangle
\\
&&\qquad\quad{}+ \bigl\langle\bigl(X^i\bigr)^{1/2}\bullet
W^{X^i}(\phi_1),\1 _{(X=0)}Y^{1/2}\bullet
W(\phi_2) \bigr\rangle,\nonumber
\end{eqnarray}
and consider the two terms on the right-hand side separately. For the
first one,
we turn the space--time white noise $W$ into functionals of $ \{
(X^i,W^{X^i});i\in\mathbb N \}$ over the space--time subset $(X>0)$ by writing
\begin{eqnarray*}
\1_{(X>0)}\bullet W &=& \1_{(X>0)}\frac{X^{1/2}}{X^{1/2}}\bullet W =\sum
_{j=1}^\infty\1_{(X>0)} \biggl(
\frac{X^j}{X} \biggr)^{1/2}\bullet W^{X^j},
\end{eqnarray*}
where the last equality follows from the compatibility condition
$X^{1/2}\bullet W=\sum_{j=1}^\infty(X^j)^{1/2}\bullet W^{X^j}$
(compare the stochastic integral terms of the SPDEs for $X$ and $\sum
_{j=1}^\infty X^j$). Note\vspace*{1pt} that the infinite series in the foregoing
display is well defined since there are only finite many immigration
events in a bounded time interval. The independence of the noises
$W^{X^j}$, $j\in\mathbb N$, and
the foregoing equality imply that
%
%
%e2.62 #&#
\begin{eqnarray}\label{covar1}
&& \bigl\langle\bigl(X^i\bigr)^{1/2}\bullet W^{X^i}(
\phi_1),\1 _{(X>0)}Y^{1/2}\bullet W(
\phi_2) \bigr\rangle
\nonumber\\[-8pt]\\[-8pt]\nonumber
&&\qquad =\1_{(X>0)}\frac{X^i}{X^{1/2}}Y^{1/2}
\bullet\lambda(\phi_1\phi_2).
\end{eqnarray}

Next, we consider the second term on the right-hand side of (\ref
{dec}). A standard property of one-parameter stochastic integrals
implies that
\[
\bigl\langle\bigl(X^i\bigr)^{1/2}\bullet W^{X^i}(
\phi_1),\1_{[s_j,s_{j+1})}\1 _{(X=0)}\bullet W(
\phi_2) \bigr\rangle\equiv0,\qquad0\leq j<i,
\]
since $X^i$ does not arrive before time $s_{j+1}$ for these pairs of
indices $(i,j)$. For $j\in\mathbb N$ with $j\geq i$, the corresponding
substitution identity in (\ref{sbstX}) [recall (vi) of Theorem~\ref
{propequiallocation}] gives $\1_{[s_j,s_{j+1})}\1_{(X=0)}\bullet W=\1
_{[s_j,s_{j+1})}\1_{(X=0)}\bullet V^{X,j}$. Hence, for all~$t$,
%
%
%e2.63 #&#
\begin{eqnarray}\label{covar2}
\qquad&& \bigl\langle\bigl(X^i\bigr)^{1/2}\bullet
W^{X^i}(\phi_1),\1 _{(X=0)}Y^{1/2}\bullet
W(\phi_2) \bigr\rangle_t
\nonumber
\\
&&\qquad =\sum_{j=i}^\infty\bigl\langle
\bigl(X^i\bigr)^{1/2}\bullet W^{X^i}(\phi
_1),Y^{1/2}\1_{[s_j,s_{j+1})}\1_{(X=0)}\bullet
V^{X,j}(\phi_2) \bigr\rangle_t
\nonumber\\[-8pt]\\[-8pt]\nonumber
&&\qquad =\sum_{j=i}^\infty\bigl\langle
\bigl(X^i\bigr)^{1/2}\bullet W^{X^i}(\phi
_1),Y^{1/2}\1_{[s_j,s_{j+1})}\1_{(X=0)}\bullet
V^{X,j}(\phi_2) \bigr\rangle_{t\wedge(s_{j+1})}
\\
&&\qquad =0,\nonumber
\end{eqnarray}
where the last equality follows from the independence (\ref{WVind}).
Applying (\ref{covar1}) and~(\ref{covar2}) to the right-hand side of
(\ref{dec}), we see that $(X^i,Y)$ has a normal covariation and get
the first equality of (\ref{eqcovardisc}).
The proof is complete.
\end{pf}

%s2.4 #&#
\subsection{Immigrant processes obeying a system of SPDEs}\label{seccappro}
Equations (\ref{eqcovardisc}) in Proposition~\ref{propcovar} give
partial information for the covariations between $X^i$ and $Y^j$. The
symmetry of these equations in $X$ and $Y$ suggests that if we consider
the case in which all of the pairs $(X^i,Y^j)$ have normal
covariations, then one possibility for $\{\llang X^i,Y^j\rrang;i,j\in
\mathbb N\}$ should be that the \textit{coexistence condition} is satisfied:
%
%
%e2.64 #&#
\begin{eqnarray}
\label{eqXiYj}
\bigl\langle\bigl\langle
X^i,Y^j{\bigr\rangle\bigr\rangle}&=&\1 _{(X>0,Y>0)}\bigl(X^i\bigr)^{1/2}
\bigl(Y^j\bigr)^{1/2} \biggl(\frac{X^i}{X}
\biggr)^{1/2} \biggl(\frac{Y^j}{Y} \biggr)^{1/2}\bullet
\lambda
\nonumber\\[-8pt]\\[-8pt]\nonumber
&=&\1_{(X>0,Y>0)}\frac{X^iY^j}{X^{1/2}Y^{1/2}}\bullet\lambda,\qquad
i,j\in\mathbb N.
\end{eqnarray}
Equations such as (\ref{eqXiYj}), if valid, would complement the fact
that all $(X^i,X^j)$ and $(Y^i,Y^j)$ have normal covariations and
%
%
%e2.65 #&#
\begin{equation}
\label{eqXiYj+1} \llanga X^i,X^j\rranga=
\delta_{ij}X^i\bullet\lambda\quad\mbox{and}\quad
\llanga
Y^i,Y^j\rranga=\delta_{ij}Y^i
\bullet\lambda,
\end{equation}
where $\delta_{ij}$ denote Kronecker's deltas
[recall that $\{X^i\}$ and $\{Y^i\}$ are families of independent
super-Brownian motions obeying the SPDE (\ref{eqmainSPDE0})]. In
terms of stochastic calculus,
(\ref{eqXiYj}) and (\ref{eqXiYj+1}) would completely characterize
the immigrant processes.

We do not pursue the question whether every interlacing pair of
${\varepsilon}$-\break approximating solutions admits immigrant processes
subject to the coexistence condition (\ref{eqXiYj}). For our purpose
to study pathwise nonuniqueness in the SPDE (\ref{eqmainSPDE}),
it is enough to turn to the converse point of view and study whether
there exist such immigrant processes so that they define an interlacing
pair of ${\varepsilon}$-approximating solutions (subject to the same
white noise) as in (\ref{eqXvep}) and (\ref{eqYvep}). More
precisely, our plan is to construct, for every ${\varepsilon}\in
(0,1]$, immigrant processes $\{X^i\}$ and $\{Y^i\}$ satisfying
conditions (i)--(v) of Theorem~\ref{propequiallocation} [we do not
require (vi)], and in addition, (\ref{eqXiYj}) so that they are
nonnegative solutions to a system of SPDEs of the form
%
%
%e2.66 #&#
\begin{equation}
\label{eqmp} \cases{ \displaystyle X^i_t(\phi)= \psi(
\1)J^{x_i}_{\varepsilon}(\phi)\1 _{t\geq s_i}+\int
_{s_i}^{s_i\vee t}X^i_s \biggl(
\frac{\Delta\phi
}{2} \biggr)\,ds
\vspace*{2pt}\cr
\displaystyle \hphantom{X^i_t(\phi)=}{}+\sum_{j=1}^\infty
\int_{s_i}^{s_i\vee t}\!\int_\R
\sigma_{2i-1,j}\bigl(X^1,Y^1,X^2,Y^2,
\ldots,s\bigr)\phi(x) \,dW^j(x,s),
\vspace*{5pt}\cr
\displaystyle Y^i_t(
\phi)= \psi(\1)J^{y_i}_{\varepsilon}(\phi)\1 _{t\geq t_i}+\int
_{t_i}^{t_i\vee t}Y^i_s \biggl(
\frac{\Delta\phi
}{2} \biggr)\,ds
\vspace*{2pt}\cr
\displaystyle\hphantom{Y^i_t(\phi)=}{}+\sum_{j=1}^\infty
\int_{t_i}^{t_i\vee t}\!\int_\R
\sigma_{2i,j}\bigl(X^1,Y^1,X^2,Y^2,\ldots,s\bigr)\phi(x)\,dW^j(x,s),}\hspace*{-38pt}
\end{equation}
for $\phi\in\C_c^\infty(\R)$ and some infinite-dimensional
deterministic diffusion coefficient matrix $\sigma
(x^1,y^1,x^2,y^2,\ldots,s)$ depending on space variables
$x^1,y^1,x^2,y^2,\ldots$ (in contrast, recall that $x_i$ and $y_i$
denote the landing targets of $X^i$ and $Y^i$, resp.) and time
variable $s$.
In (\ref{eqmp}), (1) $x_i$ and $y_i$ are i.i.d. with distribution
(\ref{eqxiyilaw}) as before, (2) $ (\sigma_{i,j}(\cdot,\cdot,s) )_{s\in[0,t]}$ are zero for all but finitely many $j$ for
every fixed $i$ and finite $t$ so that the infinite series in (\ref
{eqmp}) reduce to finite sums, and (3) $\{W^j\}$ is a family of i.i.d.
space--time white noises.
We remark that the various restrictions on~$t$ in the equations for
$X^i$ and $Y^i$ in (\ref{eqmp}) (namely, $\1_{t\geq s_i}$, $s_i\vee
t$, $\1_{t\geq t_i}$ and $t_i\vee t$) imply that $X^i$ and $Y^i$ are
nonzero only in $[s_i,\infty)$ and $[t_i,\infty)$, respectively, and
in writing $x^1,y^1,x^2,y^2,\ldots$ for the arguments of $\sigma$, we
keep track of the order in which immigrants land.
Finding an appropriate grand coefficient matrix $\sigma$ is the major
task of this section. Below we make a series of observations, and
the conclusion will be stated in Theorem~\ref{teoapprox-v2} by the
end of this section.

%
%
%pr2.9 #&#
\begin{prop}\label{propred1}
Fix ${\varepsilon}\in(0,1]$. Let $\{X^i\}$ and $\{Y^i\}$ be adapted
super-Brownian motions
defined on a filtered probability space $ (\Omega,\F,(\G_t),\P
)$, so that $\{X^i\}$ is subject to \textup{(i)} and \textup{(iii)} of Theorem~\ref
{propequiallocation} and $\{Y^i\}$ is subject to the analogous
conditions. As before, the landing targets $x_i$ and $y_i$ here are
i.i.d. with distribution~(\ref{eqxiyilaw}) and satisfy (\ref
{coneilc}). Suppose that all pairs $(X^i,Y^j)$ have normal
covariations and the coexistence condition (\ref{eqXiYj}) is
satisfied. Then
$X=\sum_iX^i$ and $Y=\sum_iY^i$ define an interlacing pair of
${\varepsilon}$-approximating solutions.
\end{prop}
\begin{pf}
As explained in Section~\ref{seccovar-imm}, we may assume the
existence of two families of independent space--time white noises $
\{W^{X^i} \}$ and $ \{W^{Y^i} \}$ so that $
(X^i,W^{X^i} )$ and $ (Y^i,W^{Y^i} )$ solve the SPDE (\ref
{eqmainSPDE0}) of super-Brownian motion.

We have to show that $X$ and $Y$ are subject to the SPDEs (\ref
{eqXvep}) and (\ref{eqYvep}), respectively, both with respect to the
same $(\G_t)$-space--time white noise $W$.
We start with the definition of $W$. Let $V$ be a $(\G_t)$-space--time
white noise
independent of $\{(X^i,W^{X^i})\}$ and $\{(Y^i,W^{Y^i})\}$.
By our assumptions and
L\'evy's theorem (cf. Theorem IV.3.6 of \cite{RYCMB}), we deduce that
%
%
%e2.67 #&#
\begin{eqnarray}
\label{WMG} W&\triangleq&\sum_{i=1}^\infty
\1_{(Y>0)} \biggl(\frac{Y^i}{Y} \biggr)^{1/2}\bullet
W^{Y^i}+\sum_{i=1}^\infty
\1_{(X>0,Y=0)} \biggl(\frac
{X^i}{X} \biggr)^{1/2}\bullet
W^{X^i}
\nonumber\\[-8pt]\\[-8pt]\nonumber
&&{}+\1_{(X=0,Y=0)}\bullet V
\end{eqnarray}
defines a $(\G_t)$-space--time white noise.
Then $Y$ is subject to the SPDE (\ref{eqYvep}) with respect to $W$
since the compatibility condition $Y^{1/2}\bullet W=\sum_{i=1}^\infty
(Y^i)^{1/2}\bullet W^{Y^i}$ holds. Indeed, we have
%
%
%e2.68 #&#
\begin{eqnarray}
\label{YMG} Y^{1/2}\bullet W&=&Y^{1/2}\1_{(Y>0)}
\bullet W=\sum_{i=1}^\infty\1
_{(Y>0)}\bigl(Y^{i}\bigr)^{1/2}\bullet
W^{Y^i}
\nonumber\\[-8pt]\\[-8pt]\nonumber
&=&\sum_{i=1}^\infty
\bigl(Y^i\bigr)^{1/2}\bullet W^{Y^i},
\end{eqnarray}
where the last equality follows from the nonnegativity of $Y^i$'s.

To prove that $X$ is also subject to $W$, one may wish that the roles
of $(\{X^i\},X)$ and $(\{Y^i\},Y)$ on the right-hand side of (\ref
{WMG}) can be exchanged, that is $W$ can also be rewritten as
%
%
%e2.69 #&#
\begin{eqnarray}
\label{WMG1} \qquad W&=&\sum_{i=1}^\infty
\1_{(X>0)} \biggl(\frac{X^i}{X} \biggr)^{1/2}\bullet
W^{X^i}+\sum_{i=1}^\infty
\1_{(Y>0,X=0)} \biggl(\frac
{Y^i}{Y} \biggr)^{1/2}\bullet
W^{Y^i}
\nonumber\\[-8pt]\\[-8pt]\nonumber
&&{}+\1_{(X=0,Y=0)}\bullet V,
\end{eqnarray}
and so the argument in (\ref{YMG}) applies to $X$. In this direction,
it is enough to claim that
%
%
%e2.70 #&#
\begin{equation}
\label{eqXYpos} \qquad\quad\sum_{i=1}^\infty
\1_{(X>0,Y>0)} \biggl(\frac{X^i}{X} \biggr)^{1/2}\bullet
W^{X^i}-\sum_{i=1}^\infty
\1_{(X>0,Y>0)} \biggl(\frac
{Y^i}{Y} \biggr)^{1/2}\bullet
W^{Y^i}=0.
\end{equation}
Recall the measure $\lambda$ defined in (\ref{deflambda}).
Since all of the pairs $(X^i,Y^j)$ have normal covariations, the
left-hand side of the foregoing equality also has a normal covariation
and we have
\begin{eqnarray*}
&& \Biggl\langle\!\Biggl\langle\sum_{i=1}^\infty
\1_{(X>0,Y>0)} \biggl(\frac{X^i}{X} \biggr)^{1/2}\bullet
W^{X^i}-\sum_{i=1}^\infty\1
_{(X>0,Y>0)} \biggl(\frac{Y^i}{Y} \biggr)^{1/2}\bullet
W^{Y^i},
\\
&&\qquad \sum_{j=1}^\infty
\1_{(X>0,Y>0)} \biggl(\frac
{X^j}{X} \biggr)^{1/2}\bullet
W^{X^j}-\sum_{j=1}^\infty\1
_{(X>0,Y>0)} \biggl(\frac{Y^j}{Y} \biggr)^{1/2}\bullet
W^{Y^j} \Biggr\rangle\!\Biggr\rangle
\\
&&\qquad =\sum_{i=1}^\infty\1_{(X>0,Y>0)}
\frac{X^i}{X}\bullet\lambda-2\sum_{i,j=1}^\infty
\1_{(X>0,Y>0)}\frac{1}{X^{1/2}Y^{1/2}}\frac
{X^iY^j}{X^{1/2}Y^{1/2}}\bullet\lambda
\\
&&\quad\qquad{}+\sum
_{i=1}^\infty\1 _{(X>0,Y>0)}\frac{Y^i}{Y}
\bullet\lambda
\\
&&\qquad =\1_{(X>0,Y>0)}\bullet\lambda-2\1_{(X>0,Y>0)}\bullet\lambda+\1
_{(X>0,Y>0)}\bullet\lambda=0,
\end{eqnarray*}
where the second equality follows from (\ref{eqXiYj}) and (\ref
{eqXiYj+1}) [the stochastic integral terms of the SPDEs for $X^i$ and
$Y^j$ are characterized by $(X^i)^{1/2}\bullet W^{X^i}$ and
$(Y^j)^{1/2}\bullet W^{Y^j}$].
We deduce our claim (\ref{eqXYpos}) from the last equality and the
fact that the action of the left-hand side of (\ref{eqXYpos}) on
every function in $\C_c^\infty(\R)$ induces a continuous martingale.
We have proved the alternative expression (\ref{WMG1}) of $W$, and
the proof is complete.
\end{pf}

%
%
%le2.10 #&#
\begin{lem}\label{lemred2}
Suppose that for $N\in\mathbb N$, $\xi^1,\ldots,\xi^N$ are continuous
nonnegative $\C_\rap(\R)$-valued solutions to the SPDE (\ref
{eqmainSPDE0}) with respect to the same filtration and independent
initial conditions, and all pairs $(\xi^i,\xi^j)$ have normal
covariations with
$\llang\xi^i,\xi^j\rrang=\delta_{ij}\xi^i\bullet\lambda$ for
$\lambda$ given by (\ref{deflambda}).
Then $\xi^1,\ldots,\xi^N$ are independent super-Brownian motions.
\end{lem}

\begin{pf*}{Sketch of proof}
The proof is to generalize the exponential duality argument for
super-Brownian motion. For each $i$, let $\phi^i$ be a nonnegative $\C
_c^\infty(\R)$-function and $u^i$ be the unique nonnegative solution
of the PDE
\[
\partial_r u_r^i=\frac{\Delta u_r^i}{2}-
\frac{1}{2}\bigl(u_r^i\bigr)^2 \mbox{ in }\R\times(0,\infty)\qquad\mbox{with } u_0^i=
\phi^i
\]
(cf. Lemma~4 in Section~II.2 of \cite{LGSBM} or pages 167--169 of \cite{PDW}).
Then for every fixed $t\in(0,\infty)$, the continuous semimartingale
$\exp\{-\sum_{i=1}^N\xi^i_{s} (u_{t-s}^i ) \}$,
$0\leq s\leq t$,
has zero finite variation by It\^{o}'s lemma and the assumption that
$\llang\xi^i,\xi^j\rrang=\delta_{ij}\xi^i\bullet\lambda$ (cf.
Proposition~II.5.7 of \cite{PDW}), and hence has constant mean. It
follows that one-dimensional marginals of $(\xi^1,\ldots,\xi^N)$ are
uniquely determined as those of independent super-Brownian motions.
A standard argument for martingale problems (cf. Section~4.4 in \cite{EK}) implies the desired result.
\end{pf*}

Thanks to Lemma~\ref{lemred2}, the main assumptions of
Proposition~\ref{propred1} are reduced to the covariation equations
(\ref{eqXiYj}) and (\ref{eqXiYj+1}) for $\{X^i\}$ and $\{Y^i\}$, as
well as other minor conditions. Then as in the standard construction of
solutions to systems of stochastic differential equations,
the issue is whether these covariation equations (\ref{eqXiYj}) and
(\ref{eqXiYj+1}) are induced by the nonnegative definite matrix
$\sigma\sigma^\top$
for some diffusion coefficient matrix $\sigma$ as in (\ref{eqmp}).

Below we write $\bs x=(x^1,x^2,\ldots)$, $\bs y=(y^1,y^2,\ldots)$,
and $\mathbf0=(0,0,\ldots)$ for which the dimensions may vary from
line to line but will be clear from the context.

%
%
%le2.11 #&#
\begin{lem}\label{lema}
Fix $n,m\in\mathbb N$, and consider the matrix-valued function
%
%
%e2.71 #&#
\begin{equation}
\label{defa} (\bs x,\bs y)\lmt a^{(n,m)}(\bs x,\bs y)=
\bigl[a^{(n,m)}_{k,\ell}(\bs x,\bs y) \bigr]_{1\leq k,\ell\leq m+n}
\end{equation}
defined on $ (\R^n\setminus\{\bs0\})\times(\mathbb R^m\setminus\{
\mathbf0\})$ as follows. For $\bs x=(x^1,x^2,\ldots,x^n)$ and $\bs
y=(y^1,y^2,\ldots,y^m)$ with $x^i,y^j\geq0$ and $\sum_{i'}x^{i'},\sum
_{j'}y^{j'}>0$, we set
\begin{eqnarray*}
&&\cases{ \displaystyle a^{(n,m)}_{i,j}(\bs x,\bs
y)=x^i\delta_{ij}, \hspace*{35pt}\mbox{$1\leq i,j\leq n$,}
\vspace*{5pt}\cr
\displaystyle a^{(n,m)}_{n+i,n+j}(\bs x,\bs y)=y^j
\delta_{ij},\qquad\mbox{$1\leq i,j\leq m$,}
\vspace*{5pt}\cr
\displaystyle
a^{(n,m)}_{i,n+j}(\bs x,\bs y)=a^{(n,m)}_{n+j,i}(
\bs x,\bs y)
\vspace*{5pt}\cr
\hspace*{54pt}\displaystyle=\bigl(x^i\bigr)^{1/2}
\bigl(y^j\bigr)^{1/2} \biggl(\frac{x^i}{\sum_{i'}x^{i'}}
\biggr)^{1/2} \biggl(\frac{y^j}{\sum_{j'}y^{j'}} \biggr)^{1/2},
\vspace*{5pt}\cr
\hspace*{126pt}\mbox{$1\leq i\leq n, 1\leq j\leq m$.}}
\end{eqnarray*}
For other $(\bs x,\bs y)\in(\R^n\setminus\{\bs0\})\times(\mathbb
R^m\setminus\{\mathbf0\})$, we set
%
%
%e2.72 #&#
\begin{equation}
\label{asym} a^{(n,m)}(\bs x,\bs y)=a^{(n,m)} \bigl(\bigl\llvert
x^1\bigr\rrvert,\bigl\llvert x^2\bigr\rrvert,\ldots,
\bigl\llvert x^n\bigr\rrvert,\bigl\llvert y^1\bigr
\rrvert,\bigl\llvert y^2\bigr\rrvert,\ldots,\bigl\llvert
y^m\bigr\rrvert\bigr).
\end{equation}
Then $a^{(n,m)}$ extends continuously to the entire space $\R^n\times
\R^m$, and the extension, still denoted by $a^{(n,m)}$, takes values
in $(m+n)$-by-$ (m+n)$ nonnegative definite matrices.
\end{lem}
\begin{pf}
Our assertion that $a^{(n,m)}$ extends continuously to $\R^n\times\R
^m$ follows plainly from the fact that
\[
\frac{\llvert x^i\rrvert }{\sum_{i'}\llvert x^{i'}\rrvert },\frac
{\llvert y^j\rrvert }{\sum_{j'}\llvert y^{j'}\rrvert }\in
[0,1]\qquad\forall\bs x\in
\R^n\setminus\{\mathbf0\}, \bs y\in\R^m\setminus\{
\mathbf0\}.
\]

We turn to the nonnegative definiteness of $a^{(n,m)}$. By continuity
and (\ref{asym}), we only need to show that $a^{(n+m)}(\bs x,\bs y)$
is nonnegative definite for $\bs x\in\R^n$ and $\bs y\in\R^m$
satisfying $x_i,y_j>0$ for all $i,j$.
Write
\[
a^{(n+m)}(\bs x,\bs y)=\lleft[\matrix{D^X &A
\vspace*{2pt}\cr
A^\top&D^Y} \rright]
\]
for an $n$-by-$n$ diagonal matrix $D^X$ and an $m$-by-$m$ diagonal
matrix $D^Y$. For any $(u,v)\in\R^{n}\times\R^{m}$, we regard $u$
and $v$ as column vectors and compute
%
%
%e2.73 #&#
\begin{eqnarray}\label{equv}
\qquad&&\bigl[\matrix{u^\top&v^\top} \bigr]
a^{(n+m)}(\bs x,\bs y) \lleft[\matrix{u
\vspace*{2pt}\cr v} \rright] \nonumber
\\
&&\qquad = \bigl[\matrix{u^\top&v^\top} \bigr] \lleft[\matrix{
D^X &A \vspace*{2pt}\cr
A^\top&D^Y} \rright] \lleft[
\matrix{u\vspace*{2pt}\cr v} \rright] \nonumber
\\
&&\qquad =u^\top D^Xu+2u^\top Av+v^\top
D^Yv
\\
&&\qquad =\sum_{i} \bigl(u^i
\bigr)^2x^i+2\sum_{i,j}u^i
\bigl(x^i\bigr)^{1/2} \biggl(\frac
{x^i}{\sum_{i'} x^{i'}}
\biggr)^{1/2}v^j\bigl(y^j\bigr)^{1/2}
\biggl(\frac{y^j}{\sum_{j'}y^{j'}} \biggr)^{1/2}\nonumber
\\
&&\quad\qquad{}+\sum_j \bigl(v^j\bigr)^2y^j\nonumber
\end{eqnarray}
by the definition of $a^{(n,m)}$. Notice that for all $\alpha^1,\ldots,\alpha^n,\beta^1,\ldots,\beta^m\in\R$,
\begin{eqnarray*}
2\sum_{i,j}\alpha^i
\beta^j&=& \biggl(\sum_i
\alpha^i+\sum_j\beta^j
\biggr)^2-\sum_i \bigl(
\alpha^i\bigr)^2-2\sum_{i_1<i_2}
\alpha^{i_1}\alpha^{i_2}
\\
&&{}-\sum_j
\bigl(\beta^j\bigr)^2-2\sum_{j_1<j_2}
\beta^{j_1}\beta^{j_2}.
\end{eqnarray*}
Applying the foregoing equality to the second term on the right-hand
side of (\ref{equv}) with the choice
\[
\alpha^i=u^i\bigl(x^i\bigr)^{1/2}
\biggl(\frac{x^i}{\sum_{i'}x^{i'}} \biggr)^{1/2}
\]
and
\[
\beta^j=v^j\bigl(y^j\bigr)^{1/2}
\biggl(\frac
{y^j}{\sum_{j'}y^{j'}} \biggr)^{1/2},
\]
we obtain
%
%
%e2.74 #&#
\begin{eqnarray}\label{eqcovar}
&& \bigl[\matrix{u^\top&v^\top} \bigr]
a^{(n+m)}(\bs x,\bs y)\lleft[\matrix{ u
\vspace*{2pt}\cr
v} \rright]
\nonumber
\\
&&\qquad =\sum_i \bigl(u^i\bigr)^2x^i\nonumber
\\
&&\quad\qquad{} + \biggl[\sum
_i u^i\bigl(x^i\bigr)^{1/2}
\biggl(\frac
{x^i}{\sum_{i'}x^{i'}} \biggr)^{1/2}+\sum
_j v^j\bigl(y^j\bigr)^{1/2}
\biggl(\frac
{y^j}{\sum_{j'}y^{j'}} \biggr)^{1/2} \biggr]^2
\nonumber\\[-8pt]\\[-8pt]\nonumber
&&\quad\qquad{}-\sum_i \bigl(u^i
\bigr)^2\frac{(x^i)^2}{\sum_{i'}x^{i'}}-2\sum_{i_1<i_2}
u^{i_1}u^{i_2}\frac{x^{i_1}}{ (\sum_{i'}x^{i'}
)^{1/2}}\frac{x^{i_2}}{ (\sum_{i'}x^{i'} )^{1/2}}
\\
&&\quad\qquad{}-\sum
_j \bigl(v^j\bigr)^2
\frac{(y^j)^2}{\sum_{j'}y^{j'}}-2\sum_{j_1<j_2} v^{j_1}v^{j_2}
\frac{y^{j_1}}{ (\sum_{j'}y^{j'} )^{1/2}}\frac{y^{j_2}}{ (\sum_{j'}y^{j'}
)^{1/2}}\nonumber
\\
&&\quad\qquad{}+\sum_j\bigl(v^j\bigr)^2y^j.\nonumber
\end{eqnarray}
The first, third and fourth terms on the right-hand side of the above
equality (with their signs taken into account as well) sum to
%
%
%e2.75 #&#
\begin{eqnarray}\label{eqcovarx}
\quad &&\sum_i \bigl(u^i
\bigr)^2x^i-\sum_i
\bigl(u^i\bigr)^2\frac{(x^i)^2}{\sum_{i'}x^{i'}}-2\sum
_{i_1<i_2} u^{i_1}u^{i_2}\frac{x^{i_1}}{ (\sum_{i'}x^{i'} )^{1/2}}
\frac{x^{i_2}}{ (\sum_{i'}x^{i'}
)^{1/2}}
\nonumber\\[-8pt]\\[-8pt]\nonumber
&&\qquad =\sum_i \bigl(u^i
\bigr)^2x^i- \biggl[\sum_i
u^i\frac{x^i}{
(\sum_{i'}x^{i'} )^{1/2}} \biggr]^2\geq0,
\end{eqnarray}
since the Cauchy--Schwarz inequality implies that
\begin{eqnarray*}
\biggl[\sum_i u^i\frac{x^i}{ (\sum_{i'}x^{i'}
)^{1/2}}
\biggr]^2&=& \biggl[\sum_i
u^i\bigl(x^i\bigr)^{1/2}\times\biggl(
\frac
{x^i}{\sum_{i'}x^{i'}} \biggr)^{1/2} \biggr]^2
\\
&\leq& \biggl(\sum
_i\bigl(u^i\bigr)^2x^i
\biggr) \biggl(\sum_i\frac{x^i}{\sum_{i'}x^{i'}} \biggr)
\\
&=&\sum_i\bigl(u^i
\bigr)^2x^i.
\end{eqnarray*}
Similarly, the last three terms on the right-hand side of (\ref
{eqcovar}) sum to
%
%
%e2.76 #&#
\begin{eqnarray}
\label{eqcovary} %
\qquad&&-\sum_j
\bigl(v^j\bigr)^2\frac{(y^j)^2}{\sum_{j'}y^{j'}}-2\sum
_{j_1<j_2} v^{j_1}v^{j_2}\frac{y^j}{ (\sum_{j'}y^{j'} )^{1/2}}
\frac
{y^j}{ (\sum_{j'}y^{j'} )^{1/2}}+\sum_j \bigl(v^j
\bigr)^2y^j
\nonumber\\[-8pt]\\[-8pt]\nonumber
&&\qquad =\sum_j \bigl(v^j
\bigr)^2y^j- \biggl[\sum_j
v^j\frac{y^j}{
(\sum_{j'}y^{j'} )^{1/2}} \biggr]^2\geq0.
\end{eqnarray}
Apply (\ref{eqcovarx}) and (\ref{eqcovary}) to the right-hand side
of (\ref{eqcovar}), and we obtain
%
%e2.77 #&#
\begin{eqnarray}
&& \bigl[\matrix{u^\top & v^\top} \bigr]
a^{(n+m)}(\bs x,\bs y)\lleft[\matrix{ u
\vspace*{2pt}\cr
v} \rright]\nonumber
\\
&&\qquad
\geq \biggl[\sum
_i u^i\bigl(x^i
\bigr)^{1/2} \biggl(\frac{x^i}{\sum_{i'}x^{i'}} \biggr)^{1/2}+\sum
_j v^j\bigl(y^j
\bigr)^{1/2} \biggl(\frac{y^j}{\sum_{j'}y^{j'}} \biggr)^{1/2}
\biggr]^2
\geq 0\nonumber
\\
\eqntext{\forall(u,v)\in\R^{n}\times\R^{m},}
\end{eqnarray}
that is, $a^{(n+m)}(\bs x,\bs y)$ is nonnegative definite. The proof is
complete.
\end{pf}

We are now ready to define the sought-after diffusion coefficient
matrix $\sigma$. For convenience, we reorder the arguments and entries
of the matrix-valued function $a^{(n,m)}(\bs x,\bs y)$ in Lemma~\ref
{lema} in accordance with the order in which immigrants land, for
$(n,m)$ equal to $(n,n)$ or $(n,n-1)$. This results in the
matrix-valued functions $A^{(n,n)}$ and $A^{(n,n-1)}$ defined by
%
%
%e2.78 #&#
\begin{eqnarray}\label{eqalphaa}
\qquad && A^{(n,n)}\bigl(x^1,y^1,x^2,y^2,\ldots,x^n,y^n\bigr)\triangleq  \Pi_{n,n}
a^{(n,n)}(\bs x,\bs y)\Pi^\top_{n,n},\qquad n\geq1,
\nonumber
\\
&& A^{(n,n-1)}\bigl(x^1,y^1,x^2,y^2,
\ldots,x^{n-1},y^{n-1},x^{n}\bigr)
\\
&&\qquad \triangleq \Pi_{n,n-1} a^{(n,n-1)}(\bs x,\bs y)\Pi^\top_{n,n-1},
\qquad n\geq2.\nonumber
\end{eqnarray}
Here, $\Pi_{n,n}$ and $\Pi_{n,n-1}$ are the permutation matrices
defined by
%
%
%e2.79 #&#
\begin{eqnarray}\label{eqPnm} %
&& \Pi_{n,n}\bigl[x^1,x^2,
\ldots,x^n,y^1,y^2,\ldots,y^n
\bigr]^\top\nonumber
\\
&&\qquad = \bigl[x^1,y^1,x^2,y^2,
\ldots,x^n,y^n\bigr]^\top,
\nonumber\\[-8pt]\\[-8pt]\nonumber
&& \Pi_{n,n-1}\bigl[x^1,x^2,\ldots,x^n,y^1,y^2,
\ldots,y^{n-1}\bigr]^\top
\\
&&\qquad = \bigl[x^1,y^1,x^2,y^2,
\ldots,x^{n-1},y^{n-1},x^{n}\bigr]^\top.\nonumber
\end{eqnarray}
Then we choose a continuous square root, denoted by $\sigma^{(n,m)}$,
of the square matrix $A^{(n,m)}$ (cf. Theorem 1.1 of \cite{CHROOT} or
\cite{FFact} for its existence) for $(n,m)=(n,n)$ or $(n,n-1)$, and
so $\sigma^{(n,m)}$ satisfies
%
%
%e2.80 #&#
\begin{equation}
\label{sigmaa} \sigma^{(n,m)} \bigl[\sigma^{(n,m)}
\bigr]^\top\equiv A^{(n,m)}.
\end{equation}
Let $\sigma^{(1,0)}(x^1)$ be the $1$-by-$1$ matrix $[\llvert
x^1\rrvert ^{1/2}]$.
Then we define $\sigma$ as follows. We set $\sigma\equiv0$ on
$[0,s_1)$, and for all $n,i,j\in\mathbb N$,
%
%
%e2.81 #&#
\begin{eqnarray}
\label{defsigmadiff} %
&& \sigma_{i,j}\bigl(x^1,y^1,x^2,y^2,
\ldots,s\bigr)
\nonumber\\[14pt]\\[-30pt]\nonumber
&&\qquad \equiv\cases{ \displaystyle\sigma_{i,j}^{(n,n-1)}
\bigl(x^1,y^1,x^2,y^2,\ldots,x^{n-1},y^{n-1},x^n\bigr),
\vspace*{3pt}\cr
\hspace*{33pt}\mbox{if $i,j \leq2n-1$ and $s\in[s_n,t_n)$,}
\vspace*{5pt}\cr
\displaystyle
\sigma_{i,j}^{(n,n)}\bigl(x^1,y^1,x^2,y^2,
\ldots,x^{n},y^{n}\bigr),
\vspace*{3pt}\cr
\hspace*{33pt}\mbox{if $i,j\leq2n$ and $s \in[t_n,s_{n+1})$,}
\vspace*{3pt}\cr
0,\qquad\mbox{otherwise}.}
\end{eqnarray}
In other words, solutions to the system (\ref{eqmp}) of SPDE's
subject to the above choice of diffusion coefficient matrix $\sigma$,
if any, can be described as follows: over $[s_n,t_n)$ for $n\in\mathbb N$,
\[
\bigl(X^1,Y^1,X^2,Y^2,
\ldots,X^{n-1},Y^{n-1},X^{n}\bigr)
\]
is subject to the diffusion coefficient $\sigma^{(n,n-1)}$ and the
independent noises $W^1,\ldots,W^{2n-1}$, and over $[t_n,s_{n+1})$ for
$n\in\mathbb N$,
\[
\bigl(X^1,Y^1,X^2,Y^2,
\ldots,X^n,Y^n\bigr)
\]
is subject to the diffusion coefficient $\sigma^{(n,n)}$ and the
independent noises $W^1,\ldots,W^{2n}$. Note that $\sigma$ depends
only on space variables between two consecutive immigration times.

%
%
%th2.12 #&#
\begin{teo}\label{teoapprox-v2}
Fix an immigration function $\psi\in\C_c^+(\R)\setminus\{0\}$. For
any ${\varepsilon}\in(0,1]$, we can construct a filtered probability
space $(\Omega,\F,(\G_t),\P_{\varepsilon})$, with $(\G_t)$
satisfying the usual conditions, on which there exist random elements
$\{x_i\}$, $\{y_i\}$, $\{X^i\}$, $\{Y^i\}$, $\{W^i\}$ and $W$ with the
following properties:
\begin{longlist}[(iii)]
\item[(i)] $x_i$ and $y_i$ are i.i.d. with law (\ref{eqxiyilaw})
and take values in the topological support of $\psi$.
\item[(ii)] $W^i$ and $W$ are $(\G_t)$-space--time white noises,
and the noises $\{W^i\}$ are independent.
\item[(iii)] For each $i\in\mathbb N$, $ (X^i_{t} )_{t\in
[s_i,\infty)}$ and $ (Y^i_{t} )_{t\in[t_i,\infty)}$ are
nonnegative processes with sample paths in $C ([s_i,\infty),\C
_\rap(\R) )$ and $C ([t_i,\infty),\C_\rap(\R) )$,
respectively.
\item[(iv)] $\{X^i\}$ and $ \{Y^i\}$ obey the system of SPDE's
(\ref{eqmp}) with respect to $\{W^i\}$ for~$\sigma$ defined by (\ref
{defsigmadiff}).
\item[(v)] The independent landing property (\ref{coneilc}) for
immigrants holds.
\item[(vi)] The sums $X=\sum_iX^i$ and $Y=\sum_iY^i$ define an
interlacing pair of \mbox{${\varepsilon}$-}approximating solutions with
respect to $W$ (see Definition~\ref{defstd}).
\end{longlist}
In particular, $\{X^i\}$ and $\{Y^i\}$ are two families of independent
super-Brownian motions for which the covariation equations (\ref
{eqXiYj}) hold.
\end{teo}

The proof of Theorem~\ref{teoapprox-v2} follows similarly as the
existence of interlacing pairs of ${\varepsilon}$-approximating
solutions (see Definition~\ref{defstd}). We introduce i.i.d. landing
targets $x_i$ and $y_i$ which are independent of a family of i.i.d.
noises $\{W^i\}$, and then
solve (\ref{eqmp}) over $[s_1,t_1),[t_1,s_2),[s_2,t_2),\ldots$
sequentially (see \cite{BMP} for a similar construction). More
precisely, over any of these intervals, (\ref{eqmp}) reduces to a
finite-dimensional system of SPDEs to which the classical Peano
approximation argument as in the proof of Theorem 2.6 of \cite{SSPDE}
applies (cf. Section~6 of \cite{MMP} as well). Indeed, by comparing
diagonal entries on both sides of (\ref{sigmaa}), we deduce that
every entry $\sigma_{i,j}^{(n,m)}$ is bouded by $z\lmt\llvert
z\rrvert ^{1/2}$,
where $z=x^k$ if $i=2k-1$, or $z=y^k$ if $i=2k$.
We omit other details.

As an immediate consequence of (\ref{eqXiYj}), we obtain the
following (see Remark~\ref{rmkmistake}).

%
%
%pr2.13 #&#
\begin{prop}\label{propcovarbdd}
Let $\{X^i\}$ and $\{Y^i\}$ be as in Theorem~\ref{teoapprox-v2}.
Then for any $i,j_1,\ldots,j_n\in\mathbb N$ for $n\in\mathbb N$ with
$j_1<j_2<\cdots<j_n$, except outside a null event, the inequality
%
%
%e2.82 #&#
%e2.83 #&#
\begin{eqnarray}\label{eqcovarSBM}
&&\Biggl\llvert\int_0^t
H_s\,d \Biggl\langle X^i(\1),\sum
_{\ell=1}^nY^{j_\ell
}(\1) \Biggr \rangle_s\Biggr\rrvert\nonumber
\\
&&\qquad \leq\int_{s_i\vee t_{j_1}}^t \llvert
H_s\rrvert\int_\R\Biggl(X^i(x,s)
\cdot\sum_{\ell=1}^n Y^{j_\ell}(x,s)
\Biggr)^{1/2}\,dx\,ds
\\
\eqntext{\forall t\in[s_i\vee t_{j_1},\infty),}
\end{eqnarray}
holds
for any locally bounded Borel measurable function $H$ on $\R_+$.
\end{prop}

\begin{cas*}
%\subsubsection*{Choice of approximating solutions}
From now on,
we \emph{only} work with $\{X^i\}$ and $\{Y^i\}$ as in Theorem~\ref
{teoapprox-v2}, and the corresponding interlacing pairs of
\mbox{${\varepsilon}$-}approxi\-mating solutions $X=\sum_iX^i$ and $Y=\sum
_iY^i$, except in Section~\ref{secpropwc}. %\qed
\end{cas*}

%s3 #&#
\section{Conditional separation of approximating solutions}\label{sec2ndsep}

%s3.1 #&#
\subsection{Basic results}\label{sec2ndsep-basic}
The theme of Section~\ref{sec2ndsep} is conditional separation of the
approximating solutions defined by
the immigrant processes $\{X^i;i\in\mathbb N\}$ and $\{Y^i;i\in\mathbb
N\}$
chosen in Theorem~\ref{teoapprox-v2}.
For any ${\varepsilon}\in(0, [8\psi(\1)]^{-1}\wedge1
]$, we
condition on the event that the total mass of a generic cluster $X^i$
hits $1$, and then the conditional separation refers to the separation
of the approximating solutions under
%
%
%e3.1 #&#
\begin{equation}
\label{defQ} \Q^i_{\varepsilon}(A)\equiv\P_{\varepsilon}
\bigl(A \mid T^{X^i}_1<\infty\bigr).
\end{equation}
Here, the restriction $[8\psi(\1)]^{-1}$ for ${\varepsilon}$ is just
to make sure that $X^i(\1)$ stays in $(0,1)$ initially, and we set
%
%
%e3.2 #&#
\begin{equation}
\label{eqTZ1} T^{H}_x\triangleq\inf\bigl\{t
\geq0;H_t(\1)=x\bigr\}
\end{equation}
for any nonnegative two-parameter process $H=(H(x,t);(x,t)\in\R\times
\R_+)$. Our specific goal is to study the differences in the growth
rates of local masses of $X$ and $Y$ over the ``initial part'' of the
space--time support of $X^i$. In the following, we prove a few basic
results concerning $\Q^i_{\varepsilon}$.

First, let us represent the Radon--Nikodym derivative process of $\Q
^i_{\varepsilon}$ relative to~$\P_{\varepsilon}$ (cf. Section VIII.1
of \cite{RYCMB} for its role in Girsanov's theorem).

%
%
%le3.1 #&#
\begin{lem}\label{lemQi}
For any $i\in\mathbb N$ and ${\varepsilon}\in(0,[8\psi(\1
)]^{-1}\wedge1 ]$,
%
%
%e3.3 #&#
\begin{eqnarray}
\P_{\varepsilon} \bigl(T^{X^i}_1<T^{X^i}_0
\bigr)&=&\psi(\1 ){\varepsilon},\label{eqhitprob}
\end{eqnarray}
and the Radon--Nikodym derivative process $\mathbb E^{\P^{\varepsilon
}}[d\mathbb Q^i_{\varepsilon}/d\P_{\varepsilon}\mid\G_t]$, $t\in
[s_i,\infty)$, of $\mathbb Q^i_{\varepsilon}$ relative to $\P
_{\varepsilon}$ is given by the stopped $ ((\G_t)_{t\geq s_i},\P
_{\varepsilon} )$-martingale $X^i(\1)^{T_1^{X^i}}/\break \psi(\1
){\varepsilon}$, that is,
%
%
%e3.4 #&#
\begin{equation}
\label{eqQirep} \mathbb Q^i_{\varepsilon}(A)=\int
_A\frac{X^i_t(\1)^{ T^{X^i}_1}}{\psi
(\1){\varepsilon}}\,d\P_{\varepsilon}\qquad\forall A\in
\G_{t} \mbox{ with }t\in[s_i,\infty).
\end{equation}
Here, $X^i(\1)^{ T^{X^i}_1}$ denotes the total mass process $X^i(\1)$
stopped at $T_1^{X^i}$ [see (\ref{eqTZ1}) for~$T_1^{X^i}$].
\end{lem}
\begin{pf}
The proof is a standard application of Doob's $h$-transforms (cf.
Section VII.3 of \cite{RYCMB}). Recall that $X^i(\1)$ under $\P
_{\varepsilon}$ is a Feller diffusion with initial condition $\psi(\1
){\varepsilon}$ and plainly the scale function of Feller diffusion is
given by $x\lmt x$. Hence, (\ref{eqhitprob}) follows from
Proposition VII.3.2 of \cite{RYCMB}. To see the second assertion, we
recall the definition (\ref{defQ}) of $\mathbb Q^i_{\varepsilon}$, and
then apply (\ref{eqhitprob}), Proposition VII.3.2 in~\cite{RYCMB}
again and the Markov property of $X^i(\1)$.
\end{pf}

Some basic properties of the total mass processes $X^i(\1)$ and $Y^j(\1
)$ for $t_j>s_i$ under $\Q^i_{\varepsilon}$ are stated in the
following lemma.

%
%
%le3.2 #&#
\begin{lem}\label{lemsep1}
Fix $i\in\mathbb N$ and ${\varepsilon}\in(0, [8\psi(\1
)]^{-1}\wedge1 ]$. Then we have the following.
\begin{longlist}[(2)]
\item[(1)] $X^i(\1)^{T_1^{X^i}}$ under $\Q^i_{\varepsilon}$ is a
copy of $\frac{1}{4}\BES Q^4(4\psi(\1){\varepsilon})$ started at
$s_i$ and stopped upon hitting $1$.
\item[(2)] For any $j\in\mathbb N$ with $t_j>s_i$, the process
$(Y^j(\1)_t)_{t\geq t_j}$ is a continuous $(\G_t)_{t\geq
t_j}$-semimartingale under $\mathbb Q^i_{\varepsilon}$ with canonical
decomposition
%
%
%e3.5 #&#
\begin{equation}
\label{eqsemiMGY} Y^j_t(\1)=\psi(\1){\varepsilon}+I^j_t+M^j_t,
\qquad t\in[ t_j,\infty),
\end{equation}
where the finite variation process $I^j$ satisfies
%
%
%e3.6 #&#
%e3.7 #&#
\begin{eqnarray}
I^j_t&=&\int_{t_j}^t
\frac{1}{X^i_s(\1)^{ T^{X^i}_1}}\,d \bigl\langle X^i(\1)^{ T^{X^i}_1},Y^j(
\1) \bigr\rangle_s,\label{defIj}
\\
\label{ineqYj1} 0&\leq& I^j_t\leq\int
_{t_j}^{t}\1_{[0,T_1^{X^i}]}(s) \frac{1}{X^i_s(\1
)}
\int_\R X^i(x,s)^{1/2}Y^j(x,s)^{1/2}\,dx\,ds,
\end{eqnarray}
for $t\in[t_j,\infty)$,
and $M^j$ is a true $(\G_t)_{t\geq t_j}$-martingale under $\mathbb
Q^i_{\varepsilon}$.
\item[(3)] For any $j\in\mathbb N$ with $t_j>s_i$,
%
%
%e3.8 #&#
\begin{equation}
\label{impilp} \qquad x_i, X^i(\1)\rest[s_i,t_j],
y_j\quad\mbox{and}\quad Y^j(\1)\rest[t_j,\infty)
\qquad\mbox{are $\P_{\varepsilon}$-independent.}
\end{equation}
\item[(4)] For any $j\in\mathbb N$,
%
%
%e3.9 #&#
\begin{eqnarray}
\label{dom} \Q^i_{\varepsilon}\bigl(\llvert
y_j-x_i\rrvert\in dx\bigr)&=&\P_{\varepsilon}\bigl(
\llvert y_j-x_i\rrvert\in dx\bigr),\qquad x\in\R,
\nonumber\\[-8pt]\\[-8pt]\nonumber
\P_{\varepsilon}(y_j\in dx)&\leq&\frac{\llVert \psi\rrVert_\infty
}{\psi(\1
)} \,dx,\qquad x \in \R.
\end{eqnarray}
\end{longlist}
\end{lem}
\begin{pf}
(1)~The proof is omitted since it is a straightforward application of
Girsanov's theorem and Lemma~\ref{lemQi}, and can be found in the
proof of Lemma 4.1 of \cite{MMP}.

(2)~Under $\P_{\varepsilon}$, the total mass process $
(Y^j_t(\1) )_{t\geq t_j}$ for any $j\in\mathbb N$ with $t_j>s_i$ is a
$(\G_t)_{t\geq t_j}$-Feller diffusion and hence a $(\G_t)_{t\geq
t_j}$-martingale.
By Lemma~\ref{lemQi} and Girsanov's theorem (cf. Theorem VIII.1.4 of
\cite{RYCMB}), $(Y^j_t(\1) )_{t\geq t_j}$ for any $j\in\mathbb
N$ with $t_j>s_i$ is a continuous $(\G_t)_{t\geq t_j}$-semimartingale
under $\mathbb Q^i_{\varepsilon}$ with canonical decomposition given by
(\ref{eqsemiMGY}).
Here, $(M^j_t )_{t\geq t_j}$ is a continuous\break $(\G_t)_{t\geq
t_j}$-local martingale under $\mathbb Q^i_{\varepsilon}$ with quadratic
variation
%
%
%e3.10 #&#
\begin{equation}
\label{eqMjQV} \bigl\langle M^j\bigr\rangle_t=\int
_{t_j}^tY^j_s(\1)\,ds,
\qquad t\in[t_j,\infty),
\end{equation}
and by Lemma~\ref{lemQi} the finite variation process $
(I^j_t )_{t\geq t_j}$ is given by (\ref{defIj}).
Applying (\ref{eqXiYj}) and (\ref{eqcovarSBM}) to (\ref{defIj}),
we obtain (\ref{ineqYj1}) at once.

For the martingale property of $M^j$ under $\mathbb Q^i_{\varepsilon}$,
we note that the one-dimensional marginals of $Y^j(\1)$ have $p$th
moments which are locally bounded on compacts, for any $p\in(0,\infty
)$. [$Y^j(\1)$ under $\P_{\varepsilon}$ is a Feller diffusion.]
Applying this to (\ref{eqMjQV}) shows that
$\E^{\mathbb Q^i_{\varepsilon}}[\langle M^j\rangle_t]<\infty$ for
every $t\in[t_j,\infty)$, and hence $M^j$ is a true martingale under
$\mathbb Q^i_{\varepsilon}$.

(3)~The assertion (\ref{impilp}) is an immediate
consequence of the independent landing property (\ref{coneilc}) and
the Markov properties of $X^i(\1)$ and $Y^j(\1)$ (cf. Theorem~\ref{teoapprox-v2}).

(4)~We consider (\ref{dom}). Recall that $x_i\in\G_{s_i}$ and $y_j\in\G
_{t_j}$ by (\ref{impilc}). If $t_j>s_i$, then we obtain from (\ref
{eqQirep}) that
%
%
%e3.11 #&#
\begin{eqnarray}\label{eqyjxi}
\Q^i_{\varepsilon}\bigl(\llvert y_j-x_i
\rrvert\in dx\bigr)&=&\frac{1}{\psi(\1){\varepsilon
}}\E^{\P_{\varepsilon}} \bigl[X_{t_j}^i(
\1)^{T_1^{X^i}};\llvert y_j-x_i\rrvert\in dx \bigr]
\nonumber\\[-8pt]\\[-8pt]\nonumber
&=& \P_{\varepsilon}\bigl(\llvert y_j-x_i\rrvert\in dx
\bigr),
\end{eqnarray}
where the last equality follows from (\ref{impilp}). If $t_j<s_i$,
then a similar argument applies [without using (\ref{eqQirep})] since
$X^i_{s_i}(\1)=\psi(\1){\varepsilon}$. Hence, the equality in (\ref
{dom}) holds. The inequality in (\ref{dom}) is obvious.
The proof is complete.
\end{pf}

%s3.2 #&#
\subsection{Setup}\label{sec2ndsep-setup}
In order to state precisely our quantifications of the local growth
rates of $X$ and $Y$, we need several preliminary results which have
similar counterparts in \cite{MMP}. First, we choose in
Proposition~\ref{propseptime} below a $(\G_t)$-stopping time $\tau
^i$ satisfying $\tau^i> s_i$, so that within $[s_i,\tau^i]$ we can
explicitly bound from below the growth rate of $X^i(\1)$.
Since $X\geq X^i$, this gives a lower bound for the size of $X$ over
the initial part of the space--time support of $X^i$. Our objective is
to study the local growth rate of $Y$ within this part.

%
%
%pr3.3 #&#
\begin{prop}\label{propseptime}
For any ${\varepsilon}\in(0,[8\psi(\1)]^{-1}\wedge1 ]$,
parameter vector
$(\eta,\break \alpha,L)\in(1,\infty)\times(0,\frac{1}{2})\times
(0,\infty)$, and $i\in\mathbb N$, we define four $(\G_t)$-stopping
times by
\begin{eqnarray*}
\tau^{i,(1)}&\triangleq& \inf\biggl\{t\geq s_i;
X^i_t(\1)^{
T^{X^i}_1}< \frac{(t-s_i)^\eta}{4} \biggr\}
\wedge T^{X^i}_1,
\\
\tau^{i,(2)}&\triangleq& \inf\biggl\{t\geq s_i; \bigl\llvert
X^i_t(\1 )^{T^{X^i}_1}-\psi(\1){\varepsilon}-
(t-s_i)\bigr\rrvert>L \biggl(\int_{s_i}^t
X^i_s(\1)^{T^{X^i}_1}\,ds \biggr)^\alpha
\biggr\}
\\
&&{}\wedge T^{X^i}_1,
\\
\tau^{i,(3)}&\triangleq& \inf\biggl\{t\geq s_i;\sum
_{j\dvtx s_i<t_j\leq
t}Y^j_t(\1)> 1 \biggr\},
\\
\tau^i&\triangleq& \tau^{i,(1)}\wedge\tau^{i,(2)}
\wedge\tau^{i,(3)}\wedge(s_i+1).
\end{eqnarray*}
Then
%
%
%e3.12 #&#
\begin{eqnarray}\label{eqseptimebdd}
&& \forall\rho>0\ \exists\delta>0\mbox{ such that}
\nonumber\\[-8pt]\\[-8pt]\nonumber
&&\qquad \sup\biggl\{\mathbb Q^i_{\varepsilon}\bigl(\tau^i
\leq s_i+\delta\bigr);i\in\mathbb N, {\varepsilon}\in\biggl(0,
\frac{1}{8\psi(\1)}\wedge1 \biggr] \biggr\}\leq\rho.
\end{eqnarray}
\end{prop}

See Section~\ref{sectauvep} for the proof of Proposition~\ref{propseptime}.

Let us explain the meanings of the parameters $\eta,\alpha,L$ in this
proposition. Since
$X^i(\1)$ is a Feller diffusion under $\P_{\varepsilon}$, a
straightforward application of Girsanov's theorem (cf. Theorem VIII.1.4
of \cite{RYCMB}) shows that
$X^i(\1)^{T^{X^i}_1}$ under $\mathbb Q_{\varepsilon}^i$ is a $\frac
{1}{4}\BES Q^4 (4\psi(\1){\varepsilon} )$ stopped upon
hitting $1$; see Lemma 4.1 of \cite{MMP} for details. As a result, by
the lower escape rate of $\BES Q^4$ (cf. Theorem 5.4.6 of \cite{KEBM}),
the time $\tau^{i,(1)}$ is strictly positive $\mathbb Q^i_{\varepsilon
}$-a.s. for any $\eta\in(1,\infty)$. In particular, we may
take the parameter $\eta$ close to $1$.

The definition of $\tau^{i,(2)}$ involves the notion of improved
modulus of continuity.
We will take the parameter $\alpha$ in the definition of $\tau
^{i,(2)}$ close to $\frac{1}{2}$ and consider
the local H\"older exponent of the martingale part of $\BES Q^4$ in terms
of its quadratic variation. The parameter $L$ bounds the associated
local H\"older coefficient.
Hence, we have the integral inequality
%
%
%e3.13 #&#
%e3.14 #&#
\begin{eqnarray}
\label{ineqintineqX} %
&& \bigl\llvert X^i_t(
\1)^{T_1^{X^i}}-\psi(\1){\varepsilon}\bigr\rrvert\leq(t-s_i)+L
\biggl(\int_{s_i}^t X^i_s(
\1)^{T_1^{X^i}}\,ds \biggr)^\alpha
\nonumber\\[-8pt]\\[-8pt]
\eqntext{\forall t\in\bigl[s_i,
\tau^i \bigr],
\mathbb Q^i_{\varepsilon}\mbox{-a.s., }\forall i
\in\mathbb N, {\varepsilon}\in\bigl(0,\bigl[8\psi(\1)\bigr]^{-1}\wedge1 \bigr],}
\end{eqnarray}
by the choice of $\tau^{i,(2)}$ in Proposition~\ref{propseptime}.
The integral inequality (\ref{ineqintineqX}) is reminiscent of the
integral inequalities to which Gronwall's lemma applies, and hence
suggests an iteration argument if we wish to bound more explicitly the
difference
$\llvert X^i_t(\1)^{T_1^{X^i}}-\psi(\1){\varepsilon}\rrvert $.
A general result for this is given by Corollary~\ref{corimc}.
Applying Corollary~\ref{corimc} to the random function
\[
t\lmt X^i_t(\1)^{T_1^{X^i}}\dvtx \bigl[s_i,
\tau^i \bigr]\lra\R,
\]
we obtain from (\ref{ineqintineqX}) that whenever $\xi\in(0,1)$ and
$N_0\in\mathbb N$ satisfies
%
%
%e3.15 #&#
\begin{equation}
\label{alphaxi} \sum_{j=1}^{N_0}
\alpha^j\leq\xi<\sum_{j=1}^{N_0+1}
\alpha^j,
\end{equation}
the following inequality holds:
%
%
%e3.16 #&#
%e3.17 #&#
\begin{eqnarray}\label{ineqXQmc}
&& \bigl\llvert X^i_t(
\1)^{T^{X^i}_1}-\psi(\1){\varepsilon}\bigr\rrvert\leq K_1^X
\bigl[\psi(\1){\varepsilon} \bigr]^{\alpha
^{N_0}}(t-s_i)^\alpha+K_2^X
(t-s_i)^\xi
\nonumber\\[-8pt]\\[-8pt]
\eqntext{\forall t\in\bigl[s_i,\tau^i \bigr],
\Q^i_{\varepsilon}\mbox{-a.s., }\forall i\in\mathbb N, {
\varepsilon}\in\bigl(0,\bigl[8\psi(\1)\bigr]^{-1}\wedge1 \bigr],}
\end{eqnarray}
where the constants $K_1^X,K^X_2\geq1$ depend only on $(\alpha,L,\xi,N_0)$.
Moreover, since $\alpha$ is close to $\frac{1}{2}$, we can choose
$N_0$ large in (\ref{alphaxi})
to make $\xi$ close to $1$, as is our intention in the sequel.
Informally, we
can interpret the foregoing inequality as the statement:
\[
\mbox{$t\lmt X^{i}_t(\1)^{T_1^{X^i}}$\qquad is H\"older-$1$
continuous at $s_i$ from the right.}
\]
A similar derivation of the improved modulus of continuity of $Y^j(\1)$
will appear in the proof of Lemma~\ref{lemFd2} below.

To use the support of $X^i$ within which we study the local growth rate
of $Y$,
we take a parameter $\beta\in(0,\frac{1}{2})$, which is now close to
$\frac{1}{2}$. We use this parameter to get a better control
of the supports of $X^i$ and $Y^j$, and this means we use the parabola
%
%
%e3.18 #&#
\begin{equation}
\label{eqPXi0} \mathcal P^{X^i}_\beta(t)\triangleq\bigl
\{(x,s)\in\R\times[s_i,t];\llvert x-x_i\rrvert\leq
\bigl({\varepsilon}^{1/2}+(s-s_i)^\beta\bigr)
\bigr\}
\end{equation}
to \emph{envelop} the space--time support of $X^i\rest[s_i,t]$, for
$t\in(s_i,\infty)$,
with a similar practice applied to other clusters $Y^j$.
(See the speed of support propagation of super-Brownian motions in
Theorem III.1.3 of \cite{PDW}.) More precisely, we can use the $(\G
_t)$-stopping time
%
%
%e3.19 #&#
\begin{eqnarray}
\label{defsigma} \qquad \sigma^{X^i}_\beta&\triangleq&\inf\bigl\{s\geq
s_i;
\nonumber\\[-8pt]\\[-8pt]\nonumber
&&\hspace*{17pt}\supp\bigl(X^i_{s}\bigr)\nsubseteq
\bigl[x_i-{\varepsilon}^{1/2}-(s-s_i)^{\beta
},x_i+{
\varepsilon}^{1/2}+(s-s_i)^{\beta} \bigr] \bigr\}
\end{eqnarray}
as well as the analogous stopping times $\sigma^{Y^j}_\beta$ for
$Y^j$ to identify the duration of the foregoing enveloping.

We now specify the clusters $Y^j$ selected for computing the local
growth rate of~$Y$.
Suppose that at time $t$ with $t>s_i$, we can still envelop the support
of $X^i$ by $\mathcal P^{X^i}_\beta(t)$ and the analogous enveloping
for the support of $Y^j$ holds for any $j\in\mathbb N$\vadjust{\goodbreak} satisfying
$t_j\in
(s_i,t]$. Informally, we can ignore the clusters $Y^j$ landing before
$X^i$, because the probability that they can invade the initial part of
the support of $X^i$ is small for small $t$ (cf. Lemma~\ref{lemoutclus+si}).
Under such circumstances, simple geometric arguments show that
only the $Y^j$ clusters landing inside the space--time rectangle
%
%
%e3.20 #&#
\begin{equation}
\label{eqPXi1} \qquad\mathcal R^{X^i}_\beta(t)\triangleq
\bigl[x_i-2 \bigl({\varepsilon}^{1/2}+(t-s_i)^\beta
\bigr),x_i+2 \bigl({\varepsilon}^{1/2}+(t-s_i)^\beta
\bigr) \bigr]\times[s_i,t]
\end{equation}
can invade the initial part of the support of $X^i$ by time $t$ (see
Lemma~\ref{lemoutclus+si}).
We remark that this choice of clusters $Y^j$ for $(y_j,t_j)\in\mathcal
R^{X^i}_\beta(t)$ is also used in \cite{MMP}.

%
%
%f2 #&#
\begin{figure}%[b]%[t]

\includegraphics{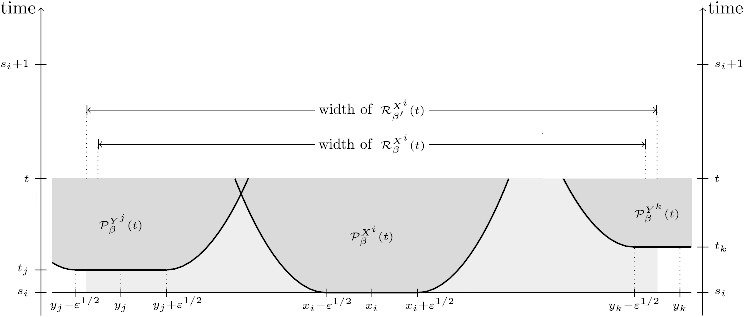}

\caption{Parabolas $\mathcal P^{X^i}_\beta(t),\mathcal P^{Y_j}_\beta
(t),\mathcal P^{Y_k}_\beta(t)$ and rectangles $\mathcal R^{X^i}_{\beta
}(t)$ and $\mathcal R^{X^i}_{\beta'}(t)$,
for $0<\beta'<\beta$ and $t\in[s_i,s_i+1)$.}\label{Fig1}
\end{figure}

For technical reasons (cf. Section~\ref{seccondsep-2} below),
however, we will consider the super-Brownian motions $Y^j$ landing
inside the slightly larger rectangle $\mathcal R^{X^i}_{\beta'}(t)$
for $t\in(s_i,s_i+1]$, where $\beta'$ is another value close to
$\frac{1}{2}$, has the same meaning as $\beta$, and satisfies $\beta
'<\beta$. See Figure~\ref{Fig1} for these rectangles as well as an
example for three parabolas $\mathcal P^{X^i}_\beta(t)$, $\mathcal
P^{Y^j}_\beta(t)$, and $\mathcal P^{Y^k}_\beta(t)$ where
$(y_j,t_j)\in\mathcal R^{X^i}_\beta(t)$ and $(y_k,t_k)\notin\mathcal
R^{X^i}_\beta(t)$. The labels $j\in\mathbb N$ of the clusters $Y^j$
landing inside $\mathcal R^{X^i}_{\beta'}(t)$ constitute the random
index set
%
%
%e3.21 #&#
\begin{equation}
\label{eqdefIitt} \mathcal J_{\beta'}^i(t)\equiv\mathcal
J^i_{\beta'}(t,t),
\end{equation}
where
%
%
%e3.22 #&#
%e3.23 #&#
\begin{eqnarray}
\label{eqdefIi} \mathcal J_{\beta'}^i
\bigl(t,t'\bigr)&\triangleq& \bigl\{j\in\mathbb N;\llvert
y_j-x_i\rrvert\leq2 \bigl({\varepsilon}^{1/2}+(t-s_i)^{\beta'}
\bigr), s_i<t_j\leq t' \bigr\}
\nonumber\\[-8pt]\\[-8pt]
\eqntext{\forall t,t'\in(s_i,\infty).}
\end{eqnarray}

%
%
%as3.4 #&#
\begin{ass}[(Choice of auxiliary parameters)]\label{ass} Throughout the
remainder of this section and Section~\ref{secsepSOL}, we fix
a parameter vector
%
%
%e3.24 #&#
\begin{eqnarray}
\label{param}
&& \bigl(\eta,\alpha,L,\beta,\beta',
\xi,N_0\bigr)
\nonumber\\[-8pt]\\[-8pt]\nonumber
&&\qquad \in (1,\infty)\times\bigl(0,\tfrac{1}{2}
\bigr)\times(0,\infty)\times\bigl[\tfrac{1}{3},\tfrac{1}{2} \bigr) \times\bigl[
\tfrac{1}{3},\tfrac{1}{2} \bigr)\times(0,1 ) \times\mathbb N
\end{eqnarray}
satisfying
%
%
%e3.25 #&#
\begin{equation}
\label{defkappa} \cases{ \displaystyle \mbox{(a)}\quad \sum_{j=1}^{N_0}
\alpha^j\leq\xi<\sum_{j=1}^{N_0+1}
\alpha^j,
\vspace*{2pt}\cr
\displaystyle\mbox{(b)}\quad \alpha<\frac{\beta'}{\beta}<1,
\vspace*{2pt}\cr
\displaystyle\mbox{(c)}\quad \beta'-\frac{\eta}{2}+\frac{3}{2}
\alpha>0,
\vspace*{2pt}\cr
\displaystyle\mbox{(d)}\quad \bigl(\beta'+1\bigr)\wedge
\biggl(\beta'-\frac{\eta
}{2}+\frac{3\xi}{2} \biggr)>\eta.}
\end{equation}
[Note that we restate (\ref{alphaxi}) in (a).]
We insist that the parameter vector in (\ref{param}) is chosen to be
\emph{independent} of $i\in\mathbb N$ and ${\varepsilon}\in
(0,[8\psi(\1)]^{-1}\wedge1 ]$. For example, we can choose these
parameters in the following order:
first choose $\eta,\alpha,\beta',\xi$ according to (c) and (d),
choose $\beta$ according to (b), and finally choose $N_0$ according to
(a) by enlarging $\xi$ if necessary; the parameter $L$, however, can
be chosen arbitrarily. %\qed
\end{ass}

The following theorem gives our quantification of the local growth
rates of $Y$ under $\Q^i_{\varepsilon}$.

%
%
%th3.5 #&#
\begin{teo}\label{teogr}
Under Assumption~\ref{ass}, set three strictly positive constants by
%
%
%e3.26 #&#
\begin{eqnarray}\label{eqkappa}
\kappa_1&=&\bigl(\beta'+1\bigr)\wedge\biggl(
\beta'-\frac{\eta}{2}+\frac
{3\xi}{2} \biggr),\qquad
\kappa_2=\frac{\alpha^{N_0}}{4},
\nonumber\\[-8pt]\\[-8pt]\nonumber
\kappa_3 &=&
\beta'-\frac{\eta}{2}+\frac{3\alpha}{2}.
\end{eqnarray}
Then there exists a constant $K^*\in(0,\infty)$, depending only on
the parameter vector in (\ref{param}) and the immigration function
$\psi$, such that
for any $\delta\in(0,\kappa_1\wedge\kappa_3)$, the following
uniform bound holds:
%
%
%e3.27 #&#
%e3.28 #&#
\begin{eqnarray}
\label{ineqgrbdd} %
&& \Q^i_{\varepsilon} \biggl(\exists s\in(s_i,t],\sum_{j\in\mc
J^i_{\beta'}(s\wedge\tau^i\wedge\sigma^{X^i}_\beta)}Y^j_s(
\1 )^{\tau^i\wedge\sigma^{X^i}_\beta\wedge\sigma_\beta
^{Y^j}}\nonumber
\\
&&\hspace*{41pt} >K^* \bigl[(s-s_i)^{\kappa_1-\delta}+{
\varepsilon}^{\kappa
_2}\cdot(s-s_i)^{\kappa_3-\delta} \bigr]
\biggr)
\nonumber\\[-8pt]\\[-8pt]\nonumber
&&\qquad  \leq\frac{2\cdot2^{\kappa_1\vee\kappa
_3}}{2^{(N+1)\delta}(1-2^{-\delta})}
\\
\eqntext{\displaystyle \forall t\in\bigl[s_i+2^{-(N+1)},s_i+2^{-N}
\bigr], N\in\mathbb Z_+, i\in\mathbb N, {\varepsilon}\in\biggl(0,\frac
{1}{8\psi(\1)}
\wedge1 \biggr],}
\end{eqnarray}
where the $(\G_t)$-stopping times $\tau^i$ are defined in
Proposition~\ref{propseptime}.
\end{teo}

%
%
%re3.6 #&#
\begin{rmk}
\label{rmkcondsep}
If we follow the aforementioned interpretation of the parameter vector
in (\ref{param}) that
$(\eta,\beta',\xi)$ is close to $(1,\frac{1}{2},1)$,
then
$\kappa_1$ in (\ref{eqkappa}) is close to $ \frac{3}{2}$.
Informally, if we regard the stopping times $\tau^i$, $\sigma
^{X^i}_\beta$, and $\sigma_\beta^{Y^j}$ as being \emph{bounded
away} from $s_i$, then by the above reason for choosing the random
index sets $\mathcal J^i_{\beta'}(\cdot)$ in (\ref{eqdefIitt}), we
can regard Theorem~\ref{teogr} as a formalization of the statement in
(\ref{statcondsep+}). %\qed
\end{rmk}

In fact, the proof of Theorem~\ref{teogr} is reduced to a study of
some nonnegative $(\G_t)_{t\geq s_i}$-submartingale dominating the process
%
%
%e3.29 #&#
\begin{eqnarray}
\label{dissumYj} \sum_{j\in\mathcal J^i_{\beta'}(t\wedge\tau^i\wedge
\sigma
^{X^i}_\beta)}Y^j_t(
\1)^{\tau^i\wedge\sigma^{X^i}_\beta\wedge
\sigma^{Y^j}_\beta},\qquad t\in[s_i,\infty),
\end{eqnarray}
in (\ref{ineqgrbdd}), and the main task will be to prove Theorem~\ref
{teogr1} below. We explain the reductions as follows.

We observe that by Lemma~\ref{lemsep1}(2), the process in (\ref{dissumYj})
is dominated by the nonnegative process
%
%
%e3.30 #&#
%e3.31 #&#
\begin{eqnarray}
\label{dissubMG}
&&\sum_{j\in\mathcal J_{\beta'}^i(t,t\wedge\tau^i\wedge
\sigma ^{X^i}_\beta)} \biggl(\psi(
\1){\varepsilon}\nonumber
\\
&&\hspace*{76pt}{}+\int_{t_j}^{t\wedge\tau
^i\wedge\sigma^{X^i}_\beta\wedge\sigma^{Y^j}_\beta} \frac
{1}{X^i_s(\1)}
\int_\R X^i(x,s)^{1/2}Y^j(x,s)^{1/2}\,dx\,ds
\nonumber\\[-8pt]\\[-8pt]\nonumber
&&\hspace*{228pt}\qquad{} +M^j_{t\wedge\tau^i\wedge\sigma^{X^i}_\beta\wedge
\sigma^{Y^j}_\beta} \biggr),
\\
\eqntext{t \in[s_i,\infty),}
\end{eqnarray}
under $\Q^i_{\varepsilon}$ for any $i\in\mathbb N$ and ${\varepsilon
}\in(0,[8\psi(\1)]^{-1}\wedge1 ]$.
The process in (\ref{dissubMG}) is in fact a nonnegative $(\G
_t)_{t\geq s_i}$-submartingale under $\Q^i_{\varepsilon}$, since for
any $j\in\mathbb N$ with $s_i<t_j$, $j\in\mathcal J_{\beta
'}^i(t,t\wedge\tau^i\wedge\sigma^{X^i}_\beta)$ if and only if the
following $\G_{t_j}$-event occurs:
\[
\bigl\{\mbox{$\llvert y_j-x_i\rrvert\leq2 \bigl({
\varepsilon}^{1/2}+(t-s_i)^{\beta'} \bigr)$ and
$t_j\leq t\wedge\tau^i\wedge\sigma^{X^i}_\beta$}
\bigr\}.
\]
(Recall that $y_j\in\G_{t_j}$ and $x_i\in\G_{s_i}$ by Theorem~\ref{teoapprox-v2}.)
It suffices to prove the bound~(\ref{ineqgrbdd}) of Theorem~\ref{teogr} with the involved process in (\ref{dissumYj}) replaced by
the nonnegative submartingale in (\ref{dissubMG}).
To further reduce the problem, we resort to the following simple
corollary of Doob's maximal inequality.

%
%
%le3.7 #&#
\begin{lem}\label{lemDoob}
Let $F$ be a nonnegative function on $[0,1]$ such that $F\rest(0,1]>0$ and
$\sup_{s,t\dvtx1\leq {t}/{s}\leq2}\frac{F(t)}{F(s)}<\infty$. In
addition, assume that
%
%
%e3.32 #&#
\begin{equation}
\mbox{for some $\delta>0$},\qquad t\lmt\frac{F(t)}{t^\delta}\mbox{ is
increasing}.\label{assDoob}
\end{equation}
Suppose that $Z$ is a nonnegative submartingale with c\`adl\`ag sample
paths such that
$\E[Z_t]\leq F(t)$ for any $t\in[0,1]$.
Then for every $N\in\mathbb Z_+$,
%
%
%e3.33 #&#
\begin{eqnarray}
\label{ineqbadprob} %
&&\sup_{t\in[2^{-(N+1)},2^{-N} ]}\P\biggl(\exists s\in
(0,t ], Z_s> \frac{F(s)}{s^\delta} \biggr)
\nonumber\\[-8pt]\\[-8pt]\nonumber
&&\qquad \leq\biggl(\sup_{s,t\dvtx1\leq {t}/{s}\leq2}\frac
{F(t)}{F(s)} \biggr)
\times\frac{1}{2^{(N+1)\delta}(1-2^{-\delta})}.
\end{eqnarray}
\end{lem}
\begin{pf}
For each $m\in\mathbb Z_+$,
\begin{eqnarray*}
&&\P\biggl(\exists s\in\bigl[2^{-(m+1)},2^{-m} \bigr],
Z_s\geq\frac
{F(s)}{s^\delta} \biggr)
\\
&&\qquad \leq \P\biggl(\sup_{2^{-(m+1)}\leq s\leq2^{-m}}Z_s\geq F \biggl(
\frac{1}{2^{(m+1)}} \biggr) \Big/ \frac{1}{2^{(m+1)\delta}} \biggr)
\\
&&\qquad \leq\frac{\E[Z_{1/2^m}]}{F ({1}/{2^{(m+1)}} ) /( {1}/{2^{(m+1)\delta}})}
\\
&&\qquad \leq\frac
{F ({1}/{2^m} )}{F ({1}/{2^{(m+1)}}
) /({1}/{2^{(m+1)\delta}})}
\\
&&\qquad =\sup_{s,t\dvtx1\leq {t}/{s}\leq2}
\frac{F(t)}{F(s)}\times\frac
{1}{2^{(m+1)\delta}},
\end{eqnarray*}
where the first inequality follows from (\ref{assDoob}) and the
second inequality follows from Doob's maximal inequality. Hence,
whenever $t\in[2^{-(N+1)},2^{-N} ]$ for $N\in\mathbb Z_+$,
the last inequality gives
\begin{eqnarray*}
&& \P\biggl(\exists s\in(0,t ], Z_s> \frac{F(s)}{s^\delta
} \biggr)
\\
&&\qquad \leq \sum_{m=N}^\infty\P\biggl(\exists s\in
\bigl[2^{-(m+1)},2^{-m} \bigr], Z_s\geq
\frac{F(s)}{s^\delta} \biggr)
\\
&&\qquad \leq \biggl(\sup_{s,t\dvtx1\leq {t}/{s}\leq2}\frac
{F(t)}{F(s)} \biggr)\sum
_{m=N}^\infty\frac{1}{2^{(m+1)\delta}}
\\
&&\qquad = \biggl(\sup_{s,t\dvtx1\leq {t}/{s}\leq2}\frac{F(t)}{F(s)} \biggr
)\times
\frac{1}{2^{(N+1)\delta}(1-2^{-\delta})}.
\end{eqnarray*}
This completes the proof.
\end{pf}

%
%
%th3.8 #&#
\begin{teo}\label{teogr1}
Under Assumption~\ref{ass}, take the same constants $\kappa_j$ as in
Theorem~\ref{teogr}. Then we can choose a constant
$K^*\in(0,\infty)$ as stated in Theorem~\ref{teogr}, such that the
following uniform bound holds:
%
%e3.34 #&#
\begin{eqnarray}
\hspace*{-4pt}&&\E^{\mathbb Q^i_{\varepsilon}} \biggl[ \sum_{j\in\mathcal J_{\beta
'}^i(t,t\wedge\tau^i\wedge\sigma^{X^i}_\beta)}\! \biggl(\psi(
\1 ){\varepsilon}\nonumber
\\[-2pt]
\hspace*{-4pt}&&\hspace*{98pt}{} +\int_{t_j}^{t\wedge\tau^i\wedge\sigma^{X^i}_\beta
\wedge\sigma^{Y^j}_\beta}\!\!
\frac{1}{X^i_s(\1)}\int_\R X^i(x,s)^{1/2}Y^j(x,s)^{1/2}\,dx\,ds
\biggr) \biggr]\nonumber
\\
\hspace*{-4pt}&&\qquad \leq K^* \bigl[(t-s_i)^{\kappa_1}+{\varepsilon
}^{\kappa_2}\cdot(t-s_i)^{\kappa_3} \bigr]\nonumber
\\[-2pt]
\hspace*{-4pt}\eqntext{\displaystyle \forall t
\in(s_i,s_i+1], i\in\mathbb N, {\varepsilon}\in\biggl(0,
\frac{1}{8\psi(\1
)}\wedge1 \biggr].}
\end{eqnarray}
\end{teo}

Now, we prove the main result of this section, that is Theorem~\ref
{teogr}, assuming Theorem~\ref{teogr1}.

\begin{pf*}{Proof of Theorem~\ref{teogr}}
In, and only in, this proof, we denote by $Z^{(0)}$ the submartingale
defined in (\ref{dissubMG}).

Since $ [j\in\mathcal J_{\beta'}^i(t,t\wedge\tau^i\wedge\sigma
^{X^i}_\beta) ]\in\G_{t_j}$, we obtain immediately from
Lem\-ma~\ref{lemsep1}(2) that the part
\[
\sum_{j\in\mathcal J_{\beta'}^i(t,t\wedge\tau^i\wedge\sigma
^{X^i}_\beta)}M^j_{t\wedge\tau^i\wedge\sigma^{X^i}_\beta\wedge
\sigma^{Y^j}_\beta},\qquad t
\in[s_i,\infty),
\]
in the definition of $Z^{(0)}$
is a true $\Q^i_{\varepsilon}$-martingale with mean zero, for any
$i\in\mathbb N$ and ${\varepsilon}\in(0,[8\psi(\1)]^{-1}\wedge
1 ]$. Hence, setting
\[
F^{(0)}(s)=K^* \bigl(s^{\kappa_1}+{\varepsilon}^{\kappa_2}
\cdot s^{\kappa_3} \bigr),\qquad s\in[0,1],
\]
we see from Theorem~\ref{teogr1} that
\[
\E_{\varepsilon} \bigl[Z^{(0)}_{t} \bigr]\leq
F^{(0)}(t-s_i)
\]
for any $t\in(s_i,s_i+1]$, $i\in\mathbb N$ and ${\varepsilon}\in
(0,[8\psi(\1)]^{-1}\wedge1 ]$. Note that
\begin{eqnarray*}
\sup_{s,t\dvtx1\leq{t}/{s}\leq2}\frac{F^{(0)}(t)}{F^{(0)}(s)}&\leq&\sup_{s,t\dvtx1\leq{t}/{s}\leq2}
\biggl(\frac{t^{\kappa
_1}}{s^{\kappa_1}}+\frac{t^{\kappa_3}}{s^{\kappa_3}} \biggr)\leq2\cdot
2^{\kappa_1\vee\kappa_3}.
\end{eqnarray*}
Hence, applying Lemma~\ref{lemDoob} with $(Z,F)$ taken to be
$(Z^{(0)},F^{(0)})$, we see that (\ref{ineqgrbdd}) with the involved
process in (\ref{dissumYj}) replaced by $Z^{(0)}$ holds. The proof is
complete.
\end{pf*}

\begin{figure}%[t]

\includegraphics{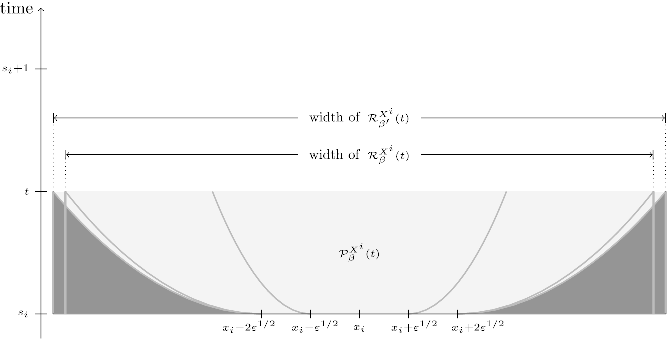}

\caption{$\mathcal P^{X^i}_\beta(t)$, $\mathcal R^{X^i}_{\beta}(t)$,
and $\mathcal R^{X^i}_{\beta'}(t)$ for $0<\beta'<\beta$ and $t\in
[s_i,s_i+1]$.}
\label{Fig2}
\end{figure}

The remainder of this section is to prove Theorem~\ref{teogr1}. For
this purpose, we need to classify the clusters $Y^j$ for $j\in\mathcal
J^i_{\beta'}(t,t\wedge\tau^i\wedge\sigma^{X^i}_\beta)$. Set
\begin{eqnarray*}
\mathcal C^i_{\beta'}(t)&\triangleq& \bigl\{j\in\mathbb N;
\llvert y_j-x_i\rrvert< 2 \bigl({\varepsilon}^{1/2}+(t_j-s_i)^{\beta'}
\bigr),s_i<t_j\leq t \bigr\},
\\[-1pt]
\mathcal L^i_{\beta'}\bigl(t,t'\bigr)
&\triangleq& \bigl\{j\in\mathbb N;2 \bigl({\varepsilon
}^{1/2}+(t_j-s_i)^{\beta'}
\bigr)\leq\llvert y_j-x_i\rrvert\leq2 \bigl({
\varepsilon}^{1/2}+(t-s_i)^{\beta'} \bigr),
\\[-1pt]
&&\hspace*{229pt} s_i<t_j\leq t' \bigr\},
\end{eqnarray*}
for $t', t\in(s_i,\infty)$ with $t\geq t'$. Hence, as far as the
clusters $Y^j$ landing inside the rectangle $\mathcal R^{X^i}_{\beta
'}(t)$ are concerned, the clusters $Y^j$, $j\in\mathcal C^i_{\beta
'}(t)$, are those landing inside the double parabola
\[
\bigl\{(x,s)\in\R\times[s_i,t];\llvert x-x_i\rrvert< 2
\bigl({\varepsilon}^{1/2}+(s-s_i)^{\beta'} \bigr)
\bigr\}
\]
(the light grey area in Figure~\ref{Fig2}), and the clusters $Y^j$,
$j\in\mathcal L^i_{\beta'}(t,t)$, are those landing outside (the dark
grey area in Figure~\ref{Fig2}).
For any $i\in\mathbb N$, we say a cluster $Y^j$ is a \textit{critical
cluster} if $j\in\mathcal C^i_{\beta'}(t)$ and a \textit{lateral
cluster} if $j\in\mathcal L^i(t,t')$ for some $t,t'$.

%
%
%f3 #&#

Since $ \{\mathcal C^i_{\beta'}(t),\mathcal L^i_{\beta
'}(t,t') \}$ is a cover of $\mathcal J^i_{\beta'}(t,t')$ by
disjoint sets,
Theorem~\ref{teogr1} can be obtained by the following two lemmas.

%
%
%le3.9 #&#
\begin{lem}\label{lemgr1-1}
Let $\kappa_j$ be as in Theorem~\ref{teogr}.
We can choose a constant $K^*\in(0,\infty)$ as in Theorem~\ref{teogr} such that the following uniform bound holds:
%
%
%e3.35 #&#
%e3.36 #&#
\begin{eqnarray}\label{sumcent}
\hspace*{-2pt}&& \E^{\Q^i_{\varepsilon}} \biggl[ \sum_{j\in\mathcal C^i_{\beta
'}(t\wedge\tau^i\wedge\sigma^{X^i}_\beta)} \biggl(\psi(
\1 ){\varepsilon}\nonumber
\\
\hspace*{-2pt}&&\hspace*{92pt}{}+\int_{t_j}^{t\wedge\tau^i\wedge\sigma^{X^i}_\beta
\wedge\sigma^{Y^j}_\beta}
\frac{1}{X^i_s(\1)}\int_\R X^i(x,s)^{1/2}Y^j(x,s)^{1/2}\,dx\,ds
\biggr) \biggr]
\nonumber\\[-8pt]\hspace*{-2pt}\\[-8pt]\nonumber
\hspace*{-2pt} &&\qquad \leq\frac{K^*}{2}
\bigl[(t-s_i)^{\kappa
_1}+{\varepsilon}^{\kappa_2}
\cdot(t-s_i)^{\kappa_3} \bigr]
\\
\hspace*{-2pt}\eqntext{\displaystyle\forall t
\in(s_i,s_i+1], i\in\mathbb N, {\varepsilon}\in\biggl(0,
\frac{1}{8\psi(\1)}\wedge1 \biggr].}
\end{eqnarray}
\end{lem}

%
%
%le3.10 #&#
\begin{lem}\label{lemgr1-2}
Let $\kappa_j$ be as in Theorem~\ref{teogr}. By enlarging the
constant $K^*$ in Lemma~\ref{lemgr1-1} if necessary, the following
uniform bound holds:
%
%
%e3.37 #&#
%e3.38 #&#
\begin{eqnarray}\label{sumlate}
&&\E^{\Q^i_{\varepsilon}} \biggl[ \sum_{j\in\mathcal L^i_{\beta
'}(t,t\wedge\tau^i\wedge\sigma^{X^i}_\beta)} \biggl(\psi(
\1 ){\varepsilon}\nonumber
\\[-4pt]
&&\hspace*{99pt}{} +\int_{t_j}^{t\wedge\tau^i\wedge\sigma^{X^i}_\beta
\wedge\sigma^{Y^j}_\beta}
\frac{1}{X^i_s(\1)}\int_\R X^i(x,s)^{1/2}\nonumber
\\[-4pt]
&&\hspace*{222pt}{}\times Y^j(x,s)^{1/2}\,dx\,ds
\biggr) \biggr]
\\[-4pt]
&&\qquad  \leq\frac{K^*}{2}
\bigl[(t-s_i)^{\kappa
_1}+{\varepsilon}^{\kappa_2}
\cdot(t-s_i)^{\kappa_3} \bigr]\nonumber
\\[-3pt]
\eqntext{\displaystyle\forall t \in(s_i,s_i+1], i\in\mathbb N, {\varepsilon}\in\biggl(0,
\frac{1}{8\psi(\1)}\wedge1 \biggr].}
\end{eqnarray}
\end{lem}

Despite some technical details, the methods of proof for 
Lemmas~\ref{lemgr1-1}~and\break \ref{lemgr1-2} are very similar. For clarity,
they are given in Sections~\ref{seccondsep-1}~and~\ref{seccondsep-2} separately, with some preliminaries set in Section~\ref
{secpre} below.

%s3.3 #&#
\subsection{Auxiliary results and notation}\label{secpre}
For each $z,\delta\in\R_+$, let $ (Z,\mathbf P_z^\delta)$
denote a copy of $\frac{1}{4}\BES Q^{4\delta}(4z)$. We assume that
$ (Z,\mathbf P^\delta_z )$ is defined by a $(\ms
H_t)$-Brownian motion $B$, where $(\ms H_t)$ satisfies the usual
conditions. This means that
\[
Z_t=z+\delta t+\int_0^t
\sqrt{Z_s}\,dB_s,\quad\mathbf P^\delta_z
\mbox{-a.s.}
\]
(Cf. Section XI.1 of \cite{RYCMB} for Bessel squared processes.) As
we will often investigate $Z$ before it hits a constant level,
we set the following notation similar to (\ref{eqTZ1}): for any
real-valued process $H=(H_t)$
\[
T^H_x=\inf\{t\geq0;H_t=x\},\qquad x\in\R.
\]

For $\delta=0$, $ (Z,\mathbf P^0_z )$ gives a Feller diffusion
and its marginals are characterized by
\[
\E^{\mathbf P^0_z} \bigl[\exp(-\lambda Z_t ) \bigr]=\exp\biggl(
\frac{-2\lambda z}{2+\lambda t} \biggr),\qquad\lambda,t\in\R_+.
\]
In particular, the survival probability of $ (Z,\mathbf P^0_z
)$ is given by
%
%
%e3.39 #&#
%e3.40 #&#
\begin{eqnarray}
\label{eqsurvprob} %
\mathbf P^0_z(Z_t>0)&=&
\lim_{\lambda\to\infty} \bigl(1-\E^{\mathbf
P^0_z} \bigl[\exp(-\lambda
Z_t ) \bigr] \bigr) =1-\exp\biggl(-\frac{2z}{t} \biggr),
\nonumber\\[-9pt]\\[-9pt]
\eqntext{z,t\in(0,\infty).}
\end{eqnarray}
Using the elementary inequality $1-e^{-x}\leq x$ for $x\in\R_+$, we
obtain from the last inequality that
%
%
%e3.41 #&#
\begin{equation}
\mathbf P^0_z(Z_t>0)\leq
\frac{2z}{t},\qquad z,t\in(0,\infty).\label
{ineqsurvprob}
\end{equation}

To save notation in the following Sections~\ref{seccondsep-1}~and~\ref{seccondsep-2}, we write $A\propleq B$ if
$A\leq CB$ for some constant $C\in(0,\infty)$ which may vary from
line to line but depends only on $\psi$ \emph{and} the parameter
vector chosen in Assumption~\ref{ass}.

%s3.4 #&#
\subsection{Proof of Lemma~\texorpdfstring{\protect\ref{lemgr1-1}}{3.9}}\label{seccondsep-1}
Fix $i\in\mathbb N$ and ${\varepsilon}\in(0,[8\psi(\1
)]^{-1}\wedge1 ]$, and henceforth we drop the subscripts
${\varepsilon}$ of $\P_{\varepsilon}$ and $\Q^i_{\varepsilon}$. In
addition, we may only consider
$t\in[s_i+\frac{{\varepsilon}}{2},s_i+1 ]$
as there are no immigrants for $Y$ arriving in $[s_i,s_i+\frac
{{\varepsilon}}{2})$.
Our analysis proceeds with the following steps.

\begin{longlist}[\textit{Step} 1.]
\item[\textit{Step} 1.]
We start with the simplification:
%
%
%e3.42 #&#
\begin{eqnarray}\label{ineqcent-bdd}
\qquad&&\sum_{j\in\mathcal C^i_{\beta'}(t\wedge\tau^i\wedge\sigma
^{X^i}_\beta)} \biggl(\psi(\1){\varepsilon}\nonumber
\\
&&\hspace*{68pt}{} +\int
_{t_j}^{t\wedge\tau
^i\wedge\sigma^{X^i}_\beta\wedge\sigma^{Y^j}_\beta} \frac{1}{X^i_s(\1
)}\int
_\R X^i(x,s)^{1/2}Y^j(x,s)^{1/2}\,dx\,ds
\biggr)
\nonumber\\[-8pt]\\[-8pt]\nonumber
&&\qquad \leq\sum_{j\in\mathcal C^i_{\beta'}(t\wedge\tau
^i)} \biggl(\psi(\1){
\varepsilon}+\int_{t_j}^{t\wedge\tau^i} \frac{1}{ [X^i_s(\1) ]^{1/2}}
\bigl[Y^j_s(\1) \bigr]^{1/2}\,ds \biggr)
\\
&&\qquad \leq\sum_{j\in\mathcal C^i_{\beta'}(t\wedge\tau
^i)} \biggl(\psi(\1){
\varepsilon}+\int_{t_j}^{t\wedge\tau^i} \frac{2}{(s-s_i)^{\eta/2}}
\bigl[Y^j_s(\1) \bigr]^{1/2}\,ds
\biggr),\nonumber
\end{eqnarray}
where the first inequality follows from the Cauchy--Schwarz inequality
and the second one follows by using the component $\tau^{i,(1)}$ of
$\tau^i$ in Proposition~\ref{propseptime}.

We claim that
%
%
%e3.43 #&#
\begin{eqnarray}\label{ineqsumpre}
&&\E^{\mathbb Q^i} \biggl[ \sum_{j\in\mathcal C^i_{\beta'}(t\wedge\tau
^i\wedge\sigma^{X^i}_\beta)} \biggl(\psi(
\1){\varepsilon}\nonumber
\\
&&\hspace*{91pt}
{}+\int_{t_j}^{t\wedge\tau^i\wedge\sigma^{X^i}_\beta\wedge
\sigma_\beta^{Y^j}} \frac{1}{X^i_s(\1)}
\int_\R X^i(x,s)^{1/2}\nonumber
\\
&&\hspace*{214pt} {}\times Y^j(x,s)^{1/2}\,dx\,ds
\biggr) \biggr]
\nonumber\\[-8pt]\\[-8pt]\nonumber
&&\qquad \propleq\sum_{j\dvtx s_i<t_j\leq t}
(t_j-s_i)^{\beta'}{\varepsilon}
\\
&&\quad\qquad{} + \sum
_{j\dvtx s_i<t_j\leq t}\int_{t_j}^t \,ds
\frac{1}{(s-s_i)^{\eta/2}} \E^{\mathbb Q^i} \bigl[ \bigl[Y^j_s(
\1) \bigr]^{1/2};
s<\tau^i,\nonumber
\\
&&\hspace*{183pt} \llvert y_j-x_i\rrvert< 2 \bigl({\varepsilon}^{1/2}+(t_j-s_i)^{\beta'}
\bigr) \bigr].\nonumber
\end{eqnarray}
Note that
%
%
%e3.44 #&#
\begin{eqnarray}\label{ineqQbdd-c}
&& \E^{\Q^i} \bigl[\psi(\1){\varepsilon}\#\mathcal C^i_{\beta
'} \bigl(t\wedge\tau^i\bigr) \bigr]\nonumber
\\
&&\qquad \propleq {\varepsilon}
\E^{\mathbb Q^i} \bigl[\#\mathcal C^i_{\beta'}(t) \bigr]\nonumber
\\
&&\qquad =\varepsilon\sum_{j\dvtx s_i<t_j\leq
t}\Q^i \bigl(\llvert
y_j-x_i\rrvert< 2 \bigl({\varepsilon}^{1/2}+(t_j-s_i)^{\beta
'}\bigr) \bigr)
\\
&&\qquad \propleq\sum_{j\dvtx s_i<t_j\leq t}4 \bigl({\varepsilon
}^{1/2}+(t_j-s_i)^{\beta'} \bigr){
\varepsilon}\nonumber
\\
&&\qquad \propleq\sum_{j\dvtx s_i<t_j\leq t}(t_j-s_i)^{\beta'}
{\varepsilon},\nonumber
\end{eqnarray}
where the second $\propleq$-inequality follows from Lemma~\ref
{lemsep1}(4),
and the last $\propleq$-inequality follows since
%
%
%e3.45 #&#
\begin{equation}
\label{ineqtjvep} {\varepsilon}^{1/2}\leq{\varepsilon}^{\beta'}
\leq2^{\beta
'}(t_j-s_i)^{\beta'} \qquad
\forall j\in\mathbb N\mbox{ with }t_j>s_i.
\end{equation}
Our claim (\ref{ineqsumpre}) follows from (\ref{ineqcent-bdd}) and
(\ref{ineqQbdd-c}).

From the display (\ref{ineqsumpre}), we see the necessity to obtain
the order of
%
%
%e3.46 #&#
%e3.47 #&#
\begin{eqnarray}
\label{defS2INT} \E^{\Q^i} \bigl[ \bigl[Y^j_s(
\1) \bigr]^{1/2};s<\tau^i,\llvert y_j-x_i
\rrvert< 2 \bigl({\varepsilon}^{1/2}+(t_j-s_i)^{\beta'}
\bigr) \bigr],
\nonumber\\[-8pt]\\[-8pt]
\eqntext{s\in(t_j,t], s_i<t_j<t,}
\end{eqnarray}
in $s_i,t_j,s,t$.

We subdivide our analysis of a generic term in (\ref{defS2INT}) into
the following \mbox{steps 2-1--2-4}, with a summary given in step~2-5.
\end{longlist}

\begin{longlist}[\textit{Step} 2-1.]
\item[\textit{Step} 2-1.] We convert the $\Q^i$-expectations in
(\ref{defS2INT}) to $\P$-expectations. Recalling that $x_i,y_j\in\G
_{t_j}$ by (\ref{impilc}), we can use Lemma~\ref{lemQi} to get
%
%
%e3.48 #&#
\begin{eqnarray}\label{eqsepINT1}
&&\E^{\Q^i} \bigl[ \bigl[Y^j_s(\1)
\bigr]^{1/2};s< \tau^i,\llvert y_j-x_i
\rrvert< 2 \bigl({\varepsilon}^{1/2}+(t_j-s_i)^{\beta'}
\bigr) \bigr]
\nonumber
\\
&&\qquad =\frac{1}{\psi(\1){\varepsilon}}\E^{\P}
\bigl[X^i_s(\1 )^{T^{X^i}_1} \bigl[Y^j_s(
\1) \bigr]^{1/2};s< \tau^i,
\\
&&\hspace*{80pt} \llvert y_j-x_i
\rrvert< 2 \bigl({\varepsilon}^{1/2}+(t_j-s_i)^{\beta'}
\bigr) \bigr].\nonumber
\end{eqnarray}

We break the $\P$-expectation in (\ref{eqsepINT1}) into finer pieces
by considering the following.
For $s>t_j$, $X^i(\1)_s^{T^{X^i}_1}$ is nonzero on the union of the two
disjoint events:
%
%
%e3.49 #&#
\begin{equation}
\label{event1} \bigl\{X^i_s(\1)^{T^{X^i}_1}>0,T_0^{X^i}
\leq t_j \bigr\}= \bigl\{ T_1^{X^i}<T_0^{X^i}
\leq t_j \bigr\}
\end{equation}
and
%
%
%e3.50 #&#
\begin{equation}
\label{event2} \bigl\{X^i_s(\1)^{ T^{X^i}_1}>0,t_j<T^{X^i}_0
\bigr\}.
\end{equation}
Here, the equality in (\ref{event1}) holds $\P$-a.s. since $0$ is an
absorbing state of $X^i(\1)$ under $\P$. In fact, $X^i(\1
)_s^{T_1^{X^i}}=1$ on the event in (\ref{event1}).
To invoke the additional order provided by the improved modulus of
continuity of $X^i(\1)$ at its starting point $s_i$, we use the trivial
inequality
\[
X^i_s(\1)^{ T^{X^i}_1}\leq\bigl\llvert
X^i_s(\1)^{T^{X^i}_1}-\psi(\1 ){\varepsilon}\bigr
\rrvert+\psi(\1){\varepsilon}
\]
on the event (\ref{event2}).

Putting things together, we see from (\ref{eqsepINT1}) that
%
%
%e3.51 #&#
%e3.52 #&#
\begin{eqnarray}\label{ineqsepINT1-2}
&& \E^{\Q^i} \bigl[ \bigl[Y^j_s(\1)
\bigr]^{1/2};s< \tau^i,\llvert y_j-x_i
\rrvert<2 \bigl({\varepsilon}^{1/2}+(t_j-s_i)^{\beta'}
\bigr) \bigr]\nonumber
\\
&&\qquad \leq\frac{1}{\psi(\1){\varepsilon}}\E^{\P} \bigl[
\bigl[Y^j_s(\1 ) \bigr]^{1/2};s\leq
T_1^{Y^j},\nonumber
\\
&&\hspace*{80pt}  \llvert y_j-x_i\rrvert<2 \bigl({\varepsilon}^{1/2}+(t_j-s_i)^{\beta
'}\bigr),T_1^{X^i}<T_0^{X^i}\leq t_j \bigr]\nonumber
\\
&&\quad\qquad{}+\frac{1}{\psi(\1){\varepsilon}}\E^{\P} \bigl[\bigl\llvert X^i_s(
\1 )^{T^{X^i}_1}-\psi(\1){\varepsilon}\bigr\rrvert\bigl[Y^j_s(
\1) \bigr]^{1/2};s<\tau^i,
\nonumber\\[-8pt]\\[-8pt]\nonumber
&&\hspace*{93pt} \llvert y_j-x_i\rrvert<2 \bigl({
\varepsilon}^{1/2}+(t_j-s_i)^{\beta '}
\bigr),X^i_s(\1)^{T^{X^i}_1}>0,t_j<T_0^{X^i}\bigr]
\\
&&\quad\qquad{}+\frac{1}{\psi(\1){\varepsilon}}\cdot\psi(\1){\varepsilon}\E^{\P}
\bigl[\bigl[Y^j_s(\1) \bigr]^{1/2};s\leq T^{Y^j}_1,\nonumber
\\
&&\hspace*{129pt} \llvert y_j-x_i\rrvert<2 \bigl({
\varepsilon}^{1/2}+(t_j-s_i)^{\beta
'}\bigr),t_j<T_0^{X^i} \bigr]\nonumber
\\
\eqntext{\forall s \in(t_j,t], s_i<t_j<t,}
\end{eqnarray}
where for the first and the third terms on the right-hand side, it is
legitimate to replace the event $ \{s<\tau^i \}$ by the
larger one $ \{s\leq T^{Y^j}_1 \}$ since, in Proposition~\ref
{propseptime}, $\tau^{i,(3)}$ is a component of $\tau^i$, and for the
third term we replace the event in (\ref{event2}) by the larger one
$ \{t_j<T^{X^i}_0 \}$.

In steps 2-2--2-4 below, we derive a bound for each of the
three terms in (\ref{ineqsepINT1-2}) which involves only Feller's
diffusion. We use the notation in Section~\ref{secpre}.
\end{longlist}

\begin{longlist}[\textit{Step} 2-2.]
\item[\textit{Step} 2-2.]
Consider the first term on the right-hand side of (\ref
{ineqsepINT1-2}), and recall the notation in Section~\ref{secpre}.
It follows from (\ref{impilp}) and (\ref{dom}) that
%
%
%e3.53 #&#
%e3.54 #&#
\begin{eqnarray}\label{ineqYjhalf1}
&&\frac{1}{\psi(\1){\varepsilon}}\E^{\P} \bigl[ \bigl[Y^j_s(
\1 ) \bigr]^{1/2};s\leq T_1^{Y^j},\nonumber
\\
&&\hspace*{47pt} \llvert
y_j-x_i\rrvert<2 \bigl({\varepsilon}^{1/2}+(t_j-s_i)^{\beta'}
\bigr),T_1^{X^i}<T_0^{X^i}\leq
t_j \bigr]\nonumber
\\
&&\qquad \propleq\frac{1}{{\varepsilon}}\P\bigl(T_1^{X^i}<T_0^{X^i}
\leq t_j \bigr) \bigl({\varepsilon}^{1/2}+(t_j-s_i)^{\beta'}
\bigr)\nonumber
\\
&&\quad\qquad{}\times  \E^{\mathbf P^0_{\psi(\1){\varepsilon}}} \bigl[ (Z_{s-t_j} )^{1/2};s-t_j
\leq T^Z_1 \bigr]
\\
&&\qquad \leq\frac{1}{{\varepsilon}}\P\bigl(T_1^{X^i}<T_0^{X^i}
\bigr) \bigl({\varepsilon}^{1/2}+(t_j-s_i)^{\beta'}
\bigr) \nonumber
\\
&&\quad\qquad{}\times \E^{\mathbf P^0_{\psi(\1){\varepsilon}}} \bigl[ (Z_{s-t_j} )^{1/2};s-t_j
\leq T^Z_1 \bigr]
\nonumber
\\
&&\qquad \propleq(t_j-s_i)^{\beta'}
\E^{\mathbf P^0_{\psi(\1){\varepsilon}}} \bigl[ (Z_{s-t_j} )^{1/2};s-t_j
\leq T^Z_1 \bigr]\nonumber
\\
\eqntext{\forall s\in(t_j,t],
s_i<t_j<t,}
\end{eqnarray}
where the last inequality follows from (\ref{ineqtjvep}) and
Lemma~\ref{lemQi}.
\end{longlist}

\begin{longlist}[\textit{Step} 2-3.]
\item[\textit{Step} 2-3.]
Let us deal with the second term in (\ref{ineqsepINT1-2}).
We claim that%
%
%e3.55 #&#
\begin{eqnarray} \label{ineqYjhalf2}
\qquad && \frac{1}{\psi(\1){\varepsilon}}\E^{\P} \bigl[\bigl\llvert X^i_s(
\1)^{T^{X^i}_1}-\psi(\1){\varepsilon}\bigr\rrvert\bigl[Y^j_s(
\1) \bigr]^{1/2};s< \tau^i,
\nonumber
\\
&&\hspace*{48pt}  \llvert y_j-x_i\rrvert< 2 \bigl({
\varepsilon}^{1/2}+(t_j-s_i)^{\beta
'}\bigr),X^i_s(\1)^{T^{X^i}_1}>0,t_j<T_0^{X^i}
\bigr]
\nonumber\\[-8pt]\\[-8pt]\nonumber
&&\qquad \propleq \bigl({\varepsilon}^{\alpha^{N_0}}(s-s_i)^\alpha
+(s-s_i)^\xi\bigr) (t_j-s_i)^{\beta'-1}
\\
&&\quad\qquad{}\times \E^{\mathbf P^0_{\psi(\1){\varepsilon}}} \bigl[(Z_{s-t_j})^{1/2};s-t_j
\leq T^Z_1 \bigr]\qquad\forall s\in(t_j,t],
s_i<t_j<t.\nonumber
\end{eqnarray}
Fix such $s$ throughout step 2-3.

First, let us transfer the improved modulus of $X^i(\1)$ under $\Q^i$
to one under $\P$. It follows from (\ref{ineqXQmc}) that on $ \{
s< \tau^i,X^i(\1)_s^{T_1^{X^i}}>0 \}\in\G_s$, we have
\[
\bigl\llvert X^i_s(\1)^{T^{X^i}_1}-\psi(\1){
\varepsilon}\bigr\rrvert\leq K_1^X\bigl[\psi(\1){
\varepsilon}\bigr]^{\alpha^{N_0}}(s-s_i)^\alpha+K_2^X
(s-s_i)^\xi\qquad\Q^i\mbox{-a.s.}
\]
and hence
%
%
%e3.56 #&#
\begin{eqnarray}\label{eqmodaux}
\qquad 0&=&\Q^i \bigl(\bigl\llvert X^i_s(
\1)^{T^{X^i}_1}-\psi(\1){\varepsilon}\bigr\rrvert>K_1^X
\bigl[\psi(\1){\varepsilon}\bigr]^{\alpha^{N_0}}(s-s_i)^\alpha\nonumber
\\
&&\hspace*{59pt} {}
+K_2^X(s-s_i)^\xi, s< \tau^i,X^i_s(
\1)^{T_1^{X^i}}>0 \bigr)
\nonumber\\[-8pt]\\[-8pt]\nonumber
 &=&\frac{1}{\psi(\1){\varepsilon}}\E^{\P}
\bigl[X^i_s(\1 )^{T_1^{X^i}};
\bigl\llvert X^i_s(\1)^{ T^{X^i}_1}-\psi(\1){\varepsilon}\bigr
\rrvert > K_1^X\bigl[\psi(\1){\varepsilon}
\bigr]^{\alpha
^{N_0}}(s-s_i)^\alpha
\\
&&\hspace*{138pt}{}+K_2^X(s-s_i)^\xi,
s<\tau^i,X^i_s(\1 )^{T_1^{X^i}}>0 \bigr],\nonumber
\end{eqnarray}
where the last equality follows from Lemma~\ref{lemQi} since the
event evaluated under $\Q^i$ is a $\G_s$-event.
Using the restriction $X^i_s(\1)^{T_1^{X^i}}>0$, we see that the
equality (\ref{eqmodaux}) implies
%
%
%e3.57 #&#
%e3.58 #&#
\begin{eqnarray}
\label{ineqimcP} &&\bigl\llvert X_s^i(
\1)^{ T^{X^i}_1}-\psi(\1){\varepsilon}\bigr\rrvert\leq K_1^X
\bigl[\psi(\1){\varepsilon}\bigr]^{\alpha^{N_0}}(s-s_i)^\alpha
+K_2^X(s-s_i)^\xi
\nonumber\\[-8pt]\\[-8pt]
\eqntext{\displaystyle \mbox{$\mathbb P$-a.s. on $ \bigl[s<\tau^i,
X^i_s(\1)^{T_1^{X^i}}>0 \bigr]$}.}
\end{eqnarray}

Using (\ref{ineqimcP}), we obtain
%
%
%e3.59 #&#
\begin{eqnarray}\label{ineqYjhalf2-0}
&&\frac{1}{\psi(\1){\varepsilon}}\E^{\P} \bigl[\bigl\llvert X^i_s(
\1 )^{T^{X^i}_1}-\psi(\1){\varepsilon}\bigr\rrvert\bigl[Y^j_s(
\1) \bigr]^{1/2};s< \tau^i,
\nonumber
\\
&&\hspace*{48pt} \llvert y_j-x_i\rrvert<2 \bigl({
\varepsilon}^{1/2}+(t_j-s_i)^{\beta
'}\bigr),X^i_s(\1)^{T_1^{X^i}}>0,t_j<T_0^{X^i}
\bigr]
\nonumber\\[-8pt]\\[-8pt]\nonumber
&&\qquad \propleq\frac{{\varepsilon}^{\alpha
^{N_0}}(s-s_i)^\alpha
+(s-s_i)^{\xi}}{{\varepsilon}}
\\
&&\quad\qquad{}\times \E^{\P} \bigl[
\bigl[Y^j_s(\1)\bigr]^{1/2};s\leq
T^{Y^j}_1,
\llvert y_j-x_i\rrvert< 2 \bigl({
\varepsilon}^{1/2}+(t_j-s_i)^{\beta
'}\bigr),t_j<T_0^{X^i} \bigr],\hspace*{-5pt}\nonumber
\end{eqnarray}
where in the last inequality we use the component $\tau^{i,(3)}$ of
$\tau^i$ in Proposition~\ref{propseptime} and discard the event
$ \{X^i_s(\1)^{T_1^{X^i}}>0 \}$.
Applying (\ref{impilp}) and (\ref{dom}) to (\ref{ineqYjhalf2-0}) gives
%
%
%e3.60 #&#
\begin{eqnarray}\label{ineqYjhalf2-1}
&&\frac{1}{\psi(\1){\varepsilon}}\E^{\P} \bigl[\bigl\llvert X^i_s(
\1)^{T^{X^i}_1}-\psi(\1){\varepsilon}\bigr\rrvert\bigl[Y^j_s(
\1) \bigr]^{1/2};s< \tau^i,
\nonumber
\\
&&\hspace*{48pt} \llvert y_j-x_i\rrvert\leq2 \bigl({
\varepsilon}^{1/2}+(t_j-s_i)^{\beta'}
\bigr), X^i_s(\1)^{T_1^{X^i}}>0, t_j<T_0^{X^i}
\bigr]
\nonumber\\[-8pt]\\[-8pt]\nonumber
&&\qquad \propleq\frac{{\varepsilon}^{\alpha^{N_0}}(s-s_i)^\alpha
+(s-s_i)^{\xi}}{{\varepsilon}}\cdot\bigl({\varepsilon
}^{1/2}+(t_j-s_i)^{\beta'} \bigr) \P
\bigl(t_j<T_0^{X^i} \bigr)
\\
&&\quad\qquad{}\times \E^{\mathbf P^0_{\psi(\1){\varepsilon}}} \bigl[(Z_{s-t_j})^{1/2};s-t_j
\leq T^Z_1 \bigr].\nonumber
\end{eqnarray}
We have
%
%
%e3.61 #&#
\begin{equation}
\label{ineqsurvprobp} \P\bigl(t_j<T^{X^i}_0 \bigr)
\leq\frac{2\psi(\1){\varepsilon}}{t_j-s_i}
\end{equation}
by (\ref{ineqsurvprob}). Applying the last display and (\ref
{ineqtjvep}) to the right-hand side of (\ref{ineqYjhalf2-1}) then
gives the desired inequality (\ref{ineqYjhalf2}).
\end{longlist}

\begin{longlist}[\textit{Step} 2-4.]
\item[\textit{Step} 2-4.]
For the third term in (\ref
{ineqsepINT1-2}), the arguments step 2-3 [cf. (\ref
{ineqYjhalf2-0}) and~(\ref{ineqYjhalf2-1})] readily give
%
%
%e3.62 #&#
%e3.63 #&#
\begin{eqnarray}\label{ineqYjhalf4}
&&\frac{1}{\psi(\1){\varepsilon}}\cdot\psi(\1){\varepsilon}\E^{\P
} \bigl[ \bigl[Y^j_s(\1) \bigr]^{1/2};s\leq
T^{Y_j}_1,\nonumber
\\
&&\hspace*{84pt}
\llvert y_j-x_i\rrvert
\leq2 \bigl({\varepsilon}^{1/2}+(t_j-s_i)^{\beta'}
\bigr),t_j<T_0^{X^i} \bigr]
\nonumber\\[-8pt]\\[-8pt]\nonumber
&&\qquad \propleq (t_j-s_i)^{\beta'-1} {
\varepsilon}\E^{\mathbf P^0_{\psi(\1
){\varepsilon}}} \bigl[(Z_{s-t_j})^{1/2};s-t_j
\leq T^Z_1 \bigr]
\\
\eqntext{\forall s\in(t_j,t],
s_i<t_j<t.}
\end{eqnarray}
\end{longlist}

\begin{longlist}[\textit{Step} 2-5.]
\item[\textit{Step} 2-5.]
We note that in (\ref{ineqYjhalf1}), (\ref{ineqYjhalf2}) and (\ref
{ineqYjhalf4}), there is a common fractional moment, or more precisely
%
%
%e3.64 #&#
\begin{equation}
\label{disfracMom} \E^{\mathbf P^0_{\psi(\1){\varepsilon}}} \bigl[
(Z_{s-t_j} )^{1/2};s-t_j
\leq T^Z_1 \bigr],
\end{equation}
left to be estimated, as will be done in this step.

Recall the filtration $(\ms H_t)$ defined in Section~\ref{secpre}.

%
%
%le3.11 #&#
\begin{lem}\label{lemFd1}
Fix $z,T\in(0,\infty)$. Under the conditional probability measure
$\mathbf P^{(T)}_z$ defined by
%
%
%e3.65 #&#
\begin{equation}
\label{defQT} \mathbf P_z^{(T)}(A)\triangleq\mathbf
P_z^0(A \mid Z_T>0),\qquad A\in\ms
H_T,
\end{equation}
the process $(Z_t)_{0\leq t\leq T}$ is a continuous $(\ms
H_t)$-semimartingale with canonical decomposition
%
%
%e3.66 #&#
\begin{equation}
\label{eqQTFV} Z_t=z+\int_0^t F
\biggl(\frac{2Z_s}{T-s} \biggr)\,ds+M_t,\qquad0\leq t\leq T.
\end{equation}
Here, $F\dvtx\R_+\lra\R_+$ defined by
%
%
%e3.67 #&#
\begin{equation}
\label{defF} F(x)\triangleq\cases{ \displaystyle\frac
{e^{-x}x}{1-e^{-x}},&\quad
$x>0$,
\vspace*{5pt}\cr
1,&\quad $x=0$,}
\end{equation}
is continuous and decreasing, and $M$ is a continuous $(\ms
H_t)$-martingale under $\mathbf P^{(T)}_z$ with quadratic variation
$\langle M\rangle_t\equiv\int_0^t Z_s\,ds$.
\end{lem}

\begin{pf}
The proof of this lemma is a standard application of Girsanov's theorem
(cf. Theorem VIII.1.4 of \cite{RYCMB}), and we proceed as follows.

First, let $(D_t)_{0\leq t\leq T}$ denote the $(\ms H_t,\mathbf
P^0_z)$-martingale associated with the Radon--Nikodym derivative of
$\mathbf P^{(T)}_z$ with respect to $\mathbf P^0_z$, that is,
%
%
%e3.68 #&#
\begin{equation}
\label{defD} D_t\equiv\frac{\mathbf P^0_z (Z_T>0\mid\ms H_t )}{\mathbf
P^0_z(Z_T>0)},\qquad0\leq t\leq T.
\end{equation}
To obtain the explicit form of $D$ under $\mathbf P^0_z$, we first note that
the $(\ms H_t,\mathbf P^0_z)$-Markov property of $Z$ and (\ref
{eqsurvprob}) imply
%
%
%e3.69 #&#
%e3.70 #&#
\begin{eqnarray}\label{eqMmg}
\mathbf P^0_z (Z_T>0\mid\ms
H_t )&=&\mathbf P^0_{Z_t}(Z_{T-t}>0)=1-
\exp\biggl(-\frac{2Z_t}{T-t} \biggr),
\nonumber\\[-12pt]\\[-8pt]
\eqntext{0\leq t<T.}
\end{eqnarray}
Hence, it follows from It\^{o}'s formula and the foregoing display
that, under $\mathbf P^0_z$,
%
%
%e3.71 #&#
%e3.72 #&#
\begin{eqnarray}
\label{eqMsemiMG} D_t&=&\frac{1}{\mathbf P^0_z(Z_T>0)} \biggl[1-\exp
\biggl(-\frac
{2z}{T} \biggr) \biggr]\nonumber
\\
&&{} +\frac{1}{\mathbf P^0_z(Z_T>0)}\int_0^t \exp\biggl(-
\frac{2Z_s}{T-s} \biggr)\cdot\biggl(\frac{2}{T-s} \biggr)\sqrt
{Z_s}\,dB_s,
\\
\eqntext{0\leq t<T.}
\end{eqnarray}

We now apply Girsanov's theorem and verify that the components of the
canonical decomposition of $(Z_t)_{0\leq t\leq T}$ under $\mathbf
P^{(T)}_z$ satisfy the asserted properties.
Under $\mathbf P^{(T)}_z$, we have
\[
Z_t=z+\int_0^t
D_s^{-1}\,d\langle D,Z\rangle_s+M_t,
\qquad0\leq t\leq T.
\]
Here,
\[
M_t=\int_0^t
\sqrt{Z_s}\,dB_s-\int_0^t
D_s^{-1}\,d\langle D,Z\rangle_s,\qquad0\leq
t\leq T
\]
is a continuous $(\ms H_t,\mathbf P_z^{(T)})$-local martingale with the
asserted quadratic variation
$\langle M_t\rangle\equiv\int_0^t Z_s\,ds$, which implies that $M$ is
a true martingale under $\mathbf P^{(T)}_z$.
In addition, it follows from (\ref{eqMmg}) and (\ref{eqMsemiMG})
that the finite variation process of $Z$ under $\mathbf P^{(T)}_z$ is
given by
\begin{eqnarray}
\int_0^t D_s^{-1}\,d
\langle D,Z\rangle_s&=&\int_0^t
\frac{1}{\mathbf
P^0_z(Z_T>0\mid\ms H_s)}\,d\bigl\langle\mathbf P^0_z(Z_T>0)D,Z
\bigr\rangle_s
\nonumber
\\
&=&\int_0^t \frac{\exp(-{2Z_s}/{(T-s)} ){2Z_s}/{(T-s)}}{1-\exp(-{2Z_s}/{(T-s)} )}\,ds
\nonumber
\\
&=&\int_0^t F \biggl(\frac{2Z_s}{T-s}
\biggr)\,ds,\qquad0\leq t\leq T,
\nonumber
\end{eqnarray}
where $F$ is given by (\ref{defF}).
The proof is complete.
\end{pf}

%
%
%le3.12 #&#
\begin{lem}\label{lemFd2}
For any $p\in(0,\infty)$, there exists a constant $K_p\in(0,\infty
)$ depending only on $p$ and $(\alpha,\xi,N_0)$ such that
%
%
%e3.73 #&#
\begin{eqnarray}
\label{ineqFd2} %
\qquad&& \E^{\mathbf P^0_z} \bigl[(Z_T)^{p};T
\leq T^Z_1 \bigr]
\nonumber\\[-8pt]\\[-8pt]\nonumber
&&\qquad \leq K_p \bigl[
\bigl(z^{p\alpha^{N_0}}T^{p\alpha}+z^p \bigr)\mathbf
P^0_z(Z_T>0) +z T^{p\xi-1} \bigr],
\qquad \forall z,T\in(0,1].
\end{eqnarray}
\end{lem}

\begin{pf}
Recall the conditional probability measure $\mathbf P^{(T)}_z$ defined
in (\ref{defQT}) and write
%
%
%e3.74 #&#
\begin{eqnarray}\label{ineqZest}
\E^{\mathbf P^0_z} \bigl[(Z_T)^p;T\leq
T^Z_1 \bigr]&\leq& \mathbf P^0_z(Z_T>0)
\E^{\mathbf P^0_z} \bigl[ (Z_{T\wedge T_1^Z} )^{p} \mid
Z_T>0 \bigr]
\nonumber\\[-8pt]\\[-8pt]\nonumber
&=& \mathbf P^0_z(Z_T>0)\E^{\mathbf P_z^{(T)}}
\bigl[ (Z_{T\wedge
T_1^Z} )^{p} \bigr].
\end{eqnarray}
Henceforth, we work under the conditional probability measure $\mathbf
P^{(T)}_z$.

We turn to the improved modulus of continuity of $Z$ at its starting
time $0$ under $\mathbf P^{(T)}_z$ in order to bound the right-hand
side of (\ref{ineqZest}). We first claim that, by enlarging the
underlying probability space if necessary,
%
%
%e3.75 #&#
\begin{equation}
\label{ineqmodZ} \qquad \llvert Z_t-z\rrvert\leq t+C^Z_\alpha
\biggl(\int_0^t Z_s\,ds
\biggr)^\alpha\qquad\forall t\in\bigl[0, T\wedge T^Z_1
\bigr] \mbox{ under }\mathbf P^{(T)}_z,
\end{equation}
where the random variable $C^Z_\alpha$ under $\mathbf P^{(T)}_z$ has
distribution depending only on~$\alpha$ and finite $\mathbf
P^{(T)}_z$-moment of any finite order.
We show how to obtain (\ref{ineqmodZ}) by using the canonical
decomposition of the continuous $(\ms H_t,\mathbf
P^{(T)}_z)$-semimartingale $(Z_t)_{0\leq t\leq T}$ in (\ref{eqQTFV}).
First, since its martingale part $M$ has quadratic variation $\int
_0^\cdot Z_s\,ds$, the Dambis--Dubins--Schwarz theorem (cf. Theorem V.1.6
of \cite{RYCMB}) implies that, by enlarging of the underlying
probability space if necessary,
\[
M_{t}=\widetilde{B} \biggl(\int_0^{t}Z_s\,ds
\biggr),\qquad t\in\bigl[0,T\wedge T_1^Z \bigr],
\]
for some standard Brownian motion $\widetilde{B}$ under $\mathbf
P^{(T)}_z$. Here, the random clock
$\int_0^{t}Z_s\,ds$, $t\in[0,T\wedge T^Z_1]$, for $\widetilde{B}$ is
bounded by $1$ by the assumption that $z,T\leq1$. On the other hand,
recall that the chosen parameter $\alpha$ lies in $(0,\frac{1}{2})$
and the uniform H\"older-$\alpha$ modulus of continuity of standard
Brownian motion on compacts has moments of any finite order. (See,
e.g., the discussion preceding Theorem I.2.2 of~\cite{RYCMB} and its
proof.) Hence,
\[
\biggl\llvert\widetilde{B} \biggl(\int_0^{t}Z_s\,ds
\biggr)\biggr\rrvert\leq C^Z_\alpha\biggl(\int
_0^{t}Z_s\,ds \biggr)^\alpha,
\qquad t\in\bigl[0, T\wedge T^Z_1\bigr],
\]
where the random variable $C^Z_\alpha$ is as in (\ref{ineqmodZ}).
Second, Lemma~\ref{lemFd1} also states that the finite variation
process of $Z$ under $\mathbf P^{(T)}_z$
given by (\ref{eqQTFV}) is a time integral with integrand uniformly
bounded by $1$. This and the last two displays are enough to obtain our
claim (\ref{ineqmodZ}).

With the integral inequality (\ref{ineqmodZ}) and the distributional
properties of $C^Z_\alpha$, we obtain the following improved modulus
of continuity of $Z$ (cf. Corollary~\ref{corimc}):
%
%
%e3.76 #&#
\begin{equation}
\label{ineqZmod} \llvert Z_{T\wedge T_1^Z}-z\rrvert\leq K_1^Zz^{\alpha
^{N_0}}T^{\alpha
}+K_2^ZT^\xi
\end{equation}
for some random variables $K_1^Z,K_2^Z\in\bigcap_{q\in(0,\infty
)}L^q (\mathbf P^{(T)}_z )$ obeying a joint law under $\mathbf
P^{(T)}_z$ depending only on $(\alpha,\xi,N_0)$ by the analogous
property of $C^Z_\alpha$ and Corollary~\ref{corimc}.

We return to the calculation in (\ref{ineqZest}). Applying (\ref
{ineqZmod}), we get
%
%
%e3.77 #&#
\begin{eqnarray}\label{Zmod1}
&& \E^{\mathbf P^0_z} \bigl[(Z_T)^{p};T\leq
T_1^Z \bigr]\nonumber
\\
&&\qquad \leq\mathbf P_z^0(Z_T>0)
\E^{\mathbf P^{(T)}_z} \bigl[ (Z_{T\wedge T^Z_1} )^{p} \bigr]
\nonumber\\[-8pt]\\[-8pt]\nonumber
&&\qquad \leq\mathbf P^0_z(Z_T>0)
\bigl(2^{p-1}\vee1 \bigr)\E^{\mathbf
P^{(T)}_z} \bigl[\llvert
Z_{T\wedge T^Z_1}-z\rrvert^{p}+z^{p} \bigr]
\\
&&\qquad \leq\mathbf P^0_z(Z_T>0)
K_p' \bigl(z^{p\alpha^{N_0}}T^{p\alpha
}+T^{p\xi}+z^p
\bigr)\nonumber
\end{eqnarray}
for some constant $K'_p$ depending only on $p$ and $(\alpha,\xi,N_0)$
by (\ref{ineqZmod}) and the distributional properties of $K^Z_j$,
where the second inequality follows from the\vadjust{\goodbreak} elementary inequality
\[
(x+y)^p\leq\bigl(2^{p-1}\vee1 \bigr)\cdot
\bigl(x^p+y^p \bigr)\qquad\forall x,y\in\R_+.
\]
The desired result follows by applying (\ref{ineqsurvprob}) to (\ref
{Zmod1}). The proof is complete.
\end{pf}
\end{longlist}

\begin{longlist}[\textit{Step} 2-6.]
\item[\textit{Step} 2-6.]
At this step, we summarize our results in steps 2-1--2-4, using
Lemma~\ref{lemFd2}.
We apply (\ref{ineqYjhalf1}), (\ref{ineqYjhalf2}) and (\ref
{ineqYjhalf4}) to (\ref{ineqsepINT1-2}). This gives
%
%
%e3.78 #&#
%e3.79 #&#
\begin{eqnarray}\label{ineqYjhalf5}
\quad && \E^{\Q^i} \bigl[ \bigl[Y^j_s(\1)
\bigr]^{1/2};s< \tau^i,\llvert y_j-x_i
\rrvert<2 \bigl({\varepsilon}^{1/2}+(t_j-s_i)^{\beta'}
\bigr) \bigr]
\nonumber
\\
&&\qquad \propleq \bigl[(t_j-s_i)^{\beta'}+(t_j-s_i)^{\beta'-1}
\bigl({\varepsilon}^{\alpha^{N_0}}(s-s_i)^\alpha+(s-s_i)^\xi
+{\varepsilon} \bigr) \bigr]
\nonumber
\\
&&\qquad\quad{} \times\E^{\mathbf P^0_{\psi(\1){\varepsilon}}} \bigl[(Z_{s-t_j})^{1/2};s-t_j
\leq T^Z_1 \bigr]
\nonumber
\\
&&\qquad \propleq \bigl[(t_j-s_i)^{\beta'}+(t_j-s_i)^{\beta'-1}
\bigl({\varepsilon}^{\alpha^{N_0}}(s-s_i)^\alpha+(s-s_i)^\xi
+{\varepsilon} \bigr) \bigr]
\nonumber
\\
&&\quad\qquad{} \times\bigl[ \bigl({\varepsilon}^{{\alpha ^{N_0}}/{2}}(s-t_j)^{{\alpha }/{2}}+{
\varepsilon}^{{1}/{2}} \bigr)\mathbf P^0_{\psi(\1){\varepsilon
}}(Z_{s-t_j}>0)+{
\varepsilon}(s-t_j)^{{\xi}/{2}-1} \bigr]
\nonumber
\\
 &&\qquad \propleq (t_j-s_i)^{\beta'-1}
\times\bigl((t_j-s_i)+{\varepsilon}^{\alpha^{N_0}}(s-s_i)^\alpha+(s-s_i)^\xi+{\varepsilon} \bigr)\nonumber
\\
&&\quad\qquad{}\times{\varepsilon}^{{\alpha^{N_0}}/{2}}(s-t_j)^{{\alpha }/{2}}\mathbf
P^0_{\psi(\1){\varepsilon}}(Z_{s-t_j}>0)
\\
&&\quad\qquad{}+(t_j-s_i)^{\beta'-1}\times\bigl((t_j-s_i)+{
\varepsilon}^{\alpha
^{N_0}}(s-s_i)^\alpha+(s-s_i)^\xi+{
\varepsilon} \bigr)\nonumber
\\
&&\hspace*{2pt}\phantom{{}+}\quad\qquad{}\times{\varepsilon}^{{1}/{2}}\mathbf P^0_{\psi(\1){\varepsilon
}}(Z_{s-t_j}>0)\nonumber
\\
&&\quad\qquad{}+(t_j-s_i)^{\beta'}\times{\varepsilon}(s-t_j)^{{\xi}/{2}-1}\nonumber
\\
&&\quad\qquad{}+(t_j-s_i)^{\beta'-1}\times\bigl({
\varepsilon}^{\alpha
^{N_0}}(s-s_i)^\alpha+{\varepsilon}
\bigr)\times{\varepsilon}(s-t_j)^{{\xi}/{2}-1}\nonumber
\\
&&\quad\qquad{}+(t_j-s_i)^{\beta'-1}\times(s-s_i)^\xi
\times{\varepsilon}(s-t_j)^{{\xi}/{2}-1}\nonumber
\\
\eqntext{\forall s \in(t_j,t], s_i<t_j<t,}
\end{eqnarray}
where the last $\propleq$-inequality follows by some algebra.

We make some simplifications for the right-hand side of (\ref
{ineqYjhalf5}) before going further. Some orders in ${\varepsilon}$
and other variables will be discarded here.
We bound the survival probability in (\ref{ineqYjhalf5}) by
%
%
%e3.80 #&#
\begin{equation}
\label{ineqsp-pow} \mathbf P^0_{\psi(\1){\varepsilon}}(Z_{s-t_j}>0)\leq
\biggl(\frac
{2\psi(\1){\varepsilon}}{s-t_j} \biggr)^{1-{\alpha^{N_0}}/{4}},
\end{equation}
as follows from the elementary inequalities $x\leq x^\gamma$ for any
$x\in[0,1]$ and $\gamma\in(0,1]$, and then (\ref{ineqsurvprob}).
Assuming $s\in(t_j,t]$ for $s_i<t_j<t$, we have
the inequalities
\[
1\geq s-s_i\geq t_j-s_i\geq
\frac{{\varepsilon}}{2}, \qquad s-t_j\leq1\quad\mbox{and}\quad0<\alpha+
\alpha^{N_0}< \xi<1
\]
[cf. (\ref{defkappa})(a) for the third inequality]. These and (\ref
{ineqsp-pow}) imply that the first term
of~(\ref{ineqYjhalf5}) satisfies
%
%
%e3.81 #&#
\begin{eqnarray}
\label{summary1} && (t_j-s_i)^{\beta'-1}
\times\bigl((t_j-s_i)+{\varepsilon}^{\alpha
^{N_0}}(s-s_i)^\alpha+(s-s_i)^\xi+{
\varepsilon} \bigr)\nonumber
\\
&&\quad{}\times{\varepsilon}^{{\alpha^{N_0}}/{2}}(s-t_j)^{{\alpha }/{2}}\mathbf
P^0_{\psi(\1){\varepsilon}}(Z_{s-t_j}>0)
\\
&&\qquad \propleq(t_j-s_i)^{\beta'-1}(s-s_i)^\alpha(s-t_j)^{\alpha /2+{\alpha^{N_0}}/{4}-1}{
\varepsilon}^{1+{\alpha^{N_0}}/{4}},\nonumber
\end{eqnarray}
the second term of (\ref{ineqYjhalf5}) satisfies
%
%
%e3.82 #&#
\begin{eqnarray}
\label{summary2} %
&& (t_j-s_i)^{\beta'-1}
\times\bigl((t_j-s_i)+{\varepsilon}^{\alpha
^{N_0}}(s-s_i)^\alpha+(s-s_i)^\xi+{
\varepsilon} \bigr)\nonumber
\\
&&\quad{} \times{\varepsilon}^{{1}/{2}}\mathbf P^0_{\psi(\1){\varepsilon
}}(Z_{s-t_j}>0)\nonumber
\\
&&\qquad \propleq(t_j-s_i)^{\beta'-1}(s-s_i)^{\alpha}(s-t_j)^{{\alpha ^{N_0}}/{4}-1}{
\varepsilon}^{{3}/{2}-{\alpha^{N_0}}/{4}}
\\
&&\qquad \propleq(t_j-s_i)^{\beta'+ ({1}/{2}-{\alpha ^{N_0}}/{2} )-1}(s-s_i)^\alpha(s-t_j)^{{\alpha ^{N_0}}/{4}-1}{
\varepsilon}^{1+{\alpha^{N_0}}/{4}}\nonumber
\\
&&\qquad \propleq(t_j-s_i)^{\beta'+{\alpha}/{2}-1}(s-s_i)^\alpha
(s-t_j)^{{\alpha^{N_0}}/{4}-1}{\varepsilon}^{1+{\alpha^{N_0}}/{4}},\nonumber
\end{eqnarray}
and, finally, the fourth term of (\ref{ineqYjhalf5}) satisfies
%
%
%e3.83 #&#
\begin{eqnarray}
\label{summary3} &&(t_j-s_i)^{\beta'-1}
\times\bigl({\varepsilon}^{\alpha
^{N_0}}(s-s_i)^\alpha+{
\varepsilon} \bigr)\times{\varepsilon}(s-t_j)^{{\xi}/{2}-1}\nonumber
\\
&&\qquad \propleq (t_j-s_i)^{\beta'-1}(s-s_i)^\alpha(s-t_j)^{{\xi }/{2}-1}{
\varepsilon}^{1+\alpha^{N_0}}
\\
&&\qquad \propleq(t_j-s_i)^{\beta'-1}(s-s_i)^\alpha(s-t_j)^{\alpha /2+{\alpha^{N_0}}/{4}-1}{
\varepsilon}^{1+{\alpha^{N_0}}/{4}}.\nonumber
\end{eqnarray}
Note that the bounds in (\ref{summary1}) and (\ref{summary3}) coincide.
Using (\ref{summary1})--(\ref{summary3}) in (\ref{ineqYjhalf5}), we obtain
%
%
%e3.84 #&#
\begin{eqnarray}\label{ineqYjhalf6}
&& \E^{\Q^i} \bigl[ \bigl[Y^j_s(\1)
\bigr]^{1/2};s< \tau^i,\llvert y_j-x_i
\rrvert\leq2 \bigl({\varepsilon}^{1/2}+(t_j-s_i)^{\beta'}
\bigr) \bigr]
\nonumber
\\
&&\qquad \propleq (t_j-s_i)^{\beta'-1}(s-s_i)^{\alpha
}(s-t_j)^{\alpha /2+{\alpha^{N_0}}/{4}-1}{
\varepsilon}^{1+{\alpha ^{N_0}}/{4}}\nonumber
\\
&&\quad\qquad{}+(t_j-s_i)^{\beta'+{\alpha}/{2}-1}(s-s_i)^\alpha(s-t_j)^{{\alpha ^{N_0}}/{4}-1}{
\varepsilon}^{1+{\alpha^{N_0}}/{4}}
\\
&&\quad\qquad{}+(t_j-s_i)^{\beta'}(s-t_j)^{{\xi}/{2}-1}{
\varepsilon}\nonumber
\\
&&\quad\qquad{}+(t_j-s_i)^{\beta'-1}(s-s_i)^\xi(s-t_j)^{{\xi }/{2}-1}{
\varepsilon}
\qquad \forall s\in(t_j,t], s_i<t_j<t.\nonumber
\end{eqnarray}
\end{longlist}

\begin{longlist}[\textit{Step} 3.]
\item[\textit{Step} 3.]
We digress to a \emph{conceptual} discussion for some elementary
integrals which will play an important role in the forthcoming
calculations in step~4.
First, for $a,b,c\in\R$ and $T\in(0,\infty)$, a straightforward
application of Fubini's theorem and changes of variables shows that
%
%
%e3.85 #&#
\begin{eqnarray}
\label{integralabc} %
&&I(a,b,c)_T\triangleq\int
_0^T dr\,r^{a}\int
_r^T ds\,s^{b}(s-r)^{c}<
\infty
\nonumber\\[-8pt]\\[-8pt]\nonumber
&&\quad\Longleftrightarrow\quad a,c\in(-1,\infty)\quad\mbox{and}\quad
a+b+c>-2.
\end{eqnarray}
Furthermore, when $I(a,b,c)_T$ is finite, it can be expressed as
\[
I(a,b,c)_T= \biggl(\int_0^1
dr\,r^{a}(1-r)^{c} \biggr)\cdot\frac
{T^{a+b+c+2}}{a+b+c+2}.
\]
Given $a+b+c>-2$ with $a,c\in(-1,\infty)$,
we consider alternative ways to show that the integral $I(a,b,c)_T$ is
finite while preserving the same order $T^{a+b+c+2}$ in~$T$, according
to $b\geq0$ and $b<0$. If $b\geq0$, then
%
%
%e3.86 #&#
\begin{eqnarray}\label{ineqredist0}
I(a,b,c)_T&\leq&\int_0^T dr\,
r^a\times T^b\times\int_0^T
ds\,s^c
\nonumber\\[-7pt]\\[-7pt]\nonumber
&=& \frac
{1}{a+1}\frac{1}{c+1}T^{a+b+c+2},
\end{eqnarray}
where the first inequality follows since $s^b\leq T^b$ for any $s\in
[r,T]$. For the case that $b<0$, we decompose the function $s\lmt s^b$
in the following way. For $b_1,b_2<0$ such that $b_1+b_2=b$, we have
%
%
%e3.87 #&#
\begin{eqnarray}
\label{ineqredist} I(a,b,c)_T&\leq&\int_0^T
dr\,r^{a+b_1}\int_r^T ds (s-r)^{b_2+c}
\nonumber\\[-7pt]\\[-7pt]\nonumber
&\leq& \int_0^T dr\,r^{a+b_1}\times\int_0^T ds\,
s^{b_2+c},
\end{eqnarray}
where the first inequality follows since for $s> r$, $s^{b_1}\leq
r^{b_1}$ and $s^{b_2}\leq(s-r)^{b_2}$. Using the following elementary
lemma, we obtain from
(\ref{ineqredist}) that
\begin{eqnarray*}
I(a,b,c)_T&\leq&\frac{1}{a+b_1+1}\frac{1}{b_2+c+1}T^{a+b+c+2}.
\end{eqnarray*}

%
%
%le3.13 #&#
\begin{lem}\label{lemalgalloc}
For any reals $a,c>-1$ and $b<0$ such that $a+b+c>-2$, there exists a
pair $(b_1,b_2)\in(-\infty,0)\times(-\infty,0)$ such that
$b=b_1+b_2$ and $a+b_1>-1$ and $b_2+c>-1$.
\end{lem}

The two simple \emph{concepts} for the inequalities (\ref
{ineqredist0}) and (\ref{ineqredist}) will be applied later on in
step~4 to bound \emph{Riemann sums} by integrals of the type
$I(a,b,c)_T$.
\end{longlist}

\begin{longlist}[\textit{Step} 4.]
\item[\textit{Step} 4.]
We complete the proof of Lemma~\ref
{lemgr1-1} in this step. Apply the bound~(\ref{ineqYjhalf6}) to the
right-hand side of
the inequality (\ref{ineqsumpre}). We have
%
%
%e3.88 #&#
\begin{eqnarray}\label{ineqYjhalf7}
&&\E^{\Q^i} \biggl[ \sum_{j\in\mathcal C^i_{\beta'}(t\wedge\tau
^i\wedge\sigma^{X^i}_\beta)} \biggl(\psi(
\1){\varepsilon}\nonumber
\\
&&\hspace*{91pt}{} +\int_{t_j}^{t\wedge\tau^i\wedge\sigma^{X^i}_\beta\wedge
\sigma_\beta^{Y^j}} \frac{1}{X^i_s(\1)}
\int_\R X^i(x,s)^{1/2}Y^j(x,s)^{1/2}\,dx\,ds
\biggr) \biggr]
\\
&&\qquad \propleq\sum_{j\dvtx s_i<t_j\leq t}(t_j-s_i)^{\beta'} {\varepsilon}\nonumber
\\
&&\quad\qquad{}+\sum_{j\dvtx s_i<t_j\leq t}(t_j-s_i)^{\beta'-1}
\int_{t_j}^t (s-s_i)^{-{\eta}/{2}+\alpha}(s-t_j)^{\alpha /2+{\alpha^{N_0}}/{4}-1}
\,ds\nonumber
\\
&&\qquad\qquad{}\times {\varepsilon}^{1+{\alpha^{N_0}}/{4}}\nonumber
\\
&&\quad\qquad{}+\sum_{j\dvtx s_i<t_j\leq t}(t_j-s_i)^{\beta'+{\alpha}/{2}-1}
\int_{t_j}^t (s-s_i)^{-{\eta}/{2}+\alpha}(s-t_j)^{{\alpha ^{N_0}}/{4}-1}\,ds\nonumber
\\
&&\qquad\qquad{}\times  {\varepsilon}^{1+{\alpha^{N_0}}/{4}}\nonumber
\\
&&\quad\qquad{}+\sum_{j\dvtx s_i<t_j\leq t}(t_j-s_i)^{\beta'}
\int_{t_j}^t(s-s_i)^{-{\eta}/{2}}
(s-t_j)^{{\xi}/{2}-1}\,ds\cdot{\varepsilon}\nonumber
\\
&&\quad\qquad{}+\sum_{j\dvtx s_i<t_j\leq t}(t_j-s_i)^{\beta'-1}
\int_{t_j}^t (s-s_i)^{-{\eta}/{2}+\xi}
(s-t_j)^{{\xi}/{2}-1}\,ds\cdot{\varepsilon}.\nonumber
\end{eqnarray}
Recall the notation $I(a,b,c)$ in (\ref{integralabc}).
It should be clear that, up to a translation of time by $s_i$, the
first, the fourth, and the fifth sums are Riemann sums of
\begin{eqnarray*}
&&I\bigl(\beta',0,0\bigr)_{t-s_i},\qquad I \biggl(
\beta',-\frac{\eta}{2},\frac{\xi}{2}-1
\biggr)_{t-s_i},\qquad I \biggl(\beta'-1,-
\frac{\eta}{2}+\xi,\frac{\xi}{2}-1 \biggr)_{t-s_i},
\end{eqnarray*}
respectively, and so are the second and the third sums after a division
by ${\varepsilon}^{{\alpha^{N_0}}/{4}}$ with the corresponding
integrals equal to
\begin{eqnarray*}
&\displaystyle I \biggl(\beta'-1,-\frac{\eta}{2}+\alpha,\frac{\alpha}{2}+
\frac{\alpha^{N_0}}{4}-1 \biggr)_{t-s_i},&
\\
&\displaystyle I \biggl(\beta'+
\frac{\alpha}{2}-1,-\frac{\eta}{2}+\alpha,\frac
{\alpha^{N_0}}{4}-1
\biggr)_{t-s_i},&
\end{eqnarray*}
respectively.
It follows from (\ref{defkappa})(c) and (d) and (\ref{integralabc})
that all of the integrals in the last two displays are finite.

We now aim to bound each of the five sums in (\ref{ineqYjhalf7}) by
suitable powers of ${\varepsilon}$ and $t$, using integral
comparisons. Observe that,
whenever $\gamma\in(-1,\infty)$, the monotonicity of $r\lmt
(r-s_i)^\gamma$ over $(s_i,\infty)$ implies
%
%
%e3.89 #&#
\begin{eqnarray}\label{ineqintcomp}
\sum_{j\dvtx s_i<t_j\leq t}(t_j-s_i)^\gamma
\cdot{\varepsilon}&\leq&2\int_{s_i}^{t+{\varepsilon}}(r-s_i)^{\gamma}\,dr\nonumber
\\
&=& \frac{2}{\gamma+1}(t+{\varepsilon}-s_i)^{\gamma+1}
\\
&\leq&\frac{2\cdot3^{\gamma+1}}{\gamma+1}(t-s_i)^{\gamma+1}\nonumber
\end{eqnarray}
since $t\geq s_i+\frac{{\varepsilon}}{2}$.
(The constant $2$ is used to accommodate the case that $\gamma<0$.)
Hence, the first sum in (\ref{ineqYjhalf7}) can be bounded as
%
%
%e3.90 #&#
\begin{equation}
\label{Yexpbdd1} \sum_{j\dvtx s_i<t_j\leq t}(t_j-s_i)^{\beta'}
{\varepsilon}\propleq(t-s_i)^{\beta'+1}.
\end{equation}

Consider the other sums in (\ref{ineqYjhalf7}). Recall our discussion
of some alternative ways to bound $I(a,b,c)$ for given $a+b+c>-2$ and
$a,c,\in(-1,\infty)$ according to $b\geq0$ or $b<0$; see (\ref
{ineqredist0}) and (\ref{ineqredist}). We use Lemma~\ref
{lemalgalloc} in the following whenever necessary.
The second sum in (\ref{ineqYjhalf7}) can be bounded as
%
%
%e3.91 #&#
\begin{eqnarray}\label{Yexpbdd2}
&&\sum_{j\dvtx s_i<t_j\leq t}(t_j-s_i)^{\beta'-1}
\int_{t_j}^t (s-s_i)^{-{\eta}/{2}+\alpha}
(s-t_j)^{\alpha /2+{\alpha^{N_0}}/{4}-1}\,ds\cdot{\varepsilon}^{1+{\alpha ^{N_0}}/{4}}
\nonumber
\\
&&\qquad  =\sum_{j\dvtx s_i<t_j\leq t}(t_j-s_i)^{\beta'-1}
\int_{t_j-s_i}^{t-s_i} s^{-{\eta}/{2}+\alpha}
\bigl[s-(t_j-s_i)\bigr]^{{\alpha}/{2}+{\alpha^{N_0}}/{4}-1}\,ds\nonumber
\\
&&\quad\qquad{}\times {\varepsilon}^{1+{\alpha^{N_0}}/{4}}
\\
&&\qquad  \propleq(t-s_i)^{\beta'-{\eta}/{2}+({3\alpha})/{2}+{\alpha^{N_0}}/{4}} \cdot{
\varepsilon}^{{\alpha ^{N_0}}/{4}}
\nonumber
\\
&&\qquad  \propleq(t-s_i)^{\beta'-{\eta}/{2}+ {(3\alpha)}/{2}} \cdot{
\varepsilon}^{{\alpha^{N_0}}/{4}}.\nonumber
\end{eqnarray}
Here, in the foregoing $\propleq$-inequality, we use the integral
comparison discussed in step~3 (with Lemma~\ref{lemalgalloc} to
algebraically allocate the exponent $-\frac{\eta}{2}+\frac{\alpha
}{2}$ if necessary) and the Riemman-sum bound (\ref{ineqintcomp}).
The other sums on the right-hand side of (\ref{ineqYjhalf7}) can be
bounded similarly as follows.
The third sum satisfies
%
%
%e3.92 #&#
\begin{eqnarray}\label{Yexpbdd3}
&&\sum_{j\dvtx s_i<t_j\leq t}(t_j-s_i)^{\beta'+{\alpha}/{2}-1}
\int_{t_j}^t (s-s_i)^{-{\eta}/{2}+\alpha}(s-t_j)^{{\alpha ^{N_0}}/{4}-1}\,ds
\cdot{\varepsilon}^{1+{\alpha ^{N_0}}/{4}}\hspace*{-20pt}
\nonumber\\[-8pt]\\[-8pt]\nonumber
&&\qquad  \propleq(t-s_i)^{\beta'-{\eta}/{2}+ {(3\alpha)}/{2}}\cdot{
\varepsilon}^{{\alpha^{N_0}}/{4}}.
\end{eqnarray}
The fourth sum satisfies
%
%
%e3.93 #&#
\begin{eqnarray}\label{Yexpbdd4}
&& \sum_{j\dvtx s_i<t_j\leq t}(t_j-s_i)^{\beta'}
\int_{t_j}^t (s-s_i)^{-{\eta}/{2}}(s-t_j)^{{\xi}/{2}-1}\,ds
\cdot{\varepsilon}
\nonumber\\[-8pt]\\[-8pt]\nonumber
&&\qquad \propleq (t-s_i)^{\beta'-{\eta}/{2}+{\xi}/{2}+1}
\propleq (t-s_i)^{\beta'-{\eta}/{2}+{(3\xi)}/{2}},
\end{eqnarray}
where the last inequality applies since $\xi\in(0,1)$.
The last sum satisfies
%
%
%e3.94 #&#
\begin{eqnarray}\label{Yexpbdd5}
&&\sum_{j\dvtx s_i<t_j\leq t}(t_j-s_i)^{\beta'-1}
\int_{t_j}^t (s-s_i)^{-{\eta}/{2}+\xi}
(s-t_j)^{{\xi}/{2}-1}\,ds\cdot{\varepsilon}
\nonumber\\[-8pt]\\[-8pt]\nonumber
&&\qquad  \propleq(t-s_i)^{\beta'-{\eta}/{2}+{(3\xi)}/{2}}.
\end{eqnarray}

The proof of Lemma~\ref{lemgr1-1} is complete upon applying (\ref
{Yexpbdd1})--(\ref{Yexpbdd5}) to the right-hand side of (\ref{ineqYjhalf7}).

%s3.5 #&#
\subsection{Proof of Lemma~\texorpdfstring{\protect\ref{lemgr1-2}}{3.10}}\label{seccondsep-2}
As in Section~\ref{seccondsep-1}, we fix $t\in[s_i+\frac
{{\varepsilon}}{2},s_i+1 ]$, $i\in\mathbb N$, and ${\varepsilon
}\in(0,[8\psi(\1)]^{-1}\wedge1 ]$ and drop the
subscripts of $\P_{\varepsilon}$ and $\Q^i_{\varepsilon}$.
For the proof of Lemma~\ref{lemgr1-2}, the arguments in Section~\ref
{seccondsep-1} work essentially.
Now, we begin to use the condition (\ref{defkappa})(b) in
Assumption~\ref{ass} and the upper limit $\sigma^{X^i}_\beta\wedge
\sigma_\beta^{Y^j}$ in the time integral in (\ref{sumlate}), which
are neglected when we prove Lemma~\ref{lemgr1-1}.

To motivate our adaptation of the arguments for critical clusters in
Section~\ref{seccondsep-1}, we discuss
some parts of Section~\ref{seccondsep-1}.
First, it is straightforward to modify the proof of (\ref
{ineqQbdd-c}) and obtain
%
%
%e3.95 #&#
\begin{eqnarray}\label{ineqQbdd}
&& \Q^i \bigl(2 \bigl({\varepsilon}^{1/2}+(t_j-s_i)^{\beta'}
\bigr)\leq\llvert y_j-x_i\rrvert\leq2 \bigl({
\varepsilon}^{1/2}+(t-s_i)^{\beta'} \bigr) \bigr)
\propleq(t-s_i)^{\beta'}.\hspace*{-30pt}
\end{eqnarray}
If we proceed as in (\ref{ineqsumpre}) and use (\ref{ineqQbdd}) in
the obvious way,
then this leads to
\begin{eqnarray*}
\hspace*{-4pt}&&\E^{\mathbb Q^i} \biggl[ \sum_{j\in\mathcal L^i_{\beta'}(t,t\wedge
\tau^i\wedge\sigma^{X^i}_\beta)}\! \biggl(\psi(
\1){\varepsilon}
\\[-3pt]
\hspace*{-4pt}&&\hspace*{97pt} {}+\int_{t_j}^{t\wedge\tau^i\wedge\sigma^{X^i}_\beta\wedge
\sigma_\beta^{Y^j}}\! \frac{1}{X^i_s(\1)}
\int_\R X^i(x,s)^{1/2}Y^j(x,s)^{1/2}\,dx\,ds
\biggr) \biggr]
\\
\hspace*{-4pt}&&\qquad \propleq\sum_{j\dvtx s_i<t_j\leq t}(t-s_i)^{\beta'}{
\varepsilon}
\\
\hspace*{-4pt}&&\hspace*{33pt} +\! \sum_{j\dvtx s_i<t_j\leq t}\int_{t_j}^t ds\frac{1}{(s-s_i)^{\eta/2}} \E^{\mathbb Q^i} \bigl[ \bigl[Y^j_s(
\1) \bigr]^{1/2};s<\tau^i,2 \bigl({\varepsilon}^{1/2}+(t_j-s_i)^{\beta
'} \bigr)
\\
\hspace*{-4pt}&&\hspace*{205pt} \leq\llvert y_j-x_i\rrvert\leq2 \bigl({
\varepsilon}^{1/2}+(t-s_i)^{\beta
'} \bigr) \bigr].
\end{eqnarray*}
[Compare this with (\ref{ineqsumpre}) for critical clusters.]
If we argue by using (\ref{ineqQbdd}) repeatedly in the steps
analogous to steps 2-2--2-4 of Section~\ref{seccondsep-1},
then we obtain the following $\propleq$-inequality similar to (\ref
{ineqYjhalf7}):
%
%
%e3.96 #&#
\begin{eqnarray}\label{eqfalse}
\hspace*{-4pt}&& \E^{\Q^i} \biggl[ \sum_{j\in\mathcal L^i_{\beta'}(t,t\wedge\tau
^i\wedge\sigma^{X^i}_\beta)}\! \biggl(\psi(
\1){\varepsilon}\nonumber
\\[-3pt]
\hspace*{-4pt}&&\hspace*{97pt}{} +\int_{t_j}^{t\wedge\tau^i\wedge\sigma^{X^i}_\beta\wedge
\sigma_\beta^{Y^j}}\! \frac{1}{X^i_s(\1)}
\int_\R X^i(x,s)^{1/2}Y^j(x,s)^{1/2}\,dx\,ds
\biggr) \biggr]\nonumber
\\
\hspace*{-4pt}&&\qquad \propleq\sum_{j\dvtx s_i<t_j\leq t}(t-s_i)^{\beta'}
{\varepsilon}\nonumber
\\
\hspace*{-4pt}&&\quad\qquad{}+\sum_{j\dvtx s_i<t_j\leq t}(t-s_i)^{\beta'}(t_j-s_i)^{-1}
\int_{t_j}^t (s-s_i)^{-{\eta}/{2}+\alpha}(s-t_j)^{\alpha /2+{\alpha^{N_0}}/{4}-1}
\,ds\nonumber
\\
\hspace*{-4pt}&&\qquad\qquad{}\times {\varepsilon}^{1+{\alpha^{N_0}}/{4}}
\\
\hspace*{-4pt}&&\quad\qquad{}+\sum_{j\dvtx s_i<t_j\leq t}(t-s_i)^{\beta'}(t_j-s_i)^{\alpha /2-1}
\int_{t_j}^t (s-s_i)^{-{\eta}/{2}+\alpha}(s-t_j)^{{\alpha^{N_0}}/{4}-1}
\,ds\nonumber
\\
\hspace*{-4pt}&&\qquad\qquad{}\times {\varepsilon}^{1+{\alpha^{N_0}}/{4}}\nonumber
\\
\hspace*{-4pt}&&\quad\qquad{}+\sum_{j\dvtx s_i<t_j\leq t}(t-s_i)^{\beta'}\int
_{t_j}^t(s-s_i)^{-{\eta}/{2}} (s-t_j)^{{\xi}/{2}-1}\,ds\cdot{\varepsilon}\nonumber
\\
\hspace*{-4pt}&&\quad\qquad{}+\sum_{j\dvtx s_i<t_j\leq t}(t-s_i)^{\beta'}(t_j-s_i)^{-1}
\int_{t_j}^t (s-s_i)^{-{\eta}/{2}+\xi}
(s-t_j)^{{\xi}/{2}-1}\,ds\cdot{\varepsilon},\nonumber
\end{eqnarray}
taking into account some simplifications similar to (\ref
{summary1})--(\ref{summary3}) where some orders are discarded. (We
omit the derivation of the foregoing display, as it will not be used
for the proof of Lemma~\ref{lemgr1-2}.)
In other words, replacing the factor $(t_j-s_i)^{\beta'}$ for each of
the sums in (\ref{ineqYjhalf7}) by $(t-s_i)^{\beta'}$ gives the
bound in the foregoing display. Applying integral domination to the
second and the last sums of the foregoing display as in step~4 of
Section~\ref{seccondsep-2} results in bounds which are \emph
{divergent} integrals.

Examining the arguments in steps 2-2--2-4 of Section~\ref{seccondsep-1} shows that the problematic factor
%
%
%e3.97 #&#
\begin{equation}
\label{fixprob} (t_j-s_i)^{-1}
\end{equation}
in (\ref{eqfalse})
results from using the bound (\ref{ineqsurvprobp}) for the survival
probability $\P(t_j<T^{X^i}_0)$.
The exponent $-1$ in the foregoing display, however, is critical, and
any decrease in this value will lead to convergent integrals.
Also, we recall that (\ref{dom})
is used repeatedly in steps 2-2--2-4 of Section~\ref
{seccondsep-1}, while (\ref{dom}) is a consequence of (\ref
{impilp}) and the proof of (\ref{impilp})
uses in particular the Markov property of $Y^j(\1)$ at $t_j$.
These observations suggest that we should modify the arguments in
Section~\ref{seccondsep-1} by replacing
$t_j$ with a ``larger'' value, subject to the condition that certain
$\P$-independence, similar to (\ref{impilp}) with $t_j$ replaced by
the resulting value, still holds.

First, let us identify the value to replace $t_j$. The idea comes from
the following observation.

\begin{ob*}
%\subsubsection*{Observation}
It takes a positive amount of time before
the support process of a lateral cluster $Y^j$ intersects the support
process of $X^i$, as leads to a time $t_j^{\mathrm{c}}$ \emph{larger} than
the landing time $t_j$ of $Y^j$. Prior to $t_j^{\mathrm{c}}$, the supports
of $X^i$ and $Y^j$ are disjoint. See Figure~\ref{Fig1}. %\qed
\end{ob*}

We formalize the definition of this time $t_j^{\mathrm{c}}$ as follows.
Let $j\in\mathbb N$ with $t_j\in(s_i,s_i+1]$. Recall that the range for
the possible values $y$ of $y_j$ associated with a lateral cluster is
%
%
%e3.98 #&#
\begin{equation}
\label{eqlateral} 2 \bigl({\varepsilon}^{1/2}+(t_j-s_i)^{\beta'}
\bigr)\leq\llvert y-x_i\rrvert\leq2 \bigl({\varepsilon
}^{1/2}+(t-s_i)^{\beta'}
\bigr),
\end{equation}
and we use $\mathcal P^{X^i}_\beta(\cdot)$ and $\mathcal
P^{Y^j}_\beta(\cdot)$ to envelop the support processes of $X^i$ and
$Y^j$, respectively. Let the\vspace*{1pt} processes of parabolas $ \{\mathcal
P^{X^i}_\beta(t);t\in[s_i,\infty) \}$ and $ \{\mathcal
P^{Y^j}_\beta(t);t\in[t_j,\infty) \}$ evolve in the
deterministic way, and consider the \textit{support contact time}
$t_j^{\mathrm{c}}(y_j)$, that is,
the first time $t$ when $\mathcal P^{X^i}_\beta(t)$ and $\mathcal
P^{Y^j}_\beta(t)$ intersect. Here, for any $y$ satisfying (\ref
{eqlateral}), $t_j^{\mathrm{c}}(y)\in(t_j,\infty)$ solves
%
%
%e3.99 #&#
\begin{equation}
\label{deftjc} \qquad\quad \cases{ \displaystyle x_i+{\varepsilon}^{1/2}+
\bigl(t_j^{\mathrm{c}}(y)-s_i \bigr)^{\beta
}=y-{
\varepsilon}^{1/2}- \bigl(t_j^{\mathrm{c}}(y)-t_j
\bigr)^{\beta},&\quad if $y>x_i$,
\vspace*{2pt}\cr
\displaystyle
x_i-{\varepsilon}^{1/2}- \bigl(t_j^{\mathrm{c}}(y)-s_i
\bigr)^{\beta}=y+{\varepsilon}^{1/2}+ \bigl(t_j^{\mathrm{c}}(y)-t_j
\bigr)^{\beta}, &\quad if $y<x_i$.}
\end{equation}
By simple arithmetic, we see that the minimum of $t_j^{\mathrm{c}}(y)$ for
$y$ satisfying (\ref{eqlateral})
is attained at the boundary cases where $y$ satisfies $2
({\varepsilon}^{1/2}+(t_j-s_i)^{\beta'} )=\llvert y-x_i\rrvert $. Let us
consider the \emph{worst} case of the support contact time as
%
%
%e3.100 #&#
\begin{equation}
t_j^{\star}\triangleq\min\bigl\{t_j^{\mathrm{c}}(y);y
\mbox{ satisfies (\ref{eqlateral})} \bigr\}.\label{deftjstar}
\end{equation}
Recall that $\beta'<\beta$ by (\ref{defkappa})(b).

%
%
%le3.14 #&#
\begin{lem}\label{lemsct}
Let $j\in\mathbb N$ with $t_j\in(s_i,s_i+1]$.
\begin{longlist}
\item[(1)] The number $t_j^\star$ defined by (\ref{deftjstar}) satisfies
%
%
%e3.101 #&#
\begin{equation}
\label{eqtjstar} t_j^\star=s_i+A(t_j-s_i)
\cdot(t_j-s_i)^{{\beta'}/{\beta}},
\end{equation}
where $A(r)$ is the unique number in $ (r^{1-{\beta'}/{\beta}},\infty)$ solving
%
%
%e3.102 #&#
\begin{equation}
\label{eqAr} A(r)^\beta+ \bigl[A(r)-r^{1-{\beta'}/{\beta}}
\bigr]^\beta=2,\qquad r\in(0,1].
\end{equation}
\item[(2)] The function $A(\cdot)$ defined by (\ref{eqAr}) satisfies
%
%
%e3.103 #&#
\begin{equation}
\label{ineqAr} 1\leq A(r)\leq1+r^{1-{\beta'}/{\beta}}\qquad\forall
r\in(0,1].
\end{equation}
\end{longlist}
\end{lem}

\begin{pf}
Without loss of generality, we may assume that $t_j^{\star}=t_j^{\mathrm
{c}}(y)$ for $y$ satisfying
\[
x_i-y=2 \bigl({\varepsilon}^{1/2}+(t_j-s_i)^{\beta'}
\bigr).
\]
Using this particular value $y$ of $y_j$ in (\ref{deftjc}), we see
that $t_j^{\star}$ solves the equation
\begin{eqnarray*}
x_i-{\varepsilon}^{1/2}- \bigl(t_j^\star-s_i
\bigr)^{\beta
}&=&y+{\varepsilon}^{1/2}+ \bigl(t_j^\star-t_j
\bigr)^{\beta}
\\
&=&x_i-{\varepsilon}^{1/2}-2(t_j-s_i)^{\beta'}+
\bigl(t_j^\star-t_j \bigr)^{\beta}.
\end{eqnarray*}
Taking $t_j^\star=s_i+A\cdot(t_j-s_i)^{{\beta'}/{\beta}}$ for
some constant $A\in(0,\infty)$ left to be determined, we obtain from
the foregoing equality that
\begin{eqnarray*}
2(t_j-s_i)^{\beta'}&=&A^\beta
\cdot(t_j-s_i)^{\beta'}+ \bigl[A\cdot
(t_j-s_i)^{{\beta'}/{\beta}}-(t_j-s_i)
\bigr]^\beta
\\
&=&A^\beta\cdot(t_j-s_i)^{\beta'}+
\bigl[A-(t_j-s_i)^{1-{\beta'}/{\beta}} \bigr]^\beta
\cdot(t_j-s_i)^{\beta'},
\end{eqnarray*}
which shows that $A=A(t_j-s_i)$ for $A(\cdot)$ defined by (\ref
{eqAr}) upon cancelling $(t_j-s_i)^{\beta'}$ on both sides. We have
obtained (1).

From the definition (\ref{eqAr}) of $A(\cdot)$, we obtain
\begin{eqnarray*}
2A(r)^\beta&\geq& A(r)^\beta+ \bigl[A(r)-r^{1-{\beta'}/{\beta}}
\bigr]^\beta=2,
\\
2 \bigl[A(r)-r^{1-{\beta'}/{\beta}} \bigr]^\beta&\leq&A(r)^\beta+
\bigl[A(r)-r^{1-{\beta'}/{\beta}} \bigr]^\beta=2,
\end{eqnarray*}
and both inequalities in (\ref{ineqAr}) follow.
The proof is complete.
\end{pf}

As a result of Lemma~\ref{lemsct}, we have
%
%
%e3.104 #&#
\begin{equation}
\label{eqsprobsct} \P\bigl(t_j^\star<T_0^{X^i}
\bigr)\propleq{\varepsilon}(t_j-s_i)^{-{\beta '}/{\beta}},
\end{equation}
where the exponent $-\frac{\beta'}{\beta}$ is an improvement in
terms of our preceding discussion about the factor (\ref{fixprob}).
The value $t_j^\star$ will serve as the desired replacement of $t_j$.

Let us show how $t_j^\star$ still allows some independence similar to
(\ref{impilp}).

%
%
%le3.15 #&#
\begin{lem}[(Orthogonal continuation)]\label{lemoc}
Let $(\ms H_t)$ be a filtration satisfying the usual conditions, and
$U$ and $V$ be two $(\ms H_t)$-Feller diffusions such that $U_0\ind
V_0$ and, for some $(\ms H_t)$-stopping $\sigma^{\perp}$, $\langle
U,V\rangle^{\sigma^\perp}\equiv0$. Then by enlarging the underlying
filtered probability space if necessary and writing again $(\ms H_t)$
for the resulting filtration with a slight abuse of notation in this
case, we can find
a $(\ms H_t)$-Feller diffusion $\widehat{U}$ such that $\widehat
{U}\ind V$ and $\widehat{U}=U$ over $[0,\sigma^\perp]$.
\end{lem}

\begin{pf}
We only give a sketch of the proof here, and leave the details, calling
for standard arguments,
to the readers. Using L\'evy's theorem, we can define a Brownian motion
$\widehat{B}$ by
\[
\widehat{B}_t=\int_0^{T^U_0\wedge\sigma^\perp\wedge t}
\frac
{1}{\sqrt{U_s}}\,dU_s+\int_0^t
\1_{\{T^U_0\wedge\sigma^\perp<s\}}\,dB_s,
\]
for some Brownian motion $B$ independent of $(U,V)$.
We can use $\widehat{B}$ to solve for a Feller diffusion $\widehat
{U}$ with initial value $U_0$.
Then the proof of pathwise uniqueness for Feller diffusions (cf. \cite
{YWSDE}) gives $\widehat{U}=U$ on $[0,\sigma^\perp]$. Note that
$\langle\widehat{U},V\rangle\equiv0$, and consider the martingale
problem associated with a two-dimensional independent Feller diffusions
with initial values $U_0$ and $V_0$. By its uniqueness, $\widehat
{U}\ind V$.
Hence, $\widehat{U}$ is the desired continuation of $U$ beyond $\sigma
^\perp$.
\end{pf}

We apply Lemma~\ref{lemoc} to the total mass processes $X^i(\1)$ and
$Y^j(\1)$ under $\P$ and prove the following analogue of (\ref{impilp}).

%
%
%pr3.16 #&#
\begin{prop}\label{propocXY}
Let $i,j\in\mathbb N$ be given so that $s_i<t_j$. Suppose that $\sigma
^\perp$ is a $(\G_t)$-stopping time such that $\sigma^\perp\geq
t_j$ and $\langle X^i(\1),Y^j(\1)\rangle^{\sigma^\perp}\equiv0$. Then
for $r_2> r_1\geq t_j$ and nonnegative Borel measurable functions
$H_1,H_2$ and $h$,
%
%
%e3.105 #&#
\begin{eqnarray}\label{eqoc-1}
&&\E^\P\bigl[H_1 \bigl(Y^j_r(
\1);r\in[t_j,r_2] \bigr)H_2
\bigl(X^i_r(\1);r\in[s_i,r_1]
\bigr)h(y_j,x_i);r_1\leq\sigma^\perp
\bigr]\nonumber
\\
&&\qquad \leq\E^\P\bigl[H_1
\bigl(Y^j_r(\1);r\in[t_j,r_2]
\bigr) \bigr]\times\E^\P\bigl[H_2 \bigl(X^i_r(
\1);r\in[s_i,r_1] \bigr) \bigr]
\\
&&\quad\qquad{} \times\E^\P
\bigl[h(y_j,x_i) \bigr].\nonumber
\end{eqnarray}
\end{prop}
\begin{pf}
By the monotone class theorem, we may only consider the case that
\begin{eqnarray*}
H_1 \bigl(Y^j_r(\1);r\in[t_j,r_2]
\bigr)&=&H_{1,1} \bigl(Y^j_r(\1 );r
\in[t_j,r_1] \bigr)H_{1,2}
\bigl(Y^j_r(\1);r\in[r_1,r_2]
\bigr),
\\
H_2 \bigl(X^i_r(\1);r\in[s_i,r_1]
\bigr)&=&H_{2,1} \bigl(X^i_r(\1 );r
\in[s_i,t_j] \bigr)H_{2,2} \bigl(X^i_r(
\1);r\in[t_j,r_1] \bigr),
\end{eqnarray*}
for nonnegative Borel measurable functions $H_{k,\ell}$.

As the first step, we condition on $\G_{r_1}$ and obtain
%
%
%e3.106 #&#
\begin{eqnarray}\label{eqHstep1}
\qquad &&\E^\P\bigl[H_1 \bigl(Y^j_r(
\1);r\in[t_j,r_2] \bigr)H_2
\bigl(X^i_r(\1);r\in[s_i,r_1]
\bigr)h(y_j,x_i);r_1\leq\sigma^\perp
\bigr] \nonumber
\\
&&\qquad =\E^\P\bigl[ H_{1,1}
\bigl(Y^j_r(\1);r\in[t_j,r_1]
\bigr)\E^\P\bigl[H_{1,2} \bigl(Y^j_r(
\1);r\in[r_1,r_2] \bigr)\mid\G_{r_1} \bigr]
\\
&&\hspace*{109pt}{}\times  H_2 \bigl(X^i_r(\1);r
\in[s_i,r_1] \bigr)h(y_j,x_i);r_1
\leq\sigma^\perp\bigr].\nonumber
\end{eqnarray}
Since $Y^j(\1)$ is a $(\G_t)$-Feller diffusion, we know that
%
%
%e3.107 #&#
\begin{equation}
\label{defH12} \E^\P\bigl[H_{1,2} \bigl(Y^j_r(
\1);r\in[r_1,r_2] \bigr)\mid\G_{r_1} \bigr]=
\widehat{H}_{1,2} \bigl(Y^j_{r_1}(\1) \bigr)
\end{equation}
for some nonnegative Borel measurable function $\widehat{H}_{1,2}$.
Hence, from (\ref{eqHstep1}), we get
%
%
%e3.108 #&#
\begin{eqnarray}\label{eqoc-2}
&&\E^\P\bigl[H_1 \bigl(Y^j_r(
\1);r\in[t_j,r_2] \bigr)H_2
\bigl(X^i_r(\1);r\in[s_i,r_1]
\bigr)h(y_j,x_i);r_1\leq\sigma^\perp
\bigr]
\nonumber
\\
&&\qquad =\E^\P\bigl[H_{1,1} \bigl(Y^j_r(
\1);r\in[t_j,r_1] \bigr)\widehat{H}_{1,2}
\bigl(Y^j_{r_1}(\1) \bigr)
\\
&&\hspace*{48pt}{}\times H_2 \bigl(X^i_r(\1);r
\in[s_i,r_1] \bigr)h(y_j,x_i);r_1
\leq\sigma^\perp\bigr].\nonumber
\end{eqnarray}

Next, since $Y^j_{t_j}(\1)\equiv\psi(\1){\varepsilon}$ is obviously
$\P$-independent of $X^i_{t_j}(\1)$ and $\sigma^\perp\geq t_j$ by
assumption, we can do an orthogonal continuation of $X^i(\1)$ over
$[\sigma^\perp,\infty)$ by Lemma~\ref{lemoc}. This gives a Feller
diffusion $\widehat{X}{}^i$ such that $\widehat{X}{}^i\ind Y^j(\1)$ under
$\P$ and $\widehat{X}{}^{i,\sigma^\perp}=X^i(\1)^{\sigma^\perp}$. Hence,
\[
X^i(\1)=\widehat{X}{}^i\qquad\mbox{over
}[s_i,r_1]\mbox{ on } \bigl\{r_1\leq
\sigma^\perp\bigr\}
\]
and from (\ref{eqoc-2}) we get
%
%
%e3.109 #&#
\begin{eqnarray}\label{eqoc-3}
&&\E^\P\bigl[H_1 \bigl(Y^j_r(
\1);r\in[t_j,r_2] \bigr)H_2
\bigl(X^i_r(\1);r\in[s_i,r_1]
\bigr)h(y_j,x_i);r_1\leq\sigma^\perp
\bigr]\nonumber
\\
&&\qquad =\E^\P\bigl[H_{1,1} \bigl(Y^j_r(
\1);r\in[t_j,r_1] \bigr)\widehat{H}_{1,2}
\bigl(Y^j_{r_1}(\1) \bigr)\nonumber
\\
&&\hspace*{48pt}{}\times H_2 \bigl(\widehat{X}{}^i_r;r
\in[ s_i,r_1] \bigr)h(y_j,x_i);r_1
\leq\sigma^\perp\bigr]
\\
&&\qquad \leq\E^\P\bigl[H_{1,1} \bigl(Y^j_r(
\1);r\in[t_j,r_1] \bigr)\widehat{H}_{1,2}
\bigl(Y^j_{r_1}(\1) \bigr)\nonumber
\\
&&\hspace*{75pt}{}\times H_2 \bigl(\widehat
{X}{}^i_r;r\in[ s_i,r_1]
\bigr)h(y_j,x_i) \bigr],\nonumber
\end{eqnarray}
where the last inequality follows from the nonnegativity of $\widehat
{H}_{1,2}$, $H_{k,\ell}$ and $h$.

Next, we condition on $\G_{t_j}$. From (\ref{eqoc-3}), we get
%
%
%e3.110 #&#
\begin{eqnarray}\label{eqoc-4}
&&\E^\P\bigl[H_1 \bigl(Y^j_r(
\1);r\in[t_j,r_2] \bigr)H_2
\bigl(X^i_r(\1);r\in[s_i,r_1]
\bigr)h(y_j,x_i);r_1\leq\sigma^\perp
\bigr]\nonumber
\\
&&\qquad \leq\E^\P\bigl[ \E^\P
\bigl[H_{1,1} \bigl(Y^j_r(\1);r\in
[t_j,r_1] \bigr)\widehat{H}_{1,2}
\bigl(Y^j_{r_1}(\1) \bigr)
\nonumber\\[-8pt]\\[-8pt]\nonumber
&&\hspace*{107pt}{}\times H_{2,2} \bigl(
\widehat{X}{}^i_r;r\in[t_j,r_1]
\bigr)\mid\G_{t_j} \bigr]
\\
&&\hspace*{90pt}{}\times H_{2,1} \bigl(\widehat{X}{}^i_r;r
\in[s_i,t_j] \bigr)h(y_j,x_i)
\bigr].\nonumber
\end{eqnarray}
To evaluate the conditional expectation in the last term, we use the
independence between $\widehat{X}{}^i$ and $Y^j(\1)$ and deduce from the
martingale problem formulation and Theorem 4.4.2 of \cite{EK} that the
two-dimensional process $ (\widehat{X}{}^i,Y^j(\1) )\rest
[t_j,\infty)$ is $(\G_t)_{t\geq t_j}$-Markov with joint law
\[
\ms L \bigl(\widehat{X}{}^i\rest[t_j,\infty) \bigr)\otimes
\ms L \bigl(Y^j(\1)\rest[t_j,\infty) \bigr).
\]
Hence,
\begin{eqnarray*}
&&\E^\P\bigl[H_{1,1} \bigl(Y^j_r(
\1);r\in[t_j,r_1] \bigr)\widehat{H}_{1,2}
\bigl(Y^j_{r_1}(\1) \bigr)H_{2,2} \bigl(
\widehat{X}{}^i_r;r\in[t_j,r_1]
\bigr)\mid\G_{t_j} \bigr]
\\
&&\qquad =\E^\P\bigl[H_{1,1} \bigl(Y^j_r(
\1);r\in[t_j,r_1] \bigr)\widehat{H}_{1,2}
\bigl(Y^j_{r_1}(\1) \bigr) \bigr]
\\
&&\quad\qquad{}\times \E^{\mathbf P^0_{\widehat{X}{}^i_{t_j}}} \bigl[H_{2,2}
\bigl(Z_r;r\in[0,r_1-t_j] \bigr) \bigr],
\end{eqnarray*}
where we recall that $(Z,\mathbf P^0_z)$ denotes a copy of $\frac
{1}{4}\BES Q^0 (4z )$. [The value of $Y^j(\1)$ at $t_j$ is
$\psi(\1){\varepsilon}$.] Applying the foregoing equality to (\ref
{eqoc-4}) and using (\ref{defH12}), we obtain
%
%
%e3.111 #&#
\begin{eqnarray}\label{eqoc-5}
\quad &&\E^\P\bigl[H_1 \bigl(Y^j_r(
\1);r\in[t_j,r_2] \bigr)H_2
\bigl(X^i_r(\1);r\in[s_i,r_1]
\bigr)h(y_j,x_i);r_1\leq\sigma^\perp
\bigr]\nonumber
\\
&&\qquad\leq\E^{\P} \bigl[H_1 \bigl(Y^j_r(
\1);r\in[t_j,r_2] \bigr) \bigr]
\nonumber
\\
&&\quad\qquad{} \times\E^\P\bigl[\E^{\mathbf P^0_{\widehat{X}{}^i_{t_j}}}
\bigl[H_{2,2} \bigl(Z_r;r\in[0,r_1-t_j]
\bigr) \bigr] H_{2,1} \bigl(\widehat{X}{}^i_r;r
\in[s_i,t_j] \bigr)h(y_j,x_i)
\bigr]
\nonumber\\[-8pt]\\[-8pt]\nonumber
&&\qquad =\E^{\P} \bigl[H_1
\bigl(Y^j_r(\1);r\in[t_j,r_2]
\bigr) \bigr]
\\
&&\quad\qquad{}\times \E^\P\bigl[\E^{\mathbf P^0_{X^i(\1)_{t_j}}} \bigl[H_{2,2}
\bigl(Z_r;r\in[0,r_1-t_j] \bigr) \bigr]\nonumber
\\
&&\hspace*{60pt}{}\times H_{2,1} \bigl(X^i_r(\1);r\in[s_i,t_j]
\bigr)h(y_j,x_i) \bigr],\nonumber
\end{eqnarray}
where the last equality follows since we only redefine $X^i(\1)_t$ for
$t\geq\sigma^\perp$ to obtain~$\widehat{X}{}^i$, whereas $\sigma
^\perp\geq t_j$. The rest is easy to obtain. Using (\ref{impilp}),
we see that~(\ref{eqoc-5}) gives
\begin{eqnarray*}
&&\E^\P\bigl[H_1 \bigl(Y^j_r(
\1);r\in[t_j,r_2] \bigr)H_2
\bigl(X^i_r(\1);r\in[s_i,r_1]
\bigr)h(y_j,x_i);r_1\leq\sigma^\perp
\bigr]
\\
&&\qquad \leq\E^{\P} \bigl[H_1 \bigl(Y^j_r(
\1);r\in[t_j,r_2] \bigr) \bigr]
\\
&&\quad\qquad{}\times\E^\P\bigl[\E^{\mathbf P^0_{X^i(\1)_{t_j}}} \bigl[H_{2,2}
\bigl(Z_r;r\in[0,r_1-t_j] \bigr) \bigr]
H_{2,1} \bigl(X^i_r(\1);r\in[s_i,t_j]
\bigr) \bigr]
\\
&&\quad\qquad{} \times\E^\P\bigl[h(y_j,x_i)\bigr]
\\
&&\qquad =\E^{\P} \bigl[H_1 \bigl(Y^j_r(
\1);r\in[t_j,r_2] \bigr) \bigr]\E^\P
\bigl[H_2 \bigl(X^i_r(\1);r
\in[s_i,r_1] \bigr) \bigr]\E^\P
\bigl[h(y_j,x_i) \bigr].
\end{eqnarray*}
We have obtained the desired inequality, and the proof is complete.
\end{pf}

We are ready to prove Lemma~\ref{lemgr1-2} with arguments similar to
those in Section~\ref{seccondsep-1}. The following steps are labelled
in the same way as their counterparts in Section~\ref{seccondsep-1},
except that steps 2-5~and~3 below correspond to steps~2-6~and~4 in Section~\ref{seccondsep-1}, respectively. Due to the
similarity, we will only point out the key changes, leaving other
details to readers.

Recall that we fix $t\in[s_i+\frac{{\varepsilon}}{2},s_i+1]$, $i\in
\mathbb N$ and ${\varepsilon}\in(0,[8\psi(\1)]^{-1}\wedge1 ]$.
\end{longlist}

\begin{longlist}[\textit{Step} 1.]
\item[\textit{Step} 1.]
We begin with a simple observation for the
integral term
\[
\int_{t_j}^{t\wedge\tau^i\wedge\sigma_\beta^{X^i}\wedge\sigma
_\beta^{Y^j}}\frac{1}{X^i_s(\1)}\int
_\R X^i(x,s)^{1/2}Y^j(x,s)^{1/2}\,dx\,ds
\]
in (\ref{sumlate}), for $y_j=y$ satisfying (\ref{eqlateral}) and
$j\in\mathbb N$ with $t_j\in(s_i,s_i+1]$. For $s\in[t_j,t\wedge\tau
^i\wedge\sigma_{\beta}^{X^i}\wedge\sigma_\beta^{Y^j}]$ with
$s<t_j^{\star}$, the support processes of $X^i$ and $Y^j$ can be
enveloped by $\mathcal P^{X^i}_\beta(\cdot)$ and $\mathcal
P^{Y^j}_\beta(\cdot)$ up to time $s$, respectively, and $\mathcal
P^{X^i}_\beta(s)\cap\mathcal P^{Y^j}_\beta(s)=\varnothing$ by the
definition of $t_j^\star$ in (\ref{deftjstar}). Hence, for such $s$,
\[
\int_\R X^i(x,s)^{1/2}Y^j(x,s)^{1/2}\,dx=0.
\]
Using the bound (\ref{ineqQbdd}), we obtain as for (\ref
{ineqsumpre}) that
\begin{eqnarray}\label{ineqsumpre-l}
\hspace*{-4pt}&&\E^{\mathbb Q^i} \biggl[\sum_{j\in\mathcal L^i_{\beta'}(t,t\wedge\tau
^i\wedge\sigma^{X^i}_\beta)}\! \biggl(\psi(
\1){\varepsilon}\nonumber
\\
\hspace*{-4pt}&&\hspace*{98pt}{} +\int_{t_j}^{t\wedge\tau^i\wedge\sigma^{X^i}_{\beta
}\wedge\sigma _\beta^{Y^j}}\!\!\frac{1}{X^i_s(\1)}
\int_\R X^i(s,x)^{1/2}Y^j(s,x)^{1/2}\,dx\,ds
\biggr) \biggr]
\nonumber
\\
\hspace*{-4pt}&&\qquad \propleq\sum_{j\dvtx s_i<t_j\leq t}
(t-s_i)^{\beta'}{\varepsilon}
\nonumber\\[-8pt]\hspace*{-4pt} \\[-8pt]\nonumber
\hspace*{-4pt}&&\quad\qquad{} + \sum_{j\dvtx s_i<t_j\leq t}\int_{t_j}^t ds\,\1_{t_j^{\star}<s}\frac
{1}{(s-s_i)^{\eta/2}}\E^{\mathbb Q^i} \bigl[
\bigl[Y^j_s(\1) \bigr]^{1/2};
s<\tau^i\wedge\sigma^{X^i}_{\beta}
\wedge\sigma_\beta^{Y^j},\nonumber
\\
\hspace*{-4pt}&&\hspace*{209pt} 2 \bigl({\varepsilon}^{1/2}+(t_j-s_i)^{\beta'}
\bigr)\leq\llvert y_j-x_i\rrvert\nonumber
\\
\hspace*{-4pt}&&\hspace*{254pt}  \leq2 \bigl({
\varepsilon}^{1/2}+(t-s_i)^{\beta'} \bigr) \bigr].\nonumber
\end{eqnarray}
Hence, for lateral clusters, we consider
%
%e3.112 #&#
\begin{eqnarray}
&&\E^{\Q^i} \bigl[ \bigl[Y^j_s(\1)
\bigr]^{1/2};s< \tau^i\wedge\sigma^{X^i}_{\beta}
\wedge\sigma_\beta^{Y^j},\nonumber
\\
&&\hspace*{23pt} 2 \bigl({\varepsilon
}^{1/2}+(t_j-s_i)^{\beta'} \bigr)\leq
\llvert y_j-x_i\rrvert\leq2 \bigl({\varepsilon
}^{1/2}+(t-s_i)^{\beta'} \bigr) \bigr],\nonumber
\\
\eqntext{\displaystyle s\in\bigl(t_j^{\star},t\bigr], s_i<t_j<t, t_j^\star<t.}
\end{eqnarray}
\end{longlist}

\begin{longlist}[\textit{Step} 2-1.]
\item[\textit{Step} 2-1.]
We partition the event $ \{X^i_s(\1)^{T_1^{X^i}}>0 \}$ into the
two events in (\ref{event1}) and (\ref{event2}) with $t_j$ replaced
by $t_j^\star$.
Then as in (\ref{ineqsepINT1-2}), we write
%
%
%e3.113 #&#
%e3.114 #&#
\begin{eqnarray}\label{ineqsepINT1-2-l}
&&\E^{\Q^i} \bigl[ \bigl[Y^j_s(\1)
\bigr]^{1/2};s< \tau^i\wedge\sigma_\beta^{X^i}
\wedge\sigma_\beta^{Y^j},2 \bigl({\varepsilon}^{1/2}+(t_j-s_i)^{\beta'} \bigr)\nonumber
\\
&&\hspace*{108pt} \leq \llvert y_j-x_i\rrvert \leq2 \bigl({\varepsilon}^{1/2}+(t-s_i)^{\beta'}
\bigr) \bigr]\nonumber
\\
&&\qquad \leq\frac{1}{\psi(\1){\varepsilon}}\E^{\P}
\bigl[\bigl[Y^j_s(\1 ) \bigr]^{1/2};s<
\tau^i\wedge\sigma_{\beta}^{X^i}\wedge\sigma_\beta^{Y^j},
2 \bigl({\varepsilon}^{1/2}+(t_j-s_i)^{\beta'}
\bigr)\nonumber
\\
&&\hspace*{91pt}\leq\llvert y_j-x_i\rrvert\leq2 \bigl({
\varepsilon}^{1/2}+(t-s_i)^{\beta'} \bigr),T_1^{X^i}<T_0^{X^i}
\leq t_j^{\star} \bigr]\nonumber
\\
&&\quad\qquad{}+\frac{1}{\psi(\1){\varepsilon}}\E^{\P} \bigl[\bigl\llvert X^i(
\1)_s^{T^{X^i}_1}-\psi(\1){\varepsilon}\bigr\rrvert
\bigl[Y^j_s(\1) \bigr]^{1/2};s<
\tau^i\wedge\sigma_{\beta}^{X^i}\wedge
\sigma_\beta^{Y^j},\nonumber
\\
&&\hspace*{95pt} 2 \bigl({\varepsilon}^{1/2}+(t_j-s_i)^{\beta'}
\bigr)\leq\llvert y_j-x_i\rrvert
\\
&&\hspace*{119pt} \leq2 \bigl({\varepsilon}^{1/2}+(t-s_i)^{\beta'} \bigr),
X^i_s(\1)^{T^{X^i}_1}>0,t_j^{\star}<T_0^{X^i}\bigr]\nonumber
\\
&&\quad\qquad{}+\frac{1}{\psi(\1){\varepsilon}}\cdot\psi(\1){\varepsilon}\E^{\P}
\bigl[\bigl[Y^j_s(\1) \bigr]^{1/2};s<
\tau^i\wedge\sigma_{\beta}^{X^i}\wedge
\sigma_\beta^{Y^j},\nonumber
\\
&&\hspace*{131pt} 2 \bigl({\varepsilon}^{1/2}+(t_j-s_i)^{\beta'}
\bigr)\leq\llvert y_j-x_i\rrvert\nonumber
\\
&&\hspace*{131pt}\leq2 \bigl({
\varepsilon}^{1/2}+(t-s_i)^{\beta'}
\bigr),t_j^{\star}<T_0^{X^i} \bigr]\nonumber
\\
\eqntext{\forall s\in\bigl(t_j^{\star},t\bigr], s_i<t_j<t, t_j^\star<t,}
\end{eqnarray}
where we replace the event $ \{X^i_s(\1)^{T^{X^i}_1}>0,t_j^{\star
}<T_0^{X^i} \}$ by the larger one $ \{t_j^\star<T_0^{X^i}
\}$ for the third term.
\end{longlist}

\begin{longlist}[\textit{Step} 2-2.]
\item[\textit{Step} 2-2.]
Consider the first term on the right-hand side of (\ref
{ineqsepINT1-2-l}). We have
%
%
%e3.115 #&#
%e3.116 #&#
\begin{eqnarray}\label{ineqYjhalf1-0-l}
&&\frac{1}{\psi(\1){\varepsilon}}\E^{\P} \bigl[ \bigl[Y^j_s(
\1 ) \bigr]^{1/2};s<\tau^i\wedge\sigma_{\beta}^{X^i}
\wedge\sigma_\beta^{Y^j},\nonumber
\\
&&\hspace*{48pt} 2 \bigl({\varepsilon}^{1/2}+(t_j-s_i)^{\beta'}
\bigr)\leq\llvert y_j-x_i\rrvert\nonumber
\\
&&\hspace*{48pt}\leq2 \bigl({
\varepsilon}^{1/2}+(t-s_i)^{\beta'} \bigr),
T_1^{X^i}<T_0^{X^i}\leq t_j^{\star} \bigr]\nonumber
\\
&&\qquad \leq\frac{1}{\psi(\1){\varepsilon}}\E^{\P}
\bigl[\bigl[Y^j_s(\1 ) \bigr]^{1/2};s\leq
T_1^{Y^j},t_j^{\star}\leq\sigma_{\beta}^{X^i}\wedge\sigma_\beta^{Y^j},
\\
&&\hspace*{82pt} 2 \bigl({\varepsilon}^{1/2}+(t_j-s_i)^{\beta'}
\bigr)\leq\llvert y_j-x_i\rrvert\nonumber
\\
&&\hspace*{82pt}\leq2 \bigl({
\varepsilon}^{1/2}+(t-s_i)^{\beta'} \bigr),
T_1^{X^i}<T_0^{X^i}
\leq t_j^{\star} \bigr]\nonumber
\\
\eqntext{\forall s\in\bigl(t_j^\star,t\bigr],
s_i<t_j<t, t_j^\star<t.}
\end{eqnarray}
We then apply Proposition~\ref{propocXY}, taking
%
%
%e3.117 #&#
\begin{equation}
\label{choiceoc} \sigma^\perp= \bigl(\sigma_{\beta}^{X^i}
\wedge\sigma_\beta^{Y^j}\wedge t_j^\star
\bigr)\vee t_j,\qquad r_1=t_j^{\star},\qquad
r_2=s.
\end{equation}
Hence, from (\ref{ineqQbdd}) and (\ref{ineqYjhalf1-0-l}), we obtain
%
%
%e3.118 #&#
%e3.119 #&#
\begin{eqnarray}\label{ineqYjhalf1-l}
\qquad &&\frac{1}{\psi(\1){\varepsilon}}\E^{\P} \bigl[ \bigl[Y^j_s(
\1 ) \bigr]^{1/2};s<\tau^i\wedge\sigma_{\beta}^{X^i}
\wedge\sigma_\beta^{Y^j},2 \bigl({\varepsilon}^{1/2}+(t_j-s_i)^{\beta'}
\bigr)\nonumber
\\[-1pt]
&&\hspace*{58pt} \leq\llvert y_j-x_i\rrvert\leq2 \bigl({
\varepsilon}^{1/2}+(t-s_i)^{\beta'}
\bigr),T_1^{X^i}<T_0^{X^i}\leq
t_j^{\star} \bigr]
\nonumber
\\[-1pt]
&&\qquad \propleq\frac{1}{{\varepsilon}}\P\bigl(T_1^{X^i}<T_0^{X^i}
\leq t_j^{\star} \bigr) (t-s_i)^{\beta'}
\E^{\mathbf P^0_{\psi(\1){\varepsilon}}} \bigl[ (Z_{s-t_j} )^{1/2};s-t_j
\leq T^Z_1 \bigr]
\\[-1pt]
&&\qquad \propleq(t-s_i)^{\beta'} \cdot
\E^{\mathbf P^0_{\psi(\1){\varepsilon}}} \bigl[ (Z_{s-t_j} )^{1/2};s-t_j
\leq T^Z_1 \bigr]\nonumber
\\[-1pt]
\eqntext{\forall s\in\bigl(t_j^{\star},t\bigr],
s_i<t_j<t, t_j^\star<t.}
\end{eqnarray}
\end{longlist}

\begin{longlist}[\textit{Step} 2-3.]
\item[\textit{Step} 2-3.]
Let us consider the second term in (\ref{ineqsepINT1-2-l}).
As before, using (\ref{ineqimcP}) gives
%
%
%e3.120 #&#
\begin{eqnarray}\label{ineqYjhalf2-0-l}
&&\frac{1}{\psi(\1){\varepsilon}}\E^{\P} \bigl[\bigl\llvert X^i_s(
\1)^{T^{X^i}_1}-\psi(\1){\varepsilon}\bigr\rrvert\bigl[Y^j_s(
\1) \bigr]^{1/2};s< \tau^i\wedge\sigma_{\beta}^{X^i}
\wedge\sigma_\beta^{Y^j},\nonumber
\\[-1pt]
&&\hspace*{49pt} 2 \bigl({\varepsilon}^{1/2}+(t_j-s_i)^{\beta'} \bigr)
\leq\llvert y_j-x_i\rrvert\nonumber
\\[-1pt]
&&\hspace*{74pt} \leq2 \bigl({
\varepsilon}^{1/2}+(t-s_i)^{\beta'}\bigr),X^i_s(\1)^{T^{X^i}_1}>0,t_j^{\star}<T_0^{X^i}\bigr]
\nonumber
\\[-1pt]
&&\qquad \propleq\frac{{\varepsilon}^{\alpha
^{N_0}}(s-s_i)^\alpha +(s-s_i)^{\xi}}{{\varepsilon}}
\nonumber\\[-8pt]\\[-8pt]\nonumber
&&\quad\qquad{}\times \E^{\P} \bigl[
\bigl[Y^j_s(\1) \bigr]^{1/2};s\leq
T^{Y^j}_1, t_j^{\star}\leq
\sigma_{\beta}^{X^i}\wedge\sigma_\beta^{Y^j},\nonumber
\\[-1pt]
&&\hspace*{64pt} 2 \bigl({\varepsilon}^{1/2}+(t_j-s_i)^{\beta'}
\bigr)\leq\llvert y_j-x_i\rrvert\nonumber
\\[-1pt]
&&\hspace*{88pt}  \leq2 \bigl({
\varepsilon}^{1/2}+(t-s_i)^{\beta'}
\bigr),t_j^{\star}<T_0^{X^i} \bigr]\nonumber
\\[-1pt]
\eqntext{\forall s\in\bigl(t_j^\star,t\bigr],
s_i<t_j<t, t_j^\star<t.}
\end{eqnarray}
Taking the choice (\ref{choiceoc}) again, we obtain from
Proposition~\ref{propocXY}, (\ref{ineqQbdd}) and the last display
that
\begin{eqnarray}
&&\frac{1}{\psi(\1){\varepsilon}}\E^{\P} \bigl[\bigl\llvert X^i_s(
\1)^{T^{X^i}_1}-\psi(\1){\varepsilon}\bigr\rrvert\bigl[Y^j_s(
\1) \bigr]^{1/2};s< \tau^i\wedge\sigma_{\beta}^{X^i}
\wedge\sigma_\beta^{Y^j},\nonumber
\\
&&\hspace*{49pt} 2 \bigl({\varepsilon}^{1/2}+(t_j-s_i)^{\beta'}
\bigr)\leq\llvert y_j-x_i\rrvert\nonumber
\\
&&\hspace*{74pt}\leq2 \bigl({
\varepsilon}^{1/2}+(t-s_i)^{\beta'}
\bigr),X^i_s(\1)^{T^{X^i}_1}>0,t_j^{\star}<T_0^{X^i}
\bigr]\nonumber
\\[-1pt]
&&\qquad \propleq\frac{{\varepsilon}^{\alpha^{N_0}}(s-s_i)^\alpha
+(s-s_i)^\xi}{{\varepsilon}} (t-s_i)^{\beta'}\P
\bigl(t_j^{\star
}<T_0^{X^i} \bigr)\nonumber
\\[-1pt]
&&\quad\qquad{}\times
\E^{\mathbf P^0_{\psi(\1){\varepsilon}}} \bigl
[(Z_{s-t_j})^{1/2};s-t_j\leq
T_1^Z \bigr]\qquad \forall s\in\bigl(t_j^\star,t\bigr],
s_i<t_j<t, t_j^\star<t.\nonumber
\end{eqnarray}
Hence, by a computation similar to (\ref{ineqYjhalf2-1}) and
Lemma~\ref{lemsct}, the foregoing display gives
%
%
%e3.121 #&#
\begin{eqnarray}\label{ineqYjhalf2-l}
&&\frac{1}{\psi(\1){\varepsilon}}\E^{\P} \bigl[\bigl\llvert X^i_s(
\1 )^{T^{X^i}_1}-\psi(\1){\varepsilon}\bigr\rrvert\bigl[Y^j_s(
\1) \bigr]^{1/2};s< \tau^i\wedge\sigma_{\beta}^{X^i}
\wedge\sigma_\beta^{Y^j},
\nonumber
\\[-1pt]
&&\hspace*{49pt} 2 \bigl({\varepsilon}^{1/2}+(t_j-s_i)^{\beta'}
\bigr)\leq\llvert y_j-x_i\rrvert\nonumber
\\ [-1pt]
&&\hspace*{74pt}  \leq2 \bigl({
\varepsilon}^{1/2}+(t-s_i)^{\beta'}
\bigr),X^i_s(\1)^{T_1^{X^i}}>0,t_j^{\star}<T_0^{X^i}
\bigr]
\nonumber\\[-8pt]\\[-8pt]\nonumber
&&\qquad \propleq \bigl({\varepsilon}^{\alpha^{N_0}}(s-s_i)^\alpha
+(s-s_i)^\xi\bigr) (t-s_i)^{\beta'}(t_j-s_i)^{-{\beta'}/{\beta}}\nonumber
\\ [-1pt]
&&\quad\qquad{}\times \E^{\mathbf P^0_{\psi(\1){\varepsilon}}} \bigl[(Z_{s-t_j})^{1/2};s-t_j
\leq T^Z_1 \bigr]\nonumber
\\[-1pt]
\eqntext{\forall s\in\bigl(t_j^{\star},t\bigr],
s_i<t_j<t, t_j^\star<t.}
\end{eqnarray}
\end{longlist}

\begin{longlist}[\textit{Step} 2-4.]
\item[\textit{Step} 2-4.]
For the third term in (\ref{ineqsepINT1-2-l}), the calculation in the
foregoing step 2-3 readily shows
%
%
%e3.122 #&#
%e3.123 #&#
\begin{eqnarray}\label{ineqYjhalf4-l}
&&\frac{1}{\psi(\1){\varepsilon}}\cdot\psi(\1){\varepsilon}\E^{\P
} \bigl[
\bigl[Y^j_s(\1) \bigr]^{1/2};s<
\tau^i\wedge\sigma_{\beta
}^{X^i}\wedge
\sigma_\beta^{Y^j},2\bigl({\varepsilon}^{1/2}+(t_j-s_i)^{\beta'}\bigr)
\nonumber
\\[-1pt]
&&\hspace*{125pt} \leq\llvert y_j-x_i\rrvert\leq2 \bigl({
\varepsilon}^{1/2}+(t-s_i)^{\beta'} \bigr),t_j^{\star}<T_0^{X^i} \bigr]
\nonumber\\[-8pt]\\[-8pt]\nonumber
&&\qquad \propleq(t-s_i)^{\beta'}(t_j-s_i)^{-{\beta '}/{\beta}}
\cdot{\varepsilon}\cdot\E^{\mathbf P^0_{\psi(\1){\varepsilon}}} \bigl
[(Z_{s-t_j})^{1/2};s-t_j
\leq T^Z_1 \bigr]
\\[-1pt]
\eqntext{\forall s\in\bigl(t_j^{\star},t\bigr],
s_i<t_j<t, t_j^\star <t.}
\end{eqnarray}
\end{longlist}

\begin{longlist}[\textit{Step} 2-5.]
\item[\textit{Step} 2-5.]
At this step, we apply (\ref{ineqYjhalf1-l}), (\ref{ineqYjhalf2-l})
and (\ref{ineqYjhalf4-l}) to (\ref{ineqsepINT1-2-l}) and give a
summary as follows:
%
%
%e3.124 #&#
%e3.125 #&#
\begin{eqnarray}\label{ineqYjhalf6-0-l}
&&\E^{\Q^i} \bigl[ \bigl[Y^j_s(\1)
\bigr]^{1/2};s<\tau^i\wedge\sigma_{\beta}^{X^i}
\wedge\sigma_\beta^{Y^j},2 \bigl({\varepsilon}^{1/2}+(t_j-s_i)^{\beta'}
\bigr)\leq\llvert y_j-x_i\rrvert
\nonumber
\\
&&\hspace*{211pt} \leq2 \bigl({
\varepsilon}^{1/2}+(t-s_i)^{\beta'} \bigr) \bigr]
\nonumber
\\
&&\qquad \propleq \bigl[(t-s_i)^{\beta'}+(t-s_i)^{\beta'}(t_j-s_i)^{-{\beta '}/{\beta}}
\bigl({\varepsilon}^{\alpha^{N_0}}(s-s_i)^\alpha
+(s-s_i)^\xi+{\varepsilon} \bigr) \bigr]
\nonumber
\\
&&\quad\qquad{} \times\E^{\mathbf P^0_{\psi(\1){\varepsilon}}} \bigl[(Z_{s-t_j})^{1/2};s-t_j
\leq T_1^Z \bigr]
\nonumber
\\
&&\qquad \propleq (t-s_i)^{\beta'}(t_j-s_i)^{-{\beta '}/{\beta}}
\bigl((t_j-s_i)^{{\beta'}/{\beta}}+{\varepsilon}^{\alpha
^{N_0}}(s-s_i)^\alpha+(s-s_i)^\xi+{
\varepsilon} \bigr)\nonumber
\\
&&\quad\qquad{} \times{\varepsilon}^{{\alpha^{N_0}}/{2}}(s-t_j)^{{\alpha }/{2}} \biggl(
\frac{{\varepsilon}}{s-t_j} \biggr)^{1-{\alpha^{N_0}}/{4}}
\\
&&\quad\qquad{}+ (t-s_i)^{\beta'}(t_j-s_i)^{-{\beta '}/{\beta}}
\bigl((t_j-s_i)^{{\beta'}/{\beta}}+{\varepsilon}^{\alpha
^{N_0}}(s-s_i)^\alpha+(s-s_i)^\xi+{\varepsilon} \bigr)\nonumber
\\
&&\qquad\qquad{} \times{\varepsilon}^{{1}/{2}} \biggl(\frac{{\varepsilon
}}{s-t_j}\biggr)^{1-{\alpha^{N_0}}/{4}} +(t-s_i)^{\beta'}{\varepsilon
}(s-t_j)^{{\xi}/{2}-1}\nonumber
\\
&&\quad\qquad{}+(t-s_i)^{\beta'}(t_j-s_i)^{-{\beta '}/{\beta}}
\bigl({\varepsilon}^{\alpha^{N_0}}(s-s_i)^\alpha+{
\varepsilon} \bigr){\varepsilon}(s-t_j)^{{\xi}/{2}-1}\nonumber
\\
&&\quad\qquad{}+(t-s_i)^{\beta'}(t_j-s_i)^{-{\beta '}/{\beta}}(s-s_i)^\xi
{\varepsilon}(s-t_j)^{{\xi}/{2}-1}\nonumber
\\
\eqntext{\forall s\in\bigl(t_j^\star,t\bigr],
s_i<t_j<t, t_j^\star<t,}
\end{eqnarray}
where as in step 2-6 of Section~\ref{seccondsep-1}, the last
``$\propleq$''-inequality follows
again from Lemma~\ref{lemFd2}, some arithmetic, and an application of
(\ref{ineqsp-pow}).

For any $s\in(t_j^\star,t]$ with $s_i<t_j<t$ and $t_j^\star<t$, we have
\[
(t_j-s_i)^{{\beta'}/{\beta}}+{\varepsilon}^{\alpha
^{N_0}}(s-s_i)^\alpha+(s-s_i)^\xi+{
\varepsilon}\propleq(s-s_i)^\alpha,
\]
which results from (\ref{defkappa})(a), (\ref{defkappa})(b),
Lemma~\ref{lemsct} and $t_j-s_i\geq\frac{{\varepsilon}}{2}$.
Hence, with some simplifications similar to (\ref{summary1})--(\ref
{summary3}), we obtain
%
%
%e3.126 #&#
%e3.127 #&#
\begin{eqnarray}\label{ineqYjhalf6-l}
&&\E^{\Q^i} \bigl[ \bigl[Y^j_s(\1)
\bigr]^{1/2};s<\tau^i\wedge\sigma_{\beta}^{X^i}
\wedge\sigma_\beta^{Y^j},
\nonumber
\\
&&\hspace*{22pt} 2 \bigl({\varepsilon}^{1/2}+(t_j-s_i)^{\beta'}
\bigr)\leq\llvert y_j-x_i\rrvert\leq2 \bigl({
\varepsilon}^{1/2}+(t-s_i)^{\beta'} \bigr) \bigr]
\nonumber
\\
&&\qquad \propleq (t-s_i)^{\beta'}(t_j-s_i)^{-{\beta '}/{\beta }}(s-s_i)^{\alpha}(s-t_j)^{\alpha /2+{\alpha
^{N_0}}/{4}-1}{
\varepsilon}^{1+{\alpha^{N_0}}/{4}}\nonumber
\\
&&\quad\qquad{} +(t-s_i)^{\beta'}(t_j-s_i)^{\alpha /2-{\beta'}/{\beta }}(s-s_i)^{\alpha}(s-t_j)^{{\alpha^{N_0}}/{4}-1}{
\varepsilon}^{1+{\alpha^{N_0}}/{4}}
\\
&&\quad\qquad{}+(t-s_i)^{\beta'}(s-t_j)^{{\xi}/{2}-1}{
\varepsilon}\nonumber
\\
&&\quad\qquad{}+(t-s_i)^{\beta'}(t_j-s_i)^{-{\beta '}/{\beta}}(s-s_i)^\xi
(s-t_j)^{{\xi}/{2}-1}{\varepsilon}\nonumber
\\
\eqntext{\forall s\in\bigl(t_j^\star,t\bigr],
s_i<t_j<t, t_j^\star<t.}
\end{eqnarray}
\end{longlist}

\begin{longlist}[\textit{Step} 3.]
\item[\textit{Step} 3.]
We complete the proof of Lemma~\ref{lemgr1-2} in this step. Apply the
bound (\ref{ineqYjhalf6-l}) to the right-hand side of
the inequality (\ref{ineqsumpre-l}). We have
\begin{eqnarray*}
\hspace*{-4pt}&&\E^{\Q^i} \biggl[ \sum_{j\in\mathcal L_{\beta'}^i(t,t\wedge\tau
^i\wedge\sigma_{\beta}^{X^i})} \biggl(\psi(
\1){\varepsilon}
\\
\hspace*{-4pt}&&\hspace*{98pt}{}
+\int_{t_j}^{t\wedge\tau^i\wedge\sigma_{\beta
}^{X^i}\wedge\sigma
_\beta^{Y^j}}\!\!\frac{1}{X^i_s(\1)}
\int_\R X^i(x,s)^{1/2}
Y^j(x,s)^{1/2}
\,dx\,ds \biggr) \biggr]
\\
\hspace*{-4pt}&&\qquad \propleq(t-s_i)^{\beta'}\sum_{j\dvtx s_i<t_j\leq t}
{\varepsilon}
\\
\hspace*{-4pt}&&\quad\qquad{}+(t-s_i)^{\beta'}\sum_{j\dvtx s_i<t_j\leq t}(t_j-s_i)^{-{\beta '}/{\beta}}
\int_{t_j}^t (s-s_i)^{-{\eta}/{2}+\alpha}
\\
\hspace*{-4pt} && \hspace*{200pt}{}\times  (s-t_j)^{\alpha /2+{\alpha^{N_0}}/{4}-1}\,ds
\cdot {\varepsilon}^{1+{\alpha^{N_0}}/{4}}
\\
\hspace*{-4pt}&&\quad\qquad{}+(t-s_i)^{\beta'}\sum_{j\dvtx s_i<t_j\leq t}(t_j-s_i)^{\alpha /2-{\beta'}/{\beta}}
\int_{t_j}^t (s-s_i)^{-{\eta}/{2}+\alpha}
\\
\hspace*{-4pt}&&\hspace*{214pt}{} \times (s-t_j)^{{\alpha^{N_0}}/{4}-1}\,ds\cdot {\varepsilon}^{1+{\alpha^{N_0}}/{4}}
\\
\hspace*{-4pt}&&\quad\qquad{}+(t-s_i)^{\beta'}\sum_{j\dvtx s_i<t_j\leq t}\int
_{t_j}^t(s-s_i)^{-{\eta}/{2}}
(s-t_j)^{{\xi}/{2}-1}\,ds\cdot{\varepsilon}
\\
\hspace*{-4pt}&&\quad\qquad{}+(t-s_i)^{\beta'}\sum_{j\dvtx s_i<t_j\leq t}(t_j-s_i)^{-{\beta '}/{\beta}}
\int_{t_j}^t (s-s_i)^{-{\eta}/{2}+\xi}
(s-t_j)^{{\xi}/{2}-1}\,ds\cdot{\varepsilon}.
\end{eqnarray*}
Thanks to the second inequality in (\ref{defkappa})(b), the integral
domination outlined in step~3 of Section~\ref{seccondsep-1} can be
applied to each term on the right-hand side of the foregoing $\propleq
$-inequality, giving bounds which are convergent integrals. As in step~4 of Section~\ref{seccondsep-1}, we obtain
\begin{eqnarray*}
\hspace*{-4pt}&&\E^{\Q^i} \biggl[ \sum_{j\in\mathcal L_{\beta'}^i(t,t\wedge\tau
^i\wedge\sigma_{\beta}^{X^i})} \biggl(\psi(
\1){\varepsilon}
\\
\hspace*{-4pt}&&\hspace*{98pt}{} +\int_{t_j}^{t\wedge\tau^i\wedge\sigma_{\beta
}^{X^i}\wedge\sigma
_\beta^{Y^j}}\!\!\frac{1}{X^i_s(\1)}
\int_\R X^i(x,s)^{1/2}Y^j(x,s)^{1/2}
\,dx\,ds \biggr) \biggr]
\\
\hspace*{-4pt}&&\qquad \propleq(t-s_i)^{\beta'+1}+(t-s_i)^{\beta'-{\eta}/{2}+({3\alpha}/{2}) }
\cdot{\varepsilon}^{{\alpha^{N_0}}/{4}}+ (t-s_i)^{\beta'-{\eta}/{2}+({3\xi}/{2})},
\end{eqnarray*}
which proves Lemma~\ref{lemgr1-2}.
\end{longlist}

%s4 #&#
\section{Uniform separation of approximating solutions}\label{secsepSOL}
In this section, we prove the main result of the present paper that
there is pathwise nonuniqueness in the SPDE (\ref{eqmainSPDE}). The
result is summarized in Theorem~\ref{teomains}.
We continue to suppress the dependence on ${\varepsilon}$ of the
approximation solutions and emphasize it only through
$\P_{\varepsilon}$, unless otherwise mentioned.

Our program to obtain uniform separation of the approximating solutions
is sketched as follows (cf. the discussion in Section~\ref{secintro}).
For small $r,{\varepsilon}\in(0,1]$, we will define an event $S(r
)=S_{\varepsilon}(r)$ which keeps track of certain separation of the
\mbox{${\varepsilon}$-}approximating solutions $X$ and $Y$ over the territory
of a ``large'' immigrant process $X^i$. The immigrant processes range
over those large and arriving approximately by time $r$. The definition
of $S(r)$ is based on
the earlier results for conditional separation of the approximating
solutions. The effect is that these particular events $S(r)$ imply the
required uniform separation:
for some $\Delta( r )\in(0,\infty)$ depending only on the parameter
vector in Assumption~\ref{ass} and $r$, we have
%
%
%e4.1 #&#
\begin{equation}
\label{eqsepinc} S( r) \subseteq\Bigl\{\sup_{0\leq s\leq2r}\llVert X
_s-Y_s\rrVert_\rap\geq\Delta(r ) \Bigr\}
\end{equation}
[recall the definition of $\llVert \cdot\rrVert_\rap$ in (\ref{eqCrap-norm})]
and, for fixed $r$,
%
%
%e4.2 #&#
\begin{equation}
\label{eqsepprob} \liminf_{{\varepsilon}\searrow0}\P_{\varepsilon} \bigl(S(r)
\bigr)>0.
\end{equation}

Let us give the precise definition of the events $S(r)$ and discuss the
ingredients. First,
recall the parameter vector chosen in Assumption~\ref{ass} as well as
the constants $\kappa_j$ defined in Theorem~\ref{teogr}. We need to
use small portions of the constants $\kappa_1$~and~$\kappa_3$, and by
(\ref{defkappa})(d) we can find a constant $\wp$ satisfying
\[
\wp\in(0,\kappa_1\wedge\kappa_3 )\qquad\mbox{such that }
\kappa_1-\wp>\eta.
\]
We insist that $\wp$ depends only on the parameter vector in (\ref{param}).
For any $i\in\mathbb N$, ${\varepsilon}\in(0,[8\psi(\1
)]^{-1}\wedge1 ]$, and random time $T\geq s_i$, let
$G^{i}(T)=G^i_{\varepsilon}(T)$ be the growth event defined by
%
%
%e4.3 #&#
\begin{eqnarray}
G^i(T)=\lleft\{ %
\begin{array} {ll}
X^{i}_s\bigl(\bigl[x_i-{
\varepsilon}^{1/2}-(s-s_i)^\beta,x_i+{
\varepsilon}^{1/2}+(s-s_i)^\beta\bigr]\bigr)
\\[3pt]
\qquad \geq\displaystyle\frac{(s-s_i)^{\eta}}{4}\quad\mbox{and}
\\[9pt]
Y_s\bigl(\bigl[x_i-{\varepsilon}^{1/2}-(s-s_i)^\beta,x_i+{
\varepsilon}^{1/2}+(s-s_i)^\beta\bigr]\bigr)
\\[3pt]
\qquad\leq K^* \bigl[(s-s_i)^{\kappa_1-\wp}+{\varepsilon
}^{\kappa_2}( s-s_i)^{\kappa_3-\wp} \bigr]
\\
\multicolumn{1}{r}{\forall s \in[s_i,T]}
\end{array}
\rright\},\label{defGiT}
\end{eqnarray}
where the constant $K^*\in(0,\infty)$ is as in Theorem~\ref{teogr}.
Note that
$G^i( \cdot)$ is \textit{decreasing} in the sense that, for any random
times $T_1,T_2$ with $T_1\leq T_2$, $G^i(T_1)\supseteq G^i(T_2)$. Later
on when taking into consideration of the support propagation of the
immigrant processes,
we will explain how the description in (\ref{defGiT}) for $X$ is
related to the stopping time $\tau^{i,(1)}$ underlying (\ref
{eqseptimebdd}), and how the
description in the same display for $Y$
is related to the event in (\ref{ineqgrbdd}) for partial sums of the
total mass processes $Y^j(\1)$.

Next, we set
%
%
%e4.4 #&#
\begin{equation}
\label{eqtheta} \qquad S(r)=S_{\varepsilon}(r)\triangleq\bigcup
_{i=1}^{\lfloor
r{\varepsilon}^{-1}\rfloor}G^i_{\varepsilon}(s_i+r),
\qquad r\in(0,\infty), {\varepsilon}\in\biggl(0,\frac{1}{8\psi(\1)}\wedge1 \biggr].
\end{equation}
Whereas the event $S(r)$ depends on all immigrant processes $X^i$
arriving approximately by time $r$,
it is intended to keep track of separation over the territories of the
``large'' ones, as we have planned above.
This idea underlies the arguments below, and will become explicit
in the proof of Lemma~\ref{lemus3} where
inclusion--exclusion and conditioning come into play in applying the
result of Theorem~\ref{teogr} immigrant-by-immigrant with respect to $X^i$'s.

By the following lemma, (\ref{eqsepinc}) is a simple consequence of
the events $G^i( \cdot)$.

%
%
%le4.1 #&#
\begin{lem}\label{lemus1}
For some $r_0\in(0,1]$, we can find ${\varepsilon}_0( r )\in
(0,r\wedge[8\psi(\1)]^{-1}\wedge1 ]$ and $\Delta(r )\in
(0,\infty)$ for any $r \in(0,r_0]$ so that the inclusion (\ref
{eqsepinc}) holds almost surely for any ${\varepsilon}\in
(0,{\varepsilon}_0(r )]$. The constant $\Delta(r )$ depends only on
$r$ and the parameter vector chosen in Assumption~\ref{ass}.
\end{lem}
\begin{pf}
First, we specify the strictly positive numbers $r_0$, ${\varepsilon
}_0( r )$, and $\Delta(r)$.
Since the small portion $\wp$ taken away from $\kappa_1$ and $\kappa
_3$ satisfies $\kappa_1-\wp>\eta$, we can choose $r_0\in(0,1]$ such that
%
%
%e4.5 #&#
\begin{equation}
\label{choicet0} \frac{r^{\eta}}{4}-2K^*r^{\kappa_1-\wp}>0\qquad\forall r
\in(0,r_0].
\end{equation}
Then we choose, for every $r\in(0,r_0]$, a number ${\varepsilon}_0(r
)\in(0,r\wedge[8\psi(\1)]^{-1}\wedge1 ]$ such that
%
%
%e4.6 #&#
\begin{equation}
\label{choicevep0} 0<{\varepsilon}_0(r )^{\kappa_2} \leq
r^{\kappa_1-\kappa_3}.
\end{equation}
Finally, we set
%
%
%e4.7 #&#
\begin{equation}
\label{eqDelta} \Delta(r )\triangleq\frac{1}{2} \biggl[ \biggl(\frac{{r ^{\eta}}/{4}-2K^*r ^{\kappa_1-\kappa}}{2+2r ^\beta} \biggr)\wedge1 \biggr
]>0,\qquad r \in(0,r_0].
\end{equation}

We check that the foregoing choices give (\ref{eqsepinc}). Fix $r\in
(0,r_0]$, ${\varepsilon}\in(0,{\varepsilon}_0(r )]$, and $1\leq
i\leq\lfloor r{\varepsilon}^{-1}\rfloor$ (note that $\lfloor
r{\varepsilon}^{-1}\rfloor\geq1$ since ${\varepsilon}\leq r$). In
this paragraph, we assume that the event $G^i(s_i+r)$ occurs.
Then by definition,
%
%
%e4.8 #&#
\begin{eqnarray}
\label{ineqYgrowth} &&Y_s\bigl(\bigl[x_i-{
\varepsilon}^{1/2}-(s-s_i)^\beta,x_i+{
\varepsilon}^{1/2}+(s-s_i)^\beta\bigr]\bigr)
\nonumber\\[-8pt]\\[-8pt]\nonumber
&&\qquad  \leq K^* \bigl[(s-s_i)^{\kappa_1-\wp}+{\varepsilon
}^{\kappa_2} (s-s_i)^{\kappa_3-\wp} \bigr]\qquad\forall s
\in[s_i,s_i+r].
\end{eqnarray}
In particular, (\ref{choicevep0}) and (\ref{ineqYgrowth}) imply that
\[
Y_{s_i+r}\bigl(\bigl[x_i-{\varepsilon}^{1/2}-r^\beta,x_i+{
\varepsilon}^{1/2}+r^\beta\bigr]\bigr) \leq 2K^*r^{\kappa_1-\wp}.
\]
Since $X\geq X^i$, the last inequality and the definition of
$G^i(s_i+r)$ imply
\begin{eqnarray*}
\hspace*{-4pt} &&X_{s_i+r}\bigl(\bigl[x_i-{\varepsilon}^{1/2}-r^\beta,x_i+{
\varepsilon}^{1/2}+r^\beta\bigr]\bigr)
-Y_{s_i+r}\bigl(\bigl[x_i-{
\varepsilon}^{1/2}-r^\beta,x_i+{
\varepsilon}^{1/2}+r^\beta\bigr]\bigr)
\\
\hspace*{-4pt}&&\qquad \geq\frac{r^{\eta
}}{4}-2K^*r^{\kappa_1-\wp},
\end{eqnarray*}
where the lower bound is strictly positive by (\ref{choicet0}). To
carry this to the $\C_\rap(\R)$-norm of $X_{s_i+r}-Y_{s_i+r}$,
we make an elementary observation: if $f$ is Borel measurable,
integrable on a finite interval $I$, and satisfies $\int_I f>A$, then
there must exist some $x\in I$ such that $f(x)> A/\ell(I)$, where
$\ell(I)$ is the length of $I$. Using this, we obtain from the last
inequality that, for some $x\in[x_i-{\varepsilon}^{1/2}-r^\beta,x_i+{\varepsilon}^{1/2}+r^\beta]$,
\[
X(x,s_i+r)-Y(x,s_i+r)\geq\frac{{r^{\eta}}/{4}-2K^*r^{\kappa
_1-\wp}}{2{\varepsilon}^{1/2}+2r^\beta} \geq
\frac{{r^\eta}/{4}-2K^*r^{\kappa_1-\wp}}{2+2r^\beta},
\]
so the definition of $\llVert \cdot\rrVert_\rap$ [in (\ref{eqCrap-norm})] and
the definition (\ref{eqDelta}) of $\Delta(r )$ entail
\[
\Delta(r )\leq\llVert X_{s_i+r}-Y_{s_i+r}\rrVert
_\rap\leq\sup_{0\leq s\leq
2r}\llVert X_s-Y_s
\rrVert_\rap,
\]
where the second inequality follows since $s_i=\frac
{(2i-1)}{2}{\varepsilon}$ and $1\leq i\leq\lfloor r{\varepsilon
}^{-1}\rfloor$.

In summary, we have shown that (\ref{eqsepinc}) holds because each
component $G^i(s_i+r)$ of $S(r)$ satisfies the analogous inclusion.
The proof is complete.
\end{pf}

We proceed to the proof of (\ref{eqsepprob}) for small enough $r>0$.
From now on, we take into account the support propagation of the
immigrant processes. The major argument will be in Lemma~\ref{lemus3} below.
As the preliminary step to use Theorem~\ref{teogr}, we bring the
involved stopping times into the events $G^i( \cdot)$ and then
translate the descriptions about $Y$ in (\ref{defGiT}) into ones
about its immigrant processes (see Lemma~\ref{lemus2} below).

Recall Figure~\ref{Fig1}, and define the event $\Gamma^i(r )=\Gamma
^i_{\varepsilon}(r )$ by
%
%
%e4.9 #&#
%e4.10 #&#
\begin{eqnarray}\label{eqRit}
\Gamma^i(r )&\triangleq& \biggl\{\mc
P_\beta^{X^i}(s_i+r)\cap\biggl(\bigcup
_{j\dvtx t_j\leq s_i}\supp\bigl(Y^j\bigr) \biggr)=\varnothing
\biggr\}\nonumber
\\
&&{} \cap\bigcap_{j\dvtx t_j\leq s_i+r} \bigl\{
\sigma_\beta^{Y^j}>t_j+3r \bigr\}
\cap\bigl\{\sigma_{\beta}^{X^i}>s_i+2r
\bigr\},
\\
\eqntext{\displaystyle r\in(0,1], i\in\mathbb N, {\varepsilon}\in\biggl(0,\frac
{1}{8\psi(\1)}
\wedge1 \biggr],}
\end{eqnarray}
where $\supp(Y^j)$ denotes the topological support of the
two-parameter (random) function $(x,s)\lmt Y^j(x,s)$, and the time
durations $r$, $3r$ and $2r$ on the right-hand side in restricting
$\mathcal P^{X^i}_\beta$, $\sigma^{Y^j}_\beta$ and $\sigma
^{X^i}_\beta$, respectively, are chosen only for technical convenience
and can be replaced by suitably large constant multiples of $r$.
Through $\Gamma^i( r)$, we confine the ranges of the supports of
$Y^j$, for $j\in\mathbb N$ satisfying $t_j\leq s_i+r$, and $X^i$. It will
become clear in passing that one of the reasons for considering this
event is to make precise the informal argument of choosing $\mathcal
J^i_{\beta'}( \cdot)$, as discussed in Section~\ref{sec2ndsep-setup}.

%
%
%le4.2 #&#
\begin{lem}\label{lemus2}
Fix $r\in(0,1]$, $i\in\mathbb N$ and ${\varepsilon}\in(0,[8\psi
(\1)]^{-1}\wedge
1 ]$.
Then on the event $\Gamma^i( r)$ defined by (\ref{eqRit}),
we have
%
%
%e4.11 #&#
%e4.12 #&#
\begin{eqnarray}
\label{eqYloc} %
&&Y_s\bigl(\bigl[x_i-{
\varepsilon}^{1/2}-(s-s_i)^\beta,x_i+{
\varepsilon}^{1/2}+(s-s_i)^\beta\bigr]\bigr)\nonumber
\\
&&\qquad =\sum_{j\in\mathcal J_{\beta
'}^i(s)}Y^j_s
\bigl(\bigl[x_i-{\varepsilon}^{1/2}-(s-s_i)^\beta,x_i+{
\varepsilon}^{1/2}+(s-s_i)^\beta\bigr]\bigr)
\\
\eqntext{\forall s\in[s_i,s_i+r].}
\end{eqnarray}
In particular, on $\Gamma^i( r)$,
%
%
%e4.13 #&#
%e4.14 #&#
\begin{eqnarray}\label{eqYloc1}
Y_s\bigl(\bigl[x_i-{
\varepsilon}^{1/2}-(s-s_i)^\beta,x_i+{
\varepsilon}^{1/2}+(s-s_i)^\beta\bigr]\bigr)\leq
\sum_{j\in\mathcal J_{\beta
'}^i(s)}Y^j_s(\1)
\nonumber\\[-8pt]\\[-8pt]
\eqntext{\forall s\in[s_i,s_i+r].}
\end{eqnarray}
\end{lem}
\begin{pf}
In this proof, we argue on the event $\Gamma^i(r ) $ and call $\Theta
_s\triangleq\{x;(x,s)\in\Theta\}$ the \textit{$s$-section} of a
subset $\Theta$ of $\R\times\R_+$ for any $s\in\R_+$.

Consider\vspace*{1pt} (\ref{eqYloc}). Since the $s$-section $\supp(Y^j)_s$
contains the support of $Y^j_s( \cdot)$, it suffices to show that, for
any $s\in[s_i,s_i+r]$ and $j\in\mathbb N$ with $t_j\leq s$ and $j\notin
\mathcal J_{\beta'}^i(s)$,
%
%
%e4.15 #&#
\begin{equation}
\bigl[x_i-{\varepsilon}^{1/2}-(s-s_i)^\beta,x_i+{
\varepsilon}^{1/2}+(s-s_i)^\beta\bigr]\cap\supp
\bigl(Y^j\bigr)_s=\varnothing.\label{eqsuppsec}
\end{equation}

If $j\in\mathbb N$ satisfies $t_j\leq s_i$, then using the first item in
the definition (\ref{eqRit}) of $\Gamma^i(r )$ gives
\[
\mathcal P^{X^i}_\beta(s_i+r)\cap\supp
\bigl(Y^j\bigr) =\varnothing.
\]
Hence, taking the $s$-sections of both $\mathcal P^{X^i}_\beta(s_i+r)$
and $\supp(Y^j)$ shows that $Y^j$ satisfies (\ref{eqsuppsec}).

Next, suppose that $j\in\mathbb N$ satisfies $s_i<t_j\leq s$ but
$j\notin
\mathcal J_{\beta'}^i(s)$. On one hand, this choice of $j$ implies
\[
\llvert y_j-x_i\rrvert>2 \bigl({\varepsilon}^{1/2}+(s-s_i)^{\beta'}
\bigr)\geq2 \bigl({\varepsilon}^{1/2}+(s-s_i)^\beta
\bigr),
\]
where the second inequality follows from the assumption $r\in(0,1]$
and the choice $\beta'<\beta$ by (\ref{defkappa})(b), so
Lemma~\ref{lemoutclus+si} entails
%
%
%e4.16 #&#
\begin{equation}
\label{eqsuppXYdis} \mathcal P^{X^i}_\beta(s)\cap\mathcal
P^{Y^j}_\beta(s)=\varnothing.
\end{equation}
On the other hand, using the second item in the definition of $\Gamma
^i(r)$, we deduce that
\[
\supp\bigl(Y^j\bigr)\cap\bigl(\R\times[t_j,t_j+3r]
\bigr)\subseteq\mathcal P^{Y^j}_\beta(t_j+3r).
\]
Using $t_j+r>s_i+r\geq s$
and taking $s$-sections of $\supp(Y^j)$ and $\mathcal P_\beta
^{Y^j}(t_j+3r)$, we obtain from the foregoing inclusion that
%
%
%e4.17 #&#
\begin{eqnarray}\label{eqsuppXYdis1}
\supp\bigl(Y^j\bigr)_s&\subseteq& \bigl[y_j-{
\varepsilon}^{1/2}-(s-t_j)^\beta,y_j+{
\varepsilon}^{1/2}+(s-t_j)^{\beta} \bigr]\nonumber
\\
& =& \mathcal
P^{Y_j}_\beta(t_j+3r)_s
\\
&=&\mathcal P^{Y_j}_\beta(s)_s.\nonumber
\end{eqnarray}
Since
\[
\mathcal P_\beta^{X^i}(s)_s=
\bigl[x_i-{\varepsilon}^{1/2}-(s-s_i)^\beta,x_i+{\varepsilon}^{1/2}+(s-s_i)^\beta
\bigr],
\]
(\ref{eqsuppXYdis}) and (\ref{eqsuppXYdis1}) give our assertion
(\ref{eqsuppsec}) for $j\in\mathbb N$ satisfying $s_j<t_j\leq s$ and
$j\notin\mathcal J_{\beta'}^i(s)$.
We have considered all cases for which $j\in\mathbb N$, $t_j\leq s$, and
$j\notin\mathcal J_{\beta'}^i(s)$. The proof is complete.
\end{pf}

Recall $r_0\in(0,1]$ and ${\varepsilon}_0(r )\in(0,r\wedge
[8\psi(\1)]^{-1}\wedge1 ]$ chosen in Lemma~\ref{lemus1} and the
events $S( r)$ in (\ref{eqtheta}).
The following lemma completes the last step (\ref{eqsepprob}) to
obtain uniform separation of the approximating solutions.

%
%
%le4.3 #&#
\begin{lem}\label{lemus3}
For some $r_1\in(0,r_0]$, we can find ${\varepsilon}_1(r )\in(0,
{\varepsilon}_0(r )]$ for any $r\in(0,r_1]$ such that
%
%
%e4.18 #&#
\begin{equation}
\label{inequnifsep} %
\inf_{{\varepsilon}\in(0,{\varepsilon}_1(r )]}\P_{\varepsilon
}
\bigl(S(r ) \bigr)>0.
\end{equation}
\end{lem}
\begin{pf}
In this proof, we transfer the $\mathbb Q^i_{\varepsilon}$-probabilities
of separation in Theorem~\ref{teogr} to $\P_{\varepsilon
}$-probabilities of separation by conditioning and use
inclusion--exclusion as in \cite{MMP}. The latter makes the $\P
_{\varepsilon}$-probabilities of separation stand out among others.

For\vspace*{2pt} any $i\in\mathbb N$, ${\varepsilon}\in(0,[8\psi(\1
)]^{-1}\wedge1 ]$, and random time $T\geq s_i$, we define
$\widehat{G}{}^i(\cdot)=\widehat{G}{}^i_{\varepsilon}(\cdot)$ by
%
%
%e4.19 #&#
\begin{equation}
\qquad \widehat{G}{}^i(T)=\lleft\{ %
\begin{array} {ll}
\displaystyle X^i_s(\1)\geq\frac{(s-s_i)^{\eta}}{4}\quad \mbox{and}
\\[6pt]
\displaystyle\sum_{j\in\mathcal J_{\beta'}^i(s)}Y^j_s(
\1)\leq K^*\bigl[(s-s_i)^{\kappa_1-\wp}
+{\varepsilon}^{\kappa_2} (s-s_i)^{\kappa_3-\wp}
\bigr]
\\
\multicolumn{1}{r}{\forall s\in[s_i,T]}
\end{array}
\rright
\}.\label{defGiT1}
\end{equation}
Note that $\widehat{G}{}^i(\cdot)$ is decreasing, and its definition
about the masses of $Y$ is the same as the event considered in
Theorem~\ref{teogr} except for the restrictions from stopping times
$\tau^i$, $\sigma^{X^i}_\beta$ and $\sigma_\beta^{Y^j}$.

The connection between $\widehat{G}{}^i(\cdot)$ and $G^i(\cdot)$ is as follows.
First, we note that
by~(\ref{eqYloc1}), the statement about the masses of $Y$
in
$\widehat{G}{}^i(r )\cap\Gamma^i(r )$ implies that in $G^i(r )\cap
\Gamma^i(r )$. Also, the statements in $G^i( \cdot)$ and $\widehat
{G}{}^i( \cdot)$ concerning the masses of $X^i$ are linked by the
obvious equality:
\[
X^i_s(\1)=X^i_s\bigl(
\bigl[x_i-{\varepsilon}^{1/2}-(s-s_i)^{\beta
},x_i+{
\varepsilon}^{1/2}+(s-s_i)^{\beta}\bigr]\bigr)\qquad
\forall s\in\bigl[s_i,\sigma_{\beta}^{X^i} \bigr].
\]
Since $\sigma^{X^i}_\beta>s_i+2r$ on $\Gamma^i( r )$,
we are led to the inclusion
%
%
%e4.20 #&#
\begin{equation}
\label{ineqGGhat0} \widehat{G}{}^i \bigl( \tau^i
\wedge(s_i+r) \bigr)\cap\Gamma^i(r )\subseteq
G^i \bigl( \tau^i\wedge(s_i+r) \bigr)\cap
\Gamma^i(r)
\end{equation}
for any $r\in(0,1]$, $i\in\mathbb N$ and ${\varepsilon}\in
(0,[8\psi(\1)]^{-1}\wedge1 ]$ ($\tau^i$ is defined in
Proposition~\ref{propseptime}).
We can also write (\ref{ineqGGhat0}) as
%
%
%e4.21 #&#
\begin{equation}\label{ineqGGhat}
\qquad \widehat{G}{}^i \bigl( \widehat{\tau}{}^i(s_i+r)
\wedge(s_i+r) \bigr)\cap\Gamma^i(r )\subseteq
G^i \bigl( \widehat{\tau}{}^i(s_i+r)
\wedge(s_i+r) \bigr)\cap\Gamma^i(r ),
\end{equation}
where
%
%
%e4.22 #&#
\begin{equation}
\label{eqtauhat} \widehat{\tau}{}^i(s_i+r)\triangleq
\tau^i\wedge\sigma_{\beta
}^{X^i}\wedge\bigwedge
_{j\dvtx s_i<t_j\leq s_i+r}\sigma_\beta^{Y^j}.
\end{equation}
Here, although the restriction $\sigma^{X^i}_\beta\wedge\bigwedge
_{j\dvtx s_i<t_j\leq s_i+r}\sigma_\beta^{Y^j}$ is redundant in (\ref{ineqGGhat})
[because $\sigma_\beta^{Y^j}>t_j+3r>s_i+r$ for each $j\in\mathbb N$
with $s_i<t_j\leq s_i+r$ by the definition of $\Gamma^i(r)$], we
emphasize its role by writing it out.

We start bounding $\P_{\varepsilon} (S(r ) )$.
For any $r\in(0,r_0]$ and ${\varepsilon}\in(0, {\varepsilon
}_0(r)]$, we have
\begin{eqnarray*}
\P_{\varepsilon} \bigl(S(r) \bigr)
&\geq& \P_{\varepsilon} \Biggl(\bigcup_{i=1}^{\lfloor r{\varepsilon
}^{-1}\rfloor}G^i
\bigl( \widehat{\tau}{}^i(s_i+r)\wedge(s_i+r)
\bigr)\cap\Gamma^i(r )\cap\bigl\{T^{X^i}_1<T^{X^i}_0\bigr\}
\\
&&\hspace*{157pt}{} \cap\bigl\{ \widehat{\tau}{}^i(s_i+r)\geq
s_i+ r \bigr\} \Biggr)
\nonumber
\\
&\geq& \P_{\varepsilon} \Biggl(\bigcup_{i=1}^{\lfloor r{\varepsilon
}^{-1}\rfloor}
\widehat{G}{}^i \bigl( \widehat{\tau}{}^i(s_i+r)
\wedge(s_i+r) \bigr)\cap\Gamma^i(r )\cap\bigl
\{T^{X^i}_1<T^{X^i}_0 \bigr\}
\\
&&\hspace*{157pt}{}\cap
\bigl\{ \widehat{\tau}{}^i(s_i+r)\geq s_i+r
\bigr\} \Biggr),
\end{eqnarray*}
where the last inequality follows from the inclusion (\ref
{ineqGGhat}). We make the restrictions $ \{T_1^{X^i}<T_0^{X^i}
\}$ in order to invoke $\Q^i$-probabilities later on.
By considering separately $\widehat{\tau}{}^i(s_i+r)\geq s_i+r$ and $
\widehat{\tau}{}^i(s_i+r)<s_i+r$, we obtain from the last inequality that
%
%
%e4.23 #&#
\begin{eqnarray}
\label{ineqsepprob0} \qquad\P_{\varepsilon} \bigl(S(r ) \bigr)&\geq&
\P_{\varepsilon} \Biggl(\bigcup_{i=1}^{\lfloor r{\varepsilon
}^{-1}\rfloor}
\widehat{G}_i \bigl( \widehat{\tau}{}^i(s_i+r)
\wedge(s_i+r) \bigr)\cap\Gamma^i(r )\cap\bigl
\{T^{X^i}_1<T^{X^i}_0 \bigr\} \Biggr)
\nonumber\\[-8pt]\\[-8pt]\nonumber
&&{}-\P_{\varepsilon} \Biggl(\bigcup_{i=1}^{\lfloor r{\varepsilon
}^{-1}\rfloor}
\bigl[ \widehat{\tau}{}^i(s_i+r)<s_i+r \bigr]
\cap\bigl\{ T^{X^i}_1<T^{X^i}_0 \bigr
\} \Biggr).
\end{eqnarray}
Applying another inclusion--exclusion to the first term on the
right-hand side of~(\ref{ineqsepprob0}) gives the main inequality of
this proof:
%
%
%e4.24 #&#
%e4.25 #&#
\begin{eqnarray}\label{ineqsepprob}
\P_{\varepsilon} \bigl(S(r ) \bigr) &\geq&
\P_{\varepsilon} \Biggl(\bigcup_{i=1}^{\lfloor r{\varepsilon
}^{-1}\rfloor} \widehat{G}{}^i \bigl(\widehat{ \tau}{}^i(s_i+r)
\wedge(s_i+r) \bigr)\cap\bigl\{T^{X^i}_1<T^{X^i}_0
\bigr\} \Biggr)\nonumber
\\
&&{}-\P_{\varepsilon} \Biggl(\bigcup_{i=1}^{\lfloor r{\varepsilon
}^{-1}\rfloor}
\Gamma^i(r )^\complementtt\cap\bigl\{ T^{X^i}_1<T^{X^i}_0
\bigr\} \Biggr)
\nonumber\\[-8pt]\\[-8pt]\nonumber
&&{}-\P_{\varepsilon} \Biggl(\bigcup_{i=1}^{\lfloor r{\varepsilon
}^{-1}\rfloor}
\bigl\{ \widehat{\tau}{}^i(s_i+r)<s_i+r \bigr\}
\cap\bigl\{ T^{X^i}_1<T^{X^i}_0 \bigr
\} \Biggr)
\\
\eqntext{\forall r\in(0,r_0], {\varepsilon}\in\bigl(0, {
\varepsilon}_0(r) \bigr].}
\end{eqnarray}
In the rest of this proof, we bound each of the three terms on the
right-hand side of (\ref{ineqsepprob}) and then choose according to
these bounds the desired $r_1$ and ${\varepsilon}_1(r )$ for~(\ref{inequnifsep}).

At this stage, we use Proposition~\ref{propseptime} and Theorem~\ref{teogr} in the following way.
For any $\rho\in(0,\frac{1}{2})$, we choose $\delta_1\in(0,1]$,
independent of $i\in\mathbb N$ and ${\varepsilon}\in(0,[8\psi(\1
)]^{-1}\wedge1 ]$, such that
%
%
%e4.26 #&#
%e4.27 #&#
\begin{eqnarray}
\qquad &&\sup\biggl\{\Q_{\varepsilon}^i\bigl(\tau^i\leq
s_i+\delta_1\bigr);i\in\mathbb N, {\varepsilon}\in\biggl(0,
\frac{1}{8\psi(\1)}\wedge1 \biggr] \biggr\} \leq \rho,\label{choicedelta1-1}
\\
&&\sup\biggl\{\Q_{\varepsilon}^i \biggl(
\exists s\in(s_i,s_i+\delta_1], \sum
_{j\in\mathcal J_{\beta'}^i(s\wedge\tau^i\wedge\sigma
_{\beta}^{X^i})}Y^j_s(\1)^{ \tau^i\wedge\sigma_{\beta}^{X^i}\wedge
\sigma_\beta^{Y^j}}\nonumber
\\
\label{choicedelta1-2}  &&\hspace*{86pt} > K^*\bigl[(s-s_i)^{\kappa_1-\wp}+{
\varepsilon}^{\kappa_2}\cdot(s-s_i)^{\kappa_3-\wp}\bigr]
\biggr);
\\
&&\hspace*{157pt} i\in\mathbb N, {\varepsilon}\in\biggl(0,\frac{1}{8\psi(\1)}\wedge1
\biggr] \biggr\}\leq\rho.\nonumber
\end{eqnarray}

Consider the first probability on the right-hand side of (\ref
{ineqsepprob}). We use the elementary inequality: for any events
$A_1,\ldots,A_n$ for $n\in\mathbb N$,
\[
\P\Biggl(\bigcup_{j=1}^nA_j
\Biggr)\geq\sum_{j=1}^n
\P(A_j)-\sum_{i=1}^n\sum
_{\stackrel{\scriptstyle j\dvtx j\neq i}{\scriptstyle1\leq
j\leq n}}\P(A_i\cap A_j).
\]
Then
%
%
%e4.28 #&#
%e4.29 #&#
\begin{eqnarray}\label{ineqsepprob1}
&&\P_{\varepsilon} \Biggl(\bigcup_{i=1}^{\lfloor r{\varepsilon
}^{-1}\rfloor}
\widehat{G}{}^i \bigl(\widehat{\tau}{}^i(s_i+r)
\wedge(s_i+r) \bigr)\cap\bigl\{T^{X^i}_1<T^{X^i}_0
\bigr\} \Biggr)
\nonumber
\\
&&\qquad \geq\sum_{i=1}^{\lfloor r{\varepsilon}^{-1}\rfloor}
\P_{\varepsilon} \bigl(\widehat{G}{}^i \bigl(\widehat{\tau
}{}^i(s_i+r)\wedge(s_i+r) \bigr)\cap\bigl
\{T^{X^i}_1<T^{X^i}_0 \bigr\} \bigr)
\nonumber\\[-8pt]\\[-8pt]\nonumber
&&\quad\qquad{} -\sum_{i=1}^{\lfloor r{\varepsilon}^{-1}\rfloor}\sum
_{\stackrel{\scriptstyle j\dvtx j\neq i}{\scriptstyle1\leq j\leq\lfloor
r{\varepsilon}^{-1}\rfloor}}\P_{\varepsilon} \bigl
(T^{X^i}_1<T^{X^i}_0,T^{X^j}_1<T^{X^j}_0
\bigr)
\\
\eqntext{\forall r\in(0,r_0], {\varepsilon}\in\bigl(0, {
\varepsilon}_0(r) \bigr].}
\end{eqnarray}
The first term on the right-hand side of (\ref{ineqsepprob1}) can be
written as
%
%
%e4.30 #&#
\begin{eqnarray}
\label{ineqsepprobm1} %
&&\sum_{i=1}^{\lfloor r{\varepsilon}^{-1}\rfloor}
\P_{\varepsilon
} \bigl(\widehat{G}{}^i \bigl(\widehat{
\tau}{}^i(s_i+r)\wedge(s_i+r) \bigr)\cap\bigl
\{T^{X^i}_1<T^{X^i}_0 \bigr\} \bigr)
\nonumber\\[-8pt]\\[-8pt]\nonumber
&&\qquad = \sum_{i=1}^{\lfloor r{\varepsilon}^{-1}\rfloor}\psi(\1){
\varepsilon}\cdot\Q_{\varepsilon}^i \bigl(\widehat{G}{}^i
\bigl(\widehat{\tau}{}^i(s_i+r)\wedge(s_i+r)
\bigr) \bigr)
\end{eqnarray}
by the definition of $\Q^i_{\varepsilon}$ in (\ref{defQ}).
By inclusion--exclusion, we have
%
%
%e4.31 #&#
%e4.32 #&#
\begin{eqnarray}\label{ineqQG1}
&&\Q_{\varepsilon}^i \bigl(\widehat{G}{}^i \bigl(
\widehat{\tau}{}^i(s_i+r)\wedge(s_i+r) \bigr)
\bigr)\nonumber
\\
&&\qquad \geq\mathbb Q_{\varepsilon}^i
\biggl(X^i_s(\1)\geq\frac{(s-s_i)^{\eta
}}{4}, \forall s\in
\bigl[s_i,\widehat{\tau}{}^i(s_i+r)\wedge
(s_i+r) \bigr] \biggr)\nonumber
\\
&&\quad\qquad{}-\mathbb Q^i_{\varepsilon} \biggl(\exists s\in\bigl(s_i,
\widehat{\tau}{}^i(s_i+r)\wedge(s_i+r) \bigr],\sum_{j\in\mathcal J_{\beta'}^i(s)}Y^j(\1)_{s}
\\
&&\hspace*{92pt}>K^*
\bigl[(s-s_i)^{\kappa
_1-\wp}+{\varepsilon}^{\kappa_2}
\cdot(s-s_i)^{\kappa_3-\wp}\bigr] \biggr)\nonumber
\\
\eqntext{\displaystyle \forall i\in\mathbb N, {\varepsilon}\in\biggl(0,\frac {1}{8\psi(\1)}
\wedge1 \biggr].}
\end{eqnarray}
Recall that $\tau^{i,(1)}\leq\tau^i$ and $X^i_{s_i}(\1)=\psi(\1
){\varepsilon}>0$. Hence,
by the definition of \mbox{$\widehat{\tau}{}^i(s_i+r)$},
%
%
%e4.33 #&#
%e4.34 #&#
\begin{eqnarray} \label{ineqQG2}
&&\mathbb Q_{\varepsilon}^i \biggl(X^i_s(
\1)\geq\frac{(s-s_i)^{\eta
}}{4}, \forall s\in\bigl[s_i,\widehat{
\tau}{}^i(s_i+r)\wedge(s_i+r) \bigr] \biggr)=1
\nonumber\\[-8pt]\\[-8pt]
\eqntext{\displaystyle\forall i\in\mathbb N, {\varepsilon}\in\biggl(0,\frac
{1}{8\psi(\1)}
\wedge1 \biggr].}
\end{eqnarray}
For
$r\in(0,\delta_1], i\in\mathbb N, \mbox{and }{\varepsilon}\in
(0,[8\psi(\1)]^{-1}\wedge1 ]$,
the second probability in (\ref{ineqQG1}) can be bounded as
%
%
%e4.35 #&#
\begin{eqnarray}\label{ineqQG3}
&&\mathbb Q^i_{\varepsilon} \biggl(\exists s\in\bigl(s_i,
\widehat{\tau}{}^i(s_i+r)\wedge(s_i+r) \bigr],\sum
_{j\in \mathcal J_{\beta'}^i(s)}Y^j_s(\1)\nonumber
\\
&&\hspace*{42pt}
>K^*\bigl[(s-s_i)^{\kappa_1-\wp}+{
\varepsilon}^{\kappa
_2}\cdot(s-s_i)^{\kappa_3-\wp}\bigr] \biggr)
\nonumber
\\
&&\qquad \leq \Q^i_{\varepsilon} \biggl(\exists s\in\bigl(s_i,
\tau^i\wedge(s_i+r)\bigr],
\sum_{j\in\mathcal J_{\beta'}^i(s\wedge\tau
^i\wedge\sigma_{\beta}^{X^i})}Y^j_s(
\1)^{ \tau^i\wedge\sigma_{\beta
}^{X^i}\wedge\sigma_\beta^{Y^j}}\nonumber
\\
&&\hspace*{122pt}
>K^*\bigl[(s-s_i)^{\kappa_1-\wp}+{
\varepsilon}^{\kappa
_2}\cdot(s-s_i)^{\kappa_3-\wp}\bigr] \biggr)
\\
&&\qquad \leq \mathbb Q^i_{\varepsilon} \biggl(\exists s\in(s_i,s_i+
\delta_1], \sum_{j\in\mathcal J_{\beta'}^i(s\wedge\tau^i\wedge\sigma
_{\beta}^{X^i})}Y^j_s(\1)^{ \tau^i\wedge\sigma_{\beta}^{X^i}\wedge\sigma
_\beta^{Y^j}}
\nonumber
\\
&&\hspace*{96pt}  >K^*\bigl[(s-s_i)^{\kappa_1-\wp}+{\varepsilon}^{\kappa_2}
\cdot(s-s_i)^{\kappa_3-\wp}\bigr] \biggr)\nonumber
\\
&&\quad\qquad{} +\Q^i_{\varepsilon}
\bigl(\tau^i\leq s_i+\delta_1\bigr) \leq2\rho.\nonumber
\end{eqnarray}
Here, the first inequality follows since for $s\in(s_i,\widehat{\tau
}{}^i(s_i+r)\wedge(s_i+r)]$, we have $s\leq\tau^i\wedge\sigma_{\beta
}^{X^i}$ and
\begin{eqnarray*}
&& j\in\mathcal J_{\beta'}^i(s) \quad\Longrightarrow\quad
j\in\mathcal J_{\beta'}^i(s_i+r)\quad \Longrightarrow\quad
t_j\in(s_i,s_i+r]
\\
&& \quad \Longrightarrow\quad
\widehat{\tau}{}^i(s_i+r)\leq\sigma_\beta^{Y^j}\quad\Longrightarrow\quad
s\leq\sigma_\beta^{Y^j},
\end{eqnarray*}
where the third implication follows from the definition of $\widehat
{\tau}{}^i(s_i+r)$ in (\ref{eqtauhat}). The first term in the second
inequality follows by considering the scenario $\tau^i>s_i+\delta_1$
and using $r\in(0,\delta_1]$, and
the last inequality follows from (\ref{choicedelta1-1}) and (\ref
{choicedelta1-2}). Applying (\ref{ineqQG2}) and (\ref{ineqQG3}) to
(\ref{ineqQG1}), we get
%
%
%e4.36 #&#
%e4.37 #&#
\begin{eqnarray}\label{ineqGi*}
\Q^i_{\varepsilon} \bigl(\widehat{G}{}^i\bigl(\widehat{
\tau}{}^i(s_i+r)\wedge(s_i+r)\bigr) \bigr)
\geq1-2\rho
\nonumber\\[-8pt]\\[-8pt]
\eqntext{\displaystyle \forall r\in(0,\delta_1\wedge r_0], i
\in\mathbb N, {\varepsilon}\in\bigl(0, {\varepsilon}_0(r)\bigr].}
\end{eqnarray}
From (\ref{ineqsepprobm1}) and the last inequality, we have shown that
%
%e4.38 #&#
\begin{eqnarray}
\sum_{i=1}^{\lfloor r{\varepsilon}^{-1}\rfloor}\P_{\varepsilon
} \bigl(
\widehat{G}_i \bigl(\widehat{\tau}{}^i(s_i+r)
\wedge(s_i+r) \bigr)\cap\bigl\{T^{X^i}_1<T^{X^i}_0
\bigr\} \bigr)\geq\psi(\1) (r-{\varepsilon} ) (1-2\rho)\nonumber
\\
\eqntext{\displaystyle \forall r\in(0,\delta_1\wedge r_0],
{\varepsilon}\in\bigl(0, {\varepsilon}_0(r) \bigr].}
\end{eqnarray}
[Recall that ${\varepsilon}_0(r )\leq r$.]
The second term on the right-hand side of (\ref{ineqsepprob1}) is
relatively easy to bound. Indeed, by using the independence between the
clusters $X^i$ and Lemma~\ref{lemQi},
%
%e4.39 #&#
\begin{eqnarray}
&& \sum_{i=1}^{\lfloor r{\varepsilon}^{-1}\rfloor}\sum
_{\stackrel
{\scriptstyle j\dvtx j\neq i}{1\leq j\leq\lfloor r{\varepsilon
}^{-1}\rfloor}}\P_{\varepsilon} \bigl(T^{X^i}_1<T^{X^i}_0,T^{X^j}_1<T^{X^j}_0
\bigr)\nonumber
\\
&&\qquad  =\sum_{i=1}^{\lfloor r{\varepsilon}^{-1}\rfloor}\sum
_{\stackrel{\scriptstyle j\dvtx j\neq i}{1\leq j\leq\lfloor
r{\varepsilon
}^{-1}\rfloor}}\P_{\varepsilon} \bigl(T^{X^i}_1<T^{X^i}_0
\bigr)\P_{\varepsilon} \bigl(T^{X^j}_1<T^{X^j}_0
\bigr) \leq\psi(\1)^2r^2\nonumber
\\
\eqntext{\displaystyle\forall r\in(0,1], {\varepsilon}\in\biggl(0,\frac {1}{8\psi(\1)}\wedge1
\biggr].}
\end{eqnarray}
Recalling (\ref{ineqsepprob1}) and using the last two displays, we
have the following bound for the first term on the right-hand side of
(\ref{ineqsepprob}):
%
%
%e4.40 #&#
%e4.41 #&#
\begin{eqnarray}
\label{inequsep-1}
&&\P_{\varepsilon} \Biggl(\bigcup
_{i=1}^{\lfloor r{\varepsilon
}^{-1}\rfloor}\widehat{G}_i \bigl(
\widehat{\tau}{}^i(s_i+t)\wedge(s_i+r) \bigr)
\cap\bigl\{ T^{X^i}_1<T^{X^i}_0 \bigr
\} \Biggr)\nonumber
\\
&&\qquad \geq\psi(\1) (r-{\varepsilon}) (1-2\rho)-\psi(\1)^2r^2
\\
\eqntext{\displaystyle \forall r\in(0,\delta_1\wedge r_0], {\varepsilon}
\in\bigl(0, {\varepsilon}_0(r) \bigr].}
\end{eqnarray}

Next, we consider the second probability on the right-hand side of
(\ref{ineqsepprob}). By the definition of $\Gamma^i(r )$ in (\ref
{eqRit}) and the general inclusion
$(A_1\cap A_2\cap A_3)^\complementtt\subseteq(A_1^\complementtt\cap
A_2\cap A_3)\cup A_2^\complementtt\cup A_3^\complementtt$,
we have
%
%
%e4.42 #&#
\begin{eqnarray}
\label{eqRitprop} %
\Gamma^i(r )^\complementtt &\subseteq&
\biggl( \biggl\{\mc P^{X^i}_\beta(s_i+r)\cap
\biggl( \bigcup_{j\dvtx t_j\leq s_i}\supp\bigl(Y^j\bigr)
\biggr)\neq\varnothing\biggr\}\nonumber
\\
&&{} \cap\bigcap_{j\dvtx t_j\leq s_i}
\bigl\{\sigma_\beta^{Y^j}>t_j+3r \bigr\}
\cap\bigl\{\sigma_{\beta}^{X^i}>s_i+2r
\bigr\} \biggr)
\\
&&{} \cup\biggl(\bigcup_{j\dvtx t_j\leq s_i+r} \bigl\{
\sigma_\beta^{Y^j}\leq t_j+3r \bigr\} \biggr)\cup
\bigl\{\sigma_{\beta}^{X^i}\leq s_i+2r \bigr\},\nonumber
\end{eqnarray}
where we note that the indices $j$ in
$\bigcap_{j\dvtx t_j\leq s_i} \{\sigma_\beta^{Y^j}>t_j+3r \}$
range only over $j\in\mathbb N$ with $t_j\leq s_i$.
Hence,
%
%
%e4.43 #&#
\begin{eqnarray}\label{ineqsepprob-2}
\qquad &&\P_{\varepsilon} \Biggl(\bigcup_{i=1}^{\lfloor r{\varepsilon
}^{-1}\rfloor}
\Gamma^i(r )^\complementtt\cap\bigl\{ T^{X^i}_1<T^{X^i}_0
\bigr\} \Biggr)
\nonumber
\\
&&\qquad \leq \P_{\varepsilon} \Biggl(\bigcup
_{i=1}^{\lfloor r{\varepsilon
}^{-1}\rfloor} \biggl( \biggl\{\mc P_\beta^{X^i}(s_i+r)
\cap\biggl(\bigcup_{j\dvtx t_j\leq s_i}\supp\bigl(Y^j
\bigr) \biggr)\neq\varnothing\biggr\}
\nonumber\\[-8pt]\\[-8pt]\nonumber
&&\quad\qquad{} \cap\bigcap_{j\dvtx t_j\leq s_i} \bigl\{
\sigma_\beta^{Y^j}>t_j+3r \bigr\}\cap\bigl\{
\sigma_{\beta}^{X^i}>s_i+2r \bigr\}\cap\bigl
\{T^{X^i}_1<T^{X^i}_0 \bigr\} \biggr)
\Biggr)
\\
&&\quad\qquad{}+\P_{\varepsilon} \Biggl(\bigcup_{j=1}^{\lfloor2r{\varepsilon
}^{-1}\rfloor+1}
\bigl\{\sigma_\beta^{Y^j}\leq t_j+3r \bigr\}
\Biggr)+\P_{\varepsilon} \Biggl(\bigcup_{i=1}^{\lfloor r{\varepsilon
}^{-1}\rfloor}
\bigl\{\sigma_{\beta}^{X^i}\leq s_i+2r \bigr\}
\Biggr),\nonumber
\end{eqnarray}
where we have the second probability in the foregoing inequality since
\[
t_j\leq s_{\lfloor r{\varepsilon}^{-1}\rfloor}+r\quad\Longrightarrow\quad
t_j\leq2r\quad\Longrightarrow\quad
j\leq\bigl\lfloor2r{\varepsilon}^{-1}\bigr\rfloor+1.
\]
Resorting to the conditional probability measures $\Q^i_{\varepsilon
}$, we see that the first probability in (\ref{ineqsepprob-2}) can be
bounded as
\begin{eqnarray*}
&&\P_{\varepsilon} \Biggl(\bigcup_{i=1}^{\lfloor r{\varepsilon
}^{-1}\rfloor}
\biggl( \biggl\{\mc P^{X^i}_\beta(s_i+r)\cap
\biggl(\bigcup_{j\dvtx t_j\leq s_i}\supp\bigl(Y^j\bigr)
\biggr)\neq\varnothing\biggr\}
\\
&&\hspace*{48pt}{} \cap\bigcap_{j\dvtx t_j\leq s_i}
\bigl\{\sigma_\beta^{Y^j}>t_j+3r \bigr\}
\cap\bigl\{\sigma^{X^i}_{\beta}>s_i+2r
\bigr\}\cap\bigl\{T^{X^i}_1<T^{X^i}_0
\bigr\} \biggr) \Biggr)
\\
&&\qquad \leq \sum_{i=1}^{\lfloor r{\varepsilon}^{-1}\rfloor}\psi(\1 ){
\varepsilon}\mathbb Q^i_{\varepsilon} \biggl( \biggl\{\mc
P_\beta^{X^i}(s_i+r)\cap\biggl(\bigcup
_{j\dvtx t_j\leq s_i}\supp\bigl(Y^j\bigr) \biggr)\neq\varnothing
\biggr\}
\\
&&\hspace*{72pt}\quad\qquad{} \cap\bigcap_{j\dvtx t_j\leq s_i} \bigl\{
\sigma_\beta^{Y^j}>t_j+3r \bigr\} \cap\bigl\{
\sigma^{X^i}_{\beta}>s_i+2r \bigr\} \biggr)
\\
&&\qquad \leq\sum_{i=1}^{\lfloor r{\varepsilon}^{-1}\rfloor}\psi(\1 ){
\varepsilon}C^1_{\supp}r^{1/6}\leq\psi(
\1)C^1_{\supp}r^{7/6}\qquad\forall r
\in(0,r_0], {\varepsilon}\in\bigl(0,{\varepsilon}_0(r)\bigr],
\end{eqnarray*}
where the next to the last inequality follows from Proposition~\ref{propoutclus-si} and the constant $C^1_{\supp}\in(0,\infty)$ is
independent of $r\in(0,r_0]$ and ${\varepsilon}\in(0,r_0]$. (Here,
we use the choice $\beta\in[\frac{1}{3},\frac{1}{2} )$
to apply this proposition.)
By Proposition~\ref{propsupppre}, the second probability in (\ref
{ineqsepprob-2}) can be bounded as [recall ${\varepsilon}_0( r)\leq r$]
%
%
%e4.44 #&#
\begin{eqnarray}\label{eqsuppbdd-1}
&& \P_{\varepsilon} \Biggl(\bigcup_{j=1}^{\lfloor2r{\varepsilon
}^{-1}\rfloor+1}\bigl\{\sigma_\beta^{Y^j}\leq t_j+3r \bigr\}
\Biggr)\nonumber
\\
&&\qquad \leq C^0_{\supp}\bigl(2r{\varepsilon}^{-1}+1
\bigr)\cdot3{\varepsilon}r
\\
&&\qquad \leq 9C_{\supp}^0r^2\qquad\forall r
\in(0,r_0], {\varepsilon}\in\bigl(0,{\varepsilon}_0(r
)\bigr],\nonumber
\end{eqnarray}
where $C_{\supp}^0$ is a constant independent of $r\in(0,r_0]$ and
${\varepsilon}\in(0,{\varepsilon}_0(r)]$. Similarly,
%
%
%e4.45 #&#
\begin{equation}
\qquad\quad \P_{\varepsilon} \Biggl(\bigcup_{i=1}^{\lfloor r{\varepsilon
}^{-1}\rfloor}
\bigl\{\sigma^{X^i}_{\beta}\leq s_i+2r \bigr\}
\Biggr)\leq2C^0_{\supp} r^2\qquad\forall r
\in(0,r_0], {\varepsilon}\in\bigl(0,{\varepsilon}_0(r)\bigr].\label{eqsuppbdd-2}
\end{equation}
From (\ref{ineqsepprob-2}) and the last three displays, we have shown
that the second probability in (\ref{ineqsepprob}) satisfies the bound
%
%
%e4.46 #&#
%e4.47 #&#
\begin{eqnarray}
\label{inequsep-2} %
&&\P_{\varepsilon} \Biggl(\bigcup
_{i=1}^{\lfloor r{\varepsilon
}^{-1}\rfloor}\Gamma^i(r )^\complementtt
\cap\bigl\{ T^{X^i}_1<T^{X^i}_0 \bigr
\} \Biggr)\leq11C_{\supp}^0 r^2+\psi(\1
)C^1_{\supp}r^{7/6}
\nonumber\\[-8pt]\\[-8pt]
\eqntext{\displaystyle \forall r\in(0,r_0], {\varepsilon}\in\bigl(0,{\varepsilon
}_0(r )\bigr].}
\end{eqnarray}

It remains to bound the last probability on the right-hand side of
(\ref{ineqsepprob}). Recall the number $\delta_1$ chosen for (\ref
{choicedelta1-1}). Similar to the derivation of (\ref
{ineqsepprob-2}), we have
\begin{eqnarray}\label{inequsep-3}
&&\P_{\varepsilon} \Biggl(\bigcup_{i=1}^{\lfloor r{\varepsilon
}^{-1}\rfloor}
\bigl\{ \widehat{\tau}{}^i(s_i+r)<s_i+r \bigr\}
\cap\bigl\{ T^{X^i}_1<T^{X^i}_0 \bigr
\} \Biggr)\nonumber
\\
&&\qquad \leq\P_{\varepsilon} \Biggl(\bigcup
_{i=1}^{\lfloor
r{\varepsilon}^{-1}\rfloor} \bigl\{\sigma^{X^i}_\beta
\leq s_i+r \bigr\} \Biggr)+\P_{\varepsilon} \Biggl(\bigcup
_{i=1}^{\lfloor
2r{\varepsilon}^{-1}\rfloor+1} \bigl\{\sigma^{Y^i}_\beta
\leq t_i+r \bigr\} \Biggr)
\nonumber\\[-8pt]\\[-8pt]\nonumber
&&\quad\qquad{} +\sum_{i=1}^{\lfloor r{\varepsilon}^{-1}\rfloor}\psi(\1){
\varepsilon}\mathbb Q^i_{\varepsilon} \bigl(\tau^i<s_i+r
\bigr)\nonumber
\\
&&\qquad \leq11C_{\supp}^0 r^2+\psi(\1)r
\rho\qquad\forall r\in(0, \delta_1\wedge r_0 ], {
\varepsilon}\in\bigl(0,{\varepsilon}_0(r) \bigr],\nonumber
\end{eqnarray}
where we use (\ref{eqsuppbdd-1}) and (\ref{eqsuppbdd-2}) in the
last inequality.

We apply the three bounds (\ref{inequsep-1}), (\ref{inequsep-2})
and (\ref{inequsep-3}) to (\ref{ineqsepprob}). This shows that
for any $\rho\in(0,\frac{1}{2})$, there exist $\delta_1>0$ such
that for any $r\in(0, \delta_1\wedge r_0]$ and ${\varepsilon}\in
(0, {\varepsilon}_0(r ) ]$ [note that ${\varepsilon}_0(r
)\leq r\wedge1$],
\begin{eqnarray*}
\P_{\varepsilon} \bigl(S(r ) \bigr)&\geq& \bigl[\psi(\1 ) (r-{\varepsilon})
(1-2\rho)-\psi(\1)^2r^2 \bigr]- \bigl(11C_{\supp}
^0 r^2+\psi(\1)C^1_{\supp}r^{7/6}
\bigr)
\\
&&{}- \bigl(11C_{\supp}^0 r^2+\psi(\1)r\rho\bigr)
\\
&=& r \bigl[\psi(\1) (1-3\rho)- \bigl(\psi(\1)^2+22C_{\supp}^0
\bigr)r-\psi(\1)C^1_{\supp}r^{1/6} \bigr]
\\
&&{}  -\psi(\1){\varepsilon}(1-2\rho).
\end{eqnarray*}

Finally, to attain the uniform lower bound (\ref{inequnifsep}), we
choose $\rho\in(0,\frac{1}{2})$ and $r_1\in(0, \delta_1\wedge
r_0]$ such that
\[
\psi(\1) (1-3\rho)- \bigl(\psi(\1)^2+22C_{\supp}^0
\bigr)r-\psi(\1 )C^1_{\supp}r^{1/6}\geq
\frac{\psi(\1)}{2}\qquad\forall r\in(0,r_1 ],
\]
and then ${\varepsilon}_1(r )\in(0, {\varepsilon}_0(r )
]$ such that
\[
\psi(\1){\varepsilon}_1(r ) (1-2\rho)\leq\frac{\psi(\1)r}{4}.
\]
By the last three displays, we obtain
\[
\P_{\varepsilon} \bigl(S(r ) \bigr)\geq\frac{\psi(\1)r }{4}\qquad
\forall{
\varepsilon}\in\bigl(0,{\varepsilon}_1(r )\bigr], r\in(0,r_1],
\]
and hence (\ref{inequnifsep}) follows. The proof is complete.
\end{pf}

We use Lemma~\ref{lemus3} to give the proof for a more precise
version of
our main theorem, namely Theorem~\ref{teoSPDE}, in Theorem~\ref
{teomains} below.

%
%
%th4.4 #&#
\begin{teo}[(Separation of limiting solutions)]\label{teomains}
Let $({\varepsilon}_n)\subseteq(0,\break [8\psi(\1)]^{-1}\wedge1
]$ with ${\varepsilon}_n\searrow0$ be such that the sequence of laws
of $ ((X,Y),\P_{{\varepsilon}_n} )$ converges to the law of
$ ((X,Y),\P_0 )$ of a pair of solutions to the SPDE (\ref
{eqmainSPDE}) in the space of probability measures on the product
space $D (\R_+,\C_\rap(\R) )\times D (\R_+,\C_\rap
(\R) )$ (cf. Proposition~\ref{propwc}). Then we have
\[
\P_0 \biggl(\sup_{0\leq s\leq2r_1}\llVert X-Y\rrVert
_\rap\geq\frac{\Delta(r_1)}{2} \biggr)\geq\inf_{{\varepsilon}\in
(0,{\varepsilon}_1(r_1)]}
\P_{\varepsilon} \bigl(S(r_1) \bigr)>0,
\]
where $\Delta(r_1)>0$ is chosen in Lemma~\ref{lemus1} and
$r_1,{\varepsilon}_1(r_1)\in(0,1]$ are chosen in Lemma~\ref{lemus3}.
\end{teo}

\begin{pf}
By Skorokhod's representation theorem, we may
take $ (X^{({\varepsilon}_n)},\break Y^{({\varepsilon}_n)} )$ to
be copies of the ${\varepsilon}_n$-approximating solutions which live
on the same probability space, and assume that $ (X^{({\varepsilon
}_n)},Y^{({\varepsilon}_n)} )$ converges almost surely to $
(X^{(0)},Y^{(0)} )$ in the product (metric) space $D (\R
_+,\C_\rap(\R) )\times D (\R_+,\C_\rap(\R) )$.

It follows from Lemmas~\ref{lemus1}~and~\ref{lemus3} that
\[
\inf_{n\dvtx{\varepsilon}_n\leq{\varepsilon}_1(r_1)}\P\Bigl(\sup
_{0\leq s\leq2r_1}\bigl\llVert
X^{({\varepsilon}_n)}_s-Y^{({\varepsilon
}_n)}_s\bigr\rrVert
_\rap\geq\Delta(r_1) \Bigr)\geq\inf_{{\varepsilon}\in(0,{\varepsilon}_1(r_1)]}
\P_{\varepsilon} \bigl(S(r_1) \bigr)>0.
\]
Hence, by Fatou's lemma, we get
\begin{eqnarray}\label{ineqsepprobMAIN}
0&<&\inf_{{\varepsilon}\in(0,{\varepsilon}_1(r_1)]}\P_{\varepsilon
} \bigl(S(r_1)\bigr)\nonumber
\\
&\leq& \limsup_{n\to\infty}\P\Bigl(\sup_{0\leq
s\leq2r_1}\bigl\llVert X^{({\varepsilon}_n)}_s-Y^{({\varepsilon
}_n)}_s\bigr\rrVert_\rap\geq\Delta(r_1) \Bigr)
\nonumber\\[-8pt]\\[-8pt]\nonumber
&\leq&\P\Bigl(\limsup_{n\to\infty} \Bigl\{\sup_{0\leq s\leq
2r_1}
\bigl\llVert X^{({\varepsilon}_n)}_s-Y^{({\varepsilon}_n)}_s
\bigr\rrVert_\rap\geq\Delta(r_1) \Bigr\} \Bigr)
\nonumber
\\
&\leq&\P\biggl(\sup_{0\leq s\leq2r_1}\bigl\llVert X^{(0)}_s-Y^{(0)}_s
\bigr\rrVert_\rap\geq\frac{\Delta(r_1)}{2} \biggr),\nonumber
\end{eqnarray}
where the last inequality follows from the convergence
\[
X^{({\varepsilon}_n)}\xrightarroww{n\to\infty} {\mathrm
{a.s.}}X^{(0)}\quad
\mbox{and}\quad Y^{({\varepsilon}_n)}\xrightarroww{n\to\infty} {\mathrm
{a.s.}}Y^{(0)}
\]
in the Skorokhod space $D (\R_+,\C_\rap(\R) )$, the
continuity of $X^{(0)}$ and $Y^{(0)}$, and Proposition 3.6.5(a) of
\cite{EK}. The proof is complete.
\end{pf}

%s5 #&#
\section{Proof of Proposition~\texorpdfstring{\protect\ref{propwc}}{2.3}}\label{secpropwc}
Many arguments in this section can be modified from the proofs in
Section~6 in \cite{MMP} because of the apparent similarity of the
involved stochastic processes, and so we only give sketches whenever
necessary. Readers interested in a complete proof of Proposition~\ref
{propwc} may see Section~3.9 of~\cite{C}. Some connections between
limit theorems for $\C_\rap(\R)$-valued processes and limit theorems
for processes taking values in the space of real-valued continuous
functions over $\R$ can be found in Section~3.11 of \cite{C}.

Throughout this section, we fix a sequence $({\varepsilon}_n)\subseteq
(0,1]$ with ${\varepsilon}_n\searrow0$ and assume that the
${\varepsilon}_n$-approximating solutions live on the same probability
space. To save notation, we write
$ \{ (X^{(n)},Y^{(n)} );n\in\mathbb N \}$
for this approximating sequence and denote by $\P$ the underlying
probability measure.
We will begin with the $C$-tightness of the sequence of joint laws of
$ \{ (X^{(n)},Y^{(n)} ) \}$ in $D (\R_+,\C
_\rap(\R) )\times D (\R_+,\C_\rap(\R) )$, where
$D (\R_+,\C_\rap(\R) )$ is equipped with Skorokhod's
$J_1$-topology. Here, $C$-tightness means not only tightness but also
the property that the limiting object of any convergent subsequence is
a continuous process. We will only discuss the $C$-tightness of the
sequence of laws of $ \{X^{(n)} \}$ in $D (\R_+,\C
_\rap(\R) )$, and the argument for $ \{Y^{(n)} \}$
follows similarly. Later on in Lemma~\ref{lemwconvXY}, we will prove
that the limit of any convergent subsequence of laws of $ \{
(X^{(n)},Y^{(n)} ) \}$ is the law of a pair of solutions to
the SPDE (\ref{eqmainSPDE}) with respect to the same space--time white noise.

Consider our first objective that the sequence of laws of $ \{
X^{(n)} \}$ is $C$-tight as probability measures on $D (\R
_+,\C_\rap(\R) )$.
The proof uses the mild forms of $ \{X^{(n)} \}$ stated below. Let
\[
p_s(x)\,dx\equiv\cases{ \displaystyle\frac{1}{\sqrt{2\pi s}}\exp
\biggl(-\frac
{x^2}{2s} \biggr)\,dx,&\quad$s\in(0,\infty)$,
\vspace*{3pt}\cr
\delta_0(dx), &\quad$s=0$,
\vspace*{3pt}\cr
0,&\quad$s\in(-\infty,0)$.}
\]
Recall the random measure $A^{X^{(n)}}$ [cf. (\ref{eqAX})] associated
with $X^{(n)}$ which is
contributed by the initial masses of its immigrants, and write
%
%
%e5.1 #&#
%e5.2 #&#
\begin{eqnarray}\label{eqMXn}
M^{X^{(n)}}_t(\phi)
&\equiv&X^{(n)}_t(\phi)-\int_0^t
X^{(n)}_s \biggl(\frac{\Delta}{2}\phi\biggr)\,ds-\int
_{(0,t]}\int_\R\phi(y)\,dA^{X^{(n)}}(y,s),
\nonumber\\[-8pt]\\[-8pt]
\eqntext{\phi\in\C_c^\infty(\R),}
\end{eqnarray}
for the martingale measure of $X^{(n)}$.

By summing up the mild forms of the immigrant processes for $X^{(n)}$
which are solutions to the SPDE (\ref{eqmainSPDE0}) and have initial
conditions taking the form $\psi(\1)J^a_{\varepsilon}( \cdot)$ for
$J^a_{\varepsilon}( \cdot)$ defined by (\ref{eqJxvep}) (see Theorem
2.1 of \cite{SSPDE}), we deduce that the mild form of $X^{(n)}$ is
given by
%
%
%e5.3 #&#
\begin{equation}
\label{eqmildXn} X^{(n)}(x,t)=p\star A^{X^{(n)}}(x,t)+p\star
M^{X^{(n)}}(x,t),\qquad(x,t)\in\R\times\R_+.
\end{equation}
Here, the convolutions on the right-hand side are given by
%
%
%e5.4 #&#
%e5.5 #&#
\begin{eqnarray}
p\star A^{X^{(n)}}(x,t)&=&\int_{(0,t]} \int
_\R p_{t-s}(x-y)\,dA^{X^{(n)}}(y,s)
\nonumber\\[-8pt]\label{eqpstarAn} \\[-8pt]\nonumber
&=&\psi(\1)\sum_{i\dvtx0<s_i\leq t}\int_{\R}
p_{t-s_i}(x-y)J^{x_i}_{{\varepsilon}_n}(y)\,dy,
\\
p\star M^{X^{(n)}}(x,t)&=&\int_0^t\!\int
_\R p_{t-s}(x-y)\,dM^{X^{(n)}}(y,s)
\nonumber\\[-8pt]\label{eqpstarMn} \\[-8pt]\nonumber
&=&\int_0^t\!\int_\R
p_{t-s}(x-y)X^{(n)}(y,s)^{1/2}\,dW(y,s).
\end{eqnarray}
More precisely, in $p\star A^{X^{(n)}}$, we read
$p_0(x-y)\,dy=\delta_0(x-dy)=\delta_x(dy)$,
and hence
%
%
%e5.6 #&#
\begin{equation}
\label{eqmildXin0} \int_\R p_0(x-y)J^{x_i}_{{\varepsilon}_n}(y)\,dy
\equiv J^{x_i}_{{\varepsilon}_n}(x).
\end{equation}

The mild form (\ref{eqmildXn}) implies the $C$-tightness of the
sequence of laws of $ \{X^{(n)} \}$ in $D (\R_+,\C
_\rap(\R) )$, provided that the sequences of laws of $ \{
p\star A^{X^{(n)}};n\in\mathbb N \}$ and $ \{p\star
M^{X^{(n)}};n\in\mathbb N \}$ are both $C$-tight as probability
measures on the same space.

%
%
%le5.1 #&#
\begin{lem}\label{lempAntight}
The sequences of laws of $ \{p\star A^{X^{(n)}};n\in\mathbb N \}$
is $C$-tight and converges in probability in $D (\R_+,\C_\rap(\R
) )$ to the deterministic process
$ (\int_0^t \int_\R p_{t-s}(\cdot-y)\psi(y)\,dy\,ds )_{t\in\R_+}$.
\end{lem}

%
%
%le5.2 #&#
\begin{lem}\label{lemmomXY}
For any $q\in[1,\infty)$ and $\lambda,T\in(0,\infty)$, there
exists a constant $\check{C}\in(0,\infty)$ depending only on $(\psi, q, \lambda,T)$ such that
\[
\sup_{n\in\mathbb N}\sup_{0\leq t\leq T}\sup
_{x\in\R}e^{\lambda
\llvert x\rrvert }\E\bigl[X^{(n)}(x,t)^q+Y^{(n)}(x,t)^q
\bigr]\leq\check{C}.
\]
\end{lem}

%
%
%le5.3 #&#
\begin{lem}\label{lempMntight}
For some universal constants $q\in(0,\infty)$ and $\gamma\in
(2,\infty)$, the following inequality holds for any $\lambda,T\in
(0,\infty)$:
%
%
%e5.7 #&#
%e5.8 #&#
\begin{eqnarray}
\label{eqmodineq} &&\sup_{n\in\mathbb N} \E\bigl[\bigl\llvert
p\star M^{X^{(n)}}\bigl(x',t'\bigr)-p\star
M^{X^{(n)}}(x,t)\bigr\rrvert^q \bigr]\nonumber
\\
&&\qquad \leq\check{C} \bigl(\bigl\llvert x'-x\bigr\rrvert
^{2\gamma}+\bigl\llvert t'-t\bigr\rrvert^{\gamma
}\bigr)e^{-\lambda\llvert x\rrvert }
\\
\eqntext{\displaystyle \forall t,t'\in[0,T], \bigl\llvert
x-x'\bigr\rrvert\leq1.}
\end{eqnarray}
Here, the constants $\check{C}$ are as in Lemma~\ref{lemmomXY} and
are enlarged if necessary.
Moreover, the sequence of laws of $ \{p\star M^{X^{(n)}} \}$ is
tight as probability measures on $C (\R_+,\C_\rap(\R) )$.
\end{lem}

The proofs of Lemmas~\ref{lempAntight}, \ref{lemmomXY} and~\ref
{lempMntight} can be obtained by arguments similar to the
proofs of Lemmas~6.6, 6.1 and~6.7 in \cite{MMP}, respectively. In this
direction, the proofs
of Lemmas~\ref{lempAntight}~and~\ref{lemmomXY} use the
particular form of the distribution~(\ref{eqxiyilaw}) of $x_i$ and
$y_i$ which is dominated by a constant multiple of Lebesgue measure
over a compact interval, as well as the fact that in our case,
immigrants can land throughout time. The latter does not create
additional difficulties since \mbox{$C$-}tightness of $\C_\rap(\R)$-valued
processes and the bound in Lemma~\ref{lemmomXY} only concern
distributional properties of the corresponding processes over compact intervals.
In addition, for the proof of Lemma~\ref{lempMntight}, we need the
moment bound in Lemma~\ref{lemmomXY} for its first assertion, and the
second assertion follows from (\ref{eqmodineq}) and Lemma 6.4 of
\cite{MMP}.

By Lemmas~\ref{lempAntight}~and~\ref{lempMntight}, the
sequence of laws of $ \{X^{(n)} \}$ is $C$-tight as
probability measures on $D (\R_+,\C_\rap(\R) )$, thanks to
(\ref{eqmildXn}). By similar arguments, the same is true for the
sequence of laws of $ \{Y^{(n)} \}$.

%
%
%le5.4 #&#
\begin{lem}\label{lemwconvXY}
Suppose that, by taking a subsequence if necessary, we have
%
%
%e5.9 #&#
\begin{equation}
\label{eqjointconv} \bigl(X^{(n)},Y^{(n)} \bigr)\xrightarroww{n
\to\infty} {{\rm(d)}} \bigl(X^{(0)},Y^{(0)} \bigr)
\end{equation}
for some continuous $\C_\rap(\R)$-valued processes $X^{(0)}$ and
$Y^{(0)}$. Then $X^{(0)}$ and $Y^{(0)}$ solve the SPDE (\ref
{eqmainSPDE}) with respect to the same space--time white noise.
\end{lem}

\begin{pf*}{Sketch of proof}
The argument in the proof of Proposition 2.2 of \cite{MMP} (in
Section~6 there) still applies and implies that both $X^{(0)}$ and
$Y^{(0)}$ are solutions to the SPDE (\ref{eqmainSPDE}).
Here, the readers may use as supporting facts a reinforcement of the
convergence in (\ref{eqconvip}) to almost-sure convergence along
$({\varepsilon}_n)$ (by the strong law of large numbers), and the
moment bound in Lemma~\ref{lemmomXY}.

It remains to show that $X^{(0)}$ and $Y^{(0)}$ can be subject to the
same space--time white noise, and moreover, all of these random objects
obey their defining properties with respect to the same filtration
satisfying the usual conditions.
Observe that, by the moment bound in Lemma~\ref{lemmomXY} and the
fact that $X^{(n)}$ and $Y^{(n)}$ are subject to the same SPDE, the
covariation of $X^{(0)}$ and $Y^{(0)}$ satisfies
%
%e5.10 #&#
\begin{eqnarray}
\bigl\langle X^{(0)}(\phi_1),Y^{(0)}(
\phi_2)\bigr\rangle_t=\int_0^t\!\int_\R X^{(0)}(x,s)^{1/2}Y^{(0)}(x,s)^{1/2}
\phi_1(x)\phi_2(x)\,dx\,ds,\nonumber
\\
\eqntext{\phi_1,\phi_2\in\C_c^\infty(\R).}
\end{eqnarray}
By an enlargement of the underlying probability space, we may assume
that for some filtration $(\ms H_t)$ satisfying the usual
conditions, $X^{(0)}$ and $Y^{(0)}$ are adapted to $(\ms H_t)$ and there
exists an $(\ms H_t)$-space--time white noise $\widetilde{W}$
independent of $(X^{(0)},Y^{(0)})$. Let $M^{X^{(0)}}$ and $M^{Y^{(0)}}$
denote the martingale measures of $X^{(0)}$ and $Y^{(0)}$, respectively
[cf. the definition of $M^{X^{(n)}}$ in (\ref{eqMXn})].
Then by the foregoing display, the required space--time white noise can
be chosen to be
\begin{eqnarray*}
W_t(\phi)&\triangleq&\int_0^t\!\int
_\R\1_{(X^{(0)}>0)}(y,s)\frac
{\phi(x)}{X^{(0)}(y,s)^{1/2}}\,dM^{X^{(0)}}(y,s)
\\
&&{}+\int_0^t\!\int_\R
\1_{(X^{(0)}=0, Y^{(0)}>0)}(y,s)\frac{\phi
(x)}{Y^{(0)}(y,s)^{1/2}}\,dM^{Y^{(0)}}(y,s)
\\
&&{}+\int_0^t\!\int_\R
\1_{(X^{(0)}=0, Y^{(0)}=0)}(y,s)\phi(x)\,d\widetilde{W}(y,s),\qquad\phi
\in
\C_c^\infty(\R)
\end{eqnarray*}
[recall the notation $(Z\in\Gamma)$ in (\ref{spacetime})].
The proof is complete.
\end{pf*}

%s6 #&#
\section{Proof of Proposition~\texorpdfstring{\protect\ref{propseptime}}{3.3}}\label{sectauvep}
In this section, we prove Proposition~\ref{propseptime} by verifying
all of the following analogues of (\ref{eqseptimebdd}):
%
%
%e6.1 #&#
\begin{eqnarray}\label{distime}
&& \forall\rho>0\ \exists\delta>0\mbox{ such that}
\nonumber\\[-8pt]\\[-8pt]\nonumber
&&\qquad \sup\biggl\{\mathbb Q^i_{\varepsilon} \bigl(\tau^{i,(j)}\leq
s_i+\delta\bigr);i\in\mathbb N, {\varepsilon}\in\biggl(0,
\frac{1}{8\psi(\1)}\wedge1 \biggr] \biggr\}\leq\rho,
\end{eqnarray}
where $1\leq j\leq3$. The proofs rely on the basic fact that for any
$i\in\mathbb N$ and
${\varepsilon}\in(0,[8\psi(\1)]^{-1}\wedge1 ]$,
$X^i(\1)^{T^{X^i}_1}$ under $\Q^i_{\varepsilon}$ is a $\frac
{1}{4}\BES Q^4(4\psi(\1){\varepsilon})$ started at $s_i$ and stopped
upon hitting $1$ (see the discussion after Proposition~\ref
{propseptime}), and we will work with various couplings of $\frac
{1}{4}\BES Q^4(4z)$. We assume that the couplings are obtained from a
$(\ms H_t)$-standard Brownian motion $B$ for a filtration $(\ms H_t)$
satisfying the usual conditions, the constructions explained in detail
later on.

Recall that we write $\mathbf P_z^1$ for the law of a copy $Z$ of
$\frac{1}{4}\BES Q^4(4z)$. Throughout this section, we do \emph{not}
impose the constraints in Assumption~\ref{ass} on the auxiliary parameters.

%
%
%le6.1 #&#
\begin{lem}\label{lemtime1}
Fix $\eta\in(1,\infty)$, and let $\tau^{i,(1)}$ be the stopping
times defined in Proposition~\ref{propseptime}. Then (\ref
{distime}) holds for $j=1$.
\end{lem}
\begin{pf}
The proof is an application of the lower escape rate of $\BES Q^4(0)$:
%
%
%e6.2 #&#
\begin{equation}
\label{eqBESQesc} \mathbf P^1_0 \bigl(\exists h>0\mbox{
such that }4Z_t\geq t^{\eta}, \forall t\in[0,h] \bigr)=1
\end{equation}
(cf. Theorem 5.4.6 of \cite{KEBM}).

We will need a monotonicity of $\BES Q^4$ in initial values. For this
purpose, we construct all $\frac{1}{4}\BES Q^4(4z)$-processes $Z^z$
with initial values $z\in\R_+$ from the \mbox{$(\ms H_t)$-}standard Brownian
motion. This is implied by the pathwise uniqueness in their stochastic
differential equations (cf. Theorems IX.1.7~and~IX.3.5 of \cite
{RYCMB}), and we can characterize them by
%
%
%e6.3 #&#
\begin{equation}
\label{eqBESQ4} Z^z_t\equiv z+t+\int
_0^t \sqrt{Z^z_s}\,dB_s,
\qquad z\in\R_+.
\end{equation}
In view of the first components in $\tau^{i,(1)}_{\varepsilon}$ (cf.
Proposition~\ref{propseptime}), we consider
\[
\sigma_z\triangleq\inf\biggl\{t\geq0;Z_{T^{Z^z}_1\wedge t}^z<
\frac
{t^\eta}{4} \biggr\},\qquad z\in\biggl[0,{\frac{1}{8}}
\biggr].
\]

Let us bound the distribution function of $\sigma_z\wedge T_1^{Z^z}$.
The comparison theorem of stochastic differential equations (cf.
Theorem IX.3.7 of \cite{RYCMB}) implies that $Z^{z_1}\leq Z^{z_2}$
whenever $0\leq z_1\leq z_2<\infty$. In particular, for any $z\in
(0,\frac{1}{8}]$,
%
%
%e6.4 #&#
\begin{equation}
T^{Z^{1/8}}_1\leq T^{Z^z}_1\leq
T^{Z^0}_1\quad\mbox{and}\quad\sigma_z\geq
\sigma_0\qquad\mbox{a.s.},\label{eqtau2aux}
\end{equation}
where the second inequality follows since
\[
Z_{t\wedge T^{Z^z}_1}^z\geq Z_{t\wedge T^{Z^0}_1}^0\geq
\frac{t^{\eta
}}{4}\qquad\forall t\in[0,\sigma_0].
\]
Hence, by (\ref{eqtau2aux}), we have
\begin{eqnarray*}
\sup_{z\in(0,{1}/{8}]}\P\bigl(\sigma_z\wedge
T^{Z^z}_1\leq\delta\bigr)&\leq&\sup_{z\in(0,{1}/{8}]}
\P(\sigma_z\leq\delta)+\sup_{z\in(0,{1}/{8}]}\P
\bigl(T^{Z^{z}}_1 \leq\delta\bigr)
\\
&\leq&\P(\sigma_0\leq\delta)+\P\bigl(T^{Z^{1/8}}_1
\leq\delta\bigr)\qquad\forall\delta\in(0,\infty).
\end{eqnarray*}
Applying the lower escape rate (\ref{eqBESQesc}) to the right-hand
side of the foregoing inequality shows that
\[
\forall\rho>0\ \exists\delta>0\qquad\mbox{such that }\sup_{z\in(0,{1}/{8}]}\P
\bigl(\sigma_z\wedge T_1^{Z^z}\leq\delta
\bigr)\leq\rho.
\]
Using the foregoing display and the distributional property of $X^i(\1
)^{T^{X^i}_1}$ under $\Q^i_{\varepsilon}$, for $i\in\mathbb N$ and
${\varepsilon}\in(0,[8\psi(\1)]^{-1}\wedge1 ]$ mentioned
above, we have proved our assertion (\ref{distime}) for $j=1$.\vadjust{\goodbreak}
\end{pf}

%
%
%le6.2 #&#
\begin{lem}\label{lemtime2}
Fix $L\in(0,\infty)$ and $\alpha\in(0,\frac{1}{2})$, and let $\tau
^{i,(2)}$ be the stopping times defined in Proposition~\ref
{propseptime}. Then (\ref{distime}) holds for $j=2$.
\end{lem}
\begin{pf}
As in the proof of Lemma~\ref{lemtime1}, we need a grand coupling of
all $\frac{1}{4}\BES Q^4(4z)$, $z\in\R_+$, on the same probability space.
For the first component of $\tau^{i,(2)}$, we need to measure the
modulus of continuity of the martingale part of a $\frac
{1}{4}\BES Q^4$ in terms of its quadratic variation. Hence, it will be
convenient to extract all of the $\frac{1}{4}\BES Q^4(4z)$'s, say
$Z^z$, from a fixed copy $Z$ of $\frac{1}{4}\BES Q^4(0)$, and
we consider
$Z^z_t\equiv Z_{T^Z_{z}+t}$ for $z\in\R_+$,
where
the stopping times $T^Z_{z}$ are finite almost surely by the transience
of $\BES Q^4$ (cf. page~442 of \cite{RYCMB}).
We may further assume that $Z=Z^0$ and is defined by (\ref{eqBESQ4}).
It follows that
%
%
%e6.5 #&#
\begin{equation}
Z^z_t=z+t+\int_{T_{z}^{Z}}^{T_{z}^Z+t}
\sqrt{Z_s}\,dB_s.\label{eqcovZ}
\end{equation}

In this case, the analogues of $\tau^{i,(2)}$ are given by, for $z\in
(0,\frac{1}{8}]$,
%
%
%e6.6 #&#
\begin{eqnarray}\label{ineqsigmavep}
\sigma_z&\triangleq&\inf\biggl\{t\geq0;\bigl\llvert
Z^z_{t\wedge
T_1^{Z^z}}-z-t\bigr\rrvert>L \biggl(\int
_0^{t}Z^z_{s\wedge T_1^{Z^z}}\,ds
\biggr)^\alpha\biggr\}\wedge T_1^{Z^z}
\nonumber
\\
&=&\inf\biggl\{t\geq0;\biggl\llvert\int
_{T_z^Z}^{(T_z^Z+t)\wedge
T^Z_1}\sqrt{Z_s}\,dB_s+
\bigl(t\wedge T_1^{Z^z}-t \bigr)\biggr\rrvert
\\
&&\hspace*{32pt} >L \biggl[ \int_{T_z^Z}^{(T_z^Z+t)\wedge T^Z_1}Z_s\,ds+
\bigl(t\wedge T_1^{Z^z}-t \bigr) \biggr]^\alpha
\biggr\}\wedge T_1^{Z^z},\nonumber
\end{eqnarray}
where the last equality follows from (\ref{eqcovZ}) and
the obvious equality
$T^{Z}_z+T^{Z^z}_1=T^Z_1$.

Let us bound the distribution function of $\sigma_z$.
By the Dambis--Dubins--Schwarz theorem (cf. Theorem V.1.6 of~\cite{RYCMB}), $\sqrt{Z}\bullet B=B'_{\langle Z\rangle}$
for some standard Brownian motion $B'$, where $\langle Z\rangle=\int
_0^\cdot Z_s\,ds$. Also,
\[
0<\int_{T^Z_z}^{(T^Z_z+t)\wedge T_1^Z}Z_s\,ds\leq t\qquad
\mbox{if } t>0,
\]
where the first inequality follows since $\{0\}$ is polar for $\BES Q^4$
(cf. page~442 of \cite{RYCMB}).
Hence, from (\ref{ineqsigmavep}), we deduce that, for any $H,\delta
\in(0,\infty)$,
%
%
%e6.7 #&#
\begin{eqnarray}\label{eqmodBM}
&& \sup_{z\in(0,{1}/{8}]}\P(
\sigma_z\leq\delta)
\nonumber\\[-8pt]\\[-8pt]\nonumber
&&\qquad \leq\P\bigl(T_1^Z>H
\bigr)+\P\biggl(\sup_{\stackrel{\scriptstyle0<\llvert t-s\rrvert \leq
2\delta}{\scriptstyle0\leq s<t\leq H}}\frac
{\llvert B'_t-B'_s\rrvert }{\llvert t-s\rrvert ^\alpha}>L \biggr)+\P
\bigl(T^{Z^{1/8}}_1\leq\delta\bigr)
\end{eqnarray}
[for the third probability, recall the inequalities for hitting times
of $1$ by $\BES Q^4$ in~(\ref{eqtau2aux})].

Let us make the dependence on $\delta$ of the second probability of
(\ref{eqmodBM}) explicit.
For the fixed $\alpha\in(0,\frac{1}{2})$, we pick $\alpha'\in
(0,\frac{1}{2})$ and $p>1$ such that $\alpha<\alpha'<\frac
{p-1}{2p}$. Then applying Chebyshev's inequality to the second term on
the right-hand side of (\ref{eqmodBM}), we get
%
%
%e6.8 #&#
\begin{eqnarray}\label{eqmodBM1}
&& \sup_{z\in(0,{1}/{8}]}\P(
\sigma_z\leq\delta)\nonumber
\\
&&\qquad \leq \P\bigl(T^Z_1>H
\bigr)+\frac{(2\delta)^{2p(\alpha'-\alpha)}}{L^{2p}}\E\biggl[ \biggl
(\sup_{ 0\leq s<t\leq H}
\frac{\llvert B'_t-B'_s\rrvert }{\llvert t-s\rrvert ^{\alpha
'}} \biggr)^{2p} \biggr]
\\
&&\quad\qquad{} +\P\bigl(T^{Z^{1/8}}_1\leq\delta\bigr)
\qquad\forall H,\delta\in(0,\infty),\nonumber
\end{eqnarray}
where
%
%
%e6.9 #&#
\begin{equation}
\label{ineqmodfin} \E\biggl[ \biggl(\sup_{0\leq s<t\leq H}
\frac
{\llvert B'_t-B'_s\rrvert }{\llvert t-s\rrvert ^{\alpha'}} \biggr
)^{2p} \biggr]<\infty
\end{equation}
(cf. the discussion preceding Theorem I.2.2 of \cite{RYCMB} as well
as its Theorem I.2.1).

By the transience of $\BES Q^4$, the first probability on the right-hand
side of~(\ref{eqmodBM1}) can be made as small as possible by choosing
sufficiently large $H$.
Since $ (\sigma_{\psi(\1){\varepsilon}},\P)$ and $
(\tau^{i,(2)},\P_{\varepsilon} )$ have the same distribution
and $\psi(\1){\varepsilon}\leq\frac{1}{8}$, (\ref{eqmodBM1}) and~(\ref{ineqmodfin}) are enough to obtain (\ref{distime}) for $j=2$.
The proof is complete.
\end{pf}

It remains to prove (\ref{distime}) for $j=3$. We need a few
preliminary results.

%
%
%le6.3 #&#
\begin{lem}\label{lemtime3-1}
Fix $i\in\mathbb N$ and ${\varepsilon}\in(0,[8\psi(\1
)]^{-1}\wedge1 ]$. Then
%
%
%e6.10 #&#
\begin{equation}
\label{ineqsupYQ} \E^{\mathbb Q^i_{\varepsilon}} \biggl[ \biggl(\sup
_{r\in[0,R]}
\sum_{j\dvtx s_i<t_j\leq s_i+r}Y^j_{s_i+r}(\1)
\biggr)^p \biggr]<\infty\qquad\forall p,R\in(0,\infty).
\end{equation}
\end{lem}
\begin{pf}
Plainly, it suffices to consider $p>1$. By Lemma~\ref{lemQi}, we have
%
%
%e6.11 #&#
\begin{eqnarray}\label{ineqsupYQ0}
&&\E^{\mathbb Q^i_{\varepsilon}} \biggl[ \biggl(\sup_{r\in[0,R]}\sum
_{j\dvtx s_i<t_j\leq s_i+r}Y^j_{s_i+r}(\1)
\biggr)^p \biggr]
\nonumber
\\
&&\qquad = \frac{1}{\psi(\1){\varepsilon}}\E^{\P_{\varepsilon}} \biggl[X^i_{s_i+R}(
\1)^{T^{X^i}_1} \biggl(\sup_{r\in[0,R]}\sum
_{j\dvtx s_i<t_j\leq s_i+r}Y^j_{s_i+r}(\1)
\biggr)^p \biggr]
\nonumber
\\
&&\qquad \leq\frac{1}{\psi(\1){\varepsilon}}\E^{\P_{\varepsilon}} \biggl[ \biggl
(\sup
_{r\in[0,R]}\sum_{j\dvtx s_i<t_j\leq s_i+r}Y^j_{s_i+r}(
\1 ) \biggr)^p \biggr]\nonumber
\\
&&\qquad \leq\frac{1}{\psi(\1){\varepsilon}}\E^{\P_{\varepsilon}} \biggl[ \biggl
(\sum
_{j\dvtx s_i<t_j\leq s_i+R}\sup_{t\in[t_j, s_i+R]}Y^j_{t}(
\1 ) \biggr)^p \biggr]
\\
&&\qquad \leq\frac{1}{\psi(\1){\varepsilon}}\#\{j;s_i<t_j\leq
s_i+R\} ^{p-1}\nonumber
\\
&&\quad\qquad{}\times \sum_{j\dvtx s_i<t_j\leq s_i+R}
\E^{\P_{\varepsilon}} \Bigl[ \Bigl(\sup_{t\in[t_j, s_i+R]}Y^j_{t}(
\1) \Bigr)^p \Bigr],\nonumber
\end{eqnarray}
where the last inequality follows from H\"older's inequality.
Since each $Y^j(\1)$ under $\mathbb P_{\varepsilon}$ is a Feller
diffusion with initial value $\psi(\1){\varepsilon}$ and started at
$t_j$, the summands on the right-hand side of (\ref{ineqsupYQ0}) are finite.
This gives (\ref{ineqsupYQ}), and the proof is complete.
\end{pf}

Next, we recall the canonical decomposition of $Y^j(\1)$ for $t_j>s_i$
under $\Q^i_{\varepsilon}$ in Lemma~\ref{lemsep1}(2). Recall (\ref
{eqcovarSBM}) and the explicit form (\ref{defIj}) of the finite
variation process $I^j$ of $Y^j(\1)$ under $\Q^i_{\varepsilon}$. Then
by the Cauchy--Schwarz inequality, we deduce that%
%
%e6.12 #&#
%e6.13 #&#
\begin{eqnarray}\label{ineqsumYfv}
\qquad && \sum_{j\dvtx s_i<t_j\leq t} \bigl(\psi(\1){\varepsilon}+I^j_t
\bigr)\nonumber
\\
&&\qquad \leq\psi(\1) {\varepsilon}\#\{j;s_i<t_j\leq t\} +
\int_{s_i}^{t\wedge T^{X^i}_1} \biggl(\frac{\sum_{j\dvtx s_i<t_j\leq
t}Y^j_s(\1)}{X^i_s(\1)}
\biggr)^{1/2}\,ds
\nonumber\\[-8pt]\\[-8pt]
&&\qquad \leq 2\psi(\1) (t-s_i)+\int_{s_i}^{t\wedge T^{X^i}_1}
\biggl(\frac
{\sum_{j\dvtx s_i<t_j\leq s}Y^j_s(\1)}{X^i_s(\1)} \biggr)^{1/2}\,ds\nonumber
\\
\eqntext{\forall t \in[s_i,\infty).}
\end{eqnarray}
Here, the last inequality follows since for $t\geq s_i+\frac
{{\varepsilon}}{2}$,
$s_i+{\varepsilon} (\#\{j;s_i<t_j\leq t\}-\frac{1}{2}
)\leq t$
and the clusters $Y^j$ with $s<t_j\leq t$ have no contributions to
$\sum_{j\dvtx s_i<t_j\leq t}Y^j_s(\1)$.

Also, recall that $M^j$ denotes the martingale part of $Y^j(\1)$ under
$\Q^i_{\varepsilon}$, and the super-Brownian motions $Y^j$ are $\P
_{\varepsilon}$-independent by Theorem~\ref{teoapprox-v2}.
Hence, we deduce from Girsanov's theorem (cf. Theorem VIII.1.4 of \cite{RYCMB}) that
%
%
%e6.14 #&#
%e6.15 #&#
\begin{eqnarray}\label{eqYQV}
\biggl\langle\sum_{j\dvtx s_i<t_j\leq\cdot}M^j
\biggr\rangle_t=\int_{s_i}^t \sum
_{j\dvtx s_i<t_j\leq t}Y^j_s(\1)\,ds =\int
_{s_i}^t \sum_{j\dvtx s_i<t_j\leq s}Y^j_s(
\1)\,ds
\nonumber\\[-8pt]\\[-8pt]
\eqntext{\forall t\in[s_i,\infty),}
\end{eqnarray}
where the omission of the clusters $Y^j$ for $s<t_j\leq t$ follows from
the same reason as in (\ref{ineqsumYfv}).

%
%
%le6.4 #&#
\begin{lem}\label{lemtime3-2}
Fix $i\in\mathbb N$ and ${\varepsilon}\in(0,[8\psi(\1
)]^{-1}\wedge1 ]$. Then
%
%
%e6.16 #&#
\begin{eqnarray}
\label{ineqEinvX} %
\qquad\quad &&\E^{\mathbb Q^i_{\varepsilon}} \biggl[\frac{1}{[X^i_{s_i+r}(\1
)]^{a}};s_i+r
\leq T^{X^i}_1 \biggr]\leq\frac{1}{r^a}
\E^{\mathbf
P^1_0} \biggl[\frac{1}{(Z_1)^a} \biggr]\qquad\forall r,a\in(0,
\infty),
\end{eqnarray}
where
%
%
%e6.17 #&#
\begin{equation}\label{ineqEinvX1}
\E^{\mathbf P^1_0} \biggl[\frac{1}{(Z_1)^a} \biggr
]<\infty\quad\Longleftrightarrow\quad a\in(-\infty,2).
\end{equation}
\end{lem}
\begin{pf}
Recall the grand coupling of $\frac{1}{4}\BES Q^4(4z)$ in the proof of
Lem\-ma~\ref{lemtime1} under which $Z^{z_1}\leq Z^{z_2}$ whenever
$0\leq z_1\leq z_2$.
Then for every $r,a\in(0,\infty)$,
\begin{eqnarray*}
\E^{\mathbb Q^i_{\varepsilon}} \biggl[\frac{1}{[X^i_{s_i+r}(\1
)]^{a}};s_i+r\leq
T^{X^i}_1 \biggr]&\leq&\E^{\mathbf P^1_{\psi(\1
){\varepsilon}}} \biggl[
\frac{1}{(Z_r)^a} \biggr]
\\
&\leq&\E^{\mathbf
P^1_0} \biggl[\frac{1}{(Z_r)^a}
\biggr]
\\
&=& \frac{1}{r^a}\E^{\mathbf P^1_0} \biggl[\frac{1}{(Z_1)^a} \biggr],
\end{eqnarray*}
where the last equality follows from the scaling property of Bessel
squared processes (cf. Proposition XI.1.6 of \cite{RYCMB}).
This gives the bound (\ref{ineqEinvX}). In addition, notice that $Z$
under $\mathbf P^1_0$ has the same distribution as the image of a
$4$-dimensional standard Brownian motion under $x\lmt\llVert x\rrVert
^2$ where $\llVert \cdot\rrVert$ denotes the Euclidean norm, and so
we deduce (\ref
{ineqEinvX1}) by writing out the expectation on its left-hand side as
an elementary integral in polar coordinates.
The proof is complete.
\end{pf}

With Lemma~\ref{lemtime3-2}, we have the following improvement of
(\ref{ineqsupYQ}).

%
%
%le6.5 #&#
\begin{lem}\label{lemtime3-3}
Fix $i\in\mathbb N$ and ${\varepsilon}\in(0,[8\psi(\1
)]^{-1}\wedge1 ]$. Then we have
%
%
%e6.18 #&#
%e6.19 #&#
\begin{eqnarray}\label{ineqY1stmom}
&&\E^{\mathbb Q^i_{\varepsilon}} \biggl[\sum
_{j\dvtx s_i<t_j\leq
s_i+r}Y^j_{s_i+r}(\1) \biggr]\nonumber
\\
&&\qquad \leq\biggl(2
\psi(\1)R+\E^{\mathbf
P^1_0} \biggl[\frac{1}{Z_1} \biggr]^{1/2}2R^{1/2}
\biggr) \exp\biggl( 2\E^{\mathbf P^1_0} \biggl[\frac{1}{Z_1}
\biggr]^{1/2}\sqrt{r} \biggr)
\\
\eqntext{\forall r\in[0,R ], R\in(0,\infty),}
\end{eqnarray}
where $\E^{\mathbf P^1_0}[1/Z_1]<\infty$ by (\ref{ineqEinvX1}).
\end{lem}
\begin{pf}
Recall that the local martingale part of $Y^j(\1)$ under $\Q
^i_{\varepsilon}$ is a true martingale by Lemma~\ref{lemsep1}(2). Hence,
for any $r\in[0,R]$, we obtain from (\ref{ineqsumYfv}) that
%
%
%e6.20 #&#
%e6.21 #&#
\begin{eqnarray}
\qquad&&\E^{\mathbb Q^i_{\varepsilon}} \biggl[\sum_{j\dvtx s_i<t_j\leq
s_i+r}Y^j_{s_i+r}(
\1) \biggr]
\nonumber
\\
&&\qquad \leq2\psi(\1)r+\int_{s_i}^{s_i+r}
\E^{\mathbb Q^i_{\varepsilon
}} \biggl[ \biggl(\frac{\sum_{j\dvtx s_i<t_j\leq s}Y^j_s(\1)}{X^i_s(\1
)} \biggr)^{1/2};s\leq
T^{X^i}_1 \biggr]\,ds\label{ineqYL0}
\\
&&\qquad \leq2\psi(\1)r\nonumber
\\
&&\quad\qquad{}+\int_{s_i}^{s_i+r}
\E^{\mathbb Q^i_{\varepsilon
}} \biggl[\frac{1}{X^i_s(\1)};s\leq T^{X^i}_1
\biggr]^{1/2}\E^{\mathbb
Q^i_{\varepsilon}} \biggl[\sum
_{j\dvtx s_i<t_j\leq s}Y^j_s(\1)
\biggr]^{1/2}\,ds
\nonumber
\\
&&\qquad \leq2\psi(\1)r\nonumber
\\
&&\quad\qquad{} +\int_{s_i}^{s_i+r}
\frac{1}{\sqrt{s-s_i}}\E^{\mathbf P^1_0} \biggl[\frac{1}{Z_1}
\biggr]^{1/2} \biggl(1+\E^{\mathbb
Q^i_{\varepsilon}} \biggl[\sum
_{j\dvtx s_i<t_j\leq s}Y^j_s(\1) \biggr] \biggr)\,ds
\nonumber
\\
&&\qquad \leq \biggl(2\psi(\1)R+\E^{\mathbf P^1_0} \biggl[
\frac{1}{Z_1} \biggr]^{1/2}2R^{1/2} \biggr)
\nonumber\\[-8pt]\label{ineqYL1} \\[-8pt]\nonumber
&&\quad\qquad{} +
\E^{\mathbf P^1_0} \biggl[\frac{1}{Z_1} \biggr]^{1/2}
\int_0^r
\frac{1}{\sqrt{s}} \E^{\mathbb
Q^i_{\varepsilon}} \biggl[\sum_{j\dvtx s_i<t_j\leq s_i+s}Y^j_{s_i+s}(
\1 ) \biggr]\,ds,
\end{eqnarray}
where the third inequality follows from Lemma~\ref{lemtime3-2}.
With the change of variables $s'=\sqrt{s}$, the foregoing inequality
with $r$ replaced by
$r^2$ and $R$ by $R^2$ becomes
%
%e6.22 #&#
\begin{eqnarray}
&& \E^{\mathbb Q^i_{\varepsilon}} \biggl[\sum_{j\dvtx s_i<t_j\leq
s_i+r^2}Y^j_{s_i+r^2}(
\1) \biggr]\nonumber
\\
&&\qquad \leq  \biggl(2\psi(\1)R^2+\E^{\mathbf P^1_0} \biggl[
\frac{1}{Z_1} \biggr]^{1/2}2R \biggr)\nonumber
\\
&&\quad\qquad{} +2\E^{\mathbf P^1_0} \biggl[\frac{1}{Z_1}
\biggr]^{1/2}\int_0^r
\E^{\mathbb Q^i_{\varepsilon}} \biggl[\sum_{j\dvtx s_i<t_j\leq
s_i+(s')^2}Y^j_{s_i+(s')^2}(
\1) \biggr]\,ds'\nonumber
\\
\eqntext{\forall r\in[0,R],}
\end{eqnarray}
so by Lemma~\ref{lemtime3-1} and Gronwall's lemma
%
%e6.23 #&#
\begin{eqnarray}
&&\E^{\mathbb Q^i_{\varepsilon}} \biggl[\sum_{j\dvtx s_i<t_j\leq
s_i+r^2}Y^j_{s_i+r^2}(
\1) \biggr]\nonumber
\\
&&\qquad \leq\biggl(2\psi(\1)R^2+\E^{\mathbf
P^1_0} \biggl[\frac{1}{Z_1} \biggr]^{1/2}2R \biggr)
\exp\biggl( 2\E^{\mathbf P^1_0} \biggl[\frac{1}{Z_1}
\biggr]^{1/2}r \biggr)\nonumber
\\
\eqntext{\forall r\in[0,R].}
\end{eqnarray}
With another change of time scales by $r'=r^2$, the foregoing gives the
desired inequality (\ref{ineqY1stmom}). The proof is complete.
\end{pf}

We are ready to prove (\ref{distime}) for $j=3$.

%
%
%le6.6 #&#
\begin{lem}
Let $\tau^{i,(3)}$ be the stopping times defined in Proposition~\ref
{propseptime}. Then~(\ref{distime}) holds for $j=3$.
\end{lem}
\begin{pf}
Fix $i\in\mathbb N$ and ${\varepsilon}\in(0,[8\psi(\1
)]^{-1}\wedge1 ]$.
It follows from (\ref{ineqsumYfv}) that, for any $R>0$ with
$\frac{1}{3}\geq2\psi(\1)R$,
we have
%
%
%e6.24 #&#
\begin{eqnarray}\label{ineqYcc}
&&\Q^i_{\varepsilon} \biggl(\sup_{r\in[0,R]}\sum
_{j\dvtx s_i<t_j\leq
s_i+r}Y^j_{s_i+r}(\1)> 1 \biggr)
\nonumber
\\
&&\qquad \leq \Q^i_{\varepsilon} \biggl(\int_{s_i}^{(s_i+R)\wedge T^{X^i}_1}
\biggl(\frac{\sum_{j\dvtx s_i<t_j\leq s}Y^j_s(\1)}{X^i_s(\1)} \biggr
)^{1/2}\,ds>\frac{1}{3} \biggr)
\nonumber
\\
&&\quad\qquad{} +\mathbb Q^i_{\varepsilon} \biggl(\sup
_{r\in[0,R]}\biggl\llvert\sum_{j\dvtx s_i<t_j\leq s_i+r}M^j_{s_i+r}
\biggr\rrvert>\frac{1}{3} \biggr)
\\
&&\qquad \leq3\E^{\mathbb Q^i_{\varepsilon}} \biggl[\int_{s_i}^{(s_i+R)\wedge
T^{X^i}_1}
\biggl(\frac{\sum_{j\dvtx s_i<t_j\leq s}Y^j_s(\1)}{X^i_s(\1
)} \biggr)^{1/2}\,ds \biggr]\nonumber
\\
&&\quad\qquad{} +9\sup_{r\in[0,R]}\E^{\mathbb Q^i_{\varepsilon}} \biggl[ \biggl(
\sum_{j\dvtx s_i<t_j\leq s_i+r}M^j_{s_i+r}
\biggr)^2 \biggr],\nonumber
\end{eqnarray}
where the first term of the last inequality follows from Chebyshev's
inequality, and the second term follows by applying Doob's
$L^2$-inequality to the $\Q^i_{\varepsilon}$-martingale $\sum_{j\dvtx
s_i<t_j\leq\cdot}M^j$.\vspace*{2pt}

We claim that the right-hand side of (\ref{ineqYcc}) converges to
zero uniformly in $i\in\mathbb N$ and ${\varepsilon}\in(0,[8\psi
(\1){\varepsilon}]^{-1}\wedge1 ]$ as $R\lra0+$.
Inspecting the arguments from (\ref{ineqYL0}) to (\ref{ineqYL1})
shows that the first term in (\ref{ineqYcc}) satisfies
\begin{eqnarray*}
&&3\E^{\mathbb Q^i_{\varepsilon}} \biggl[\int_{s_i}^{(s_i+R)\wedge
T^{X^i}_1}
\biggl(\frac{\sum_{j\dvtx s_i<t_j\leq s}Y^j_s(\1)}{X^i_s(\1
)} \biggr)^{1/2}\,ds \biggr]\nonumber
\\
&&\qquad \leq 3 \biggl(2\psi(\1)R+\E^{\mathbf P^1_0} \biggl[\frac
{1}{Z_1}
\biggr]^{1/2}2R^{1/2} \biggr)
\\
&&\quad\qquad{} +3\E^{\mathbf P^1_0} \biggl[
\frac{1}{Z_1} \biggr]^{1/2}\int_0^R
\frac{1}{\sqrt{s}} \E^{\mathbb
Q^i_{\varepsilon}} \biggl[\sum_{j\dvtx s_i<t_j\leq s_i+s}Y^j_{s_i+s}(
\1 ) \biggr]\,ds.\nonumber
\end{eqnarray*}
For the second term on the right-hand side of (\ref{ineqYcc}), we use
(\ref{eqYQV}) and obtain
\begin{eqnarray*}
&& 9\sup_{r\in[0,R]}\E^{\mathbb Q^i_{\varepsilon}} \biggl[ \biggl(\sum
_{j\dvtx s_i<t_j\leq s_i+r}M^j_{s_i+r} \biggr)^2
\biggr]
\\
&&\qquad \leq\int_{0}^{R}\E^{\Q^i_{\varepsilon}}
\biggl[\sum_{j\dvtx s_i<t_j\leq s_i+s}Y^j_{s_i+s}(\1 )
\biggr]\,ds.
\end{eqnarray*}
Applying the uniform bound (\ref{ineqY1stmom}) to the right-hand
sides of the last two displays shows the existence of a constant $C\in
(0,\infty)$ depending only on $\psi$ such that
%
%e6.25 #&#
\begin{eqnarray}
\Q^i_{\varepsilon} \biggl(\sup_{r\in[0,R]}\sum
_{j\dvtx s_i<t_j\leq
s_i+r}Y^j_{s_i+r}(\1)> 1 \biggr)
\leq CR^{1/2}\nonumber
\\
\eqntext{\displaystyle \forall R\in\biggl(0,\frac{1}{6\psi(\1)} \biggr], i\in\mathbb N, {
\varepsilon}\in\biggl(0,\frac{1}{8\psi(\1)}\wedge1 \biggr],}
\end{eqnarray}
where the restriction on $R$ follows since $\frac{1}{3}\geq2\psi(\1)R$.
The foregoing inequality proves our claim and is enough for our
assertion of the present lemma.
\end{pf}

%s7 #&#
\section{Some properties of support processes}\label{secsupp}
We study the supports of the immigrant processes $X^i,Y^j$ in this
section. Recall that Assumption~\ref{ass} does not apply to the
present section.

%
%
%pr7.1 #&#
\begin{prop}\label{propsupppre}
There is a constant $C_{\supp}^0\in(0,\infty)$ depending only on the
immigration function $\psi$ and the parameter $\beta\in[\frac
{1}{4},\frac{1}{2})$ such that
%
%
%e7.1 #&#
%e7.2 #&#
\begin{eqnarray}\label{ineqsuppprob}
&&\P_{\varepsilon} \bigl(\sigma^{X^i}_\beta-
s_i\leq r \bigr)+\P_{\varepsilon} \bigl(\sigma^{Y^i}_\beta-
t_i\leq r \bigr)\leq C_{\supp}^0 {\varepsilon}(r
\vee{\varepsilon})
\nonumber\\[-8pt]\\[-8pt]
\eqntext{\forall{\varepsilon},r\in(0,1], i\in\mathbb N.}
\end{eqnarray}
\end{prop}

\begin{pf}
The immigrant processes satisfy the SPDE (\ref{eqmainSPDE0}) with
initial condition taking the form $\psi(\1)J^a_{\varepsilon}$ [recall
that $J^a_{\varepsilon}$ is defined by (\ref{eqJxvep})]. Hence,
Corollary~7.2 of \cite{MMP} applies to the normalized processes
$X^i/\psi(\1)$ and $Y^i/\psi(\1)$ with the parameter $a$ in equation~(7.1)
set to be $\psi(\1)^{-1/2}$. Our assertion follows.
\end{pf}

In the remainder of this section, we consider, under $\Q
^i_{\varepsilon}$, the supports of the immigrant processes $Y^j$
landing by time $s_i+r\in(s_i,\infty)$ and with space--time locations
$(y_j,t_j)$ lying \emph{outside} the rectangle $\mathcal R_{\beta
}^{X^i}(s_i+r)$ defined by (\ref{eqPXi1}).
We start with the immigrants $Y^j$ landing before time $s_i$.

%
%
%pr7.2 #&#
\begin{prop}\label{propoutclus-si}
There exists a constant $C_{\supp}^1\in(0,\infty)$ depending only on
the immigration function $\psi$ such that whenever $\beta\in[\frac
{1}{3},\frac{1}{2})$,
\begin{eqnarray*}
&&\mathbb Q_{\varepsilon}^i \biggl(\mathcal P_\beta^{X^i}(s_i+r)
\cap\biggl(\bigcup_{j\dvtx t_j\leq s_i}\supp\bigl(Y^j
\bigr) \biggr)\neq\varnothing,
\\
&&\hspace*{29pt} \min_{j\dvtx t_j\leq s_i} \bigl(
\sigma_\beta^{Y^j}-t_j \bigr)>3r,
\sigma^{X^i}_\beta-s_i>2r \biggr)
\\
&&\qquad  \leq C_{\supp}^1r^{1/6}\qquad\forall
i\in\mathbb N\mbox{ with $s_i\leq1$}, r\in[s_i,1], {
\varepsilon}\in(0,r].
\end{eqnarray*}
\end{prop}

The proof of Proposition~\ref{propoutclus-si} is similar to the proof
of Lemma 8.4 in \cite{MMP} for $\gamma=1/2$ (note that our notation
$\beta$ is denoted by $\alpha$ there instead), except that in \cite
{MMP} the immigrant processes are subject to i.i.d. space--time white
noises, but in our case they are not. For this reason, we need a
slightly different argument
whenever covariations between the involved immigrants may be nonzero.
Roughly speaking, we will handle the $Y^j$-immigrants which land a bit
``far away'' from the support of $X^i$ in both space and time.
Since these immigrants do not interfere with $X^i$ immediately, we can
apply orthogonal continuation (Lemma~\ref{lemoc}) to $X^i(\1)$ and
then argue as in \cite{MMP} accordingly.

\begin{pf*}{Sketch of proof of Proposition~\ref{propoutclus-si}}
We give the details to handle the $Y^j$-immigrants mentioned above and
sketch the rest of the proof. A complete proof can be found in
Section~3.12 of \cite{C}.

Fix $(\beta,i,r,{\varepsilon})$ as described in the statement of
Proposition~\ref{propoutclus-si}.
We will argue throughout this proof on the event that
%
%
%e7.3 #&#
\begin{equation}
\label{minsigmaY} \min_{j\dvtx t_j\leq s_i} \bigl(\sigma_\beta^{Y^j}-t_j
\bigr)>3r\quad\mbox{and}\quad\sigma^{X^i}_\beta-s_i>2r.
\end{equation}
Let $n_0$ and $n_1$ be nonnegative integers chosen as equation (8.5) in
\cite{MMP}, that is,
%
%
%e7.4 #&#
\begin{equation}
2^{-n_0-1}<r\leq2^{-n_0}\quad\mbox{and}\quad2^{-n_1-1}<{
\varepsilon}\leq2^{-n_1}.\label{eqn0n1}
\end{equation}
Then as in the proof of Lemma 8.4 in \cite{MMP} [cf. equation (8.8) there],
we have
%
%
%e7.5 #&#
\begin{eqnarray}
\label{ineqoutside} %
&& \bigl\{(y_j,t_j);t_j
\in(0,s_i),\mc P_\beta^{X^i}(s_i+r)\cap
\supp\bigl(Y^j\bigr)\neq\varnothing\bigr\}
\nonumber\\[-8pt]\\[-8pt]\nonumber
&&\qquad  \subseteq\bigl[x_i-7\cdot2^{-n_0\beta},x_i+7
\cdot2^{-n_0\beta} \bigr]\times[0,s_i).
\end{eqnarray}
The inclusion in (\ref{ineqoutside}) rules out a number of clusters
$Y^j$ landing before $s_i$ whose space--time supports can intersect $\mc
P^{X^i}_\beta(s_i+r)$ by time $s_i+r$. In the following, we handle the
remaining immigrant processes $Y^j$ for $j\in\mathbb N$ with $t_j<s_i$.

As in the proof of Lemma 8.4 of \cite{MMP}, we classify the clusters
$Y^j$ for $j\in\mathbb N$ satisfying $t_j\in(0,s_i)$
and
$y_j\notin[x_i-7\cdot2^{-n_0\beta},x_i+7\cdot2^{-n_0\beta
} ]$
according to the space--time landing locations $(y_j,t_j)$. Define the
following \emph{random} rectangles
\begin{eqnarray*}
\mathcal R^{\mathtt0}_n&=& \bigl[x_i-7
\cdot2^{-n\beta},x_i+7\cdot2^{-n\beta} \bigr]\times
\bigl[s_i-2^{-n+1},s_i-2^{-n} \bigr],
\\
\mc R^{\mathtt L}_n&=& \bigl[x_i-7
\cdot2^{-n\beta},x_i-7\cdot2^{-(n+1)\beta} \bigr]\times
\bigl[s_i-2^{-n},s_i \bigr],
\\
\mc R^{\mathtt R}_n&=& \bigl[x_i+7
\cdot2^{-(n+1)\beta},x_i+7\cdot2^{-n\beta} \bigr]\times
\bigl[s_i-2^{-n},s_i \bigr],
\end{eqnarray*}
for nonnegative integers $n\geq n_0$, and we group the clusters $Y^j$
according to these rectangles by
\[
Y^{(n),q}\triangleq\sum_{j\dvtx t_j\leq s_i}
\1_{\mathcal
R^q_n}(y_j,t_j)Y^j,\qquad q={\mathtt
L},{\mathtt0},{\mathtt R}, n\geq n_0.
\]
Then as in equation~(8.11) of \cite{MMP}, the probability under
consideration can be bounded as
%
%
%e7.6 #&#
\begin{eqnarray}\label{ineqoutside1}
&&\Q^i_{\varepsilon} \biggl(\mathcal P_\beta^{X^i}(s_i+r)
\cap\biggl(\bigcup_{j\dvtx t_j\leq s_i}\supp\bigl(Y^j
\bigr) \biggr)\neq\varnothing,\nonumber
\\
&&\hspace*{28pt} \min_{j\dvtx t_j\leq s_i} \bigl(
\sigma_\beta^{Y^j}-t_j \bigr)>3r,\sigma
^{X^i}_\beta-s_i>2r \biggr)\nonumber
\\
&&\qquad \leq\mathbb Q^i_{\varepsilon} \Biggl(
\bigcup_{n=n_1+1}^\infty\bigcup
_{q={\mathtt L},{\mathtt0},{\mathtt R}} \bigl\{\mc P_\beta^{X^i}(s_i+r)
\cap\supp\bigl(Y^{(n),q} \bigr)\neq\varnothing\bigr\} \Biggr)
\\
&&\quad\qquad{}+\sum_{n=n_0}^{n_1}\sum
_{q={\mathtt L}, {\mathtt0},{\mathtt R}}\mathbb Q^i_{\varepsilon} \Bigl(\mc
P_\beta^{X^i}(s_i+r)\cap\supp
\bigl(Y^{(n),q} \bigr)\neq\varnothing,\nonumber
\\
&&\hspace*{119pt}\min_{j\dvtx t_j\leq s_i} \bigl(
\sigma_\beta^{Y^j}-t_j \bigr)>3r,
\sigma^{X^i}_\beta-s_i>2r \Bigr).\nonumber
\end{eqnarray}
Recall that the landing locations $x_i$ and $y_i$ have distributions
given by (\ref{eqxiyilaw}).
Then following the arguments between equations (8.11)~and~(8.22) in \cite{MMP}, we deduce that
%
%
%e7.7 #&#
%e7.8 #&#
\begin{eqnarray}
&&\mathbb Q^i_{\varepsilon} \Biggl(\bigcup
_{n=n_1+1}^\infty\bigcup_{q={\mathtt L},{\mathtt0},{\mathtt R}}
\bigl\{\mc P_\beta^{X^i}(s_i+t)\cap\supp
\bigl(Y^{(n),q} \bigr)\neq\varnothing\bigr\} \Biggr) \leq
C_\psi{\varepsilon}^\beta,\label{eqR0n-1}
\\
&&\Q^i_{\varepsilon} \Bigl(\mathcal
P_\beta^{X^i}(s_i+r)\cap\supp
\bigl(Y^{(n),{\mathtt0}} \bigr)\neq\varnothing,\nonumber
\\
&&\hspace*{20pt} \min_{j\dvtx t_j\leq
s_i} \bigl(\sigma_\beta^{Y^j}-t_j \bigr)>3r,
\sigma^{X^i}_\beta-s_i>2r \Bigr)\label{ineqYn0prob}
\\
&&\qquad \leq C_\psi2^{-{n}/{6}} \nonumber
\end{eqnarray}
for some constant $C_\psi$ depending only on the immigration function
$\psi$.

It remains to deal with the summands on the right-hand side of (\ref
{ineqoutside1}) associated with $Y^{(n),\mathtt R}$ for $n_0\leq n\leq
n_1$ (the\vspace*{1pt} probability bounds for $Y^{(n),\mathtt L}$ follow similarly).
In this case, the $Y^j$ summands in $Y^{(n),{\mathtt R}}$ can arrive up to
$s_i-\frac{{\varepsilon}}{2}$, and hence, $Y^{(n),{\mathtt R}}$ can
survive beyond $s_i$ when the covariation between $Y^{(n),\mathtt R}$
and $X^i$ may become nonzero.

Fix $n$ such that $n_0\leq n\leq n_1$. Following the argument from equation (8.24) to equation (8.26) in \cite{MMP}, we deduce that
%
%
%e7.9 #&#
\begin{eqnarray}\label{ineqYnrprob}
&&\mathbb Q^i_{\varepsilon} \Bigl(\mc P^{X^i}_\beta(s_i+r)
\cap\supp\bigl(Y^{(n),\mathtt R} \bigr)\neq\varnothing,\nonumber
\\
&&\hspace*{19pt} \min_{j\dvtx t_j\leq
s_i} \bigl(\sigma^{Y^j}_\beta-t_j \bigr)>3r,
\sigma^{X^i}_\beta-s_i>2r \Bigr)
\\
&&\qquad \leq\Q^i_{\varepsilon} \Bigl(Y^{(n),\mathtt
R}_{s_i+2^{-n}}(
\1)>0,\min_{j\dvtx t_j\leq s_i} \bigl(\sigma^{Y^j}_\beta
-t_j \bigr)>3r,\sigma^{X^i}_\beta-s_i>2r
\Bigr).\nonumber
\end{eqnarray}
We can use a calculation of Feller diffusions to bound the right-hand
side of (\ref{ineqYnrprob}).
Let us start with the inequality:
%
%
%e7.10 #&#
\begin{eqnarray}\label{eqYnr-1}
&&\Q^i_{\varepsilon} \Bigl(Y^{(n),\mathtt R}_{s_i+2^{-n}}(\1)>0,
\min_{j\dvtx t_j\leq s_i} \bigl(\sigma^{Y^j}_\beta-t_j
\bigr)>3r,\sigma^{X^i}_\beta-s_i>2r \Bigr)
\nonumber
\\
&&\qquad\leq\frac{1}{\psi(\1){\varepsilon}}\E^{\P_{\varepsilon}} \bigl
[X^i_{s_i+2^{-n}}(
\1)^{T^{X^i}_1};Y^{(n),\mathtt R}_{s_i+2^{-n}}(\1 )>0,\sigma^{X^i}_\beta-s_i>2^{-n},
\\
&&\hspace*{158pt}
\sigma^{Y^j}_\beta-s_i>2^{-n},\forall
(y_j,t_j)\in\mc R^{\mathtt R}_n
\bigr],\nonumber
\end{eqnarray}
where the restriction for $\sigma^{Y^j}_\beta$ applies since $ r\geq
s_i$ and
$2r\geq2^{-n_0}\geq2^{-n}$.

To evaluate the right-hand side of (\ref{eqYnr-1}), we apply
orthogonal continuation in the following way.
First, note that under $\P_{\varepsilon}$, $X^i(\1)\rest[s_i,\infty
)$ and $Y^{(n),\mathtt R}(\1)\rest[s_i,\infty)$ are $(\G_s)_{s\geq
s_i}$-Feller diffusions with independent starting values by the
independent landing property (\ref{coneilc}).
Define a $(\G_{s})_{s\geq s_i}$-stopping time $\sigma^\perp$ by
\[
\sigma^\perp= \biggl(\sigma^{X^i}_\beta\wedge
\bigwedge_{j\dvtx t_j\leq
s_i}\widehat{\sigma}{}^{Y^j}_\beta
\wedge\bigl(s_i+2^{-n}\bigr) \biggr)\vee s_i,
\]
where the $(\G_{s})_{s\geq s_i}$-stopping times $\widehat{\sigma
}{}^{Y^j}_\beta$ are given by
\begin{eqnarray*}
\widehat{\sigma}{}^{Y^j}_\beta= \cases{ \displaystyle
\sigma^{Y^j}_\beta, &\quad$(y_j,t_j)
\in\mc R^{\mathtt R}_n$,
\vspace*{2pt}\cr
\displaystyle\infty,&\quad
otherwise.}
\end{eqnarray*}
Through $\sigma^\perp$, we control the support propagation of $X^i$
and $Y^j$ for $(y_j,t_j)\in\mc R^{\mathtt R}_n$.
Note that%
%
%e7.11 #&#
\begin{equation}
\label{ineqYjRrn} \mc P_\beta^{X^i}
\bigl(s_i+2^{-n} \bigr)\cap\mc P_\beta^{Y^j}
\bigl(s_i+2^{-n} \bigr) =\varnothing
\end{equation}
for any $j\in\mathbb N$ with $(y_j,t_j)\in\mc R^{\mathtt R}_n$, since the
distance between $\mathcal P^{X^i}(s_i+2^{-n})$ and $\mc P^{Y^j}_\beta
(s_i+2^{-n})$ is given by
\begin{eqnarray*}
&& \bigl(y_j-\bigl(s_i+2^{-n}-t_j
\bigr)^{\beta}-{\varepsilon}^{1/2} \bigr)- \bigl(x_i-2^{-n\beta}-{
\varepsilon}^{1/2} \bigr)
\nonumber
\\
&&\qquad \geq7\cdot2^{-(n+1)\beta}-\bigl[s_i+2^{-n}-
\bigl(s_i-2^{-n}\bigr)\bigr]^{\beta
}-2^{-n\beta}-2
\cdot2^{-n\beta}
\\
&&\qquad  \geq\bigl(7\cdot2^{-\beta}-2^\beta-3 \bigr)\cdot2^{-n\beta}>0,
\end{eqnarray*}
where for the first inequality, we recall (\ref{eqn0n1}) and $\beta
\in[\frac{1}{3},\frac{1}{2})$.
The equality (\ref{ineqYjRrn}) implies $ \langle X^i(\1
),Y^{(n),\mathtt R}(\1) \rangle^{\sigma^\perp}=0$. This allows for
orthogonal continuation of $X^i(\1)$ beyond $\sigma^\perp$ (cf.
Lemma~\ref{lemoc}), and thereby we get a $(\G_s)_{s\geq s_i}$-Feller
diffusion $\widehat{X}{}^i$ independent of $Y^{(n),\mathtt R}(\1)\rest
[s_i,\infty)$ and satisfying $\widehat{X}{}^i=X^i(\1)$ over
$[s_i,\sigma^\perp]$, under~$\P_{\varepsilon}$.\vadjust{\goodbreak}

We use $\widehat{X}{}^i$ to compute the right-hand side of (\ref
{eqYnr-1}) and get
%
%
%e7.12 #&#
\begin{eqnarray}\label{eqYnr-2}
&&\Q^i_{\varepsilon} \Bigl(Y^{(n),\mathtt R}_{s_i+2^{-n}}(\1)>0,
\min_{j\dvtx t_j\leq s_i} \bigl(\sigma^{Y^j}_\beta-t_j
\bigr)>3r,\sigma^{X^i}_\beta-s_i>2r \Bigr)
\nonumber
\\
&&\qquad \leq\frac{1}{\psi(\1){\varepsilon}}\E^{\P
_{\varepsilon}} \bigl[X^i_{s_i+2^{-n}}(
\1)^{T^{X^i}_1};Y^{(n),\mathtt
R}_{s_i+2^{-n}}(\1)>0,\sigma^\perp=s_i+2^{-n}
\bigr]
\nonumber
\\
&&\qquad =\frac{1}{\psi(\1){\varepsilon}}\E^{\P_{\varepsilon
}} \bigl[ \bigl(
\widehat{X}_{s_i+2^{-n}}^{i} \bigr)^{T^{\widehat
{X}{}^i}_1};Y^{(n),\mathtt R}_{s_i+2^{-n}}(
\1)>0,\sigma^\perp=s_i+2^{-n} \bigr]
\\
&&\qquad \leq\frac{1}{\psi(\1){\varepsilon}}\E^{\P
_{\varepsilon}} \bigl[ \bigl(
\widehat{X}_{s_i+2^{-n}}^{i} \bigr)^{T^{\widehat{X}{}^i}_1};Y^{(n),\mathtt
R}_{s_i+2^{-n}}(
\1)>0 \bigr]
\nonumber
\\
&&\qquad =\P_{\varepsilon} \bigl(Y^{(n),\mathtt R}_{s_i+2^{-n}}(\1 )>0
\bigr),\nonumber
\end{eqnarray}
where the last quantity follows from the independence of $\widehat
{X}{}^i$ and $Y^{(n),\mathtt R}$ and the martingale property of
$\widehat
{X}{}^i$, both under $\P_{\varepsilon}$. With an argument similar to
equation~(8.22) in \cite{MMP}, we have
%
%
%e7.13 #&#
\begin{eqnarray}\label{eqYnr-3}
\qquad && \P_{\varepsilon} \bigl(Y^{(n),\mathtt R}_{s_i+2^{-n}}(\1)>0 \bigr)\nonumber
\\
&&\qquad \leq
\P_{\varepsilon} \bigl(Y^{(n),\mathtt R}_{s_i}\geq2^{-n
(1+\beta-{1}/{6} )}
\bigr)+2\cdot2^{n}\cdot2^{-n
(1+\beta-{1}/{6} )}
\nonumber
\\
&&\qquad \leq \psi(\1){\varepsilon}2^{n (1+\beta-{1}/{6} )}\# \bigl\{j\in
\mathbb N;s_i-2^{-n}\leq t_j<s_i\bigr\}
\cdot\frac{\llVert \psi\rrVert_\infty
}{\psi(\1)}14\cdot2^{-n\beta}
\\
&&\quad\qquad {} +2\cdot2^{-n (\beta-{1}/{6} )}
\nonumber
\\
&&\qquad \leq \bigl(28\llVert\psi\rrVert_\infty+2 \bigr)2^{-{n}/{6}},\nonumber
\end{eqnarray}
where the last inequality follows since $\beta\geq\frac{1}{3}$ and
$\#\{j\in\mathbb N;s_i-2^{-n}\leq t_j<s_i\}\leq{\varepsilon}^{-1}2^{-n}+1$.
Then we apply (\ref{eqYnr-2}) and (\ref{eqYnr-3}) to bound the
probability on the right-hand side of (\ref{ineqYnrprob}). By
symmetry, the resulting bound also holds when $Y^{(n),\mathtt R}$ is
replaced by $Y^{(n),\mathtt L}$. We have shown that, by enlarging the
constant $C_\psi$ for (\ref{eqR0n-1}) and (\ref{ineqYn0prob}) if necessary,
%
%
%e7.14 #&#
\begin{eqnarray}\label{ineqYnrprob-4}
&&\mathbb Q^i_{\varepsilon} \Bigl(\mc P^{X^i}_\beta(s_i+r)
\cap\supp\bigl(Y^{(n),q} \bigr)\neq\varnothing,\nonumber
\\
&&\hspace*{18pt} \min_{j\dvtx t_j\leq s_i}
\bigl(\sigma^{Y^j}_\beta-t_j \bigr)>3r, \sigma^{X^i}_\beta-s_i>2r \Bigr)
\\
&&\qquad \leq C_\psi2^{-{n}/{6}},\qquad q={\mathtt L},{\mathtt R},
\forall n_0\leq n\leq n_1. \nonumber
\end{eqnarray}

We apply (\ref{eqR0n-1}), (\ref{ineqYn0prob}) and (\ref
{ineqYnrprob-4}) to (\ref{ineqoutside1}). This gives the conclusion that
\begin{eqnarray*}
&&\mathbb Q_{\varepsilon}^i \biggl(\mc P^{X^i}_\beta(s_i+r)
\cap\biggl(\bigcup_{j\dvtx t_j\leq s_i}\supp\bigl(Y^j
\bigr) \biggr)\neq\varnothing,\nonumber
\\
&&\hspace*{7pt}\qquad  \min_{j\dvtx t_j\leq s_i} \bigl(
\sigma^{Y^j}_\beta-t_j \bigr)>3r,\sigma
^{X^i}_\beta-s_i>2r \biggr)
\\
&&\qquad \leq C_\psi{\varepsilon}^\beta+\sum
_{n=n_0}^{n_1}(3C_\psi)2^{-{n}/{6}}
\\
&&\qquad \leq C_\psi{\varepsilon}^\beta+ \Biggl[ \Biggl(\sum
_{n=0}^{\infty }(3C_\psi)2^{-{n}/{6}}
\Biggr)\cdot2^{{1}/{6}} \Biggr]\cdot2^{{(-n_0-1)}/{6}}.
\end{eqnarray*}
Since $2^{-n_0-1}\leq r$ by (\ref{eqn0n1}) an ${\varepsilon}^\beta
\leq r^\beta\leq r^{{1}/{6}}$, our assertion follows from the
last inequality.
\end{pf*}

Finally, we deal with the simple case where the clusters land after the
landing time $s_i$ of $X^i$ but outside the rectangle $\mathcal
R^{X^i}_\beta(s_i+r)$ defined by (\ref{eqPXi1}).

%
%
%le7.3 #&#
\begin{lem}\label{lemoutclus+si}
Let $r\in(0,\infty)$. Then for any $j\in\mathbb N$ with $t_j\in
(s_i,s_i+r]$ and $\llvert y_j-x_i\rrvert >2 ({\varepsilon
}^{1/2}+r^{\beta
} )$,
$\mc P^{X^i}_\beta(s_i+r)\cap\mc P_\beta^{Y^j}(s_i+r)=\varnothing$.
\end{lem}
\begin{pf}
We only consider the case that $x_i<y_j$, as the other case follows by
symmetry. Note that the distance between
$\mc P^{X^i}_\beta(s_i+r)$ and $\mc P^{Y^j}_\beta(s_i+r)$ is strictly
positive since
\begin{eqnarray*}
&& \bigl(y_j-(s_i+r-t_j)^\beta-{\varepsilon}^{1/2} \bigr)- \bigl(x_i+r^\beta+{
\varepsilon}^{1/2} \bigr)
\\
&&\qquad > 2{\varepsilon}^{1/2}+2r^\beta
-(s_i+r-t_j)^\beta-r^\beta-2{
\varepsilon}^{1/2}
\\
&&\qquad \geq 2r^\beta-r^\beta-r^\beta= 0.
\end{eqnarray*}
It follows that $\mc P^{X^i}_\beta(s_i+r)$
and $\mc P^{Y^j}_\beta(s_i+r)$ are disjoint.
\end{pf}

%s8 #&#
\section{Improved modulus of continuity}\label{secimprovmod}
In this section, we study pointwise modulus of continuity for bounded
Borel-measurable functions satisfying certain Gronwall-type integral
inequalities.

%
%
%th8.1 #&#
\begin{teo}\label{teoimc}
Let $T\in(0,\infty)$.
Suppose that $(f_t)_{t\in[0,T]}$ is a real-valued bounded
Borel-measurable function such that for some $b,a\in(0,\infty)$ and
$C,B,\break A\in\R_+$ which are all independent of $t\in[0,T]$, we have
%
%
%e8.1 #&#
\begin{eqnarray}
\llvert f_t-f_0\rrvert&\leq&C +B t^b+A
\biggl(\int_0^t \llvert f_s
\rrvert \,ds \biggr)^a\qquad\forall t\in[0,T].\label{ineqmodX0}
\end{eqnarray}
Set $\llVert f\rrVert_\infty\triangleq\sup_{s\in[0,T]}\llvert
f_t\rrvert $ and
$D_a\triangleq2^{a-1}\vee1$.
Then for any $n\in\mathbb N$,
%
%
%e8.2 #&#
\begin{eqnarray}
\label{ineqseries}
\llvert f_t-f_0 \rrvert &\leq& C+B t^b\nonumber
\\
&&{}+ (D_a)^{2n} \sum_{j=1}^n
\biggl[\frac{\prod_{k=1}^j (A)^{a^{k-1}}\cdot\prod
_{k=1}^{j-1}(D_{a})^{2(n-k)a^{k}}\cdot(\llvert f_0\rrvert
+C)^{a^j}}{\prod_{k=1}^{j-1}(a_k+1)^{a^{j-k}}} \biggr]t^{a_j}
\\
&&{}+ (D_a)^{2n}\sum_{j=1}^n
\biggl[\frac{\prod_{k=1}^j(A)^{a^{k-1}}\cdot\prod
_{k=1}^{j-1}(D_a)^{2(n-k)a^{k}}\cdot ({B}/{(b+1)} )^{a^j}}{\prod_{k=1}^{j-1}(b_k+1)^{a^{j-k}}} \biggr]t^{b_j}\nonumber
\\
&&{}+ (D_a)^{2n} \biggl[\frac{
(A)^{c_{n}}\cdot\prod_{k=1}^n(D_a)^{2(n-k)a^{k}}\cdot\llVert f\rrVert
_\infty
^{a^{n+1}}}{\prod_{k=1}^{n}(a_k+1)^{a^{n-k+1}}}
\biggr]t^{a_{n+1}}\qquad\forall t\in[0,T],\nonumber
\end{eqnarray}
with the convention that $\prod_{k=1}^0\equiv1$,
where the sequences $\{a_k\}$, $\{b_k\}$, and $\{c_k\}$ are given by
%
%
%e8.3 #&#
\begin{equation}
\label{defabc} a_{k}=\sum_{j=1}^{k}a^{j},
\qquad b_{k}=\sum_{j=1}^{k-1}a^j+(b+1)a^{k}
\quad\mbox{and}\quad c_{k}=\sum_{j=0}^{k}a^j
\end{equation}
with the convention that $\sum_{j=1}^0\equiv0$.
\end{teo}

\begin{pf}
By (\ref{defabc}), we can characterize the sequences $\{a_k\}$, $\{
b_k\}$ and $\{c_k\}$ alternatively by the equations:
%
%
%e8.4 #&#
\begin{eqnarray}\label{eqabc}
a_1&=&a,\qquad a_{k+1}=a(a_k+1),\nonumber
\\
b_1&=&(b+1)a,\qquad b_{k+1}=a(b_k+1),
\\
c_1&=&a+1,\qquad c_{k+1}=a c_k+1.\nonumber
\end{eqnarray}
We use these identifications in the following argument.

We prove the theorem by an induction on $n\in\mathbb N$. We will need the
following elementary inequality: for any $n\in\mathbb N$ with $n\geq2$,
%
%
%e8.5 #&#
\begin{equation}
\label{ineqsubadd} \Biggl(\sum_{j=1}^n
x_j \Biggr)^a\leq(D_a)^{n-1} \Biggl(
\sum_{j=1}^n x_j^a
\Biggr)\qquad\forall x_1,\ldots,x_n\in\R_+.
\end{equation}
Consider (\ref{ineqseries}) for $n=1$. Note that (\ref{ineqmodX0}) implies
%
%
%e8.6 #&#
\begin{equation}
\label{ineqmodX1} \llvert f_t-f_0\rrvert\leq C+B
t^b+A \llVert f\rrVert_\infty^a
t^{a}\qquad\forall t\in[0,T].
\end{equation}
Apply (\ref{ineqmodX1}) to (\ref{ineqmodX0}), and we obtain
\begin{eqnarray*}
&&\llvert f_t-f_0\rrvert
\\
&&\qquad \leq C+B t^b+A \biggl(\llvert f_0\rrvert t+\int
_0^t \bigl(C+ B s^{b}+A\llVert f
\rrVert^a_\infty s^a \bigr)\,ds
\biggr)^a
\\
&&\qquad =  C+B t^b+A \biggl(\bigl(\llvert f_0\rrvert+C\bigr)
t+\frac{B}{b+1}t^{b+1}+\frac{A\llVert f\rrVert
_\infty^a}{a+1}t^{a+1}
\biggr)^a
\\
&&\qquad \leq C+B t^b
\\
&&\quad\qquad{} +A \cdot(D_a )^2 \biggl[\bigl(
\llvert f_0\rrvert+C\bigr)^a t^a + \biggl(
\frac{B}{b+1} \biggr)^a t^{(b+1)a}+ \biggl(
\frac{A\llVert f\rrVert_\infty
^a}{a+1} \biggr)^a t^{(a+1)a} \biggr]
\\
&&\qquad = C+B t^b
\\
&&\quad\qquad{} + (D_a )^2 \biggl[A\bigl(\llvert
f_0\rrvert+C\bigr)^a t^{a_1} +A \biggl(
\frac{B}{b+1} \biggr)^a t^{b_1}+\frac{(A)^{a+1}\llVert f\rrVert_\infty
^{a^2}}{(a+1)^a}t^{a_2}
\biggr],
\end{eqnarray*}
where the second inequality follows from (\ref{ineqsubadd}).
This proves (\ref{ineqseries}) for $n=1$.

Suppose that (\ref{ineqseries}) holds for some $n\in\mathbb N$. Then
for any $t\in[0,T]$, we have
\begin{eqnarray*}
\hspace*{-3pt} &&\int_0^t \llvert f_s\rrvert \,ds
\\
\hspace*{-3pt} &&\qquad \leq\llvert f_0\rrvert t+\int_0^t
\llvert f_s-f_0\rrvert \,ds
\\
\hspace*{-3pt} &&\qquad \leq\bigl(\llvert f_0\rrvert+C\bigr) t+\frac{B}{b+1}t^{b+1}
\\
\hspace*{-3pt} &&\quad\qquad{}+(D_a)^{2n}
\\
\hspace*{-3pt} &&\qquad\qquad{}\times \sum_{j=1}^n
\biggl[\frac{\prod_{k=1}^j(A)^{a^{k-1}}\cdot
\prod_{k=1}^{j-1}(D_a)^{2(n-k)a^{k}}\cdot(\llvert f_0\rrvert
+C)^{a^j}}{\prod_{k=1}^{j-1}(a_k+1)^{a^{j-k}}} \biggr]\frac{1}{(a_j+1)}t^{a_j+1}
\\
\hspace*{-3pt} &&\quad\qquad{}+(D_a)^{2n}
\\
\hspace*{-3pt} &&\qquad\qquad{}\times \sum_{j=1}^n
\biggl[\frac{\prod_{k=1}^j(A)^{a^{k-1}}\cdot\prod
_{k=1}^{j-1}(D_a)^{2(n-k)a^{k}}\cdot ({B}/{(b+1)} )^{a^j}}{\prod_{k=1}^{j-1}(b_k+1)^{a^{j-k}}} \biggr
]\frac{1}{b_j+1}t^{b_j+1}
\\
\hspace*{-3pt} &&\quad\qquad{}+(D_a)^{2n} \biggl[\frac{(A)^{c_{n}}
\cdot
\prod_{k=1}^n(D_a)^{2(n-k)a^{k}}\cdot\llVert f\rrVert_\infty
^{a^{n+1}}}{\prod_{k=1}^{n}(a_k+1)^{a^{n-k+1}}} \biggr]
\frac{1}{a_{n+1}+1}t^{a_{n+1}+1},
\end{eqnarray*}
where the right-hand side is a sum of $2n+3$ many terms (there are $2n$
terms in total under the two summation signs).
Recall the recursive equations in (\ref{eqabc}). Applying (\ref
{ineqmodX0}) and (\ref{ineqsubadd}) for $n$ replaced by $2n+3$ to
the foregoing inequality, we obtain,
for every $t\in[0,T]$,
\begin{eqnarray*}
&&\llvert f_t-f_0\rrvert
\\
&&\qquad \leq C+B t^b+A \cdot(D_a)^{2n+2}\bigl(\llvert
f_0\rrvert+C\bigr)^a t^a+A \cdot
(D_a)^{2n+2} \biggl(\frac{B}{b+1}
\biggr)^a t^{(b+1)a}
\\
&&\quad\qquad{} +A \cdot(D_a)^{2n+2}\cdot(D_a)^{2na}
\\
&&\qquad\qquad{}\times \sum_{j=1}^n \biggl[\frac{\prod_{k=2}^{j+1}(A)^{a^{k-1}}\cdot\prod
_{k=2}^{j}(D_a)^{2[(n+1)-k]a^{k}}\cdot(\llvert f_0\rrvert
+C)^{a^{j+1}}}{\prod_{k=1}^{j-1}(a_k+1)^{a^{(j+1)-k}}}
\biggr]
\\
&&\qquad\qquad{} \times\frac{1}{(a_j+1)^a}t^{a_{j+1}}
\\
&&\quad\qquad{} +A\cdot(D_a)^{2n+2}\cdot(D_a)^{2na}
\\
&&\qquad\qquad{}\times \sum_{j=1}^n \biggl[\frac{\prod_{k=2}^{j+1}(A)^{a^{k-1}}\cdot\prod
_{k=2}^j(D_a)^{2[(n+1)-k]a^{k}}\cdot({B}/{(b+1)}
)^{a^{j+1}}}{\prod_{k=1}^{j-1}(b_k+1)^{a^{(j+1)-k}}}
\biggr]
\\
&&\qquad\qquad{} \times\frac{1}{(b_j+1)^a}t^{b_{j+1}}
\\
&&\quad\qquad{} +A\cdot(D_a)^{2n+2}\cdot(D_a)^{2na}
\biggl[\frac
{(A)^{c_{n}a}\cdot\prod_{k=2}^{n+1}(D_a)^{2[(n+1)-k]a^{k}}\cdot\llVert
f\rrVert
_\infty^{a^{n+2}}}{\prod_{k=1}^{n}(a_k+1)^{a^{(n+1)-k+1}}} \biggr]
\\
&&\qquad\qquad{}\times \frac
{1}{(a_{n+1}+1)^a}t^{a_{n+2}}.
\end{eqnarray*}
Then the rest follows by writing the right-hand side of the foregoing
inequality into the desired form:
\begin{eqnarray*}
&& \llvert f_t-f_0\rrvert
\\[-1pt]
&&\qquad \leq C+B t^b+A
\cdot(D_a)^{2n+2}\bigl(\llvert f_0\rrvert+C
\bigr)^a t^a+A \cdot(D_a)^{2n+2}
\biggl(\frac{B}{b+1} \biggr)^a t^{(b+1)a}
\\[-1pt]
&&\qquad\quad{} + (D_a)^{2n+2}
\\[-1pt]
&&\qquad\qquad{}\times \sum
_{j=1}^n \biggl[\frac{\prod_{k=1}^{j+1}(A)^{a^{k-1}}\cdot\prod
_{k=1}^j(D_a)^{2[(n+1)-k]a^{k}}\cdot(\llvert f_0\rrvert +C)^{a^{j+1}}}{\prod_{k=1}^{j}(a_k+1)^{a^{(j+1)-k}}}\biggr]t^{a_{j+1}}
\\[-1pt]
&&\qquad\quad{} + (D_a)^{2n+2}
\\[-1pt]
&&\qquad\qquad{}\times \sum
_{j=1}^n \biggl[\frac{\prod_{k=1}^{j+1}(A)^{a^{k-1}}\cdot\prod
_{k=1}^j(D_a)^{2[(n+1)-k]a^{k}}\cdot({B}/{(b+1)}
)^{a^{j+1}}}{\prod_{k=1}^{j}(b_k+1)^{a^{(j+1)-k}}} \biggr]
t^{b_{j+1}}
\\[-1pt]
&&\qquad\quad{} + (D_a)^{2n+2} \biggl[\frac{
(A)^{c_{n+1}}\cdot\prod_{k=1}^{n+1}(D_a)^{2[(n+1)-k]a^{k}}\cdot\llVert
f\rrVert
_\infty^{a^{n+2}}}{\prod_{k=1}^{n+1}(a_k+1)^{a^{(n+1)-k+1}}}
\biggr]t^{a_{n+2}}
\\[-1pt]
&&\qquad = C+B t^{b}
\\[-1pt]
&&\qquad\quad{}+ (D_a)^{2n+2}
\\[-1pt]
&&\qquad\qquad{}\times \sum
_{j=1}^{n+1} \biggl[\frac{\prod_{k=1}^{j}(A)^{a^{k-1}}\cdot\prod
_{k=1}^{j-1}(D_a)^{2[(n+1)-k]a^{k}}\cdot(\llvert f_0\rrvert
+C)^{a^{j}}}{\prod_{k=1}^{j-1}(a_k+1)^{a^{j-k}}}
\biggr]t^{a_{j}}
\\[-1pt]
&&\qquad\quad{} + (D_a)^{2n+2}
\\[-1pt]
&&\qquad\qquad{}\times \sum
_{j=1}^{n+1} \biggl[\frac{\prod_{k=1}^{j}(A)^{a^{k-1}}\cdot\prod
_{k=1}^{j-1}(D_a)^{2[(n+1)-k]a^{k}}\cdot({B}/{(b+1)}
)^{a^{j}}}{\prod_{k=1}^{j-1}(b_k+1)^{a^{j-k}}} \biggr]
t^{b_{j}}
\\[-1pt]
&&\qquad\quad{} + (D_a)^{2n+2} \biggl[\frac{
(A)^{c_{n+1}}\cdot\prod_{k=1}^{n+1}(D_a)^{2[(n+1)-k]a^{k}}\cdot\llVert
f\rrVert
_\infty^{a^{n+2}}}{\prod_{k=1}^{n+1}(a_k+1)^{a^{(n+1)-k+1}}}
\biggr]t^{a_{n+2}}.
\end{eqnarray*}
This proves our assertion for (\ref{ineqseries}) when $n$ is replaced
by $n+1$, and the proof is complete by mathematical induction.
\end{pf}

%
%
%co8.2 #&#
\begin{cor}[(Improved modulus of continuity)]\label{corimc}
Let $T\in(0,1]$, $a\in(0,\frac{1}{2})$, and $B,A\in\R_+$. If
$(f_t)_{t\in[0,T]}$ is a real-valued Borel-measurable function
uniformly bounded by $1$ and satisfies (\ref{ineqmodX0}) with $C=0$
and $b=1$,
then for $\xi'\in(0,1)$ and $N'\in\mathbb N$ satisfying
$\sum_{j=1}^{N'}a^j\leq\xi'<\sum_{j=1}^{N'+1}a^j$,
we have
%
%
%e8.7 #&#
%e8.8 #&#
\begin{eqnarray}\label{eqegimcA0}
\llvert f_t-f_0\rrvert
&\leq& \Biggl[ \bigl(A^{{1}/{(1-a)}}+1 \bigr)\sum_{j=1}^{N'}
\llvert f_0\rrvert^{a^j} \Biggr]t^{a}\nonumber
\\
&&{} + \Biggl[B+ \bigl(A^{{1}/{(1-a)}}+1 \bigr)\sum
_{j=1}^{N'} \biggl(\frac{B}{2}\biggr)^{a^j} + A^{{1}/{(1-a)}}+1 \Biggr]t^{\xi'}
\\
\eqntext{\forall t\in[0,T].}
\end{eqnarray}
\end{cor}
\begin{pf}
We simplify the right-hand side of (\ref{ineqseries}) with elementary
algebra, using the present assumptions.
First, $D_a=1$ since $a\in(0,\frac{1}{2})$. Next, since
$\sum_{k=1}^{j}a^{k-1}\leq\frac{1}{1-a}$ for all $j\in\mathbb N$,
we have
%
%
%e8.9 #&#
\begin{equation}
\label{eqcccccbdd} %
\qquad \prod_{k=1}^jA^{a^{k-1}}
\leq A^{{1}/{(1-a)}}+1,\qquad1\leq j\leq n\quad\mbox{and}\quad A^{c_n}
\leq A^{{1}/{(1-a)}}+1.
\end{equation}
Finally, let us handle the exponents $b_j$ in the second sum in (\ref
{ineqseries}).
Using $b=1$ and the definition of $\{b_k\}$ in (\ref{defabc}), we obtain
\begin{eqnarray*}
b_k&=&\frac{a(1-a^{k-1})}{1-a}+2a^{k} =\frac{a-a^k+2a^k-2a^{k+1}}{1-a}
\\
&=& \frac{a+a^k(1-2a)}{1-a}\searrow\frac{a}{1-a}
= \sum_{j=1}^\infty a^j
\end{eqnarray*}
as $k$ tends to infinity since $a\in(0,\frac{1}{2})$. The inequality
(\ref{eqegimcA0}) follows by applying the above observations to (\ref
{ineqseries}).
The proof is complete.
\end{pf}

%\begin{appendix}
%\section{}
%\end{appendix}

%\section*{Acknowledgments}
\section*{Acknowledgements}
The results of the present paper appear in my Ph.D. thesis \cite{C}. I
am grateful to my advisor Professor Ed Perkins for many enlightening
discussions, and this work is based on an important idea of his.
My special thanks go to Leonid Mytnik as the problem considered here
was originally in his joint research program with
Ed Perkins, and he kindly agreed to having it be part of my thesis
research. I wish to thank the anonymous referee for pointing out a
mathematical gap in an early version of this paper as well as giving
several suggestions to improve readability.

%\begin{supplement}[id=suppA]
%\slink[doi]{10.1214/00-AOPXXXXSUPP} %[doi,text={...}] - jei reikia
%suskaldyti doi
%\sfilename{aopXXXX\_supp.pdf}
%\end{supplement}

% imsref loaded by linak, 2014-09-17 14:47:01
%
% imsref loaded by linak, 2015-09-08 14:07:44

\printaddresses
\end{document}